\documentclass[11pt,leqno]{article}
\usepackage[all]{xy}
\usepackage{epsfig}

\usepackage{amsfonts,amssymb,amsthm,amsmath,amscd}
\usepackage[mathscr]{eucal}
\usepackage{mathrsfs}
\usepackage{hyperref}
\usepackage{url}
\usepackage{enumitem}
\usepackage{blindtext}
\usepackage{tensor}
\usepackage{authblk}

      \makeatletter

\usepackage[OT2,T1]{fontenc}
\DeclareSymbolFont{cyrletters}{OT2}{wncyr}{m}{n}
\DeclareMathSymbol{\Sha}{\mathalpha}{cyrletters}{"58}
\usepackage[utf8]{inputenc}

\usepackage{titlesec}
\titleformat{\subsection}[runin]
{\normalfont\large\bfseries}{\thesubsection}{1em}{}
\titleformat{\subsubsection}[runin]
{\normalfont\large\bfseries}{\thesubsubsection}{1em}{}

\theoremstyle{plain}
\newtheorem{thm}[subsubsection]{Theorem}
\newtheorem{thm*}{Theorem}
\newtheorem{cor}[subsubsection]{Corollary}
\newtheorem{lem}[subsubsection]{Lemma}
\newtheorem{prop}[subsubsection]{Proposition}

\theoremstyle{definition}
\newtheorem{defn}[subsubsection]{Definition}
\theoremstyle{remark}
\newtheorem{rem}[subsubsection]{Remark}

\numberwithin{equation}{subsubsection}

\newcommand{\N}{\mathbb N}
\newcommand{\Z}{\mathbb Z}
\newcommand{\Q}{\mathbb Q}
\newcommand{\R}{\mathbb R}
\newcommand{\C}{\mathbb C}
\newcommand{\A}{\mathbb A}

\newcommand{\F}{\mathbb F}
\newcommand{\Fp}{\mathbb{F}_{p}}
\newcommand{\Fpb}{\bar{\mathbb{F}}_{p}}
\newcommand{\Zp}{{\mathbb Z}_p}
\newcommand{\Qp}{{\mathbb Q}_p}
\newcommand{\Ql}{{\mathbb Q}_l}
\newcommand{\Qv}{{\mathbb Q}_v}
\newcommand{\Qb}{\overline{\mathbb Q}}
\newcommand{\Qpb}{\overline{\mathbb{Q}}_p}
\newcommand{\Qlb}{\overline{\mathbb{Q}}_l}
\newcommand{\Qvb}{\overline{\mathbb{Q}}_v}
\newcommand{\Qpnr}{\Q_p^{\mathrm{ur}}} %maximal unramified extn of Qp

\newcommand{\nr}{\mathrm{ur}} %(maxiam) unramified
\newcommand{\Gal}{\mathrm{Gal}} 

\newcommand{\Hom}{\mathrm{Hom}}
\newcommand{\Aut}{\mathrm{Aut}}
\newcommand{\im}{\mathrm{im}}
\newcommand{\Nm}{\mathrm{N}} %norm
\newcommand{\et}{\text{\'et}} %etale
\newcommand{\Int }{\mathrm{Int}} %interior automorphism
\newcommand{\lisom}{\stackrel{\sim}{\longrightarrow}} %long isomorphism
\newcommand{\isom}{\stackrel{\sim}{\rightarrow}} %isomorphism

\newcommand{\tors}{\mathrm{tors}} %torsion
\newcommand{\ra}{\rightarrow}

\newcommand{\hra}{\hookrightarrow}

\newcommand{\Res}{\mathrm{Res}} %Restriction of Scalars
\newcommand{\cO}{\mathcal{O}} %mathcal O
\newcommand{\sS}{\mathscr{S}} %integral canonical model of Shimura variety
\newcommand{\dS}{\mathbb{S}} %Deligne torus
\newcommand{\Sh}{\mathrm{Sh}} %Shimura variety
\newcommand{\sF}{\mathscr{F}} %automorphic sheaf
\newcommand{\cris}{\mathrm{cris}} 
\newcommand{\dR}{\mathrm{dR}}
\newcommand{\Betti}{\mathrm{B}}
\newcommand{\Fr}{\mathrm{Fr}}

\newcommand{\GL}{\mathrm{GL}}

\newcommand{\Gm}{\mathbb{G}_{\mathrm{m}}} 
\newcommand{\ab}{\mathrm{ab}} %maximal abelian quotient
\newcommand{\ad}{\mathrm{ad}} %adjoint
\newcommand{\der}{\mathrm{der}} %derived
\newcommand{\uc}{\mathrm{sc}} %universal covering

\mathchardef\mhyphen="2D

\newcommand{\mfk}{\mathfrak{k}} 
\newcommand{\mfkb}{\bar{\mathfrak{k}}}

\begin{document}

\title{Lefschetz number formula for Shimura varieties of Hodge type}
\author{Dong Uk Lee \thanks{Institute for Industrial and Applied Mathematics, Chungbuk National University, (28644), 1, Chungdae-ro, Seowon-gu, Cheongju-si, Chungcheongbuk-do, Korea \texttt{machhama@gmail.com}}}
\date{}

\maketitle

\begin{abstract}
For any Shimura variety of Hodge type with hyperspecial level at a prime $p$ and automorphic lisse sheaf on it, we prove a formula, conjectured by Kottwitz \cite{Kottwitz90}, for the Lefschetz numbers of Frobenius-twisted Hecke correspondences acting on the compactly supported \'etale cohomology.
Our proof is an adaptation of the arguments of Langlands and Rapoport \cite{LR87} of deriving the Kottwitz's formula from their conjectural description of the set of mod-$p$ points of Shimura variety (Langlands-Rapoport conjecture), which replaces their Galois gerb theoretic arguments by more geometric ones.
We also prove a generalization of Honda-Tate theorem in the context of Shimura varieties and  fix an error in Kisin's work \cite{Kisin17}.
We do not assume that the derived group is simply connected.
\end{abstract}


\

%\begin{keyword} Shimura varieties, Lefschetz number formula \end{keyword}
%\subjclass[2010]{11G18 (primary); 11G25, 11S40 (secondary)}

\tableofcontents

%%%%%%%%%%%%%%%%%%%%%%%%%%%%%%%%%%%%%%%%
%%%%%%%%%%%%%%%%%%%%%%%%%%%%%%%%%%%%%%%%
\section{Introduction}

\subsection{Statement of the main results}

The main result of this work is the following description of the Lefschetz numbers of Frobenius-twisted Hecke correspondences on the compactly supported cohomology of Shimura varieties of Hodge type with hyperspecial level, which was conjectured by Kottwitz \cite[(3.1)]{Kottwitz90} in the case that the derived group is simply connected, and proved in the PEL-type cases of simple Lie type $A$ or $C$ \cite[(19.6)]{Kottwitz92}.

Fix two distinct primes $p$, $l$ of $\Q$.
Let $(G,X)$ be a Shimura datum of Hodge type, $\mathbf{K}_p$ a hyperspecial subgroup of $G(\Qp)$, $\mathbf{K}^p$ a (sufficiently small) compact open subgroup of $G(\A_f^p)$; put $\mathbf{K}=\mathbf{K}_p\mathbf{K}^p$. Let $Sh_{\mathbf{K}}(G,X)$ be (the canonical model over the reflex field $E(G,X)$ of) the associated Shimura variety, and $\sF_{\mathbf{K}}$ a $\lambda$-adic lisse sheaf on it defined by a finite-dimensional algebraic representation $\xi$ of $G$ over a number field $L$ and choice of a place $\lambda$ of $L$ over $l$. For $i\in\N_{\geq0}$, let $H_c^i(Sh_{\mathbf{K}}(G,X)_{\Qb},\sF_{\mathbf{K}})$ be the compactly supported cohomology of $\sF_{\mathbf{K}}$.
Let $\wp$ be a prime of $E=E(G,X)$ above $p$ and $\Phi\in \Gal(\Qpb/E_{\wp})$ be a geometric Frobenius for an embedding $E_{\wp}\hookrightarrow \Qpb$. It is known \cite{Vasiu99}, \cite{Kisin10} that $Sh_{\mathbf{K}}(G,X)$ has a canonical smooth integral model over the ring of integers $\mathcal{O}_{E_{\wp}}$ and $\sF_{\mathbf{K}}$ extends over it.

%%%%%%%%%%%%%%%%%%%%
%%%%%%%%%%%%%%%%%%%%
\begin{thm} \label{thm_intro:Kottwitz_formula} 
Assume that the center $Z(G)$ of $G$ has same ranks over $\Q$ and $\R$.%%
\footnote{In this work, for Hodge-type $(G,X)$, we always assume that $G$ is the smallest algebraic $\Q$-group such that every $h\in X$ factors through $G_{\R}$; then this assumption holds.}
For any $f^p$ in the Hecke algebra $\mathcal{H}(G(\A_f^p)/\!\!/ \mathbf{K}^p)$, there exists $m(f^p)\in\N$, depending on $f^p$, such that for each $m\geq m(f^p)$, the Lefschetz number of $\Phi^m\times f^p$ is given by
\begin{align} \label{eq_intro:Kottwitz_formula}
\sum_{i\geq0}(-1)^i\mathrm{tr}( & \Phi^m\times f^p | H^i_c(Sh_{\mathbf{K}}(G,X)_{\Qb},\sF_{\mathbf{K}})) \\
& = \sum_{(\gamma_0;\gamma,\delta)} c(\gamma_0;\gamma,\delta)\cdot \mathrm{O}_{\gamma}(f^p)\cdot \mathrm{TO}_{\delta}(\phi_p) \cdot \mathrm{tr}\xi(\gamma_0), \nonumber
\end{align}
where the sum is over a set of representatives $(\gamma_0;\gamma,\delta)$ of all equivalence classes of stable Kottwitz triples of level $n$ having trivial Kottwitz invariant and $c(\gamma_0;\gamma,\delta)$ is a constant defined in terms of the triple $(\gamma_0;\gamma,\delta)$ by means of Galois cohomology. Moreover, if the adjoint group $G^{\ad}$ is $\Q$-anisotropic or $f^p$ is the identity, we can take $m(f^p)$ to be $1$.
\end{thm}

See Theorem \ref{thm:Kottwitz_formula:Kisin} for precise definitions of the terms appearing in this formula and more details. We emphasize that we do \emph{not} assume that the derived group $G^{\der}$ of $G$ is simply connected.
%: to the best of the author's knowledge our result is the first result beyond such assumption. 
In those cases that $G^{\der}=G^{\uc}$, $G^{\uc}$ being the simply connected cover of $G^{\der}$, the expression on the right-hand side of (\ref{eq_intro:Kottwitz_formula}) specializes to the one in the original conjecture \cite[(3.1)]{Kottwitz90}.
Regarding this generality, we remark on two points. First, we take the sum only over certain Kottwitz triples which we call \emph{stable} (Definition \ref{defn:stable_Kottwitz_triple}). This notion of stableness has the same meaning as in the stable conjugacy, so is a natural condition in this context; for example,
a Kottwitz triple (e.g. in the sense of \cite{Kottwitz92}) is always stable in our sense if $G^{\der}=G^{\uc}$.
But there is also some subtle point at $p$: our condition is weaker than the condition (previously adopted when $G^{\der}=G^{\uc}$) that the stable norm of $\delta$ equals the stable conjugacy class of $\gamma_0$ (cf. Remark \ref{rem:Kottwitz_triples}).
Secondly, the condition of ``\textit{having trivial Kottwitz invariant}'' should be also taken with a grain of salt: unless $G^{\der}=G^{\uc}$, Kottwitz invariant is not a notion well-defined by a stable Kottwitz triple alone, although its vanishing is so by our definition (cf. \ref{subsubsec:Kottwitz_invariant}).

In another work \cite{Lee18c} (see \cite[Thm.7.2.2]{Lee18a}), we carried out the stabilization of the right-hand side of (\ref{eq_intro:Kottwitz_formula}), which in the case that $G^{\der}=G^{\uc}$ was done by Kottwitz \cite{Kottwitz90}:

%%%%%%%%%%%%%%%%%%%
%%%%%%%%%%%%%%%%%%%%
\begin{thm} \label{thm_intro:Kottwitz_conj}
For every $f^p$ in the Hecke algebra $\mathcal{H}(G(\A_f^p)/\!\!/ \mathbf{K}^p)$, there exists $m(f^p)\in\N$, depending on $f^p$, such that for each $m\geq m(f^p)$, the Lefschetz number of $\Phi^m\times f^p$ is given by
\[\sum_i (-1)^i \mathrm{tr}(\Phi^m\times f^p | H_c^i(Sh_{\mathbf{K}}(G,X)_{\Qb},\sF_{\mathbf{K}}))=\sum_{\underline{H}\in \mathscr{E}_{\mathrm{ell}}(G)} \iota(G,\underline{H}) \mathrm{ST}_{\mathrm{ell}}^{H_1}(f^{H_1}),\]
where the right sum is over a set $\mathscr{E}_{\mathrm{ell}}(G)$ of representatives of the isomorphism classes of elliptic endoscopic data of $G$ with a choice of a $z$-pair $(H_1,\xi_1)$ for each $\underline{H}\in \mathscr{E}_{\mathrm{ell}}(G)$,  and $\mathrm{ST}_{\mathrm{ell}}^{H_1}(f^{H_1})$ is the elliptic part of the geometric side of the stable trace formula for a suitable function $f^{H_1}$ on $H_1(\A)$, $\iota(G,\underline{H})$ is a certain constant depending only on the pair $(G,\underline{H})$.

Moreover, if $G^{\ad}$ is $\Q$-anisotropic or $f^p$ is the identity, we can take $m(f^p)$ to be $1$.
\end{thm}

This result is the first step in the Langlands's program for establishing the celebrated conjecture that the Hasse-Weil zeta function of an arbitrary Shimura variety is a product of automorphic $L$-functions.

%%%%%%%%%%%%%%%%%%%%
%%%%%%%%%%%%%%%%%%%%
\subsection{Strategy and ideas of proof} \label{subsec:sketch_of_proof}

The strategy of our proof of the formula (\ref{eq_intro:Kottwitz_formula}) (and its stabilization) originates from R. Langlands who conceived a systematic approach (after previous successes with examples by Eichler, Shimura, Kuga, Sato, and Ihara) based on Lefschetz-Grothendieck trace formula. Langlands \cite{Langlands72}, \cite{Langlands76}, \cite{Langlands79}, together with others, especially Kottwitz \cite{Kottwitz84b}, developed a set of concrete ideas and techniques. This Langlands-Kottwitz method however needed further refinement and completion. To that end, later Langlands and Rapoport \cite{LR87} gave a conjectural explicit description of the set of $\Fpb$-points of Shimura variety in terms of the theory of Galois gerbs (Langlands-Rapoport conjecture) and provided Galois gerb theoretic arguments of deriving the Kottwitz formula (\ref{eq_intro:Kottwitz_formula}) from that description.

Our proof of the formula (\ref{eq_intro:Kottwitz_formula}) emulate these deduction arguments in \cite{LR87}, but avoids Galois gerb theory.%%
\footnote{except at one place: Theorem \ref{thm:LR_special_admissible_morphism}. But, if one assumes that $G^{\der}=G^{\uc}$, even this theorem is unnecessary.}
In PEL-type cases, it is different from the Kottwitz's arguments \cite{Kottwitz92} (for which Honda-Tate theory \cite{Tate71} is a crucial ingredient).

In the rest of this subsection, we discuss the ideas of our proof in detail.

%%%%%%%%%%%%%%%%%%%%
\subsubsection{Description of the mod-$p$ points of Shimura variety and Langlands-Kottwitz method}  \label{subsubsec:fromLRtoKF}
The Langlands-Kottwitz method is based on certain specific geometric description of the set of $\Fpb$-points of the mod-$p$ reduction of Shimura variety as being equipped with Hecke operators and Frobenius automorphism.
Over the reflex field $E(G,X)$, a Hodge-type Shimura variety $\Sh_{\mathbf{K}_p}(G,X)$ is a moduli space of abelian varieties endowed with a certain prescribed set of additional structures defined by $G$ (consisting of a polarization, a set of absolute Hodge cycles fixed under $G$, and a level structure): we call it $G$-structure, for short. 
One expects that there exists a (smooth) integral model $\sS_{\mathbf{K}_p}(G,X)$ of $\Sh_{\mathbf{K}_p}(G,X)$ over $\cO_{E_{\wp}}$
whose reduction affords a similar moduli description (at least over $\Fpb$), more precisely such that there is a bijection 
\begin{equation} \label{eqn:LRconj-ver1}
\sS_{\mathbf{K}_p}(G,X)(\Fpb)\isom \bigsqcup_{\mathscr{I}}S(\mathscr{I})
\end{equation}
where the disjoint union is over distinct isogeny classes of abelian varieties with $G$-structure and each $S(\mathscr{I})$ parameterizes isomorphism classes in the corresponding isogeny class:
\[S(\mathscr{I})=\varprojlim_{\mathbf{K}^p} I_{\mathscr{I}}(\Q)\backslash X_p(\mathscr{I})\times X^p(\mathscr{I})/\mathbf{K}^p,\]
where $X_p(\mathscr{I})$ and $X^p(\mathscr{I})$ consist respectively of the isogenies of $p$-power order and prime-to-$p$ order (say, leaving from a fixed member in the isogeny class $\mathscr{I}$) preserving $G$-structure, and $I_{\mathscr{I}}(\Q)$ is the automorphism group (in the isogeny category) of the chosen abelian variety with $G$-structure, which thus acts naturally on $X_p(\mathscr{I})$ and $X^p(\mathscr{I})$.
Each of the sets $S(\mathscr{I})$ carries compatible actions of the group $G(\A_f^p)$ (such that $S(\mathscr{I})/\mathbf{K}^p=I_{\mathscr{I}}(\Q)\backslash X_p(\mathscr{I})\times X^p(\mathscr{I})/\mathbf{K}^p$) and the Frobenius automorphism $\Phi$ (an element of $\Gal(\Fpb/\kappa(\wp))$), and the bijection (\ref{eqn:LRconj-ver1}) should be compatible with these actions.
This geometric description is well-known in PEL-type cases due to the existence of a fine moduli scheme over mixed characteristic base, and for Hodge types (or even for abelian types) is now proved by Kisin \cite{Kisin17}.

The Langlands-Kottwitz method requires \emph{group-theoretic} parametrizations of both the isogeny classes $S(\mathscr{I})$ and the index set $\{\mathscr{I}\}$ of isogeny classes. 
First, regarding parametrization of isogeny classes, it is immediate from the geometric contents of (\ref{eqn:LRconj-ver1}) that to each isogeny class $\mathscr{I}$, one can attach a set of group-theoretic data:
\begin{itemize}  [noitemsep]
\item[(i)] a connected reductive $\Q$-group $I_{\mathscr{I}}$;
\item[(ii)] an embedding of $\A_f^p$-group schemes $i^p:(I_{\mathscr{I}})_{\A_f^p} \hookrightarrow G_{\A_f^p}$;
\item[(iii)] an element of $b\in G(\mfk)$ together with an embedding of $\Qp$-group schemes $i_p:(I_{\mathscr{I}})_{\Qp} \hookrightarrow J_b$ ($\mfk:=W(\F)[1/p]$, $J_b$ is the $\sigma$-centralizer of $b$ (\ref{eq:J_b})),
\end{itemize}
which is canonical up to conjugation under $g=(g_p,g^p)\in G(\mfk)\times G(\A_f^p)$: $(b\sigma,i_p,i^p)\sim g(b\sigma,i_p,i^p)g^{-1}$ (here, $b\sigma$ is an element of the semi-direct product $G(\mfk)\rtimes \langle\sigma\rangle$ wtih $\langle\sigma\rangle$ being the free abelian group generated by $\sigma$, and the conjugation of $b\sigma$ is taken inside this group).
With a fixed choice of these data, one obtains a (partial) group-theoretic description of $S(\mathscr{I})$. The sets $X^p(\mathscr{I})$, $X_p(\mathscr{I})$ can be identified respectively with $G(\A_f^p)$ and (the underlying set of) a suitable affine Deligne-Lusztig variety $X(\{\mu_X^{-1}\},b)_{\mathbf{K}_p}$:
\[ X^p(\mathscr{I})=G(\A_f^p),\quad X_p(\mathscr{I})=X(\{\mu_X^{-1}\},b)_{\mathbf{K}_p},\]
in such a way that $I_{\mathscr{I}}(\Q)$ acts on $X_p(\mathscr{I})\times X^p(\mathscr{I})$ diagonally via $i_p\times i^p$ while $\Phi$ acts on $X_p(\mathscr{I})$ by $(b\sigma)^{[\kappa(\wp):\Fp]}$. 

At this point, the first step of the deduction arguments (as explained in \cite{Kottwitz84b}) runs as follows; for simplicity, here we only consider the trivial correspondence case (so aim at a formula of the local zeta function of $\sS_{\mathbf{K}}(G,X)_{\kappa(\wp)}$) and also assume that the derived group $G^{\der}$ of $G$ is simply connected (so, in the rest of this subsection \ref{subsec:sketch_of_proof}, we ignore the issue of stable conjugacy versus (geometric) conjugacy for Kottwitz triples).
An elementary, formal argument (ibid., $\S$1.4) shows that for each isogeny class $\mathscr{I}$, the corresponding fixed-point set $(S(\mathscr{I})/\mathbf{K}^p)^{\Phi^m=\mathrm{Id}}$ breaks up further as a disjoint union of subsets $S(\mathscr{I},\epsilon)_{\mathbf{K}}$ indexed by (the equivalence classes of) pairs $(\mathscr{I},\epsilon)$ (which we call \textit{K-pairs}), where $\epsilon\in I_{\mathscr{I}}(\Q)$ (it is to be considered as (the inverse of) the relative $p^n$-Frobenius endomorphism of $\mathscr{I}$ when $\mathscr{I}$ is defined over $\F_{p^n}, \ n=m[\kappa(\wp):\Fp]$): 
\begin{equation*}
\sS_{\mathbf{K}}(G,X)(\F_{q^m})\isom \bigsqcup_{[\mathscr{I}]}(S(\mathscr{I})/\mathbf{K}^p)^{\Phi^m=\mathrm{id}} = \bigsqcup_{[\mathscr{I},\epsilon]} S(\mathscr{I},\epsilon)_{\mathbf{K}}.
\end{equation*}
Here, we restrict ourselves only to the K-pairs, called \emph{admissible}, satisfying certain natural condition that is necessary (but, not sufficient in general!) for the corresponding set $S(\mathscr{I},\epsilon)_{\mathbf{K}}$ to be non-empty (Definition \ref{defn:admissible_pair2}); this is inspired from an idea of Langlands and Rapoport \cite{LR87}.
Then, with any admissible K-pair $(\mathscr{I},\epsilon)$, one expects to be able to associate a triple of group elements (\textit{Kottwitz triple})
\[(\gamma_0;\gamma=(\gamma_l)_{l\neq p},\delta)\in G(\Q)\times G(\A_f^p)\times G(L_n)\] 
such that $\gamma$ and $\delta$ are determined directly from $\epsilon$ and the data (ii), (iii) by an explicit recipe:
\begin{equation} \label{eq:AP->K-triple}
 (\gamma,\delta):=(i^p(\epsilon),cb\sigma(c^{-1})),
 \end{equation}
where $c$ is any element of $G(\mfk)$ satisfying $c(i_p(\epsilon)^{-1}(b\sigma)^n)c^{-1}=\sigma^n$ (the existence of such $c$ will be guaranteed by the admissibility condition),
and that $\gamma_0$ is an (a priori, arbitrary) $\R$-elliptic, rational element which is conjugate to both $\gamma_l$ (under $G(\Qlb)$) and $\Nm_n\delta:=\delta\sigma(\delta)\cdots\sigma^{n-1}(\delta)$ (under $G(\Qpb)$), where $L_n$ is the unramified extension of $\Qp$ of degree $n$ and $\sigma$ is the Frobenius automorphism of $L_n$: so, the problem here is to find such rational $\gamma_0$. 
When $\epsilon$ is (the inverse of)  the relative $p^n$-Frobenius endomorphism of an abelian variety in an isogeny class $\mathscr{I}$ defined over $\F_{p^n}$, $\gamma_l$ is its realization on \'etale homology group and $\delta$ is the absolute Frobenius acting on crystalline homology group.
Kottwitz triples are considered up to the equivalence relation that $\gamma_0$ is defined up to geometric(=stable) conjugacy, $\gamma$ up to $G(\A_f^p)$-conjugacy, and $\delta$ up to $\sigma$-conjugacy under $G(L_n)$.
For a fixed $\epsilon\in I_{\mathscr{I}}(\Q)$, equivalent data $(b,i_p,i^p)$ produce equivalent Kottwitz triples.

To proceed with deduction, we have to assume some further properties on the data $(b,i_p,i^p)$: 
for each admissible $(\mathscr{I},\epsilon)$, $i^p$ and $i_p$ must induce isomorphisms
\begin{equation} \label{eq:Tate_thm_SV} \tag{\textrm{Ta}}
 i^p:(I_{\mathscr{I},\epsilon})_{\A_f^p}\isom G_{i^p(\epsilon)},\quad i_p:(I_{\mathscr{I},\epsilon})_{\Qp}\isom Z_{J_b}(i_p(\epsilon))
\end{equation}
where $I_{\mathscr{I},\epsilon}$ is the centralizer of $\epsilon$ in $I_{\mathscr{I}}$ (and $G_{i^p(\epsilon)}:=Z_{G}(i^p(\epsilon))$, $Z_{J_b}(i_p(\epsilon))$ are similar centralizers); this property is a version in the context of Shimura varieties of the Tate's theorem on endomorphisms of abelian varieties over finite fields.
It then follows by another easy argument \cite[$\S$1.4-$\S$1.5]{Kottwitz84b} that for any sufficiently small $\mathbf{K}^p$, the cardinality of the set $S(\mathscr{I},\epsilon)_{\mathbf{K}}$ is expressed only in terms of the associated Kottwitz triple, as 
\begin{equation} \label{eq:|S(sI,epsilon)_K|}
|S(\mathscr{I},\epsilon)_{\mathbf{K}}|=\mathrm{vol}(G_{\gamma_0}(\Q)\backslash G_{\gamma_0}(\A_f))\cdot O_{\gamma}(f^p)\cdot TO_{\delta}(\phi_p), 
\end{equation}
where $G_{\gamma_0}$ is the centralizer of $\gamma_0$ in $G$, $O_{\gamma}(f^p)$ and $TO_{\delta}(\mathscr{I}_p)$ are the orbital and twisted orbital integrals of some functions $f^p$, $\phi_p$ (which are determined respectively by $\mathbf{K}^p$ and the pair $(\mathbf{K}_p,X)$ alone).  
Thus, the final formula for $|\sS_{\mathbf{K}}(G,X)(\F_{q^m})|$ (for any sufficiently small $\mathbf{K}^p$) takes the form of a sum, indexed by (the equivalence classes of) Kottwitz triples $(\gamma_0;\gamma,\delta)$, of a product of quantities that can be defined purely group theoretically:
\begin{equation} \label{eqn:formulra_for_number_of_pts}
|\sS_{\mathbf{K}}(G,X)(\F_{q^m})|=\sum_{(\gamma_0;\gamma,\delta)}\iota(\gamma_0;\gamma,\delta)\cdot \mathrm{vol}(G_{\gamma_0}(\Q)\backslash G_{\gamma_0}(\A_f))\cdot O_{\gamma}(f^p)\cdot TO_{\delta}(\phi_p),
\end{equation}
where $\iota(\gamma_0;\gamma,\delta)$ is \emph{by definition} the number of equivalence classes of admissible K-pairs giving rise to a fixed Kottwitz triple $(\gamma_0;\gamma,\delta)$. 

With the explicit cohomological expression $\iota(\gamma_0;\gamma,\delta)=|\ker[\ker^1(\Q,G_{\gamma_0})\rightarrow \ker^1(\Q,G)]|$, this is the formula conjectured (and proved in PEL-type cases) by Kottwitz \cite[(3.1)]{Kottwitz90}, except that here the sum is only over the \emph{effective} Kottwitz triples, namely those Kottwitz triples arising via the recipe (\ref{eq:AP->K-triple}) from some admissible K-pair $(\mathscr{I},\epsilon)$, while in the original formula, one takes \emph{all} Kottwitz triples with trivial Kottwitz invariant (it is known that every effective Kottwitz triple has trivial Kottwitz invariant).
We bring the reader's attention to this usuage of the terminology \emph{effectivity} and the possible confusion that a Kottwitz triple $(\gamma_0;\gamma,\delta)$ which is effective in this sense may not appear as a summation index ``effectively'' in the sum (\ref{eqn:formulra_for_number_of_pts}), i.e. the corresponding summand could be zero (this is because admissibility is only a \emph{necessary} condition for $S(\mathscr{I},\epsilon)_{\mathbf{K}}$ to be non-empty).

Therefore, when one assumes two things: existence of a suitable association $(\mathscr{I},\epsilon)\rightsquigarrow (\gamma_0;\gamma,\delta)$ (especially, of $\gamma_0$) and validity of the isomorphisms (\ref{eq:Tate_thm_SV}),
in order to reconcile the formula (\ref{eqn:formulra_for_number_of_pts}) with the Kottwitz's formula (\ref{eq_intro:Kottwitz_formula}) (in the case $G^{\der}=G^{\uc}$), it suffices to establish the following two statements:
\begin{align} 
\label{eq:E} \tag{E} &\begin{minipage}{.92\textwidth} \emph{Effectivity criterion of Kottwitz triple}: a Kottwitz triple with trivial Kottwitz invariant is effective (in the sense just defined) if the corresponding summand $O_{\gamma}(f^p)\cdot TO_{\delta}(\phi_p)$ is non-zero. \end{minipage} \\ 
\label{eq:C} \tag{C} &\begin{minipage}{.92\textwidth} for effective Kottwitz triple $(\gamma_0;\gamma,\delta)$, one has $\iota(\gamma_0;\gamma,\delta)=|\ker[\ker^1(\Q,G_{\gamma_0})\rightarrow \ker^1(\Q,G)]|$.
\end{minipage}
\end{align}

In fact, the two assumptions: the existence of the association $(\mathscr{I},\epsilon)\rightsquigarrow (\gamma_0;\gamma,\delta)$ and the property (\ref{eq:Tate_thm_SV}), are relatively easy to establish in PEL-type cases (cf. \cite{Kottwitz92}), and for general Hodge-types follow without much difficulty from Kisin's results \cite{Kisin17} (see (\ref{eq:K-triple_attached_to_admissible_K-pair}) and Corollary \ref{cor:Tate_thm2} in this article). But, the two statements (\ref{eq:E}) and (\ref{eq:C}) are much harder and subtle problems.

%%%%%%%%%%%%%%%%%%%%
\subsubsection{The works of Langlands and Rapoport \cite{LR87}  and of Kisin \cite{Kisin17}}

To prove statements like (\ref{eq:E}) and (\ref{eq:C}), naturally one looks for some group-theoretic parameterization of the index set of isogeny classes. Indeed, Langlands \cite{Langlands76}, \cite{Langlands79} already came up with the idea of parametrizing the isogeny classes $\{\mathscr{I}\}$ over $\Fpb$ in terms of a special point lifting (a point in) $\mathscr{I}$ and a relative Frobenius endomorphism attached to $\mathscr{I}$ (relative to any fixed finite field over which $\mathscr{I}$ is defined), more group-theoretically, in terms of a pair $(h,\epsilon\in T(\Q))$ for some special Shimura subdatum $(T,h)$ and a rational element $\epsilon$ of $T$ which satisfy certain condition (see Definition \ref{defn:Frobenius_pair}). But, the equivalence relation that Langlands put on these ``Frobenius pairs'' was defined by everywhere local equivalences, so it could not distinguish two different isogeny classes which are locally isomorphic everywhere (which amounts to having the same associated equivalence classes of Kottwitz triples). Although this notion of Frobenius pair and the very idea of studying isogeny classes in terms of their CM liftings have been and are also here of fundamental importance, this defect (i.e. weak equivalence relation) made the notion inadequate for the purpose of proving (\ref{eq:E}) and (\ref{eq:C}), so it needed some refinement.

Then, Langlands and Rapoport \cite{LR87} proposed a conjectural description of $\sS_{\mathbf{K}_p}(G,X)(\Fpb)$ resolving such problem, a description which is \emph{precise} enough to allow one to derive the Kottwitz formula (\ref{eq_intro:Kottwitz_formula}) from it. The desired precision was provided by Galois-gerb theory, most importantly by a Galois-gerb theoretic version of (\ref{eqn:LRconj-ver1}) (i.e. Langlands-Rapoport conjecture), where each isogeny class $\mathscr{I}$ and the corresponding set $S(\mathscr{I})$ are replaced respectively by a certain object in Galois-gerb theory called ``admissible morphism'' $\phi$ and a set with similar structure $S(\phi)=I_{\phi}(\Q)\backslash X_p(\phi)\times X^p(\phi)$ (cf. \cite{LR87}, \cite{Lee18a}).
We remark that such description is expected from the belief that a large class of Shimura varieties (including abelian types satisfying a mild condition) are to be moduli spaces of Grothendieck motives with additional structures, together with the general fact that Galois gerb theory is a group theoretic language describing Tannakian categories, in general non-neutral, such as the Tannakian cateogry of Grothendieck motives over $\Fpb$ (cf. \cite{DM82}, \cite{Milne92}, \cite{Milne94}).

Meanwhile, recently Kisin \cite{Kisin17} obtained a description of $\sS_{\mathbf{K}_p}(G,X)(\Fpb)$ which is formulated in Galois-gerb language and quite similar to Langlands-Rapoport conjecture. But, in his version of the conjecture (i.e. ibid., Theorem 0.3), the action of the group $I_{\phi}(\Q)$ on $X_p(\phi)\times X^p(\phi)$ (for admissible morphism $\phi$) is not the \emph{natural} one specified in the original conjecture.
Unfortunately, for such imprecise description, the Langlands-Rapoport arguments in \cite{LR87} for establishing (Galois-gerb theoretic versions of) (\ref{eq:E}) and (\ref{eq:C}) do not work by themselves.%%
\footnote{In fact, this is an issue which was well-known to Langlands and as such was the very motivation for the work of Langlands and Rapoport: see their comment  \cite[p.116, line+15]{LR87}. Recently, Kisin-Shin-Zhu announced a proof of Theorem \ref{thm_intro:Kottwitz_formula} for general abelian-type cases (\url{https://arxiv.org/abs/2110.05381}) which seemed to get around to this difficulty.}

%%%%%%%%%%%%%%%%%%%%
\subsubsection{Our proof of (\ref{eq:E}) and (\ref{eq:C})}
For this reason, we cannot simply invoke the deduction arguments of \cite{LR87} which depend on the Langlands-Rapoport conjecture in full.
Instead, we translate their Galois-gerb theoretic arguments into arguments in more familiar notions in algebraic geometry (abelian varieties, torsors, etc) and algebraic groups (reductive groups over global and local fields and their Galois cohomology). 
We use mostly standard results in these areas and some specific geometric (but not any Galois-gerb theoretic) results in \cite{Kisin17} (see the introduction to Section \ref{sec:Pf_Kottwitz_formula} for an exact list of these geometric results).

But, depending on the coarser geometric description (\ref{eqn:LRconj-ver1}) of $\sS_{\mathbf{K}_p}(G,X)(\Fpb)$ afforded by \cite{Kisin17} rather than the much more precise, Galois-gerb theoretic description in \cite{LR87}, such imitation in return entails some new problems which are often subtle.
For example, the key method in \cite{LR87} proving the Galois-gerb theoretic statements corresponding to (\ref{eq:E}) and (\ref{eq:C}) is that of ``twisting of an admissible morphism/pair'' by certain cohomology classes. For admissible K-pairs $(\mathscr{I},\epsilon)$, adapting a similar technique in \cite{Kisin17}, we seek to find (\ref{eq:E}) or to count (\ref{eq:C}) appropriate cohomology classes, the twists by which of a relevant $(\mathscr{I},\epsilon)$ produce a given Kottwitz triple $(\gamma_0;\gamma,\delta)$ by the recipe (\ref{eq:AP->K-triple}). One novel problem we have to take care of, whereas it has natural solution in Galois gerb theory, is that in twisting an admissible K-pair $(\mathscr{I},\epsilon)$, one also needs to keep track of its effect on the group-theoretic data $(b,i_p,i^p)$ that is used in attaching the triple $(\gamma_0;\gamma,\delta)$: in fact, there is one more such datum of vital importance which is necessary for producing $\gamma_0$ out of $\epsilon$ (a maximal torus $T$ of $I_{\mathscr{I}}$ containing $\epsilon$).
A related point which is critical in our faithful imitation of \cite{LR87} is that inner twistings of algebraic groups should be constructed not just group-theoretically (i.e. from cohomology classes only), but via isomorphisms of underlying motives (i.e. abelian varieties) or vector spaces.
See Remark \ref{rem:issue_with_twisting_method} and Remark \ref{rem:Kisin17's_error_Prop.4.4.13} for illustrations of these points. 

The geometric inputs from \cite{Kisin17} are mostly generalizations in Hodge-type setting of well-known results in PEL-type cases, except one strong result which is new even in PEL-type cases and as such allows for our approach to the formula different from the one by Kottwitz \cite{Kottwitz92}, especially our solution to  (\ref{eq:E}) and (\ref{eq:C}): \cite[Cor.2.2.5]{Kisin17}.
We have dubbed it ``strong CM lifting theorem'' to distinguish it from the weaker version \cite[Thm.2.2.3]{Kisin17} (traditionally called ``CM lifting theorem'', cf. \cite{CCO14}): the latter theorem simply says that every isogeny class $\mathscr{I}$ in a Shimura variety of Hodge type over $\Fpb$ contains a point which lifts to a CM point, while the former one makes the more sophisticated claim that such liftings can be constructed by any cocharacter satisfying certain conditions, in which case this cocharacter becomes the Hodge cocharacter of the lifted CM point. It is known \cite[5.11]{LR87} that \emph{every} maximal torus $T$ of $I_{\mathscr{I}}$ admits such a cocharacter into it, so embeds in the Mumford-Tate group of the lifting CM abelian variety, that is, $T$ always embeds in $G$ (albeit in a non-canonical way). 
In particular, this embedding provides a wanted auxiliary means of attaching a rational component $\gamma_0\in G(\Q)$ (in a Kottwitz triple) to any admissible pair $(\mathscr{I},\epsilon)$ (or to any $\F_{p^n}$-rational point in $\sS_{\mathbf{K}_p}$). 
But we need a fine control on such embeddings and we provide results to that effect.

To prove (\ref{eq:E}), we need the fact (Theorem \ref{thm:LR-Satz5.21} (2)) that non-vanishing of $TO_{\delta}(\phi_p)$ implies a condition introduced by Langlands and Rapoport (condition $\ast(\epsilon)$); although this fact is an important (missing) ingredient in the original arguments of \cite{LR87}, it is a purely group theoretic statement about algebraic groups over local fields and is independent of Galois-gerb theory (cf. Remark \ref{rem:criterion_for_ast(epsilon)}).

%%%%%%%%%%%%%%%%%%%%
\subsection{Further results: Honda-Tate theorem for Hodge-type Shimura varieties}
Our direct geometric approach and methods of proof give further geometric results which in the Galois-gerb theoretic setting of \cite{LR87} are attainable only with the complete description of the Langlands-Rapoport conjecture.

When granting the results leading to (\ref{eq:|S(sI,epsilon)_K|}), the effectivity criterion (\ref{eq:E}) becomes actually a stronger ``geometric effectivity criterion'' which tells when a level-$n$ Kottwitz triple $(\gamma_0;\gamma,\delta)$ with trivial Kottwitz invariant arises from an $\F_{p^n}$-rational point in $\sS_{\mathbf{K}}$. In PEL-type cases, this geometric effectivity criterion was a key part in the Kottwitz's proof \cite{Kottwitz92} of the formula (\ref{thm_intro:Kottwitz_formula}) who derived it from Honda-Tate theorem \cite{Tate71}, while we deduce the criterion (\ref{eq:E}) first and as its consequence what can be regarded as Honda-Tate theorem in the context of Shimura varieties (Corollary \ref{cor:geom_effectivity_of_K-triple}):  
an $\R$-elliptic, stable conjugacy class of an element $\gamma_0$ of $G(\Q)$ is geometrically effective (i.e. is the class of an element $\gamma_0$ attached to some $\F_{p^n}$-point of $\sS_{\mathbf{K}}$) if and only if some $G(\A_f^p)$-conjugate of $\gamma_0$ lies in $\mathbf{K}^p$ and there exists $\delta\in G(L_n)$ such that $\gamma_0$ is conjugate to $\Nm_n\delta$ and $\mathrm{TO}_{\delta}(\phi_p)\neq0$.

%%%%%%%%%%%%%%%%%%%%
\subsection{Comments on our previous works \cite{Lee18a}, \cite{Lee18c}, and Kisin's \cite{Kisin17}}
This is a minor revision of (the first part of) our previous work with the same title (submitted in 2018, July) \cite{Lee18c}, which is itself a revision of another earlier work (\url{http://arxiv.org/abs/1801.03057v2}) \cite{Lee18a}.
%: the stabilization theorem (Theorem \ref{thm_intro:Kottwitz_conj}) is included only in \cite{Lee18a}, \cite{Lee18c}, but not here, and also from \cite{Lee18c} the level subgroup is restricted to hyperspecial subgroup from more general special parahoric subgroups. In addition to these, the major differences among these works are as follows: 
The proof of \cite{Lee18a} is Galois-gerb theoretic, following \cite{LR87} closely, whereas in \cite{Lee18c} we replaced the only argument in \cite{Lee18a} (proof of  (\ref{eq:E})) that used Galois gerb theory by another one which does not, thereby making the entire work completely free of Galois gerb theory. 

However, \cite{Lee18a} contained two errors which also slipped in to \cite{Lee18c}, and the current work fixes these errors.
The first error is \cite[Prop.5.2.7]{Lee18a}=\cite[Prop.3.24]{Lee18c}. This proposition is likely to be false in full generality although it is visibly true for the groups whose derived groups are simply connected and further can be shown to hold in much more generality. It was used to attach \emph{stable} Kottwitz triples, not just Kottwitz triples, to (equivalence classes of) admissible (K-)pairs. In the original Galois-gerb theory setting \cite{Lee18a}, this problem (and certain related minor issues) can be solved quite easily (by using only ``special'' admissible pairs among their equivalence classes), so the error causes no trouble. But in the proof explained here (and in \cite{Lee18c}) which is based on Kisin's (incomplete) description (\ref{eqn:LRconj-ver1}) and only imitates \cite{LR87}, such simple solution is not feasible, and here we solve it by proving some nontrivial results (Theorem \ref{thm:stable_isogeny_diagram}, Corollary \ref{cor:stable_conjugacy_from_nice_tori}:
the former theorem is the only part in this work that draws on Galois-gerb theory results in \cite{LR87}.%%
\footnote{More precisely, all the necessary Galois-gerb theoretic results from \cite{LR87} are summarized in Theorem \ref{thm:LR_special_admissible_morphism}. The statement and proof of this theorem will be given in plain group-theoretic terms though. Moreover in the proof, no actual understanding of Galois gerb theory will be required as one just quotes results in \cite{LR87}.})

The second error already appears in Kisin's work: the proof of \cite{Kisin17}, Proposition 4.4.13 is wrong (see Remark \ref{rem:Kisin17's_error_Prop.4.4.13} for an explanation). In \cite[Thm.6.2.11]{Lee18a} and \cite[Thm.4.19]{Lee18c}, we had claimed to have found a correct proof of it (in fact, some generalization of it), but they also suffered from a gap. Here, we fill this gap (Theorem \ref{thm:equiv_K-triples}). 

Finally, we remark that except for these two errors, the group-theoretic ingredients in the current work as well as in \cite{Lee18c} are all taken from \cite{Lee18a}. Besides, the results in \cite{Lee18a} hold for more general special parahoric level subgroups. See Remark \ref{rem:Zhou} for possible generalizations of the main results in such setting and their potential applications to semisimple local zeta functions of Shimura varieties at bad reductions.

This article is organized as follows. 

Section 2 is devoted to group theoretic preparations. 
In subsections 2.1 to 2.3, we review Kottwitz and Newton maps, some results in the $\sigma$-stabilizer groups $J_b$, $G_{\delta\theta}$ (also denoted by $G_{\delta\sigma}$), and fix notations.
In 2.4 and 2.5, we recall the notion of (stable) Kottwitz triple and Kottwitz invariant, discussing certain subtle points arising from working in the general setup that the derived group is not necessarily simply connected. 
In 2.6 (which contains the only (essentially) new results in the second section), we study the important condition \ref{itm:ast(epsilon)} introduced by Langlands and Rapoport. In particular, we prove the statement that for a (stable) Kottwitz triple with trivial Kottwitz invariant, the non-vanishing of $TO_{\delta}(\phi_p)$ implies condition \ref{itm:ast(epsilon)}, which is a key fact in our proof of the effectivity criterion of stable Kottwitz triple.
In 2.7, we prove some Galois cohomological results.
Except in 2.6, the discussions in this section apply to general groups which are not necessarily unramified.

In Section 3, we prove Theorem \ref{thm_intro:Kottwitz_formula} from the theorems and corollaries in Subsection 3.2 which are structured as follows: 
(\ref{eq:Tate_thm_SV}) is proved in Corollary \ref{cor:Tate_thm2}, and (\ref{eq:E}) and (\ref{eq:C}) in Corollary \ref{cor:LR-Satz5.25}. We also need a result (Corollary \ref{cor:stable_conjugacy_from_nice_tori}) ensuring that any admissible K-pair produces a stable Kottwitz triple, not just a Kottwitz triple, in general (i.e. when $G^{\uc}\neq G^{\der}$). These corollaries are deduced by simple purely group-theoretic arguments from five theorems and various lemmas and propositions in the same subsection. The theorems are also all group-theoretic statements about properties of group-theoretic objects (automorphism groups, Kottwitz triples) attached to geometric objects (abelian varieties and isogeny classes) and describing some geometric constructions (notably, twisting), while the lemmas and propositions are results mainly in Galois cohomology theory of reductive groups over global and local fields. 
The proof of Theorem \ref{thm_intro:Kottwitz_formula} is completed in the last subsection (Theorem \ref{thm:Kottwitz_formula:Kisin}) by the same argument as used in \cite{Kottwitz92}. 

In Appendix \ref{sec:existence_of_aniostropic_and_unramified_torus}, we present a proof of a certain result in Bruhat-Tits theory which is perhaps well-known (as it is already mentioned in \cite[p.172]{LR87}). One novelty of this work is the generalization of numerous well-known arguments and statements beyond the setup that $G^{\der}=G^{\uc}$ which are sometimes nontrivial (e.g. the statement of the formula (\ref{eq_intro:Kottwitz_formula}) itself). For that, it is necessary to work with abelianized cohomology groups which are cohomology groups of complexes of tori. In appendix \ref{sec:abelianization_complex}, we collect basic facts about complexes of tori attached to connected reductive groups.

Our sign convention is the same as that of \cite{Kottwitz92}, so opposite to that of \cite{LR87} (and \cite{Kisin17}); the difference is due to our working with homology groups (of abelian varieties) rather than cohomology groups.

%%%%%%%%%%%%%%%%%%%%
\textbf{Acknowledgement}
This work was supported by the National Research Foundation of Korea(NRF) grant funded by the Korea government (NRF-2019R1I1A1A01062321).
The author would like to thank M. Rapoport and C.-L. Chai for their interests in this work and encouragement.

%%%%%%%%%%%%%%%%%%%%
\textbf{Notations}

Throughout this paper, $\Qb$ denotes the algebraic closure of $\Q$ inside $\C$. 

%By our convention, a reductive group is not necessarily connected, even though it is so in most of the cases. When we discuss a non-necessarily connected group (e.g. the centralizer of an element in a (connected) reductive group), we make it explicit.

For a connected reductive group $G$ over a field, we let $G^{\uc}$ be the universal covering of its derived group $G^{\der}$, and for a (linear algebraic) group $G$, $Z(G)$, and $G^{\ad}$ denote its center, and the adjoint group $G/Z(G)$, respectively. 

For a group $I$ and an $I$-module $A$, we let $A_I$ denote the quotient group of $I$-coinvariants: $A_I=A/\langle ia-a\ |\ i\in I, a\in A\rangle$. For an element $a\in A$, we write $\underline{a}$ for the image of $a$ in $A_I$. In case of need for distinction, sometimes we write $\underline{a}_A$.

For a finitely generated abelian group $A$, we denote by $A_{\mathrm{tors}}$ its subgroup of torsion elements. For a locally compact abelian group $A$, we let $X^{\ast}(A):=\Hom_{\mathrm{cont}}(A,\C^{\times})$ (continuous character group) and $A^D:=\Hom_{\mathrm{cont}}(A,S^1)$ (Pontryagin dual).
For a (commutative) algebraic group $A$ over a field $F$, $X_{\ast}(A):=\Hom_{\mathrm{alg}}(\Gm,A)$, $X^{\ast}(A):=\Hom_{\mathrm{alg}}(A,\Gm)$. So, for a diagonalizable $\C$-group $A$, we have $\pi_0(A)^D=X^{\ast}(A)_{\mathrm{tors}}$ (with the embeddings $\Q/\Z\subset \R/Z=S^1\subset \C^{\times}$ understood).

In this article, the german letter $\mfk$ denotes the completion of the maximal unramified extension (in a fixed algebraic closure $\Qpb$) of $\Qp$, and for $n\in\N$, $L_n$ will denote $\mathrm{Frac}(W(\F_{p^n}))$.

%%%%%%%%%%%%%%%%%%%%%%%%%%%%%%%%%%%%%%%%
%%%%%%%%%%%%%%%%%%%%%%%%%%%%%%%%%%%%%%%%

\section{Group theoretic preparations}

%%%%%%%%%%%%%%%%%%%%
%%%%%%%%%%%%%%%%%%%%
\subsection{Kottwitz maps and Newton map}

In this section, we briefly recall the definitions of the Kottwitz maps and the Newton map. We refer to \cite{Kottwitz97}, \cite{Kottwitz85}, \cite{RR96}, and references therein for further details.

%%%%%%%%%%%%%%%%%%%%
\subsubsection{The Kottwitz maps $w_G$, $v_G$, $\kappa_{G}$} \label{subsubsec:Kottwitz_hom}

Let $L$ be a complete discrete valued field with algebraically closed residue field and set $I:=\Gal(\overline{L}/L)$. 
For any connected reductive group $G$ over $L$, Kottwitz \cite[$\S$7]{Kottwitz97} constructs a group homomorphism
\[w_G:G(L)\rightarrow X^{\ast}(Z(\widehat{G})^{I})=\pi_1(G)_{I}.\]
Here, $\widehat{G}$ denotes the Langlands dual group of $G$, $\pi_1(G)=X_{\ast}(T)/\Sigma_{\alpha\in R^{\ast}}\Z\alpha^{\vee}$ is the fundamental group of $G$ (\`a la Borovoi) (i.e. the quotient of $X_{\ast}(T)$ for a maximal torus $T$ over $F$ of $G$ by the coroot lattice), and $\pi_1(G)_I$ is the (quotient) group of coinvariants of the $I$-module $\pi_1(G)$. This map $w_G$ is sometimes denoted by $\widetilde{\kappa}_G$, e.g. in \cite{Rapoport05}. When $G^{\der}$ is simply connected (so that $\pi_1(G)=X_{\ast}(G^{\ab})$ for $G^{\ab}=G/G^{\der}$), $w_G$ factors through $G^{\ab}$: $w_G=w_{G^{\ab}}\circ p_G$, where $p_G:G\rightarrow G^{\ab}$ is the natural projection \cite[7.4]{Kottwitz97}.

There is also a homomorphism 
\[v_G:G(L)\rightarrow \Hom(X_{\ast}(Z(\widehat{G}))^I,\Z)\] 
sending $g\in G(L)$ to the homomorphism $\chi\mapsto \mathrm{val}(\chi(g))$ from $X_{\ast}(Z(\widehat{G}))^I=\Hom_L(G,\Gm)$ to $\Z$, where $\mathrm{val}$ is the usual valuation on $L$, normalized so that uniformizing elements have valuation $1$. It is clear from this definition that $v_G=v_{G^{\ab}}\circ p_G$ for \emph{any} $G$ (i.e. not necessarily having the property $G^{\der}=G^{\uc}$).

There is the relation: 
\[v_G=q_G\circ w_G,\] 
where $q_G$ is the natural surjective map 
\[q_G:X^{\ast}(Z(\widehat{G})^{I})=X^{\ast}(Z(\widehat{G}))_I \rightarrow \Hom(X_{\ast}(Z(\widehat{G}))^I,\Z).\]
The kernel of $q_G$ is the torsion subgroup of $X^{\ast}(Z(\widehat{G}))_I$, i.e. $\Hom(X_{\ast}(Z(\widehat{G}))^I,\Z)\cong \pi_1(G)_I/\text{torsions}$; in particular, $q_G$ is an isomorphism if the coinvariant group $X^{\ast}(Z(\widehat{G}))_I$ is free (e.g. the $I$-module $X^{\ast}(Z(\widehat{G}))$ is trivial or more generally \textit{induced}, i.e. has a $\Z$-basis permuted by $I$).

For example, when $G$ is a torus $T$, we have $\langle \chi,w_T(t)\rangle= \mathrm{val}(\chi(t))$ for $ t\in T(L)$, $\chi\in X^{\ast}(T)^I$, where $\langle\ ,\ \rangle$ is the canonical pairing between $X^{\ast}(T)^I$ and $X_{\ast}(T)_I$.

Now suppose that $G$ is defined over a local field $F$, i.e. a finite extension of $\Qp$ (in a fixed algebraic closure $\Qpb$) with residue field $\F_q$. Let $L$%%
\footnote{In this case that the residue field is $\Fpb$, we will write $\mfk$ for $L$ more often.}
be the completion of the maximal unramified extension $F^{\nr}$ of $F$ in $\Qpb$ 
and let $\sigma$ denote the Frobenius automorphism on $L$ which fixes $F$ and induces $x\mapsto x^q$ on the residue field of $L$ ($\cong\Fpb$). In this situation, the maps $v_{G_L}$, $w_{G_L}$ each induce notable maps.

First, as $w_{G_L}$ (and $v_{G_L}$ too) commutes with the action of $\Gal(F^{\nr}/F)$, by taking $H^{0}(\Gal(F^{\nr}/F),-)$ on both sides of $w_{G_L}$, we obtain a homomorphism
\[\lambda_G:G(F)\rightarrow X^{\ast}(Z(\widehat{G})^I)^{\langle\sigma\rangle},\]
where $I\cong\Gal(\overline{F}/F^{\nr})$. This map is introduced in \cite[$\S$3]{Kottwitz84b} (cf. \cite[7.7]{Kottwitz97}) when $G$ is \emph{unramified} over $F$ (in which case the canonical action of $I$ on $Z(\widehat{G})$ is trivial) and used in \cite{LR87} (with the same notation) under the additional assumption $G^{\der}=G^{\uc}$ so that $w_{G_L}=v_{G_L}$. We remark that in our general setup that $G_{\Qp}$ is not necessarily simply connected, to achieve what $\lambda_G$ did in \cite{LR87}, we use $v_G$, or $w_G$ depending on the situation. 

Next, let $B(G)$ denote the set of $\sigma$-conjugacy classes:
\[B(G):=G(L)/\stackrel{\sigma}{\sim},\]
where two elements $b_1$, $b_2$ of $G(L)$ are said to be \textit{$\sigma$-conjugate}, denoted $b_1\stackrel{\sigma}{\sim} b_2$, if there exists $g\in G(L)$ such that $b_2=gb_1\sigma(g)^{-1}$. 
Then, $w_{G_L}$ induces a map of sets 
\begin{equation} \label{eq:kappa_G}
\kappa_G:B(G)\rightarrow X^{\ast}(Z(\widehat{G})^{\Gamma_F})=\pi_1(G)_{\Gamma_F}: \kappa_G([b])=\underline{w_{G_L}(b)}_F.
\end{equation}
Here, for $b\in G(L)$, $[b]=[b]_G$ denotes its $\sigma$-conjugacy class, and for $x\in \pi(G)_{I}$, $\underline{x}_F$ denotes its image under the natural quotient map $\pi(G)_{I}\rightarrow \pi(G)_{\Gamma_F}$. For further details, see \cite[7.5]{Kottwitz97}.

All these maps are functorial in $G$ (i.e. for group homomorphisms).

%%%%%%%%%%%%%%%%%%%%
\subsubsection{The Newton map $\nu_G$} \label{subsubsec:Newton_map}

Let $\mathbb{D}$ denote the protorus $\varprojlim\Gm$ with the character group $\Q=\varinjlim\Z$.
For an algebraic group $G$ over a $p$-adic local field $F$, we put
\[\mathcal{N}(G):=(\Hom_{L}(\mathbb{D},G)/\Int(G(L)))^{\sigma}\]
(the subset of $\sigma$-invariants in the set of $G(L)$-conjugacy classes of $L$-rational quasi-cocharacters into $G_L$). We will use the notation $\overline{\nu}$ for the the conjugacy class of $\nu\in\Hom_{L}(\mathbb{D},G)$. 

For every $b\in G(L)$, Kottwitz \cite[$\S$4.3]{Kottwitz85} constructs an element $\nu=\nu_G(b)=\nu_b\in\Hom_L(\mathbb{D},G)$%%
\footnote{We interchangeably write $\nu_G(b)$ or $\nu_b$.}
uniquely characterized by the property that there are an integer $s>0$, an element $c\in G(L)$ and a uniformizing element $\pi$ of $F$ such that:
\begin{itemize} \addtolength{\itemsep}{-4pt}
\item[(i)] $s\nu\in\Hom_L(\Gm,G)$.
\item[(ii)] $\Int (c)\circ s\nu$ is defined over the fixed field of $\sigma^s$ in $L$.
\item[(iii)] $c\cdot (b\sigma)^s\cdot c^{-1}=c\cdot (s\nu)(\pi)\cdot c^{-1}\cdot \sigma^{s}$.
\end{itemize}
In (iii), the product (and the equality as well) take place in the semi-direct product group $G(L)\rtimes\langle\sigma\rangle$. We call $\nu_b$ the \textit{Newton homomorphism} attached to $b\in G(L)$

When $G$ is a torus $T$, $\nu_b=\mathrm{av}\circ w_{T_L}(b)$, where $\mathrm{av}:X_{\ast}(T)_I\rightarrow X_{\ast}(T)_{\Q}^{\Gamma_F}$ is ``the average map'' $X_{\ast}(T)_I\rightarrow X_{\ast}(T)_{\Gamma_F}\rightarrow X_{\ast}(T)_{\Q}^{\Gamma_F}$ sending $\underline{\mu}\ (\mu\in X_{\ast}(T))$ to $|\Gamma_F\cdot\mu|^{-1} \sum_{\mu'\in \Gamma_F\cdot\mu}\mu'$ (cf. \cite[Thm. 1.15, (iii)]{RR96}). Hence, it follows that if $T$ is split by a finite Galois extension $K\supset F$, for $b\in T(L)$, $[K:F]\nu_b\in X_{\ast}(T)$ and that $\langle \chi,\nu_b\rangle= \mathrm{val}(\chi(b))$ (especially $\in\Z$) for every $F$-rational character $\chi$ of $T$.

The map $b\mapsto \nu_b$ has the following properties.
\begin{itemize}\addtolength{\itemsep}{-4pt}
\item[(a)] $\nu_{\sigma(b)}=\sigma(\nu_b)$.
\item[(b)] $gb\sigma(g)^{-1}\mapsto \Int (g)\circ \nu,\ g\in G(L)$.
\item[(c)] $\nu_b=\Int (b)\circ\sigma(\nu_b)$.
\end{itemize}
It follows from (b) and (c) that $\nu_G:G(L)\rightarrow \Hom_L(\mathbb{D},G)$ gives rise to a map $\overline{\nu}_G:B(G)\rightarrow \mathcal{N}(G)$, which we call the \textit{Newton map}. This can be also regarded as a functor from the category of connected reductive groups to the category of sets (endowed with partial orders defined as below):
\begin{equation} \label{eq:nubar}
\overline{\nu}:B(\cdot)\rightarrow \mathcal{N}(\cdot)\ ;\ \overline{\nu}_{G}([b])=\overline{\nu}_{b},\quad b\in[b].
\end{equation}

%%%%%%%%%%%%%%%%%%%%
\subsubsection{The set $B(G,\{\mu\})$} \label{subsubsec:B(G,{mu})}  
For a connected reductive group $G$ over an arbitrary (i.e. not necessarily $p$-adic) field $F$, let $\mathcal{BR}(G)=(X^{\ast},R^{\ast},X_{\ast},R_{\ast},\Delta)$ be the based root datum of $G$: we may take $X^{\ast}=X^{\ast}(T)$, $X_{\ast}=X_{\ast}(T)$ for a maximal $F$-torus $T$ of $G$ and $R^{\ast}\subset X^{\ast}(T)$, $R_{\ast}\subset X_{\ast}(T)$ are respectively the roots and the coroots for the pair $(G,T)$ with a choice of basis $\Delta$ of $R^{\ast}$ (whose choice corresponds to that of a Borel subgroup $B\subset G_{\overline{F}}$ containing $T_{\overline{F}}$).
Let $\overline{C}\subset (X_{\ast})_{\Q}$ denote the closed Weyl chamber associated with the root base $\Delta$. It comes with a canonical action of $\Gamma_F:=\Gal(\overline{F}/F)$ on $\overline{C}$. 

For a cocharacter $\mu\in\Hom_{\overline{F}}(\Gm,G)$ lying in $\overline{C}$, we set 
\begin{equation} \label{eq:mubar}
\overline{\mu}:=|\Gamma_F\cdot\mu|^{-1} \sum_{\mu'\in\Gamma_F\cdot\mu}\mu'\quad\in\overline{C}.
\end{equation}
Here, the orbit $\Gamma_F\cdot\mu$ is obtained using the canonical Galois action on $\overline{C}$. Once a Weyl chamber $\overline{C}$ (equivalently, a Borel subgroup $B$ or a root base $\Delta$) is chosen, $\overline{\mu}$ depends only on the $G(\overline{F})$-conjugacy class $\{\mu\}$ of $\mu$.

%%%%%%%%%%%%%%%
Suppose that $\mu\in X_{\ast}(T)\cap\overline{C}$. 
As $X_{\ast}(T)=X^{\ast}(\widehat{T})$ for the dual torus $\widehat{T}$ of $T$, regarded as a character on $\widehat{T}$, we can restrict $\mu$ to the subgroup $Z(\widehat{G})^{\Gamma_F}$ of $\widehat{T}$, obtaining an element 
\begin{equation} \label{eqn:mu_natural}
\mu^{\natural}\in X^{\ast}(Z(\widehat{G})^{\Gamma_F})=\pi_1(G)_{\Gamma_F}.
\end{equation}
Again, $\mu^{\natural}$ depends only on the $G(\overline{F})$-conjugacy class $\{\mu\}$ of $\mu$.
Alternatively, $\mu^{\natural}$ equals the image (sometimes, also denoted by $\underline{\mu}$) of $\mu\in X_{\ast}(T)$ under the canonical map $X_{\ast}(T)\rightarrow \pi_1(G)_{\Gamma_F}$.

%%%%%%%%%%%%%%%%%%%%
Again, let us return to a $p$-adic field $F$. We fix a closed Weyl chamber $\overline{C}$ (equiv. a Borel subgroup $B$ over $\overline{F}$). 
Suppose given a $G(\overline{F})$-conjugacy class $\{\mu\}$ of cocharacters into $G_{\overline{F}}$.
Let $\mu$ be the representative of $\{\mu\}$ in $\overline{C}$; so we have $\overline{\mu}=\overline{\mu}(G,\{\mu\})\in\overline{C}$ and $\mu^{\natural}\in X^{\ast}(Z(\widehat{G})^{\Gamma_F})$. 
We define a finite subset $B(G,\{\mu\})$ of $B(G)$ (cf. \cite[Sec.6]{Kottwitz97}, \cite[Sec.4]{Rapoport05}): 
\[B(G,\{\mu\}):=\left\{\ [b]\in B(G)\ |\quad \kappa_{G}([b])=-\mu^{\natural},\quad \overline{\nu}_{G}([b])\preceq -\overline{\mu}\ \right\},\]
where $\preceq$ is the natural partial order on the closed Weyl chamber $\overline{C}$ defined by that $\nu\preceq \nu'$ if $\nu'-\nu$ is a nonnegative linear combination (with \emph{rational} coefficients) of simple coroots in $R_{\ast}(T)$ \cite{RR96}, Lemma 2.2). 
The two minus sign are due to our sign convention (which results from our decision to work with homology instead of cohomology).
One knows \cite[4.13]{Kottwitz97} that the map 
\begin{equation}  \label{eq:kappa_times_nu} 
(\overline{\nu},\kappa):B(G)\rightarrow \mathcal{N}(G)\times X^{\ast}(Z(\widehat{G})^{\Gamma_F})
\end{equation}
is injective, hence $B(G,\{\mu\})$ can be identified with a subset of $\mathcal{N}(G)$; we will give further information of this map later.

%%%%%%%%%%%%%%%%%%%%
\subsection{Galois hypercohomology of crossed modules}
In order to extend various results in Galois cohomology used previously by Langlands-Rapoport and Kottwitz for the problem at hand
%(especially, Satz 5.25 of \cite{LR87}, namely Corollary \ref{cor:LR-Satz5.25} here)
in the case $G^{\der}$ is simply connected to a general case, we need to consider Galois cohomology groups of crossed modules and (length-$2$) complexes of tori quasi-isomorphic to them (for example, $H^1_{\ab}(\Q,G)$ instead of the cohomology of the quotient $G^{\ab}:=G/G^{\der}$). 
Here we give a (very) brief review of the theory of Galois hypercohomology of crossed modules; for details, see \cite{Borovoi98}, \cite[Ch.1]{Labesse99}. 
For a connected reductive group $H$ over a field $k$, we denote by 
\[\rho_H:H^{\uc}\rightarrow H\] 
the canonical map from the simply connected cover $H^{\uc}$ of $H^{\der}$ to $H$. Unless stated otherwise, every bounded complex of groups considered in this paper will be concentrated in non-negative homological degrees.

For a connected reductive group $G$ over a field $k$ (of characteristic zero), if $T$ is a maximal $k$-torus of $G$, the complex $(\rho^{-1}(T)\rightarrow T)$ of $k$-tori, where $\rho^{-1}(T)$ and $T$ are placed in degree $-1$ and $0$ respectively, is quasi-isomorphic to its sub-complex $(\rho^{-1}(Z(G))=Z(G^{\uc})\rightarrow Z(G))$ (cf. \cite[$\S$2, $\S$3]{Borovoi98}), hence, as an object in the derived category of complexes of commutative algebraic $k$-group schemes, depends only on $G$. We denote it by $G_{\mathbf{ab}}$:
\begin{equation} \label{eq:abelianized_cx}
G_{\mathbf{ab}}:=\rho^{-1}(T)\rightarrow T.
\end{equation}
Then, following Borovoi \cite{Borovoi98}, we define the abelianized Galois cohomology group $H^i_{\ab}(k,G)\ (i\in\Z)$ of $G$ to be the hypercohomology group of $G_{\mathbf{ab}}(\bar{k}):=(\rho^{-1}(T)(\bar{k})\rightarrow T(\bar{k}))$, two-term complex of discrete $\Gal(\bar{k}/k)$-modules:
\[H^i_{\ab}(k,G):=\mathbb{H}^i(k,G_{\mathbf{ab}}).\] 
(We will also often write $H^i(k,G_{\mathbf{ab}})$ for $\mathbb{H}^i(k,G_{\mathbf{ab}})$.)
For $-1\leq i\leq 1$, this group is also equal to the cohomology (group) 
\[\mathbb{H}^i(k,G^{\uc}\stackrel{\rho}{\rightarrow} G)\] 
of the crossed module $G^{\uc}(\bar{k})\stackrel{\rho}{\rightarrow} G(\bar{k})$ of $\Gal(\bar{k}/k)$-groups (ibid., (3.3.2)): for a proof, see the proof of the next lemma.

%%%%%%%%%%%%%%%%%%%%
\begin{lem} \label{lem:abelianization_exact_seq}
Let $H\subset G$ be a (not necessarily connected) $k$-subgroup containing a maximal $k$-torus of $G$ and $(a_{\tau})_{\tau}$ a cochain on $\Gal(\bar{k}/k)$ valued in $H$ whose coboundary $a_{\tau_1}\tau_1(a_{\tau_2})a_{\tau_1\tau_2}^{-1}$ belongs to $Z(H)^{\mathrm{o}}$. Then, for the (simultaneous) inner twist $\rho_1:\tilde{H}_1\rightarrow H_1$ of the canonical map $\rho:\tilde{H}:=\rho^{-1}(H)\rightarrow H$ via the cocyle $a_{\tau}^{\ad}\in Z^1(\Q,H^{\ad})$ (the image of $a_{\tau}$ in $H^{\ad}$), there exists a natural exact sequence 
\begin{equation} \label{eq:abelianization_from_Levi}
H^1(k,\tilde{H}_1) \stackrel{\rho_{1\ast}}{\longrightarrow} H^1(k,H_1) \stackrel{ab_1}{\longrightarrow} H^1_{\ab}(k,G).
\end{equation}
The abelianization complexes $\tilde{H}_{\mathbf{ab}}$, $\tilde{H}_{1\mathbf{ab}}$ are isomorphic.
\end{lem}

Note that the condition on (the coboundary of) $a_{\tau}$ allows us to twist $\tilde{H}$ via  $a_{\tau}^{\ad}$.
In the applications, $H$ will be the centralizer of a semisimple element of $G(k)$.

\begin{proof}
The first map is the obvious one (induced by $\rho_1$) and the second map is the composite of the natural map 
\begin{equation} \label{eq:H^1(ab)}
H^1(k,H_1)\rightarrow \mathbb{H}^1(k,\tilde{H}_1\rightarrow H_1)
\end{equation}
resulting from the map $(1\rightarrow H_1)\rightarrow (\tilde{H}_1\rightarrow H_1)$ of crossed modules of $k$-groups, and the isomorphisms 
\begin{equation*} \label{eq:isom_of_abelianized_coh}
\mathbb{H}^1(k,\tilde{H}_1\rightarrow H_1)\isom \mathbb{H}^1(k,\tilde{Z}(G)\rightarrow Z(G))\isom \mathbb{H}^1(k,G^{\uc}\rightarrow G)
\end{equation*}
($\tilde{Z}(G):=\rho^{-1}(Z(G))=Z(G^{\uc})$) resulting from the quasi-isomorphisms of crossed modules of $k$-groups 
\begin{equation*}
(\tilde{H}_1\rightarrow H_1)\leftarrow (\tilde{Z}(G)\rightarrow Z(G)) \rightarrow (G^{\uc}\rightarrow G)
\end{equation*}
(cf. \cite{Borovoi98}, Lem. 2.4.1 and its proof: the key point is that $H_1$ contains $Z(G)$). 
Now, the exactness follows from \cite[Cor. 3.4.3]{Borovoi98}.
\end{proof}

When $k$ is a non-archimedean local field and $G$ is a connected reductive group, the natural map $ab:H^1(k,G)\ra H^1(k,G_{\mathbf{ab}})$ (\ref{eq:H^1(ab)}) is a bijection, giving canonical and functorial (in $G$) bijections 
\begin{equation} \label{eq:H^1(Qv,G)=H^1(Qv,G_{ab})}
H^1(k,G)=H^1(k,G_{\mathbf{ab}})=(\pi_1(G)_{\Gamma_k})_{\mathrm{tor}}
\end{equation}
\cite[Cor.5.4.1, Cor.5.5]{Borovoi98}. In particular, if $G_1$ is an inner twist of $G$, there are canonical bijections
\begin{equation} \label{eq:H^1(Qv,G)=H^1(Qv,G_1)}
 H^1(k,G)=(\pi_1(G)_{\Gamma_k})_{\mathrm{tor}} =(\pi_1(G_1)_{\Gamma_k})_{\mathrm{tor}} = H^1(k,G_1) : 
\end{equation} 
the middle isomorphism is induced by any inner twist $G_{\bar{k}}\isom (G_1)_{\bar{k}}$, but is also independent of the choice of such one.

%%%%%%%%%%%%%%%%%%%%%%%%%%%%%%%%%%%%%%%%
\subsection{The groups $G_{\delta\theta}$, $J_b$}

%%%%%%%%%%%%%%%%%%%%
\subsubsection{The group $G_{\delta\theta}$}

For the next discussion, we follow \cite[$\S$5]{Kottwitz82} (but will use slightly different notations and conventions).
Let $E/F$ be a cyclic extension of degee $n$ contained in $\bar{F}$ (of characteristic zero), and let $\sigma$ be a generator of $\Gal(E/F)$. Let $G$ be a connected reductive group over $F$, and $R:=\Res_{E/F}G_E$ the Weil restriction of the base-change of $G$ (this group was denoted by $I$ in \textit{loc. cit.}). 
There exists a natural isomorphism $R_E\isom G_E\times \cdots G_E$, where the factors are ordered such that the $i$-th factor corresponds to $\sigma^i\in \Gal(E/F)\ (i=1,\cdots,l)$ (note the convention is different from that of \textit{loc. cit.}).
The element $\sigma\in \Gal(E/F)$ determines an automorphism $\theta\in\Aut_F(R)$, which on $R(E)$ takes the form
\[(x_1,\cdots,x_{n-1},x_n)\mapsto (x_2,\cdots,x_n,x_1)\]
The inclusion
\[ \Delta_{\sigma}:G(E)=R(F)\rightarrow R(E)=G(E)\times \cdots\times G(E) \]
is given by
\[x\mapsto (x,x^{\sigma},\cdots, x^{\sigma^{n-1}}).\]
We define an $F$-morphism $N:R\rightarrow R$ by $Nx=x\cdot x^{\theta}\cdots x^{\theta^{n-1}}$; clearly, one has $N\circ \Delta_{\sigma}=\Delta_{\sigma}\circ\Nm$ on $G(E)=R(F)$ ($\Nm$ denotes the map $\Nm_n$ on $G(E)$).

For $x\in G(E)=R(F)$, the \emph{$\sigma$-centralizer} of $x$ is by definition the $F$-subgroup of $R$:
\begin{equation} \label{eq:G_{xtheta}}
G_{x\theta}=\{g\in R\ :\ gxg^{-\theta}=x\}.
\end{equation}
(We sometimes write $g^{-\theta}$ for $\theta(g^{-1})$; in \cite[$\S$5]{Kottwitz82}, $\theta$ and $G_{x\theta}$ were denoted respectively by $s$ (on p. 801) and $I_{sx}$ (on p. 802)).
If $p:R_E=G_E\times\cdots G_E\rightarrow G_E$ denotes the projection onto the factor indexed by the identity element of $\Gal(E/F)$, by restriction $p$ induces an isomorphism 
\begin{equation} \label{eq:p_x}
p_x:(G_{x\theta})_E\rightarrow (G_E)_{\Nm x}
\end{equation}
(ibid., Lemma 5.4). Define $G_{x\theta}^{\ast}$ to be the inverse image under the $E$-isomorphism $p_x$ of the subgroup $G_{\Nm x}^{\ast}$ of $G_{\Nm x}$ defined in ibid., $\S$3 (this equals $G_{\Nm x}^{\mathrm{o}}$ if $\Nm x$ is semisimple). This is an $F$-subgroup of $G_{x\theta}$ (ibid., Lemma 5.5) which equals the neutral component $G_{x\theta}^{\mathrm{o}}$ of $G_{x\theta}$ if $\Nm x$ is semisimple, and equals $G_{x\theta}$ when $G^{\der}=G^{\uc}$.

People (and we) also write $G_{\delta\sigma}$ for $G_{\delta\theta}$.

%%%%%%%%%%%%%%%%%%%%
\subsubsection{The group $J_b$} \label{subsubsec:J_b}

From now on until the end of this subsection, we will be in the following $p$-adic setup.

Let $k$ be an algebraically closed field of characteristic $p>0$ and $K=W(k)[\frac{1}{p}]$.
Suppose that $F$ is a finite extension of $\Qp$ in a fixed algebraic closure $\bar{K}$ of $K$ and let $L$ be the composite of $K$ and $F$ in $\bar{K}$. For $n\in\N$, let $L_n\subset L$ denote the unramified extension of $F$ of degree $n$.

Let $W(\bar{K}/F)$ be the Weil group, i.e. the group of continuous automorphisms of $\bar{K}$ which fix $F$ pointwise and induce on the residue field $k$ an integral power of the Frobenius automorphism (cf. \cite[$\S$1]{Kottwitz85}); when $k=\Fpb$, the restriction homomorphism $W(\bar{K}/F)\rightarrow \Gal(\bar{F}/F)$ identifies $W(\bar{K}/F)$ with the usual absolute Weil group $W_F$ of $F$ (ibid., (1.4)). There exists an exact sequence $1\rightarrow \Gal(\bar{K}/L)\rightarrow W(\bar{K}/F)\rightarrow \langle\sigma\rangle \rightarrow 1$ which endows $W(\bar{K}/F)$ with a natural topology such that the injection identify $\Gal(\bar{K}/L)$ with an open subgroup. 

For a linear algebraic group $G$ over $F$ and  $b\in G(L)$, we define a cocyle 
\begin{equation} \label{eq:b_{tau}}
 b_{\tau}:=\Nm_{i(\tau)}b\ \in Z^1(\langle\sigma\rangle, G(L)),
\end{equation}
where for each $\tau\in W(\bar{K}/F)$, $i(\tau)\in\Z$ is determined by $\tau|_L=\sigma^{i(\tau)}$ and for $n<0$, $\Nm_{n}b:=\sigma^{n}(b_{\sigma^{-n}}^{-1})$. Sometimes, via inflation, we regard $b_{\tau}$ as a cocycle in $Z^1(W(\bar{K}/F),G(\bar{K}))$. We recall that when $G$ is connected, as $H^1(L,G)$ is trivial, the following inflation map on cohomology sets is a bijection: 
\[ B(G)=H^1(\langle\sigma\rangle,G(L)) \ \isom \ H^1(W(\bar{K}/F), G(\bar{K})). \] 

Next, we recall (\cite[3.3 - 3.5]{Kottwitz97}) that for any linear algebraic $F$-group $G$ and $b\in G(L)$, there is associated an $F$-group $J_b=J_b^G$ defined by that for any $F$-algebra $R$, 
\begin{align} \label{eq:J_b}
J_b(R) &:=G(\bar{K}\otimes R)^W=\{g\in G(\bar{K}\otimes R)\ |\ b_{\tau}\tau(g)b_{\tau}^{-1}=g,\ \forall\tau\in W(\bar{K}/F)\} \\
&=G(L\otimes R)^{\langle\sigma\rangle}=\{g\in G(L\otimes R)\ |\ b\sigma(g)=gb\} \nonumber
\end{align}
Here, $W:=W(\bar{K}/F)$ acts on $G(\bar{K}\otimes R)$ by the \emph{twisted} action $\tau(g)=b_{\tau}{} \tau(g) b_{\tau}^{-1}$ while $\tau(g)$ is the \emph{standard} action through the factor $\bar{K}$ (the action of $\langle\sigma\rangle$ on $G(L\otimes R)$ is similarly defined). For any $\bar{K}$-algebra $R$, the map $G(\bar{K}\otimes R)\rightarrow G(R)$ induced by the canonical $\bar{K}$-algebra homomorphism $\bar{K}\otimes R\rightarrow R:l\otimes x\mapsto lx$ restricts to an injection
\begin{equation*} \label{eq:Xi}
\Xi:J_b(R)=G(\bar{K}\otimes R)^W \hookrightarrow G(R),
\end{equation*}
and this identifies $(J_b)_{\bar{K}}$ with the centralizer subgroup scheme of $\nu_b$ in $G_{\bar{K}}$. Its restriction to $J_b(F)$ equals the inclusion $J_b(F)=G(\bar{K})^{W}\subset G(R)$ (as it is induced by the inclusion $G(\bar{K})\rightarrow G(\bar{K}\otimes R)\rightarrow G(R)$). 
When $R=\bar{K}$, the resulting map $J_b(\bar{K})\hookrightarrow G(\bar{K})$ is $W$-equivariant for the standard action on the source $J_b(\bar{K})$ (i.e. via the factor $R=\bar{K}$ in $G(\bar{K}\otimes R)$) and the twisted action defined above on the target $G(\bar{K})$. These continue to hold if we use $L$ in place of $\bar{K}$; for example, the same map gives an injection 
\begin{equation} \label{eq:Xi2}
\Xi:J_b(R)=G(L\otimes R)^{\langle\sigma\rangle}\hookrightarrow G(R)
\end{equation} 
for any $L$-algebra $R$, which is $W$-equivariant when $R=\bar{K}$. It follows that $\Xi$ gives rise to a map 
\begin{align*} 
 j_b=j_b^G\  :\  Z^1(W,J_b(\bar{K})) & \rightarrow  Z^1(W,G(\bar{K})) \\
 x_{\tau} \quad  & \mapsto \quad  \Xi(x_{\tau})\cdot b_{\tau}.
\end{align*}
which induces a map $H^1(W,J_b(\bar{K}))\ra H^1(W,G(\bar{K}))$.
When $b'=gb\sigma(g)^{-1}$ for $g\in G(L)$, $\Int(g)$ gives an $F$-isomorphism $J_b\ra J_{b'}$ and the maps $j_b$, $j_{b'}\circ \Int(g):H^1(W,J_b)\ra H^1(W,G)$ are the same.
The restriction homomorphism $W(\bar{K}/F)\ra \Gal(\bar{F}/F)$ and the inclusion $G(\bar{F})\subset G(\bar{K})$ gives an injection 
\[ H^1(F,G)\hookrightarrow H^1(W,G(\bar{K})) \] 
\cite[1.8.2]{Kottwitz85}.
Upon restriction to the subset $H^1(F,J_b)$ of $H^1(W,J_b(\bar{K}))$, $j_b^G$ induces an injection (with the same notation)
\begin{equation} \label{eq:j_b}
j_b=j_b^G:H^1(F,J_b)\hookrightarrow H^1(W,G(\bar{K})).
\end{equation}
which sends the distinguished element to the $\sigma$-conjugacy class of $b$ and whose image is $\bar{\nu}_{G}^{-1}(\bar{\nu}_G([b]))$ for the Newton map $\bar{\nu}_G$ (\ref{eq:kappa_G}) \cite[1.15]{RR96}, \cite[3.5]{Kottwitz97}.

%%%%%%%%%%%%%%%%%%%%
\begin{lem} \label{lem:Xi_{b'}_inner-twisting}
Fix a uniformizer $\pi$ of $F$. Let $H$ be a connected reductive $F$-group and $b\in H(L)$ a \emph{basic} element (i.e. the Newton homomorphism $\nu:=\nu_H(b)\in\Hom_L(\mathbb{D},H)$  (\ref{subsubsec:Newton_map}) maps to $Z(H)$).

(1) Choose $h\in H(L)$ and $n\in\N$ satisfying $\Nm_n(b')=\nu(\pi)^n$ for $b':=h^{-1}b\sigma(h)$. 
Then, the cocycle $b'^{\ad}_{\tau}$ (\ref{eq:b_{tau}}) defined by the image $b'^{\ad}$ of $b'$ in $H^{\ad}(L)$ in fact lies in $Z^1(L_n/F,H^{\ad}(L_n)$), and the $L$-isomorphism (\ref{eq:Xi2})
\begin{equation} \label{eq:Xi_{b'}}
 \Xi:(J_{b'}^{H})_{L} \isom H_L
 \end{equation}
 is the base-change to $L$ of an inner twisting defined over $F^{\nr}$ corresponding to $[b'^{\ad}_{\tau}]\in H^1(F,H^{\ad})$.
 
(2) The composite map 
\[ H^1(F, J_b) \stackrel{j_b^H}{\hookrightarrow} B(H)_{basic} \stackrel{\kappa_H}{\isom} \pi_1(H)_{\Gamma} \]
is the translation by $\kappa_H([b])$ under the canonical bijections $H^1(F,J_b)\cong (\pi_1(J_b)_{\Gamma})_{\mathrm{tor}}\cong (\pi_1(H)_{\Gamma})_{\mathrm{tor}}$ (\ref{eq:H^1(Qv,G)=H^1(Qv,G_1)}), where $B(H)_{basic}$ is the subset of $B(H)$ consisting of basic elements and $\Gamma:=\Gal(\bar{F}/F)$.

(3) Let $I_0\subset G$ be connected reductive $F$-groups and $b\in I_0(L)$ a basic element. There exist canonical bijections
\[ \ker[H^1(F,J_b^{I_0})\rightarrow H^1(F,J_b^G)]=\ker[H^1(F,I_0)\rightarrow H^1(F,M)]  =\ker[H^1(F,I_0)\rightarrow H^1(F,G)], \]
where for $H=I_0,G$, $J_b^H$ denotes the corresponding group $J_b$ and $M:=Z_G(\nu_{I_0}(b))$. 
The resulting injection
\begin{equation} \label{eq:j_[b]^{I_0}}
j_{[b]}^{I_0}:\ker[H^1(F,I_0)\rightarrow H^1(F,G)] \hookrightarrow B(I_0)
\end{equation}
depends only on the $\sigma$-conjugacy class $[b]_{I_0}$, not on the choice of its representative in $I_0(L)$ (which justifies the notation $j_{[b]}^{I_0}$).
\end{lem}

\begin{proof}
(1) That $\Xi$ (\ref{eq:Xi2}) is an isomorphism in this case follows from that $\nu$ maps to $Z(H)$. The claim is proved in \cite[5.2]{Kottwitz85}: the fact that $\Xi:J_{b'}(L)\hookrightarrow H(L)$ is $\langle\sigma\rangle$-equivariant for the standard action on $J_b(L)$ and the twisted action on $H(L)$ implies that $\Xi$ serves as an inner twisting denoted by $u$ in ibid.

(2) First, as $j_b^H(H^1(F,J_b))=\bar{\nu}_{G}^{-1}(\bar{\nu}_G([b]))$, we have $j_b^H(H^1(F,J_b))\subset B(H)_{basic}$. Also, the restriction of $\kappa_H:B(H)\ra \pi_1(H)_{\Gamma}$ to $B(H)_{basic}$ becomes an isomorphism \cite[5.6]{Kottwitz85}.

The map $j_b^H:H^1(F,J_b)\hookrightarrow H^1(W,H(\bar{K}))$ is the restriction of $j_b^H:H^1(W,J_b(\bar{K})) \ra H^1(W,H(\bar{K}))$, and 
the diagram 
\[ \xymatrix@C=4pc@R=1.5pc{ B(J_b^H) \ar[r]^{j_b} \ar[d]_{\kappa_{J_b^H}} & B(H) \ar[d]^{\kappa_H} \\ \pi_1(J_b^H)_{\Gamma} \ar[r]^{\quad \cdot \kappa_{J_b}(b)} & \pi_1(H)_{\Gamma} } \] 
is commutative (cf. \cite{Kottwitz97}, the diagram in the proof of Cor. 4.12), where the map in the bottom is the translation by $\kappa_{J_b}(b)$ on the abelian group
$\pi_1(J_b^H)_{\Gamma}= \pi_1(H)_{\Gamma}$ (as the map $\Int(h):J_b^H\isom J_{b'}^H$ is an $F$-isomorphism, $J_b$ is an inner twist of $H$ by (1)).
Since the restriction of $\kappa_{J_b^H}$ to $H^1(F,J_b^H)$ is the canonical isomorphism $H^1(F,J_b^H)\cong \pi_1(J_b^H)_{\Gamma,\mathrm{tor}}$ \cite[5.7]{Kottwitz85}, the claim follows.

(3) Let $b'\in I_0(L)$ be as in (1). The inclusion $J_{b'}^{I_0}\hookrightarrow J_{b'}^G$ is a simultaneous inner twist of the inclusion $I_0\hookrightarrow M$, via the cochain $b'_{\tau}\in C^1(W,I_0)$ defined by $b'$ using the recipe (\ref{eq:b_{tau}}) which becomes cocycles in both $I_0^{\ad}(F^{\nr})$ and $M^{\ad}(F^{\nr})$ according to (1). 
The first equality is induced by this simultaneous inner twisting $((J_b^{I_0})_{\bar{F}}\hookrightarrow (J_b^G)_{\bar{F}}) \isom ((I_0)_{\bar{F}}\hookrightarrow M_{\bar{F}})$ and follows from (\ref{eq:H^1(Qv,G)=H^1(Qv,G_{ab})}); so this equality is canonical, independent of the choice of such simultaneous inner twisting. 
The second equality is due to the fact that as $M$ is a Levi-subgroup of $G$, the map $H^1(F,M)\rightarrow H^1(F,G)$ is injective (cf. \cite[(4.13.3)]{Kottwitz97}).
\end{proof}

For a linear algebraic group $G$ over a field $F$, let $G(F)_{\mathrm{ss}}$ be the subset of $G(F)$ of semisimple elements.

%%%%%%%%%%%%%%%%%%%%
\begin{prop} \label{prop:psi_p}
Let $G$ be a connected reductive group over $F=\Qp$; we write $\mfk$ for $L=K$.
Let $(\gamma_0,\delta,c)\in G(\Qp)_{\mathrm{ss}}\times G(L_n)\times G(\mfk)$ be a triple satisfying that
\begin{equation} \label{eq:(iii)'}
c\gamma_0 c^{-1}=\Nm_n\delta,\quad \text{ and }\quad b:=c^{-1}\delta\sigma(c)\text{ lies in }I_0(\mfk)\ (I_0:=G_{\gamma_0}^{\mathrm{o}}).
\end{equation}
(a priori, one only has $b\in G_{\gamma_0}(\mfk)$: $\delta^{-1}c\gamma_0 c^{-1}\delta=\sigma(\delta)\cdots\sigma^n(\delta)=\sigma(c\gamma_0 c^{-1})=\sigma(c)\gamma_0\sigma(c^{-1})$.)

(1) The element $b=c^{-1}\delta\sigma(c)\in I_0(\mfk)$ has Newton quasi-cocharacter $\nu_{I_0}(b)=\frac{1}{n}\nu_{I_0}(\gamma_0)$. In particular, its $\sigma$-conjugacy class $[b]_{I_0}\in B(I_0)$ is basic. 

(2) The functorial map $\Int(c^{-1}):G(\mfk\otimes R)\isom G(\mfk\otimes R)$ (for $\Qp$-algebras $R$) restricts to a $\Qp$-isomorphism
\[ G_{\delta\theta}^{\mathrm{o}}\isom J_b^{I_0}.\]
Its composite with the map $\Xi:(J_b^{I_0})_{\mfk} \isom (I_0)_{\mfk}$ (\ref{eq:Xi2}):
\[ \Xi\circ \Int(c^{-1}) \ : \ (G_{\delta\theta}^{\mathrm{o}})_{\mfk} \isom (J_b^{I_0})_{\mfk} \isom (I_0)_{\mfk}\] 
equals the restriction to the neutral components of the isomorphism 
\[ \Int(c^{-1})\circ p_{\delta} \ : \ (G_{\delta\theta})_{\mfk} \isom (G_{\mfk})_{\Nm \delta} \isom (G_{\gamma_0})_{\mfk}, \]
where $p_{\delta}$ is (\ref{eq:p_x}). 

(3) For any $h\in I_0(\mfk)$ as in Lemma \ref{lem:Xi_{b'}_inner-twisting} (1), 
\begin{equation} \label{eq:psi_p}
\psi_p:=p_{\delta}^{-1} \circ \Int(ch) \ : \ (I_0)_{\mfk} \isom (G_{\delta\theta}^{\mathrm{o}})_{\mfk} 
\end{equation} 
is the base-change to $\mfk$ of an inner twisting over $\Qpnr$ with the cohomology class $[b^{\ad}]\in B(I_0^{\ad})=H^1(\langle\sigma\rangle,I_0^{\ad}(\mfk))$ which in fact lies in the subset $H^1(\Qpnr/\Qp,I_0^{\ad}(\Qpnr))$. 
\end{prop}

In (3), the element $c':=ch\in G(\mfk)$ satisfies (\ref{eq:(iii)'}) and the identity $\Nm_n(b')=\nu(\pi)^n$ for $b':=c'^{-1}\delta\sigma(c')$ and $\nu:=\nu_{I_0}(b')\in\Hom_{\Qp}(\mathbb{D},Z(I_0))$. In fact, it is easy to see that (3) holds for any such element $c'$.

\begin{proof}
(1) This is Lemma \ref{lem:equality_of_two_Newton_maps} (2) (which slightly refines \cite[Lem.5.15]{LR87}). 

(2) One has $\gamma_0=(b\sigma)^n\cdot c^{-1}\cdot \sigma^{-n}\cdot c$ (identity in $G(\mfk)\rtimes\langle\sigma\rangle$), which implies that an element $g\in G(\mfk\otimes R)$ commutes with $\gamma_0$ and $b\sigma$ if and only if $cgc^{-1}$ does with $\delta\sigma=c(b\sigma)c^{-1}$ and $\sigma^n$ ($\sigma$ acts on $G(\mfk\otimes R)$ via its action on $\mfk$). This proves the first statement.

For $\mfk$-algebra $R$, the map $\Int(c^{-1})\circ p_{\delta}$ is the restriction to $G_{\delta\theta}(R)$ of the composite
\[ G(L_n\otimes R)\cong \prod_{L_n\hookrightarrow R}G(R) \stackrel{pr_1}{\longrightarrow} G(R)\stackrel{\Int(c^{-1})}{\longrightarrow} G(R), \]
where the first decomposition is over the set of embeddings $L_n\ra R$ and $pr_1$ is the projection to the factor corresponding to the inclusion $L_n\subset \mfk \ra R$.
On the other hand, the first map is the restriction to $G_{\delta\theta}(R)$ of the composite
\[  G(\mfk\otimes R) \stackrel{\Int(c^{-1})}{\longrightarrow} G(\mfk\otimes R) \stackrel{\Xi}{\longrightarrow} G(R) \]
where $\Xi$ is induced by $\mfk\otimes R \ra R:l\otimes x\mapsto lx$.
It is obvious that these two maps become identical upon restriction to $G(L_n\otimes R)$. This proves the second statement.

(3) For $b':=h^{-1}b\sigma(h)\in I_0(\mfk)$, we have $b'^{\ad}\in I_0(\Qpnr)$ and $\Xi:(J_{b'}^{I_0})_{\mfk} \isom (I_0)_{\mfk}$ is an inner twisting with cohomology class $[b'^{\ad}]\in H^1(\Qpnr/\Qp,I_0^{\ad}(\Qpnr))\subset B(I_0^{\ad})$ (Lemma \ref{lem:Xi_{b'}_inner-twisting} (1)). If $\Xi_b:(J_b^{I_0})_{\mfk} \isom (I_0)_{\mfk}$ and $\Xi_{b'}:(J_{b'}^{I_0})_{\mfk} \isom (I_0)_{\mfk}$ denote the isomorphisms $\Xi$ (\ref{eq:Xi2}) for $b$ and $b'$ respectively, then $f:=\Int(h^{-1}):J_b\ra J_{b'}$ is a $\Qp$-isomorphism and we have $\Int(h^{-1})\circ \Xi_b=\Xi_{b'}\circ f$. 
Hence, $\Xi_{b'}\circ f\circ \Int(c^{-1})=\Int(h^{-1})\circ \Xi_b\circ \Int(c^{-1})= \Int(ch)^{-1}\circ p_{\delta}: (G_{\delta\theta}^{\mathrm{o}})_{\mfk} \isom (J_{b'}^{I_0})_{\mfk} \isom (I_0)_{\mfk}$ is an inner twisting defined over $\Qpnr$ with class $[b^{\ad}]=[b'^{\ad}]\in B(I_0^{\ad})$.
\end{proof}

%%%%%%%%%%%%%%%%%%%%%%%%%%%%%%%%%%%%%%%%
\subsection{Kottwitz triple} \label{subsubsec:pre-Kottwitz_triple}
Our main references for the material covered here are \cite[p.182-183]{LR87}, \cite[$\S2$]{Kottwitz90}, and \cite{Kottwitz92}.

%%%%%%%%%%%%%%%%%%%%
\begin{defn}  \label{defn:Kottwitz_triple}
A Kottwitz triple is a triple $(\gamma_0;\gamma=(\gamma_l)_{l\neq p},\delta)$, where
\begin{itemize} \addtolength{\itemsep}{-4pt} 
\item[(i)] $\gamma_0$ is a semisimple element of $G(\Q)$ that is elliptic in $G(\R)$, defined up to conjugacy in $G(\overline{\Q})$;
\item[(ii)] for $l\neq p$, $\gamma_l$ is a semisimple element in $G(\Q_l)$, defined up to conjugacy in $G(\Q_l)$, such that $\gamma=(\gamma_l)\in G(\A_f^p)$ and is conjugate to $\gamma_0$ under $G(\bar{\A}_f^p)$;
\item[(iii)] $\delta$ is an element of $G(L_n)$ (for some $n$), defined up to $\sigma$-conjugacy in $G(L_n)$, such that the norm $\Nm_n\delta$ of $\delta$  is conjugate to $\gamma_0$ under $G(\bar{\mfk})$, where $\Nm_n\delta:=\delta\cdot\sigma(\delta)\cdots\sigma^{n-1}(\delta)\in G(L_n)$, and the following condition holds:
\item[$\ast(\delta)$] \label{itm:ast(delta)}
the image of $[\delta]$ under the Kottwitz homomorphism $\kappa_{G_{\Q_p}}:B(G_{\Q_p})\rightarrow \pi_1(G_{\Q_p})_{\Gamma(p)}$ (\ref{subsubsec:Kottwitz_hom}) is equal to $-\mu^{\natural}$ (defined in (\ref{eqn:mu_natural})).%%{sign convention}
\end{itemize}
\end{defn}

We call the natural number $n$ in (iii) the \emph{level} of the given Kottwitz triple. 
Two Kottwitz triples $(\gamma_0;(\gamma_l)_{l\neq p},\delta)$, $(\gamma_0';(\gamma_l')_{l\neq p},\delta')$ are said to be \textit{(geometrically) equivalent}, if $\gamma_0$ is $G(\Qb)$-conjugate to $\gamma_0'$, $\gamma_l$ is $G(\Ql)$-conjugate to $\gamma_l'$ for each finite $l\neq p$, and $\delta, \delta'\in G(L_n)$ for some $n$ and $\delta$ is $\sigma$-conjugate to $\delta'$ in $G(L_n)$ (i.e. there exists $d\in G(L_n)$ such that $\delta'=d\delta\sigma(d^{-1})$). Normally, we consider Kottwitz triples having level $n$ divisible by $[\kappa(\wp):\Fp]$.

When the derived group of $G$ is not simply connected, the relevant notion of Kottwitz triple is a stable version of the above definition.
Recall \cite[$\S$3]{Kottwitz82} that for a connected reductive group $F$ over a perfect field $F$, two rational elements $x,y\in G(F)$ are said to be \emph{stably conjugate} if there exists $g\in G(\bar{F})$ such that $gxg^{-1}=y$ and $g^{-1}{}^{\tau}g\in G_s^{\mathrm{o}}$ for all $\tau\in \Gal(\bar{F}/F)$, where $s$ is the semisimple part of $x$ in its Jordan decomposition. When we just refer to the relation $gxg^{-1}=y$ for some $g\in G(\bar{F})$, we will say that $x,y$ are $\bar{F}$-(or $G(\bar{F})$-)conjugate or \emph{geometrically} conjugate. By definition, a \textit{stable conjugacy class} in $G(F)$ is an equivalence class in $G(F)$ with respect to this stable conjugation relation. Let $G_{\gamma_0}:=Z_G(\gamma_0)$ (centralizer of $\gamma_0$ in $G$).

%%%%%%%%%%%%%%%%%%%%
\begin{defn} \label{defn:stable_Kottwitz_triple}
A Kottwitz triple $(\gamma_0;\gamma=(\gamma_l)_{l\neq p},\delta)$, say of level $n\in\N$, is \emph{stable} if it satisfies the following conditions (in addition to (i) - (iii)):
\begin{itemize} \addtolength{\itemsep}{-4pt} 
\item[(i$'$)] $\gamma_0\in G(\Q)_{\mathrm{ss}}$ is defined up to stable conjugacy;
\item[(ii$'$)] for each $l\neq p$, $\gamma_0$ is stably conjugate to $\gamma_l$;
\item[(iii$'$)] there exists $c\in G(\mfk)$ such that $c\gamma_0 c^{-1}=\Nm_n\delta$ and $b:=c^{-1}\delta\sigma(c)$ lies in $I_0(\mfk)$ for $I_0:=G_{\gamma_0}^{\mathrm{o}}$ (a priori, one only has $b\in G_{\gamma_0}(\mfk)$: $\delta^{-1}c\gamma_0 c^{-1}\delta=\sigma(\delta)\cdots\sigma^n(\delta)=\sigma(c\gamma_0 c^{-1})=\sigma(c)\gamma_0\sigma(c^{-1})$).
\end{itemize}
\end{defn}

If $G^{\der}$ is simply connected, one has $G_{\gamma_0}=I_0$ by Steinberg's theorem, so these two sets of conditions are the same, and every Kottwitz triple is stable. 
It is easy to verify the following fact: let $(\gamma_0;\gamma,\delta)$ be a stable Kottwitz triple of level $n$. If $\gamma_0'\in G(\Q)$ is stably conjugate to $\gamma_0$, $\gamma$ is $G(\A_f^p)$-conjugate to $\gamma'$, and $\delta'$ is $\sigma$-conjugate to $\delta$ in $G(L_n)$, then $(\gamma_0';\gamma',\delta')$ is also a stable Kottwitz triple (of same level). We say that two stable Kottwitz triples $(\gamma_0;\gamma,\delta)$, $(\gamma_0';\gamma',\delta')$ are \emph{stably equivalent} if they are (geometrically) equivalent as Kottwitz triples and $\gamma_0$, $\gamma_0'$ are stably conjugate. A priori, two triples satisfying the conditions of Definition \ref{defn:stable_Kottwitz_triple} except for (i$'$) can be geometrically equivalent without being stably equivalent (see the next remark (3)).

%%%%%%%%%%%%%%%%%%%%
\begin{rem} \label{rem:Kottwitz_triples}
(1) There is one more condition imposed on Kottwitz triples which was taken as a part of the definition by Kottwitz \cite[p.171-172]{Kottwitz90} (cf. \cite[4.3.1(iv)]{Kisin17}). Since in our work (a stronger version of) it will be always satisfied by the Kottwitz triples that are geometrically produced (Corollary \ref{cor:Tate_thm2}), we just recall this condition briefly (with a natural modification in the general case $G_{\gamma_0}\neq I_0$). First, for each finite place $v\neq p$ of $\Q$, set \[ I(v):=\begin{cases} Z_{G_{\Q_v}}(\gamma_v)^{\mathrm{o}} & \text{ if }v\neq p, \\ G_{\delta\theta}^{\mathrm{o}} \ (\ref{eq:G_{xtheta}}) & \text{ if } v=p. \end{cases} \] At the infinite place, we choose an elliptic maximal torus $T_{\R}$ of $G_{\R}$ containing $\gamma_0$ and $h\in X\cap \Hom(\dS,T_{\R})$. We twist $G_{\gamma_0}$ using the Cartan involution $\Int (h(i))$ on $I_0/(Z(G)\cap I_0)$, and get an inner twisting $\psi_{\infty}:(I_0)_{\C}\rightarrow I(\infty)_{\C}$ with $I(\infty)/Z(G)$ being compact. Then, any elements $(g_l)_l$ in (ii$'$), $c$ in (iii$'$) determine, for every place $v$ of $\Q$, an inner twisting \[\psi_v:(I_0)_{\Qvb}\rightarrow I(v)_{\Qvb}\] over $\Qvb$ by $\psi_l:=\Int(g_l)$ for $l\neq p$ and $\psi_p$ (\ref{eq:psi_p}) (the morphism with the same notation). When $G_{\gamma_0}$ is connected, the inner class of $\psi_v$ in $H^1(\Qv,I_0^{\ad})$ depends only on $(\gamma_0;\gamma,\delta)$, but in general depends on the choice of $(g_l)_l$, $c$ and there will be no distinguished choice of such inner twisting. Kottwitz considered the following condition:

\begin{itemize} \addtolength{\itemsep}{-4pt} %\label{eq:iv'} \tag{iv$'$}
\item[(iv$'$)] Set $I(v):=H(v)^{\mathrm{o}}$ for each place $v$. There exists a quadruple $(I,\psi,(j_v),\psi_v)$ consisting of a $\Q$-group $I$, an inner twisting $\psi:I_0\rightarrow I$, for each place $v$ of $\Q$, an isomorphism $j_v:I_{\Qv}\rightarrow I(v)$ over $\Q_v$, unramified almost everywhere, and an inner twisting $\psi_v:(I_0)_{\Qvb}\isom I(v)_{\Qvb}$, such that $j_v\circ\psi$ and $\psi_v$ differ by an inner automorphism of $I_0$ over $\Qb_v$. \end{itemize}

We remark that a stable Kottwitz triple will fulfill this condition (iv$'$) if for some choice of $(g_l)_l$, $c$, the associated Kottwitz invariant (to be recalled below) is trivial.

(2) The condition (iii$'$) should be distinguished from the following stronger condition: 
\begin{itemize}
\item[(iii$''$)] there exists $c\in G(\Qpnr)$ (not just in $G(\mfk)$) fulfilling the same condition as (iii$'$) (i.e. $c\gamma_0 c^{-1}=\Nm_n\delta$ and $b:=c^{-1}\delta\sigma(c)\in I_0(\Qpnr)$). 
\end{itemize}
The two conditions (iii$'$), (iii$''$) are the same if $G_{\gamma_0}=I_0$ (by Steinberg's theorem: $H^1(\Qpnr,G_{\gamma_0})=\{1\}$). But, in general, condition (iii$''$) seems to be strictly stronger than condition (iii$'$). We explain why.

As $G_{\Qp}$ is quasi-split, there exists a norm mapping $\mathscr{N}=\mathscr{N}_n$ from $G(L_n)$ to the set of stable conjugacy classes in $G(\Qp)$ \cite[$\S$5]{Kottwitz82}: if $G^{\der}=G^{\uc}$, for any $\delta\in G(L_n)$, the $G(\Qpb)$-conjugacy class of $\Nm_n\delta$, being defined over $\Qp$, contains a rational element (i.e. lying in $G(\Qp)$) by \cite[Thm.4.1]{Kottwitz82} and $\mathscr{N}(\delta)$ is defined to be its stable (=geometric) conjugacy class. For general $G$, one uses a $z$-extension to reduce to the former situation (see \textit{loc. cit.} for a detailed argument). Then, 

($\heartsuit$) \emph{the condition (iii$''$) holds if and only if $\mathscr{N}(\delta)$ is the stable conjugacy class of $\gamma_0$ ($\gamma_0$ being regarded as an element of $G(\Qp)$)}, in which case we say that $\gamma_0$ is the \emph{stable norm} of $\delta$. 

Indeed, it follows immediately from the definition of $\mathscr{N}(\delta)$ that there exists a representative $\gamma_s\in G(\Qp)$ of $\mathscr{N}(\delta)$ which satisfies (together with $\delta$) the condition (iii$''$). As the validity of the condition (iii$''$) depends only on the stable conjugacy class of $\gamma_0$ (or $\gamma_s$), we see that condition (iii$''$) already holds for $\delta$ and any rational representative of $\mathscr{N}(\delta)$; in particular, we see that the implication $\Leftarrow$ holds.
Conversely, if the condition (iii$''$) holds for $\gamma_0\in G(\Qp)$ and $\delta$, $\gamma_0$ is stably conjugate to any (semisimple) rational representative $\gamma_s$ of $\mathscr{N}(\delta)$, because for any choice of $c$, $c_s\in G(\Qpnr)$ satisfying  (iii$''$) for $\gamma_0$, $\gamma_s$ respectively, we have $\gamma_s=\Int(c_s^{-1}c)(\gamma_0)$ and
\[ c^{-1}c_s\sigma(c_s^{-1}c)= c^{-1}c_s\cdot b_s^{-1}c_s^{-1}\delta \cdot \delta^{-1}cb=\Int(c^{-1}c_s)(b_s)\cdot b\in I_0(\Qpnr).\]

Now, given this equivalence, it does not seem that the weaker condition (iii$'$) implies that $\gamma_s$ and $\gamma_0$ are stably conjugate. Indeed, if $c\in G(\mfk)$ is as in the condition (iii$'$), there is an element $g$ in $G(\mfk)$ such that $\gamma_s=\Int(g)(\gamma_0)$ and $g^{-1}\tau(g)\in I_0$ for every $\tau\in W_{\Qp}$ (i.e. $g:=c_s^{-1}c$), but it is not clear whether one can find such $g$ in $G(\Qpb)$ (and satisfying that $g^{-1}\tau(g)\in I_0$ for all $\tau\in \Gal(\Qpb/\Qp)$).

For further results related to this subtle difference which is not used in this work, see \cite[3.5.8]{Lee18a}.

(3) Suppose that $(\gamma_0;\gamma,\delta)$, $(\gamma_0';\gamma',\delta')$ are two triples satisfying the conditions of Definition \ref{defn:stable_Kottwitz_triple} except for (i$'$) and that $\gamma$ and $\gamma'$ are $G(\A_f^p)$-conjugate and $\delta$ and $\delta'$ are $\sigma$-conjugate under $G(L_n)$; so $\gamma_0$ and $\gamma_0'$ are geometrically conjugate. Then, the obstruction to $\gamma_0$ and $\gamma_0'$ being stably conjugate lies in the abelian group
\begin{equation} \label{eq:trivaility_of_Sha^{infty}}
\Sha^{(\infty)}(\Q,\pi_0(G_{\epsilon})):=\ker[H^1(\Q,\pi_0(G_{\epsilon}))\rightarrow \prod_{v\neq\infty}H^1(\Qv,\pi_0(G_{\epsilon}))].
\end{equation}
More precisely, for $g\in G(\Qb)$ such that $\gamma_0'=g\gamma_0g^{-1}$, the cohomology class $[\zeta_{\tau}]\in H^1(\Q,\pi_0(G_{\gamma_0}))$ of the cocycle $\zeta_{\tau}:=(g^{-1}\tau(g))_{\tau}\in Z^1(\Q,G_{\gamma_0})$ lies in this group. At $l\neq p$, $\gamma_0$ and $\gamma_0'$ are stably conjugate as elements $G(\Ql)$, thus if we choose $g_l\in G(\Qlb)$ such that $g_l^{-1}\tau(g_l)\in I_0(\Qlb)$ for all $\tau\in\Gal(\Qlb/\Ql)$,
the image of $[\zeta_{\tau}]$ in $H^1(\Ql,G_{\epsilon})$ equals $[(g_l^{-1}\tau(g_l))_{\tau}]$. At $p$, choose $c,c'\in G(\mfk)$ as in Definition \ref{defn:stable_Kottwitz_triple}, (iii$'$): $c\gamma_0c^{-1}=\Nm_n\delta$, $c'\gamma_0'c'^{-1}=\Nm_n\delta'$ and we obtain cocycles $b_{\tau}:=c^{-1}\Nm_{i(\tau)}\delta\tau(c)$, $b'_{\tau}:=c'^{-1}\Nm_{i(\tau)}\delta\tau(c')$ in $Z^1(W_{\Qp},I_0(\mfk))$ (\ref{eq:b_{tau}}). 
We may assume that $\delta'=\delta$. Then, as $c'g\gamma_0g^{-1}c'^{-1}=c\gamma_0c^{-1}$, we have $h:=c^{-1}c'g\in G_{\gamma_0}(\mfkb)$ and $\zeta_{\tau}=h^{-1}c^{-1}c'\tau(c'^{-1}c)\tau(h)$ for $\tau\in W_{\Qp}$. So, it suffices to show that $c^{-1}c'\tau(c'^{-1}c)\in I_0$ for every $\tau\in W_{\Qp}$. But, since
\[ c^{-1}c'\tau(c'^{-1}c)=c^{-1}c'b_{\tau}'^{-1}c'^{-1}cb_{\tau} \]
and $\Int(c'^{-1}c)=\Int(gh^{-1})$ induces a $\mfk$-isomorphism $I_0\isom I_0'$, the claim follows.

Of course, this obstruction vanishes (so, for stable Kottwitz triples, geometric equivalence becomes the same as stable equivalence) for the groups $G$ with $G^{\uc}=G^{\der}$. But one can show that for Hodge-type Shimura datum $(G,X)$, the obstruction group (\ref{eq:trivaility_of_Sha^{infty}}) is trivial under much less strict assumptions on $G$.
\end{rem}

%%%%%%%%%%%%%%%%%%%%%%%%%%%%%%%%%%%%%%%%
\subsection{Kottwitz invariant} \label{subsec:Kottwitz_invariant}

Here, we extend the definition of Kottwitz invariant which previously (in the context of Shimura varieties, as in \cite[$\S$2]{Kottwitz90}, \cite[$\S$5]{Kottwitz92}) was used under the assumption $G^{\der}=G^{\uc}$.
In this case, the notion of Kottwitz invariant is well-defined by a Kottwitz triple only. For general groups, however, this is not the case any longer, and its natural generalization requires some auxiliary choices; in fact, as we will see, it seems better to regard Kottwitz invariant as being defined on the set $\ker[H^1(\A_f,I_0)\rightarrow H^1(\A_f,G)]$.

Recall our convention for the notation $X_{\ast}(A)$, $X^{\ast}(A)$, and $A^D$: for a locally compact abelian group $A$, they denote the (co)character groups $\Hom(\C^{\times},A)$, $\Hom(A,\C^{\times})$, and the Pontryagin dual group $\Hom(A,S^1)$, respectively. Let $\Gamma:=\Gal(\Qb/\Q)$.

%%%%%%%%%%%%%%%%%%%%
\subsubsection{Invariants refining stable Kottwitz triples} \label{subsubsec:alpha_v}

Let $(\gamma_0;\gamma,\delta)$ be a stable Kottwitz triple; put $I_0:=G_{\gamma_0}^{\mathrm{o}}$. For each place $v$ of $\Q$, we will construct invariants in $\pi_1(I_0)_{\Gamma_p}=X^{\ast}(Z(\hat{I}_0)^{\Gamma(p)})$. This will depend on some auxiliary datum, in addition to $(\gamma_0;\gamma,\delta)$.

At a finite place $l\neq p$, when we choose $g_l\in G(\Qlb)$ such that 
\begin{equation} \label{eq:stable_g_l}
g_l\gamma_0g_l^{-1}=\gamma_l,\ \text{ and }\ g_l^{-1}\tau(g_l) \in I_0(\Qlb)\quad \forall\tau\in\Gamma_l,
\end{equation}
we obtain a cocycle $(g_l^{-1}\tau(g_l))_{\tau}$ in $Z^1(\Ql,I_0)$. We define $\alpha_l(\gamma_0;\gamma_l;g_l)$ to be its cohomology class:
\begin{equation}  \label{eq:alpha_l}
\alpha_l(\gamma_0;\gamma_l;g_l):= [(g_l^{-1} \tau(g_l))_{\tau}]\ \in\ \ker[H^1(\Ql,I_0)\rightarrow H^1(\Ql,G)]. 
\end{equation}
In view of the canonical isomorphism $H^1(\Ql,I_0)\isom\pi_0(Z(\hat{I}_0)^{\Gamma_l})^D$ \cite[Thm.1.2]{Kottwitz86}, we also regard it as a character of $Z(\hat{I}_0)^{\Gamma_l}$ (of finite order). 

In general, $\alpha_l(\gamma_0;\gamma_l;g_l)$ depends on the choice of $g_l$ as well as on the pair $(\gamma_0,\gamma_l)$, while its image in $H^1(\Ql,G_{\gamma_0})$ does not. The set $\ker[H^1(\Ql,G_{\gamma_0})\rightarrow H^1(\Ql,G)]$ classifies the $G(\Ql)$-conjugacy classes in the $G(\Qlb$)-conjugacy class of $\gamma_0\in G(\Ql)$, and its subset
\begin{equation} \label{eq:C_l(gamma_0)}
\mathfrak{C}_l(\gamma_0):= \im\left(\ker[H^1(\Ql,I_0)\rightarrow H^1(\Ql,G)]\rightarrow \ker[H^1(\Ql,G_{\gamma_0})\rightarrow H^1(\Ql,G)]\right)
\end{equation}
classifies the $G(\Ql)$-conjugacy classes in the \emph{stable} conjugacy class of $\gamma_0\in G(\Ql)$.
The invariant $\alpha_l(\gamma_0;\gamma_l;g_l)$ is then a lifting in $\ker[H^1(\Ql,I_0)\rightarrow H^1(\Ql,G)]$  of the $G(\Ql)$-conjugacy class of $\gamma_l$.

At $v=p$, choosing $g_p$ to be $c\in G(\mfk)$ satisfying condition (iii$'$) of Definition \ref{defn:stable_Kottwitz_triple}, we obtain
a basic element $b=b(\gamma_0;\delta;g_p):=g_p^{-1}\delta\sigma(g_p)$ of $I_0(\mfk)$ (Proposition \ref{prop:psi_p}); let $[b]_{I_0}$ denote its $\sigma$-conjugacy class in $B(I_0)$. We define $\alpha_p(\gamma_0;\delta;g_p)$ to be the image of $[b]_{I_0}$ under the canonical map $\kappa_{I_0}:B(I_0)\rightarrow X^{\ast}(Z(\hat{I}_0)^{\Gamma(p)})\cong \pi_1(I_0)_{\Gamma(p)}$ (or, its restriction $B(I_0)_{basic}\isom \pi_1(I_0)_{\Gamma(p)}$ \cite[Prop.6.2]{Kottwitz85}):
\begin{equation}  \label{eq:alpha_p}
\alpha_p(\gamma_0;\delta;g_p):=\kappa_{I_0}([b(\gamma_0;\delta;g_p)]_{I_0})\ \in \pi_1(I_0)_{\Gamma(p)}.
\end{equation}
So, $[b(\gamma_0;\delta;g_p)]_{I_0}\in B(I_0)_{basic}$ and $\alpha_p(\gamma_0;\delta;g_p) \in \pi_1(I_0)_{\Gamma(p)}$ determine each other.

\begin{rem} \label{rem:D_p^{(n)}(gamma_0)}
In general, the $\sigma$-conjugacy class of $b=g_p^{-1}\delta\sigma(g_p)$ in $B(I_0)$ depends on the choice of $g_p$ as well as on the pair $(\gamma_0,\delta)$, while its $\sigma$-conjugacy class in $B(G_{\gamma_0})$ does not. 
Let $\mathfrak{C}_p^{(n)}(\gamma_0)$ denote the set of $\sigma$-conjugacy classes $[\delta]$ of elements in $G(L_n)$ satisfying conditions (iii$'$) and $\ast(\delta)$, and $\mathfrak{D}_p^{(n)}(\gamma_0)$ the subset of $B(I_0)$ consisting of $\sigma$-conjugacy classes $[b=c^{-1}\delta\sigma(c)]_{I_0}$ arising from the pairs $(\delta,c)$ satisfying conditions (iii$'$) and $\ast(\delta)$. 
Then, one knows \cite[Prop.3.5.13]{Lee18a} that with choice of a reference element $[b]_{I_0}$ in $\mathfrak{D}_p^{(n)}(\gamma_0)\subset B(I_0)_{basic}$, the map $\kappa_{I_0}-\kappa_{I_0}([b])$ provides an identification
\[\mathfrak{D}_p^{(n)}(\gamma_0)=\ker[H^1(\Qp,I_0)\rightarrow H^1(\Qp,G)] \]
sending $[b]_{I_0}$ to the distinguished element of the target (whose inverse is given by $j_{[b]}$) and there exists a natural identification
\begin{equation} \label{eq:C_p(gamma_0)}
\mathfrak{C}_p^{(n)}(\gamma_0)=\im\left(\ker[H^1(\Qp,G_{\delta\theta}^{\mathrm{o}})\rightarrow H^1(\Qp,G)]\rightarrow \ker[H^1(\Qp,G_{\delta\theta})\rightarrow H^1(\Qp,R)]\right)
\end{equation}
(recall that $H^1(\Qp,I_0)=H^1(\Qp,G_{\delta\theta}^{\mathrm{o}})$).
\end{rem}

For $v=\infty$, we define $\alpha_{\infty}(\gamma_0)$ to be the image of $\mu_h$ in $\pi_1(I_0)_{\Gamma(\infty)}\cong X^{\ast}(Z(\hat{I}_0)^{\Gamma(\infty)})$ for some $h\in X$ factoring through a maximal torus $T$ of $G$ containing $\gamma_0$ and elliptic over $\R$. This character is independent of the choice of $(T,h)$ \cite[Lem.5.1]{Kottwitz90}.

%%%%%%%%%%%%%%%%%%%%
\subsubsection{Kottwitz invariants} \label{subsubsec:Kottwitz_invariant}

Let $(\gamma_0;\gamma,\delta)$ be a \emph{stable} Kottwitz triple and put $I_0:=G_{\gamma_0}^{\mathrm{o}}$, $\tilde{I}_0:=\rho^{-1}(I_0)$ for the canonical homomorphism $\rho=\rho_G:G^{\uc}\rightarrow G$.
The exact sequence
\[1\rightarrow Z(\hat{G})\rightarrow Z(\hat{I}_0)\rightarrow Z(\hat{\tilde{I}}_0) \rightarrow 1\]
induces a homomorphism \cite[Cor.2.3]{Kottwitz84a} 
\begin{equation} \label{eq:boundary_map_for_center_of_dual} 
\partial: \pi_0(Z(\hat{\tilde{I}}_0)^{\Gamma})\rightarrow H^1(\Q,Z(\hat{G})). 
\end{equation}
Let $\ker^1(\Q,Z(\hat{G}))$ denote the kernel of the localization map $H^1(\Q,Z(\hat{G}))\rightarrow \prod_v H^1(\Q_v,Z(\hat{G}))$. Then we define $\mathfrak{K}(I_0/\Q)$ by
\[\mathfrak{K}(I_0/\Q) :=\{ a\in \pi_0(Z(\hat{\tilde{I}})^{\Gamma}) \ |\ \partial(a)\in \ker^1(\Q,Z(\hat{G}))\}.\]
This is known to be a finite group. Since $\gamma_0$ is elliptic, there is also an identification
\[\mathfrak{K}(I_0/\Q)=\biggl( \bigcap_v Z(\hat{I}_0)^{\Gamma_v}Z(\hat{G}) \biggl)/Z(\hat{G}).\]

The Kottwitz invariant $\alpha(\gamma_0;\gamma,\delta;(g_v)_v)$ will then be a character of $\mathfrak{K}(I_0/\Q)$ which is determined by the stable Kottwitz triple $(\gamma_0;\gamma,\delta)$ and a tuple of elements $((g_l)_{l\neq p},g_p)\in G(\bar{\A}_f^p)\times G(\mfk)$ satisfying (\ref{eq:stable_g_l}) and Definition \ref{defn:stable_Kottwitz_triple} (iii$'$).
It is defined as a product over all places of $\Q$ of the restrictions to $\mathfrak{K}(I_0/\Q)$ of some
characters $\tilde{\alpha}_v(\gamma_0;\gamma,\delta;(g_v)_{v})\in  X^{\ast}(Z(\hat{I}_0)^{\Gamma_v}Z(\hat{G}))$:
\begin{equation} \label{eq:Kottwitz_invariant}
\alpha(\gamma_0;\gamma,\delta;(g_v)_v)=\prod_v \tilde{\alpha}_v|_{\bigcap_v Z(\hat{I}_0)^{\Gamma_v}Z(\hat{G})}
\end{equation}
(this factors through the quotient $( \bigcap_v Z(\hat{I}_0)^{\Gamma_v}Z(\hat{G}) )/Z(\hat{G})$).
For each place $v$, $\tilde{\alpha}_v$ is itself defined as the unique extension of the character $\alpha_v$ on $Z(\hat{I}_0)^{\Gamma_v}$ (introduced in (\ref{subsubsec:alpha_v}) whose restriction to $Z(\hat{G})$ is %%{sign convention}
\begin{equation} \label{eq:restriction_of_alpha_to_Z(hatG)}
\tilde{\alpha}_v|_{Z(\hat{G})}=
\begin{cases}
\quad  \mu_{\natural} &\text{ if }\quad  v=\infty \\
\ -\mu_{\natural} &\text{ if }\quad  v=p \\
\text{ trivial } &\text{ if }\quad  v\neq p,\infty.
\end{cases}
\end{equation}
Such extension of $\alpha_v\in X^{\ast}(Z(\hat{I}_0)^{\Gamma_v}Z(\hat{G}))$ is possible, since $\alpha_v$'s have the same restrictions to $Z(\hat{I}_0)^{\Gamma_v}\cap Z(\hat{G})$ as these characters (\ref{eq:restriction_of_alpha_to_Z(hatG)}) on $Z(\hat{G})$: at $p$, since the image of $\alpha_p$ in $X^{\ast}(Z(\hat{G})^{\Gamma(p)})$ equals that of the $\sigma$-conjugacy class of $b$ in $B(G_{\Qp})$ under $\kappa_{G_{\Qp}}$, it is equal to $-\mu^{\natural}$ by condition $\ast(\delta)$.
It also follows that the product $\prod_v \tilde{\alpha}_v$ is trivial on $Z(\hat{G})$.

A proper notation for $\tilde{\alpha}$ will be
\[ \tilde{\alpha}_l(\gamma_0;\gamma_l;g_l)\ (l\neq p,\infty),\quad \tilde{\alpha}_p(\gamma_0;\delta;g_p),\quad \tilde{\alpha}_{\infty}(\gamma_0).\]
But, here and thereafter, to have uniform notation for all $v$'s, we will write $\alpha_v(\gamma_0;\gamma,\delta;g_v)$ and $\tilde{\alpha}_v(\gamma_0;\gamma,\delta;g_v)$, although for specific $v$, these depend only on either $\gamma_l\ (l\neq p)$ or $\delta$, and $g_v$.

Clearly, the Kottwitz invariant $\alpha(\gamma_0;\gamma,\delta;(g_v)_v)$ is determined by the adelic class 
\[ (\alpha_v(\gamma_0;\gamma,\delta;g_v))_v\in \ker[H^1(\A_f,I_0)\rightarrow H^1(\A_f,G)]. \] 
(once we fix a reference $\sigma$-conjugacy class $[b]=[g_p^{-1}\delta\sigma(g_p)]\in B(I_0)$). 

Later, we will often use the expression ``\textit{stable Kottwitz triple with trivial Kottwitz invariant}''. By this, we will mean that for the given stable Kottwitz triple $(\gamma_0;(\gamma_l)_{l\neq p},\delta)$ there exist elements $(g_v)_v\in G(\bar{\A}_f^p)\times G(\mfk)$ satisfying conditions (\ref{eq:stable_g_l}) and Definition \ref{defn:stable_Kottwitz_triple} (iii$'$) such that the associated Kottwitz invariant $\alpha(\gamma_0;\gamma,\delta;(g_v)_v)$ vanishes.

%%%%%%%%%%%%%%%%%%%%
%%%%%%%%%%%%%%%%%%%%
\subsection{Langlands-Rapoport condition \ref{itm:ast(epsilon)} and admissible Frobenius pairs} \label{subsec:LR-cond.ast(epsilon)}

The next condition on an element $\gamma_0\in G(\Q)$ generalizes a condition introduced by Langlands and Rapoport \cite[p.183]{LR87} (of the same name $\ast(\epsilon)$), especially to groups whose derived groups are not necessarily simply connected:
\begin{enumerate} [label=$\ast(\epsilon)$]
\item \label{itm:ast(epsilon)} 
Let $H$ be the centralizer in $G_{\Qp}$ of the maximal $\Qp$-split torus in the center of $(G_{\gamma_0})_{\Qp}$.
Then, there exists a cocharacter of $\mu$ of $H$ lying in the $G(\Qpb)$-conjugacy class $\{\mu_X\}$
such that
\begin{equation} \label{eq:lambda(gamma_0)}
w_H(\gamma_0)=-\sum_{i=1}^{n}\sigma^{i-1}\underline{\mu} %%{sign convention}
\end{equation}
for some $n\in\N$, where $w_H:H(\mfk)\rightarrow \pi_1(H)_I$ is the map from (\ref{subsubsec:Kottwitz_hom}) and $\underline{\mu}$ denotes the image of $\mu$ in $\pi_1(H)_I$ ($I=\Gal(\Qpb/\Qpnr)$).
\end{enumerate}

We refer to the number $n$ appearing in \ref{itm:ast(epsilon)} as the \textit{level} of this condition (e.g., we will say that condition \ref{itm:ast(epsilon)} holds for $\gamma_0$ with level $n$).
There are two minor differences between the original condition and ours: first there is the sign difference (which results from our sign convention) and secondly, unlike \cite{LR87}, we do not require $\mu$ to be defined over $L_n$ (or even over an unramified extension of $\Qp$).
When $G_{\Qp}$ is unramified, so is $H$, hence $\pi_1(H)_I=\pi_1(H)$.

%%%%%%%%%%%%%%%%%%%%
\begin{lem} \label{lem:invariance_of_(ast(gamma_0))_under_transfer_of_maximal_tori}
Let $\gamma_0\in G(\Q)$ be a semisimple element. 

(1) Suppose $\gamma_0\in T(\Q)$ for a maximal torus $T\subset G$. Then condition \ref{itm:ast(epsilon)} holds for $\gamma_0$ if and only if it holds for $\gamma_0':=\Int g(\gamma_0)$ for any transfer of maximal torus $\Int g:T\hra G$.

(2) If $\gamma_0':=\Int g(\gamma_0)\in G(\Q)$ for some $g\in G(\Qb)$, for each place $v$ of $\Q$, the image of $\gamma_0$ in $G^{\ad}(\Ql)$ lies in a compact open subgroup of $G^{\ad}(\Ql)$ if and only if $\gamma_0'$ is so.
\end{lem}

See \cite[$\S$9]{Kottwitz84a} for the notion of ``\emph{transfer} (or \emph{admissible embedding}) \emph{of maximal torus}''.

\begin{proof}
(1) Let $H$ be the centralizer of the maximal $\Qp$-split torus in the center of $(G_{\gamma_0})_{\Qp}$ and $H'$ the similarly defined group for $\gamma_0'$; one has $H'=\Int(g)(H)$ as $\Int(g)$ induces a $\Q$-isomorphism between the centers of $G_{\gamma_0}$ and $G_{\gamma_0'}$. Choose $\mu\in X_{\ast}(T)\cap \{\mu_X\}$ which is conjugate under $H$ to a cocharacter satisfying condition \ref{itm:ast(epsilon)}.
Let $K$ be a finite Galois extension of $\Qp$ over which $\mu$ is defined, $K_0\subset K$ its maximal unramified sub-extension, and choose a uniformizer $\pi$ of $K$.
In view of the equality 
\begin{equation} \label{eq:Kottwitz97-(7.3.1)}
w_{H}(\Nm_{K/K_0}(\mu(\pi)))=\underline{\mu},
\end{equation}
where $\underline{\mu}$ is the image of $\mu$ in $\pi_1(H)_{\Gal(\Qpb/\Qpnr)}$ (commutativity of diagram (7.3.1) of \cite{Kottwitz97}), condition \ref{itm:ast(epsilon)} holds for $\gamma_0$ if and only if
\begin{equation} \label{eq:ast(gamma_0)}
\gamma_0\cdot\prod_{i=1}^n\sigma^{i-1}\Nm_{K/K_0}(\mu(\pi))\in \ker(w_H)\cap T(\mfk). %%{sign convention}
\end{equation}
By the Steinberg's theorem $H^1(\Qpnr,T)=0$, we may find $g_p\in G(\Qpnr)$ such that $\Int(g)|_{T_{\Qp}}=\Int(g_p)|_{T_{\Qp}}$.
Hence, it follows from funtoriality (with respect to $\Int(g_p):H\rightarrow H'=\Int(g_p)(H)$) of the functor $w_H$ \cite[$\S$7]{Kottwitz97} that (\ref{eq:ast(gamma_0)}) holds if and only if the same condition holds for $(\gamma_0',T',\mu')=\Int(g_p)(\gamma_0,T,\mu)$.

(2) Let $P$ and $P'$ be the $\Q$-subgroups of $G$ generated by $\gamma_0$, $\gamma_0'$, respectively.
Then, the $\Qb$-isomorphism $\Int g:G_{\Qb}\isom G_{\Qb}$ restricts to a $\Q$-isomorphism $P\isom P'$, thus if the image of $\gamma_0$ in $G^{\ad}(\Ql)$ lies in a compact open subgroup of $G^{\ad}(\Ql)$, then as it also lies in a compact open subgroup of $Q(\Ql)$, where $Q$ is the image of $P$ in $G^{\ad}$, the same property holds for $\gamma_0'$.
\end{proof}

%%%%%%%%%%%%%%%%%%%%
\begin{defn} \label{defn:stable_conj_MT,SD}
(1) For a $\Q$-torus $T$ and a connected reductive $\Q$-group $G$ with $\Qb\operatorname{-}\mathrm{rk}(T)=\Qb\operatorname{-}\mathrm{rk}(G)$, a \emph{stable conjugacy class of $\Q$-embeddings} $T\hookrightarrow G$ is, by definition, an equivalence class of $\Q$-embeddings $T\hookrightarrow G$ with respect to stable conjugacy relation: two $\Q$-embeddings $i_T,i_T':T\hookrightarrow G$ are \emph{stably conjugate} if and only if there exists $g\in G(\Qb)$ such that $i_T'=\Int(g)\circ i_T$ (in particular, $\Int(g)|_{i_T(T)}$ induces a transfer of maximal torus $i_T(T)\isom i_T'(T)$). 

(2) For a Shimura datum $(G,X)$, a \emph{stable conjugacy class of special Shimura sub-data} $(T,h)$ is, by definition, an equivalence class of special Shimura subdatum $(T,h)$ with respect to the following (stable conjugacy) equivalence relation: two special Shimura sub-data $(T,h)$, $(T',h')$ are \emph{stably conjugate} if and only if there exist $g\in G(\Qb)$ inducing a $\Q$-isomorphism $\Int (g)|_T:T\isom T'$ such that $h'=\Int(g)(h)$.
\end{defn}

We note that for $g\in G(\Qb)$ as in (2), since $h'=\Int(g_{\infty})(h)$ for some $g_{\infty}\in G(\R)$, we see that 
$g_{\infty}^{-1} g\in K_{\infty}(\C)$ for the centralizer $K_{\infty}$ of $h$ (real aglebraic group which is a compact-mod-center inner form of $G_{\R}$). As the map $H^1(\R,T)\rightarrow H^1(\R,K_{\infty})$ is injective \cite[4.4.5]{Kisin17}, we see that the class $[g^{-1}\cdot{}^{\tau}g]\in H^1(\R,T)$ is trivial, and so $\Int(g)|_{T_{\R}}=\Int(g')_{T_{\R}}$ for some $g'\in G(\R)$.

We put $\tilde{\mathbf{K}}_p:=G_{\Zp}(W(\F))$, the hyperspecial subgroup of $G(\mfk)$ corresponding to $\mathbf{K}_p$ (i.e. $G_{\Zp}$ is the reductive group scheme over $\Zp$ with generic fiber $G_{\Qp}$ such that $\mathbf{K}_p=G_{\Zp}(\Zp)$), and 
\begin{equation}  \label{eq:Adm_K(mu)}
 \mathrm{Adm}_{\tilde{\mathbf{K}}_p}\bigl(\{\mu_X^{-1}\}\bigr):=\tilde{\mathbf{K}}_p \mu(p^{-1}) \tilde{\mathbf{K}}_p.
\end{equation}
the $\tilde{\mathbf{K}}_p$-double coset containing $\mu(p^{-1})$ for any unramified representative $\mu$ of the conjugacy class $\{\mu_X\}$.

%%%%%%%%%%%%%%%%%%%%
\begin{defn} \label{defn:Frobenius_pair}
(1) A \emph{Frobenius pair} (F-pair, for short) of level $n=m[\kappa(\wp):\Fp]\  (m\in\N)$ is a pair $(h,\gamma_0)$ consisting of an element $\gamma_0$ of $G(\Q)$ and $h\in X\cap\Hom(\dS,(G_{\gamma_0})_{\R})$ such that $(\gamma_0,\mu_h)$ satisfies equation (\ref{eq:lambda(gamma_0)}) of (\ref{itm:ast(epsilon)}) with this $n$, and there exists an $\R$-elliptic, maximal $\Q$-torus $T$ of $G$ which $h$ factors through and contains $\gamma_0$.

Two Frobenius pairs $(h,\gamma_0)$, $(h',\gamma_0')$ are said to be \emph{equivalent} if there exist an $\R$-elliptic, maximal $\Q$-torus $T$ of $G$ which $h$ factors through and contains $\gamma_0$, and a transfer of maximal torus $\Int(g):T\hookrightarrow G\ (g\in G(\Qb))$ (i.e. $\Int(g)|_{T_{\Qb}}$ defines a $\Q$-rational embedding $T\rightarrow G$) such that $h'=\Int(g)(h)$ and $\gamma_0'=\Int(g)(\gamma_0)$.

(2) A Frobenius pair $(h,\gamma_0)$ is said to be \emph{admissible} of level $n=m[\kappa(\wp):\Fp]\  (m\in\N)$ if for some $\R$-elliptic, maximal $\Q$-torus $T$ of $G$ which $h$ factors through and contains $\gamma_0$, there exists $x_p\in G(\mfk)/\tilde{\mathbf{K}}_p$ such that 
\[\gamma_0\cdot x_p=(\Nm_{K/K_0}(\mu_h(\pi))^{-1}\sigma)^n\cdot x_p,\]
where $K$ is a finite Galois extension of $\Qp$ splitting $T$ with a uniformizer $\pi$, $K_0$ is its maximal subfield unramified over $\Qp$: this condition depends only on the F-pair $(h,\gamma_0)$ (again, we regard $(\Nm_{K/K_0}(\mu_h(\pi))^{-1}\sigma$ as an element of $T(\Qpnr)\rtimes\langle\sigma\rangle$, which acts on $G(\mfk)/\tilde{\mathbf{K}}_p$).

(3) In (2), if one can find such $x_p$ from the set 
\begin{equation} \label{eq:X_p(b,epsilon')}
X_p(b,\gamma_0):=\{ x\in X(\{\mu_X^{-1}\},b)_{\mathbf{K}_p}\ |\ \gamma_0 x=(b\sigma)^nx \},
\end{equation}
where $b:=\Nm_{K/K_0}(\mu_h(\pi))^{-1}$ and 
\begin{equation} \label{eq:X_p(mu_X,b)}
X(\{\mu_X^{-1}\},b)_{\mathbf{K}_p}:=\left\{\ x\in G(\mfk)/\tilde{\mathbf{K}}_p\ \Big| \ \mathrm{inv}_{\tilde{\mathbf{K}}_p}(x,b\sigma x)\in \mathrm{Adm}_{\tilde{\mathbf{K}}_p}\bigl(\{\mu_X^{-1}\}\bigr) \ \right\},
\end{equation}
then we say that $(h,\gamma_0)$ is \emph{$\mathbf{K}_p$-effective}. Here, $\mathrm{inv}_{\tilde{\mathbf{K}}_p}$ is the map
\[ \mathrm{inv}_{\tilde{\mathbf{K}}_p}:G(\mfk)/\tilde{\mathbf{K}}_p \times G(\mfk)/\tilde{\mathbf{K}}_p \rightarrow \tilde{\mathbf{K}}_p\backslash G(\mfk)/\tilde{\mathbf{K}}_p:\ (x\tilde{\mathbf{K}}_p,y\tilde{\mathbf{K}}_p)\mapsto \tilde{\mathbf{K}}_p x^{-1}y \tilde{\mathbf{K}}_p.\]
\end{defn}
%Note that any rational element $\epsilon$ of $I_{\mathscr{I}}$ is semisimple).

To see that the condition in (2) is also independent of the choice of a torus $T$ as in that statement,
we note (using (\ref{eq:Kottwitz97-(7.3.1)}) and \cite[2.8]{Kottwitz85}) that the image of $\Nm_{K/K_0}(\mu_h(\pi))^{-1}$ under the injective map
$\kappa_{G}\times \bar{\nu}_{G}: B(G) \rightarrow \pi_1(G)_{\Gamma_p}\times\mathcal{N}(G)$  (\ref{eq:kappa_times_nu}) is $(-\underline{\mu_h}_G,\frac{-1}{[K:\Qp]}\sum_{\tau\in \Gal(K/\Qp)}\tau\mu_h)$, where the average is well-defined independent of the torus $T$ (containing the image of $\mu_h$) and $K$.

%%%%%%%%%%%%%%%%%%%%
\begin{rem}
The notion of Frobenius pair was first introduced by Langlands \cite{Langlands76}, \cite[p.47]{Langlands79} and later taken up again by Zink \cite[4.6]{Zink83}) for the same goal of parameterizing isogeny classes over $\Fpb$.
However, our definition of Frobenius pair differs from that of Langlands, in two respects: first, Langlands considers a slightly stronger condition \cite[p.413-414]{Langlands76} than ours \ref{itm:ast(epsilon)}, and his equivalence relation among Frobenius pairs is different from the above one (it is defined by everywhere local conditions for the Kottwitz triples that Frobenius pairs produce \cite[p.417]{Langlands76}).
\end{rem}

%%%%%%%%%%%%%%%%%%%%  
%%%%%%%%%%%%%%%%%%%%
\begin{thm} \label{thm:LR-Satz5.21}
Assume that $G_{\Qp}$ is unramified and $G$ is of classical Lie type. Let $\mathbf{K}_p$ be a hyperspecial subgroup of $G(\Qp)$ and $\gamma_0$ an element of $G(\Q)$ which is $\R$-elliptic. 

(1) If $\gamma_0$ satisfies condition \ref{itm:ast(epsilon)}  of (\ref{subsubsec:pre-Kottwitz_triple}), there exists an admissible Frobenius pair $(T,h,\epsilon)$ such that $\epsilon$ is stably conjugate to $\gamma_0$ (in fact, $\epsilon$ is the image of $\gamma_0$ under a transfer of a maximal torus $T_0$ containing $\gamma_0$ such that $(T_0)_{\Qp}$ is elliptic in $(G_{\gamma_0}^{\mathrm{o}})_{\Qp}$). For some $t\in\N$, the admissible F-pair $(h,\epsilon^t)$ is also $\mathbf{K}_p$-effective.

(2) If there exists $\delta\in G(L_n)$ (with $[\kappa(\wp):\Fp]|n$) such that $\gamma_0$ and $\Nm_n(\delta)$ are $G(\mfk)$-conjugate 
and the set 
\begin{equation} \label{eq:Y_p(delta)}
Y_p(\delta):=\left\{\ x\in G(\mfk)/\tilde{\mathbf{K}}_p\ \Big| \ \sigma^nx=x,\ \mathrm{inv}_{\tilde{\mathbf{K}}_p}(x,\delta\sigma x)\in \mathrm{Adm}_{\tilde{\mathbf{K}}_p}\bigl(\{\mu_X^{-1}\}\bigr)\ \right\}
\end{equation}
is non-empty, condition \ref{itm:ast(epsilon)}  holds for the stable conjugacy class of $\gamma_0$ with level $n$.\end{thm}

%%%%%%%%%%%%%%%%%%%%
\begin{rem} \label{rem:criterion_for_ast(epsilon)}
(1) In the setup of Langlands and Rapoport \cite{LR87} that $G^{\der}=G^{\uc}$, statement (1) is ibid., Satz 5.21, except that it is formulated in terms of Galois gerb language. The assertion (1) and its proof here is a translation (with certain generalization in our setup) of the original statement and proof using only more standard notions in Galois cohomology of algebraic groups.
The statement (2) appears for the first time in \cite[Thm.6.2.11]{Lee18a}. This is a key to our proof of the effectivity criterion (\ref{eq:E}) (of Kottwitz triple), and is also the (only) missing ingredient in the arguments in \cite{LR87} deducing (the constant coefficient case of) the Kottwitz formula (Theorem \ref{thm_intro:Kottwitz_formula}) from their conjecture \cite[$\S$5]{LR87}.%%
\footnote{More precisely, the arguments in \cite{LR87} for that are essentially those sketched in the introduction (with isogeny class replaced by admissible morphism). Langlands and Rapoport showed the cardinality statement (\ref{eq:C}), but only some weaker statement of (\ref{eq:E}) \cite[Satz 5.25]{LR87} which uses condition $\ast(\epsilon)$ instead of the more natural one of non-vanishing of $O_{\gamma}(f^p)\cdot TO_{\delta}(\phi_p)$: they did not relate these two conditions, and as such their deduction arguments in \cite{LR87} remained incomplete. See \cite[6.1.8]{Lee18a} for some misconception around this matter.}
%%
%(cf. \cite[Conj. 3.4.1]{Lee18a}).

(2) For $b\in G(\mfk)$ and $\gamma_0\in J_b(\Qp)$, if the set $X_p(b,\gamma_0)$ (\ref{eq:X_p(b,epsilon')}) contains an element $c\in G(\mfk)$, then $\delta:=cb\sigma(c^{-1})$ belongs to $G(L_n)$ and satisfies $c\gamma_0 c^{-1}=\Nm_n\delta$, and the mapping $x\mapsto cx$ gives a bijection 
\[X_p(b,\gamma_0)\rightarrow Y_p(\delta). \]
Conversely, for a pair $(\gamma_0,\delta)\in G(\mfk)\times G(L_n)$ satisfying the relation $c\gamma_0 c^{-1}=\Nm_n\delta$ for some $c\in G(\mfk)$, we have $\gamma_0\in J_b(\Qp)$ for 
$b:=c^{-1}\delta\sigma(c)$ (since $\delta\sigma$ commutes with $\Nm_n\delta$), and the mapping $x\mapsto c^{-1}x$ gives a bijection $Y_p(\delta)\rightarrow X_p(b,\gamma_0)$. 

(3). To make comparison with the original source \cite[5.21]{LR87} (which is reproduced in \cite{Lee18a}),
in this proof, we decided to go back to the original sign convention, which means that now (\ref{eq:lambda(gamma_0)}) has no minus sign, replacing $\mathrm{Adm}_{\tilde{\mathbf{K}}_p}(\{\mu_X^{-1}\})$ by $\mathrm{Adm}_{\tilde{\mathbf{K}}_p}(\{\mu_X\})$.
\end{rem}

\begin{proof}
(1) Set $I_0:=G_{\gamma_0}^{\mathrm{o}}$;
by the well-known fact that every semisimple element in a connected reductive group lies in a maximal torus, we have $\gamma_0\in I_0(\Q)$.
We take a maximal $\Q$-torus $T_0$ of $I_0$ such that $(T_0)_{\R}$ is elliptic in $G_{\R}$ and $(T_0)_{\Qp}$ is elliptic in $(I_0)_{\Qp}$. To see that such torus exists, we choose a maximal $\Q$-torus $T_{\infty}$ of $G$ that contains $\gamma_0$ (so $T_{\infty}\subset I_0$) and is elliptic in $G$ over $\R$ (which exists as $\gamma_0$ is elliptic over $\R$). We also choose an elliptic maximal $\Qp$-torus $T_p$ of $(I_0)_{\Qp}$ (which exists as $(I_0)_{\Qp}$ is reductive, \cite[Thm.6.21]{PR94}). Then, one can deduce (cf. Step 1 of the proof of Theorem 4.1.1 of \cite{Lee18b}) that there exists a maximal $\Q$-torus $T_0$ of $I_0$ which is $I_0(\Q_v)$-conjugate to $T_v$ for each $v=\infty$ and $p$. 

Let $H$ be the centralizer of the maximal $\Qp$-split torus in the center of $(G_{\gamma_0})_{\Qp}$ (so, $(G_{\gamma_0})_{\Qp}\subset H$). Let $\mu_0$ be a cocharacter of $T_0$ that is conjugate under $H(\Qpb)$ to some $\mu$ satisfying condition \ref{itm:ast(epsilon)} . Clearly, condition \ref{itm:ast(epsilon)}  still holds for $\mu_0$.
Moreover, as $(T_0)_{\Qp}$ is elliptic in $(I_0)_{\Qp}$, we see that when $K$ is a finite Galois extension of $\Qp$ splitting $T_0$, the following property holds:

($\dagger$) \emph{$\Nm_{K/\Qp}(\mu_0)$ ($\Qp$-rational cocharacter of $T_0$) maps into the maximal $\Qp$-split torus in the center of $I_0$, thus into the maximal $\Qp$-split torus in the center of $H$.}

Then, with such $(T_0,\mu_0)$, one applies the argument of proof of Lemma 5.12 \cite{LR87} (which continues to hold in our general setup without any change; or, see Step 2 of the proof of \cite[Thm. 3.2.1]{Lee18b}, where this argument is reproduced) to find a transfer of maximal torus $\Int g_0:T_0\hra G$ such that $\Int(g_0)(\mu_0)=\mu_h$ for some special Shimura subdatum $(T:=\Int(g_0)(T_0),h)$; we denote $\Int(g_0)(\gamma_0)$ by $\epsilon$. 
By Lemma \ref{lem:invariance_of_(ast(gamma_0))_under_transfer_of_maximal_tori}, conditions \ref{itm:ast(epsilon)}  and ($\dagger$) continue to hold for the new object $(T,\mu_h,\epsilon)$ which thus comes from a Frobenius pair $(T,h,\epsilon)$.
Now, we check that the Frobenius pair 
\[ (h,\epsilon) \]
is admissible and also the pair $(h,\epsilon^t)$ is $\mathbf{K}_p$-effective admissible for some $t\in\N$. We need some preparation.

Since $H$ is a semi-standard $\Qp$-Levi subgroup of quasi-split $G_{\Qp}$, there exists $g\in G(\Qp)$ such that $H(\Qp)\cap {}^g\mathbf{K}_p$ is a hyperspecial subgroup of $H(\Qp)$. Indeed, by \cite[4.2.3, (3)]{Lee18a},  it is enough to find $g\in G(\Qp)$ and a maximal split $\Qp$-torus $S$ of $G_{\Qp}$ contained in $H$ such that the apartment $\mathcal{A}(\Int g^{-1}(S),\Qp)$ contains a hyperspecial point in $\mathcal{B}(G,\Qp)$ giving $\mathbf{K}_p$: this is easy since such a hyperspecial point lies in some apartment of $G_{\Qp}$ and the apartments of $G_{\Qp}$ (which correspond to maximal split $\Qp$-torus of $G$) are permuted under $G(\Qp)$.
Then, we apply Proposition \ref{lem:existence_of_aniostropic_and_unramified_torus} to $H(\Qp)\cap {}^g\mathbf{K}_p$ to choose an unramified, elliptic maximal $\Qp$-torus $T'$ of $H$ such that the (unique) maximal bounded subgroup $T'(\mfk)_1$ of $T'(\mfk)$ is contained in $T'(\mfk)\cap {}^g\tilde{\mathbf{K}}_p$; this Proposition is also stated without proof in \cite[p.172]{LR87}.
We choose $j\in H(\Qpb)$ with $T'=\Int (j)(T_{\Qp})$, and set
\[ \mu':=\Int (j)(\mu_h)\ \ \in X_{\ast}(T').\]
We take $K$ to be big enough so that $K$ splits $T'$ (as well as $T$). Fix a uniformizer $\pi$ of $K$ and let $K_0$ be the maximal subfield of $K$ unramified over $\Qp$. We put
\begin{equation} \label{eq:xi_p'}
F:=\Nm_{K/K_0}(\mu'(\pi)) \sigma\in T'(\Qpnr)\rtimes\langle\sigma\rangle.
\end{equation}

%%%%%%%%%%%%%%%%%%%%
\begin{lem} \label{lem:Phi_for_special_morphism}
(1) Put $\nu_p':=-\Nm_{K/\Qp}\mu'\in X_{\ast}(T')^{\Gamma_p}$. For $j\in \N$ divisible by $[K:\Qp]$, we have 
\[ F^j=(p^{-\nu_p'}\cdot u_0)^{\frac{j}{[K:\Qp]}}\sigma^{j} \]
for some $u_0\in T'(\Qp)_1(:=\ker (w_{T'})\cap T'(\Qp))$. In particular, $b_j:=F^j\sigma^{-j}\in T'(\Qp)$.

(2) The following two $\Qp$-rational cocharacters of $H$ are equal and map into $Z(H)$:
\[ \Nm_{K/\Qp}\mu'=\Nm_{K/\Qp}\mu_h.\] 
There exists $v\in H(\mfk)$ such that
\[ \Nm_{K/K_0}(\mu'(\pi))=v \Nm_{K/K_0}(\mu_h(\pi))\sigma(v^{-1}) .\]

(3) Put $\epsilon':=\Int (v)(\epsilon)\in H(\mfk)$; so, $\epsilon'$ commutes with $F$. 
For $t\in\N$, define an element $k_t\in H(\mfk)$ by
\[ k_t  \sigma^{nt} := (\epsilon'^{-1}F^n)^{t}. \]
Then, for any neighborhood $V$ of $1$ in $H(\mfk)$, there exists $t\in\N$ with $k_t\in V$.
\end{lem}

\begin{proof}
(1) Let $e_K:=[K:K_0]$, $f_K:=[K_0:\Qp]$. 
Put $u:=\pi^{e_K}p^{-1}\in \cO_K^{\times}$ and $t':=j/[K:\Qp]$. We express $F^j$ in terms of $\nu_p'$:
\begin{align*}
F^{j}&=(\Nm_{K/K_0}(\mu'(\pi)) \sigma)^j=\Nm_{K/\Qp}(\mu'(\pi^{e_K}))^{t'} \sigma^{j} \\
&=(p^{\Nm_{K/\Qp}\mu'}\cdot \Nm_{K/\Qp}(\mu'(u)))^{t'} \sigma^{j} \\
&=p^{-t'\nu_p'}\cdot u_0^{t'} \sigma^{j}
\end{align*}
where $u_0:=\Nm_{K/\Qp}(\mu'(u))$. Clearly, $u_0\in T'(\Qp)_0(:=\ker(v_{T'})\cap T'(\Qp))$, the maximal compact subgroup of $T'(\Qp)$. Since $T'$ is unramified over $\Qp$, this subgroup equals  $T'(\Qp)_1(:=\ker (w_{T'})\cap T'(\Qp))$, the maximal bounded subgroup of $T'(\Qp)$; in fact, it is not difficult to see that one always has $u_0\in T'(\Qp)_1$ without the unramifiedness condition.

(2) We first see that the two cocharacters both map into the center of $H$: for $\Nm_{K/\Qp}\mu_h$, this is by definition of $H$ (property ($\dagger$)), while $\nu_p'=-\Nm_{K/\Qp}\mu'$ maps into a $\Qp$-split sub-torus of the elliptic maximal torus $T'$ of $H$, so factors through the center $Z(H)$. Since their projections into $H^{\ab}=H/H^{\der}$ are also the same (as $\mu'$, $\mu_h$ are conjugate under $H(\Qpb)$), they must be equal themselves.
These are the Newton cocharacters attached to $b:=\Nm_{K/K_0}(\mu_h(\pi))$, $b':=\Nm_{K/K_0}(\mu'(\pi))$ as $\kappa_H(b)=\underline{\mu_h}$, $\kappa_H(b')=\underline{\mu'}$ by (\ref{eq:Kottwitz97-(7.3.1)}) \cite[(2.8.1)]{Kottwitz85}. Hence their $\sigma$-conjugacy classes in $B(H)$ are both basic and are equal, as $\underline{\mu_h}=\underline{\mu'}$ \cite[5.6]{Kottwitz85}.

(3) Let $Z_{\epsilon}:=Z(G_{\epsilon})$ and $Z_{\epsilon}^{\mathrm{o}}$ be the center of $G_{\epsilon}$ and its neutral component, respectively. 
Then, we claim that
\[k_0:=\epsilon^{-[K:\Qp]}\cdot p^{n\Nm_{K/\Qp}\mu_h}\]
lies in $\ker(v_H)\cap Z_{\epsilon}(\Qp)$, and for some $a\in\N$,
\begin{equation} \label{eq:k_0_is_bounded}
k_0^{a}\in \ker(v_{Z_{\epsilon}^{\mathrm{o}}}).
\end{equation}
In particular, $k_0^a$ lies in the maximal compact subgroup of $T(\Qp)$.

Indeed, condition \ref{itm:ast(epsilon)}  implies (the first line of) the following equality:
\begin{align} \label{eq:v_H(k_0)=0}
[K:\Qp]w_{H}(\epsilon) &=\sum_{j=1}^{[K:\Qp]}\sigma^{j-1} \left(\sum_{i=1}^n\sigma^{i-1}\underline{\mu'} \right) \\
&=\sum_{j=1}^{[K:\Qp]}\sigma^{j-1}\left(\sum_{i=1}^n\sigma^{i-1}w_{H}(\Nm_{K/K_0}(\mu'(\pi))) \right) \nonumber\\
&=nw_{H}(\prod_{j=1}^{[K:\Qp]}\sigma^{j-1}(\Nm_{K/K_0}(\mu'(\pi))) ) \nonumber \\
&=w_{H}(p^{n\Nm_{K/\Qp}\mu_h}). \nonumber
\end{align}
Here, the first equality also uses that $\epsilon\in H(\Qp)$  (so $\sigma(w_H(\epsilon))=w_H(\epsilon)$).
The second equality holds by (\ref{eq:Kottwitz97-(7.3.1)}), and the next two equalities follow from Lemma \ref{lem:Phi_for_special_morphism} (1): $\prod_{j=1}^{[K:\Qp]}\sigma^{j-1}(\Nm_{K/K_0}(\mu'(\pi)))=\Nm_{K/\Qp}(\mu'(pu))\in T'(\Qp)$ ($\pi^{e_K}=pu$), and $\Nm_{K/\Qp}(\mu'(u))=u_0\in \ker (w_{T'})$. 

Now, by property ($\dagger$), $\Nm_{K/\Qp}\mu_h$ maps into $Z_{\epsilon}^{\mathrm{o}}$, so one has $k_0\in Z_{\epsilon}(\Qp)$. Also, the equality (\ref{eq:v_H(k_0)=0}) above shows that $k_0\in \ker(v_H)$. 
To see the property (\ref{eq:k_0_is_bounded}),
we first take $a_1\in\N$ with $\epsilon^{a_1}\in Z_{\epsilon}^{\mathrm{o}}(\Q)$, so that $k_0^{a_1}\in Z_{\epsilon}^{\mathrm{o}}(\Qp)$.
Let $A_{\epsilon}$ and $B_{\epsilon}$ be respectively the isotropic kernel and the anisotropic kernel of $Z_{\epsilon}^{\mathrm{o}}$. Then, there exists $a_2\in\N$ such that $k_0^{a_1a_2}=x\cdot y$ with $x\in A_{\epsilon}(\Qp)$ and $y\in B_{\epsilon}(\Qp)$. 
Since $v_H=v_{H^{\ab}}\circ p_H$ for the quotient map $p_H:H\rightarrow H^{\ab}$, there are the implications:
\[v_{B_{\epsilon}}(y)=0\ \Rightarrow\ v_{H^{\ab}}(p_H(y))=0\ \Rightarrow\ v_{H^{\ab}}(p_H(x))=0.\]
Then, since $A_{\epsilon}\subset Z(H)$ and the natural map $Z(H)\rightarrow H^{\ab}$ is an isogeny, it follows that $v_{A_{\epsilon}}(x)=0$. Therefore, one has $k_0^{a}\in \ker(v_{Z_{\epsilon}^{\mathrm{o}}})$ for $a:=a_1a_2$; in particular, $k_0^a$ lies in a compact (thus bounded) subgroup of $T(\Qp)$. We also note in passing that when $G^{\der}=G^{\uc}$ (so that $H^{\der}=H^{\uc}$ too), the same argument establishes that 
$k_0^{a}\in \ker(w_{Z_{\epsilon}^{\mathrm{o}}})$ (since in such case $w_H=w_{H^{\ab}}\circ p_H$).

Finally, the claim (3) follows from the property (\ref{eq:k_0_is_bounded}), in view of the equality (for $t\gg1$):
\begin{align*} \label{eqn:epsilon'^{-1}F^n}
k_t  &= v \epsilon^{-t} v^{-1} \cdot (p^{-\nu_p'}\cdot u_0)^{\frac{nt}{[K:\Qp]}} \nonumber \\
&= v (\epsilon^{-[K:\Qp]}\cdot p^{n\Nm_{K/\Qp}\mu_h})^{\frac{t}{[K:\Qp]}} v^{-1} \cdot u_0^{\frac{nt}{[K:\Qp]}} \nonumber \\
&= v k_0^{\frac{t}{[K:\Qp]}} v^{-1}\cdot (u_0)^{\frac{nt}{[K:\Qp]}}.\quad \qedhere \nonumber
\end{align*}
\end{proof}

It follows from this lemma that for sufficiently large $t$, $k_t$ lies in the hyperspecial subgroup $H(\mfk)\cap {}^g\tilde{\mathbf{K}}_p$ of $H(\mfk)$, which then implies (\cite[Prop.3]{Greenberg63}) the existence of $h\in H(\mfk)\cap {}^g\tilde{\mathbf{K}}_p$ such that
\begin{equation*} 
%\label{eqn:epsilon'^{-1}F^n_has_fixed_pt}
(\epsilon'^{-1}\Phi^m)^{t}=h\sigma^{tn}(h^{-1})\rtimes \sigma^{tn},
\end{equation*}
where $\Phi:=F^{[\kappa(\wp):\Fp]}$; we fix such $t\in\N$. 

We note that $\epsilon'^{-t}\Phi^{mt}$ fixes $hg x^{\mathrm{o}}=g x^{\mathrm{o}}$ ($x^{\mathrm{o}}:=1\cdot \tilde{\mathbf{K}}_p$, the base point of $G(\mfk)/\tilde{\mathbf{K}}_p$), and that by (\ref{eq:Kottwitz97-(7.3.1)}),
\[\mathrm{inv}_{T'(\mfk)_1}(T'(\mfk)_1,F T'(\mfk)_1)=\underline{\mu'}, \]
where $\underline{\mu'}$ denotes the image of $\mu'$ in $X_{\ast}(T')_{\Gamma_{\mfk}}$ ($\Gamma_{\mfk}:=\Gal(\mfkb/\mfk)$); of course, $X_{\ast}(T')_{\Gamma_{\mfk}}=X_{\ast}(T')$, but we prefer to keep this underline notation.
Let $\tilde{W}:=N_G(T')(\mfk)/T'(\mfk)_1$ be the extended affine Weyl group and $\tilde{W}^{\tilde{\mathbf{K}}_p}:=\tilde{W}\cap\tilde{\mathbf{K}}_p$; there is a canonical isomorphism $\tilde{\mathbf{K}}_p\backslash G(\mfk)/\tilde{\mathbf{K}}_p\simeq \tilde{W}^{\tilde{\mathbf{K}}_p}\backslash\tilde{W}/\tilde{W}^{\tilde{\mathbf{K}}_p}$ \cite[2.4]{Rapoport05}. For $\lambda\in X_{\ast}(T')_{\Gamma_{\mfk}}$, let $t^{\lambda}$ denote the corresponding element of $\tilde{W}$ via $X_{\ast}(T')_{\Gamma_{\mfk}}\cong T'(\mfk)/T'(\mfk)_1\subset \tilde{W}$. 
Then, further, one has $\tilde{W}^{\tilde{\mathbf{K}}_p}t^{\underline{\mu'}} \tilde{W}^{\tilde{\mathbf{K}}_p}\in\mathrm{Adm}_{\tilde{\mathbf{K}}_p}(\{\mu_X\})$, so that as $g^{-1}T'(\mfk)_1g\subset \tilde{\mathbf{K}}_p$,
\[ \mathrm{inv}_{\tilde{\mathbf{K}}_p}(gx^{\mathrm{o}},F gx^{\mathrm{o}})=\tilde{W}^{\tilde{\mathbf{K}}_p}\ t^{\underline{g^{-1}\mu'g}}\ \tilde{W}^{\tilde{\mathbf{K}}_p}\in\mathrm{Adm}_{\tilde{\mathbf{K}}_p}(\{\mu_X\}) \]
(regarding $g^{-1}\mu'g$ as a cocharacter of $g^{-1}T'g$),
which proves that $(\phi,\epsilon^t)$ is $\mathbf{K}_p$-effective admissible.%% 
\footnote{This fact will not be used in the rest of this paper though, especially in our proof of the Lefschetz number formula.}

Next, we show that there exists $e\in H(\mfk)$ such that \[e^{-1}(\epsilon')^{-1}\Phi^m e=\sigma^n,\] 
which will establish the admissibility of $(T,h,\epsilon)$. 
We first claim that there exists $c\in H(\mfk)$ such that \[c^{-1}(\epsilon')^{-1}\Phi^m c\in \rho(H^{\uc}(\mfk))\times\sigma^n,\] 
where $\rho:H^{\uc}\rightarrow H$ is the canonical morphism. 
By Lemma \ref{lem:kernel_of_w} below, it suffices to show that $w_{H}(\epsilon'^{-1}b_n)=0$ (recall that $\Phi^m=F^n=b_n \sigma^n$). By (\ref{eq:xi_p'}), $b_n=\prod_{i=1}^n\sigma^{i-1}(\Nm_{K/K_0}(\mu'(\pi)))$ so 
\[w_H(b_n)=\sum_{i=1}^n\sigma^{i-1}\underline{\mu'}=\sum_{i=1}^n\sigma^{i-1}\underline{\mu_h}=w_H(\epsilon)=w_H(\epsilon')\] 
(the first equality is (\ref{eq:Kottwitz97-(7.3.1)})).

Next, we proceed as in the proof (on p. 193) of \cite{LR87}, Satz 5.21, to find $d\in \rho(H^{\uc}(\mfk))\cap {}^g\tilde{\mathbf{K}}_p$ such that $d^{-1}c^{-1}\epsilon'^{-1}\Phi^mcd=\sigma^n$. For that, when we pick $k'\in H^{\uc}(\mfk)$ mapping to $c^{-1}\epsilon'^{-1}b_n \sigma^{n}(c)$ (i.e. $\rho(k') \sigma^n=c^{-1}(\epsilon'^{-1}\Phi^m)c$), by \cite[Prop.5.4]{Kottwitz85}, it suffices to show that $k'$ is basic (in $B(H^{\uc})$).
By definition \cite[(4.3.3)]{Kottwitz85}, this is the same as the existence of $d'\in H^{\uc}(\mfk)$ with $(k' \sigma^n)^{t}=d'(1\rtimes\sigma^{nt})d'^{-1}$ for some sufficiently large $t$.
But, we have
\[c^{-1}(\epsilon'^{-1}\Phi^m)^{t}c=c^{-1}k_t\sigma^{nt}(c) \sigma^{nt}=(c^{-1}k_tc)(c^{-1}\sigma^{nt}(c)) \sigma^{nt},\]
and for sufficiently large $t\in\N$, both $c^{-1}k_tc$ and $c^{-1}\sigma^{tn}(c)$ are contained in any neighborhood of $1$ in $H(\mfk)$, in particular, in the hyperspeical subgroup $H(\mfk)\cap {}^g\tilde{\mathbf{K}}_p$ of $H(\mfk)$. Thus, 
if $k'_t:=k'\cdot\sigma^n(k')\cdots \sigma^{n(t-1)}(k')\in H^{\uc}(\mfk)$ (i.e. 
$\rho(k_t') \sigma^{nt}=c^{-1}(\epsilon'^{-1}\Phi^m)^{t}c$), $\rho(k_t')$ lies in the hyperspeical subgroup $H^{\der}(\mfk)\cap {}^g\tilde{\mathbf{K}}_p$ of $H^{\der}(\mfk)$. In view of the canonical equality of reduced buildings $\mathcal{B}(H^{\uc},\mfk)=\mathcal{B}(H^{\der},\mfk)$, this implies that $k_t'$ also lies in the stabilizer in $H^{\uc}(\mfk)$ of the corresponding hyperspecial vertex, which is itself a hyperspecial subgroup of $H^{\uc}(\mfk)$ (as $w_{H^{\uc}}$ is trivial). Hence, again by \cite[Prop.3]{Greenberg63} there exists $d'\in H^{\uc}(\mfk)$ such that $k'_t=d'\sigma^{nt}(d'^{-1})$, as required.

(2) Suppose that $c\gamma_0c^{-1}=\Nm_n(\delta)$ for $c\in G(\mfk)$; then, $b:=c^{-1}\delta\sigma(c)\in Z_{G_{\Qp}}(\gamma_0)(\mfk)$; thus, $\gamma_0,b\in H(\mfk)$. As a matter of fact, the arguments coming next work in general for \emph{any} semi-standard $\Qp$-Levi subgroup $H$ of $G_{\Qp}$ containing $\gamma_0$ and $b$. 

We proceed in several steps:

\textbf{Step (i)}: 
First, we prepare the setup.

We choose a maximal $\Qp$-split torus in $H$ and a maximal $\Qp$-split torus $S$ in $G_{\Qp}$ containing it; as $H$ is the centralizer of a $\Qp$-split torus, we have $S\subset H$ and $S$ is also a maximal $\Qp$-split torus in $H$. There exists a $\Qp$-torus $T'$ in $H$ whose extension to $\Qpnr$ becomes a maximal $\Qpnr$-split torus of $H_{\Qpnr}$ containing $S_{\Qpnr}$ \cite[5.1.12]{BT84}; clearly, $T'$ is again such a torus for $G_{\Qp}$. As we assume that $G_{\Qp}$ is unramified, $T'$ is a maximal torus of $H$.

Let $\mathcal{A}^{G}_{\mfk}$ and $\mathcal{A}^{H}_{\mfk}$ be respectively the apartments corresponding to $T'$ in $\mathcal{B}(G_{\Qp},\mfk)$ and $\mathcal{B}(H_{\Qp},\mfk)$. By conjugation, we assume that the given hyperspecial point $\mathbf{0}$ defining $\mathbf{K}_p$ lies in $\mathcal{A}^{G}_{\mfk}$.
We fix a $\sigma$-stable alcove $\mathbf{a}$ in $\mathcal{A}^{G}_{\mfk}$ whose closure contains $\mathbf{0}$, and let $K_{\mathbf{a}}(\mfk)$ be the corresponding Iwahoric subgroup of $G(\mfk)$ (so, $K_{\mathbf{a}}(\mfk)\subset \tilde{\mathbf{K}}_p$). Then, there exists a $H(\mfk)\rtimes\Gal(\Qpnr/\Qp)$-equivariant embedding $\mathcal{B}(H,\mfk)\hra\mathcal{B}(G,\mfk)$; we may assume that $\mathbf{0}$ lies in the image of the apartment $\mathcal{A}^{H}_{\mfk}$, so that $H(\mfk)\cap \mathbf{K}$ is also a hyperspecial subgroup of $H(\mfk)$.
Also, the choice of a base alcove $\mathbf{a}$ presents the extended affine Weyl group $\tilde{W}=N_G(T')(\mfk)/T'(\mfk)_1$ as the semidirect product $W_a\rtimes \Omega_{\mathbf{a}}$ of the affine Weyl group $W_a$ (attached to $T'$) with the normalizer subgroup $\Omega_{\mathbf{a}}\subset \tilde{W}$ of $\mathbf{a}$, thereby fixes a Bruhat order $\leq$ on $\tilde{W}$ and on $\tilde{W}^{\tilde{\mathbf{K}}_p}\backslash \tilde{W} /\tilde{W}^{\tilde{\mathbf{K}}_p}\cong X_{\ast}(T')/W_0$ ($W_0$ being the absolute Weyl group of $(G,T')$).

\textbf{Step (ii)}:  
We write $\mathbf{K}$ for $\tilde{\mathbf{K}}_p$.
We claim existence of $c_1\in G(\mfk)$ such that for $\delta':=c_1b\sigma(c_1^{-1})$, one has
\begin{equation} \label{eq:stable_conjugacy_1}  
c_1\gamma_0c_1^{-1}=\Nm_n(\delta'), 
\end{equation}
and 
\begin{equation} \label{eq:Deligne-Lustizg_1}
\mathbf{K}\cdot \delta' \cdot \mathbf{K} \in \mathrm{Adm}_{\mathbf{K}}(\{\mu_X\}).
\end{equation}
Indeed, we pick $g\mathbf{K}\in Y_p(\delta)$. Then, as $g^{-1}\sigma^n(g)\in \mathbf{K}$, by \cite[Prop.3]{Greenberg63}, there exists $k_0\in \mathbf{K}$ with $g^{-1}\sigma^n(g)=k_0\sigma^n(k_0^{-1})$, i.e. $d:=(gk_0)^{-1}\in G(L_n)$, and so $c_1:=dc$ satisfies (\ref{eq:stable_conjugacy_1}). It also fulfills (\ref{eq:Deligne-Lustizg_1}), since 
\[ \mathbf{K}\cdot \delta'\cdot \mathbf{K}=\mathbf{K}\cdot g^{-1}\delta\sigma(g)\cdot \mathbf{K}\in \mathrm{Adm}_{\mathbf{K}}(\{\mu_X\})\]
Note that $\delta'=d\delta\sigma(d^{-1})$, i.e. $\delta'$ is $\sigma$-conjugate to $\delta$ under $G(L_n)$.

\textbf{Step (iii):} 
Next, we choose a $\Qp$-parabolic subgroup $Q$ of $G_{\Qp}$ of which $H$ is a Levi factor; let $N_Q$ be the unipotent radical of $Q$.
We claim that there exists $m\in H(\mfk)$ such that 
\begin{equation} \label{eq:stable_conjugacy_2}
\Nm_n(m^{-1}b\sigma(m))^{-1}\cdot m^{-1}\gamma_0m \in H(\mfk)\cap (\mathbf{K}\cdot N_Q(\mfk))
\end{equation}
and that
\begin{equation} \label{eq:Deligne-Lustizg_2}
\mathbf{K} \cdot m^{-1}b\sigma(m)n'\cdot \mathbf{K} \in \mathrm{Adm}_{\mathbf{K}}(\{\mu_X\})
\end{equation}
for some $n'\in N_Q(\mfk)$.

Indeed, we have the Iwasawa decomposition $G(\mfk)=Q(\mfk) \mathbf{K}$, which follows from the classical Iwasawa decomposition $G(\mfk)=Q(\mfk) \mathrm{Fix}(\mathbf{0})$ \cite[3.3.2]{Tits79} and that $\mathrm{Fix}(\mathbf{0})\subset T'(\mfk)\cdot \mathbf{K}$. Indeed, $\mathrm{Fix}(\mathbf{0})\subset G(\mfk)=\mathbf{K} T'(\mfk) \mathbf{K}$, so any $g\in \mathrm{Fix}(\mathbf{0})$ is written as $k_1 tk_2$ with $k_1,k_2\in \mathbf{K}$ and $t\in T'(\mfk)\cap \mathrm{Fix}(\mathbf{0})$. But, as $\mathbf{K}$ is normal in $\mathrm{Fix}(\mathbf{0})$, we see that $g\in T'(\mfk)\cdot \mathbf{K}$.
Using this Iwasawa decomposition, we write 
\[c_1^{-1}=nmk\] 
with $n\in N_Q(\mfk)$, $m\in H(\mfk)$, and $k\in \mathbf{K}$, so that
\begin{align*}
\delta' & =c_1b\sigma(c_1)^{-1}=k^{-1}m^{-1}n^{-1}b\sigma(n)\sigma(m)\sigma(k) \\
& =k^{-1}m^{-1}b\sigma(m)n'\sigma(k)
\end{align*}
for $n':=\sigma(m)^{-1}b^{-1}n^{-1}b\sigma(n)\sigma(m)$ (which belongs to $N_Q(\mfk)$, as $H$ normalizes $N_Q$). Hence, we have 
\[ \mathbf{K} \cdot m^{-1}b\sigma(m)n'\cdot \mathbf{K}=\mathbf{K}\cdot \delta'\cdot \mathbf{K},\] 
establishing (\ref{eq:Deligne-Lustizg_2}) by (\ref{eq:Deligne-Lustizg_1}).

On the other hand, the left-hand side of (\ref{eq:stable_conjugacy_1}) becomes
\[ c_1\gamma_0c_1^{-1}=k^{-1}m^{-1}n^{-1}\gamma_0 nmk =k^{-1}m^{-1}\gamma_0m n_1 k \]
for $n_1:=m^{-1}\gamma_0^{-1}n^{-1}\gamma_0nm\in N_Q(\mfk)$, so that if $b_1:=m^{-1}b\sigma(m)\in H(\mfk)$, the right-hand side of (\ref{eq:stable_conjugacy_1}) becomes (again using that $H$ normalizes $N_Q$):
\begin{align*}
\Nm_n(\delta') &=\Nm_n(k^{-1}b_1n'\sigma(k)) \\
&= k^{-1}\cdot b_1n'\cdot \sigma(b_1)\sigma(n')\cdots\sigma^{n-1}(b_1)\sigma^{n-1}(n')\cdot \sigma^n(k) \\
& =k^{-1}\Nm_n(b_1)n_2\sigma^n(k).
\end{align*}
for some $n_2\in N_Q(\mfk)$.
Therefore, (\ref{eq:stable_conjugacy_1}) becomes:
\[m^{-1}\gamma_0m\cdot n_1\cdot k= \Nm_n(b_1)\cdot n_2 \cdot \sigma^n(k),\]
which reduces to: 
\[ \Nm_n(b_1)^{-1}\cdot m^{-1}\gamma_0m \cdot n_3=\sigma^n(k)k^{-1} \]
for some $n_3 \in N(\mfk)$ (i.e. $n_3$ is defined by $hn_3=\sigma^n(k)k^{-1}=n_2^{-1}h\cdot n_1$ for $h:=\Nm_n(b_1)^{-1}\cdot m^{-1}\gamma_0m\in H(\mfk)$), establishing (\ref{eq:stable_conjugacy_2}). 

\textbf{Step (iv):} 
Next, we claim that
\begin{equation} \label{eq:vanishing_of_w_H}
 w_H(H(\mfk)\cap (\mathbf{K}\cdot N_Q(\mfk)))=0.
\end{equation}
For that, recalling that $T'_{\mfk}$ is a maximal $\mfk$-split torus  of $G_{\mfk}$ and also of $H$, we choose a Borel subgroup $B'$ of $G_{\mfk}$ containing $T'_{\mfk}$ and $N_Q$. Then, we pick a Borel subgroup $B$ of $H_{\mfk}$ containing $T'_{\mfk}$ and contained in $B'$; let $N_B$ be its unipotent radical. Thanks to the condition that $H(\mfk)\cap \mathbf{K}$ is also a hyperspecial subgroup of $H(\mfk)$, it is enough to show that 
\[H(\mfk)\cap (\mathbf{K}\cdot N_Q(\mfk)) \subseteq (H(\mfk)\cap \mathbf{K})\cdot N_B(\mfk).\]

Suppose given $m=kn$ for some $m\in H(\mfk)$, $k\in\mathbf{K}$, and $n\in N_Q(\mfk)$.
Since $H(\mfk)\cap \mathbf{K}$ is a hyperspecial subgroup of $H(\mfk)$, we use the Iwasawa decomposition for $(H(\mfk),H(\mfk)\cap \mathbf{K},T',B)$, to write
\[m=k_0 t u\]
where $k_0\in H(\mfk)\cap \mathbf{K}$, $t\in T'(\mfk)$, and $u\in N_B(\mfk)$. So, we have
$k_0^{-1}k=t(un^{-1})$. Then, since one has 
\[ B'(\mfk)\cap\mathbf{K}=(T'(\mfk)\cap\mathbf{K})\cdot  (N_{B'}(\mfk)\cap\mathbf{K}) \]
\cite{BT84}, 5.2.4, applied to $(G(\mfk),\mathbf{K}, T'_{\mfk},B')$), we must have
$t\in \mathbf{K}$, $un^{-1}\in \mathbf{K}$.

\textbf{Step (v):} 
Let $\mu_B$ be a dominant representative of $\{\mu\}$, where we choose the dominant Weyl chamber opposite to the unique Weyl chamber containing the base alcove $\mathbf{a}$ with apex at the hyperspecial vertex $\mathbf{0}$ (following the convention of \cite{HeRapoport17}).
From (\ref{eq:Deligne-Lustizg_2}), there exists $\mu'\in X_{\ast}(T')\cap W_0\cdot\mu_B$ that 
\begin{eqnarray} \label{eq:Iwahori_inequality} 
\mathbf{K}\cdot m^{-1}b\sigma(m)n'\cdot \mathbf{K} &\leq&\tilde{W}^{\mathbf{K}}\cdot t^{\underline{\mu'}}\cdot \tilde{W}^{\mathbf{K}}
\end{eqnarray} 
under the canonical isomorphism $\mathbf{K}\backslash G(\mfk)/\mathbf{K}\simeq \tilde{W}^{\mathbf{K}}\backslash\tilde{W}/\tilde{W}^{\mathbf{K}}$ (\cite[3.4]{Rapoport05}). 

Next, let $\mu''\in X_{\ast}(T')$ be the cocharacter defined by the Cartan decomposition
\begin{equation} \label{eq:mu''} 
m^{-1}b\sigma(m)\in (\mathbf{K}\cap H(\mfk))\ t^{\underline{\mu''}}\ (\mathbf{K}\cap H(\mfk))
\end{equation}
for $(m^{-1}b\sigma(m),H,\mathbf{K}\cap H(\mfk))$.

Then, one has 
\[\mu''=w\mu'\] 
for some $w\in W$; 
in particular, one has $w_H(m^{-1}b\sigma(m))=\underline{\mu''}=\underline{w\mu'}$ in $\pi_1(H)_{\Gal(\overline{\mfk}/\mfk)}(=\pi_1(H))$.
In our case that $\mathbf{K}$ is hyperspecial and $T'$ is unramified, this was proved in \cite[p. 178, line 11-12]{LR87} by using Satake transform \cite[2.3.3, 2.3.7]{Kottwitz84b}.

Therefore, by (\ref{eq:vanishing_of_w_H}) we have
\begin{align*}
w_H(\gamma_0) &=w_H(m^{-1}\gamma_0m)=w_H(\Nm_n(m^{-1}b\sigma(m))) \\
&=\sum_{i=1}^n\sigma^{i-1}w_H(m^{-1}b\sigma(m))=\sum_{i=1}^n\sigma^{i-1}\underline{w\mu'}.
\end{align*} 
As $w\mu'\in X_{\ast}(T')\cap \{\mu_X\}$, the lemma is proved. 

This completes the proof of the theorem.
\end{proof}

%%%%%%%%%%%%%%%%%%%%
\begin{lem} \label{lem:kernel_of_w}
Let $H$ be a connected reductive group over $\Qp$. For any element $h\in H(\mfk)$ with $w_H(h)=0$, there exists $c\in H(\mfk)$ such that $ch\sigma^n(c^{-1})\in \rho(H^{\uc}(\mfk))$, where $\rho:H^{\uc}\rightarrow H$ is the canonical homomorphism 
\end{lem}

\begin{proof}
Using the fact \cite[(3.3.3)]{Kottwitz84b} that for any maximal $\mfk$-split torus $S$ of $H_{\mfk}$ and its centralizer $T$ (which is a maximal torus of $H$ as $H_{\mfk}$ is automatically quasi-split), one has
\[H(\mfk)=\rho(H^{\uc}(\mfk))T(\mfk),\]
we write $h=\rho(h')t$, where $h'\in H^{\uc}(\mfk)$ and $t\in T(\mfk)$.
Since $w_H$ vanishes on $\rho(H^{\uc}(\mfk))$ (cf. diagrams (7.4.1), (7.4.2) of \cite{Kottwitz97}), we have
$w_H(t)=0$. When we put $T^{\uc}:=\rho^{-1}(T)$ (maximal torus of $H^{\uc}$), as $X_{\ast}(T^{\uc})$ is an induced module for $I=\Gal(\bar{\mfk}/\mfk)$ \cite[4.4.16]{BT84} so that
one has
\begin{equation} \label{eq:HR08_p.196}
X_{\ast}(T)_I/X_{\ast}(T^{\uc})_I=\pi_1(H)_I
\end{equation}
\cite[p.196]{HainesRapoport08}, and as $w_{T^{\uc}}$ is surjective, we can find $t'\in T^{\uc}(\mfk)$ with $w_{T}(\rho(t')t^{-1})=0$. Hence, by \cite[Prop.3]{Greenberg63} there exists $c\in T(\mfk)$ with $\rho(t')t^{-1}=c^{-1}\sigma^n(c)$, namely with $c\rho(t')\sigma^n(c^{-1})=t$.
Finally, as $\rho(H^{\uc}(\mfk))$ is a normal subgroup of $H(\mfk)$, this establishes the claim.
\end{proof}

%%%%%%%%%%%%%%%%%%%%
\subsection{Some results in Galois cohomology of algebraic groups}
For the next discussion, 
it is necessary to use Tate hypercohomology groups $\widehat{\mathbb{H}}^i\ (i\in\Z)$ of (bounded) complexes of discrete $\mathcal{G}$-modules for a finite group $\mathcal{G}$: they are defined by means of either the complete (standard) resolution of the trivial $\Z[\mathcal{G}]$-module $\Z$ or hypercochains in the usual manner (cf. \cite[$\S$1]{Koya90}).%% 
\footnote{For finite $\mathcal{G}$, the Tate hypercohomology $\widehat{\mathbb{H}}^i(\mathcal{G},-)$ factors through the stable module category $\mathcal{T}(\mathcal{G})$ of the group algebra $\Z[\mathcal{\mathcal{G}}]$, and as such equals the cohomological functor $A^{\bullet}\mapsto \Hom_{\mathcal{T}(\mathcal{G})}(\Z,A^{\bullet}[i])$ on that triangulated category.}
In this work, we will be only interested in the bounded complexes $A^{\bullet}$ of discrete $\mathcal{G}$-modules \emph{whose positive terms are zero}, in which case one has 
\[ \widehat{\mathbb{H}}^i(\mathcal{G},A^{\bullet})=\begin{cases} 
\mathbb{H}^0(\mathcal{G},A^{\bullet})/\Nm_{\mathcal{G}}\mathcal{H}^0(A^{\bullet}) &\text{ if } i=0\\
\mathbb{H}^i(\mathcal{G},A^{\bullet}) &\text{ if } i>0, \end{cases} \]
where $\mathcal{H}^0(A^{\bullet})$ denotes the $0$-th cohomology $\mathcal{G}$-module of the complex $A^{\bullet}$ and $\Nm_\mathcal{G}$ is the norm map \cite[Prop.1.2]{Koya90}. 

For a diagonalizable $\C$-group $D$ with (algebraic) action of a finite group $\mathcal{G}$,
the  canonical surjection and injection 
\begin{equation} \label{eq:H^0_H^-1}
\pi_0(D^{\mathcal{G}})\twoheadrightarrow \widehat{\mathbb{H}}^{0}(\mathcal{G},D)=\widehat{\mathbb{H}}^{1}(\mathcal{G},X_{\ast}(D)),\quad \widehat{\mathbb{H}}^{-1}(\mathcal{G},X^{\ast}(D))\hookrightarrow X^{\ast}(D)_{\mathcal{G},\tors},
\end{equation}
identify the canonical duality $X^{\ast}(D)_{\mathcal{G},\tors} \isom \pi_0(D^{\mathcal{G}})^D$ with the duality induced from the cup-product pairing
\[\widehat{\mathbb{H}}^{1}(\mathcal{G},X_{\ast}(D))\otimes\widehat{\mathbb{H}}^{-1}(\mathcal{G},X^{\ast}(D))\rightarrow \widehat{\mathbb{H}}^{0}(\mathcal{G},\Z)=\frac{1}{|\mathcal{G}|}\Z.\]

For a field $k$, either global or local, we let 
\begin{align} \label{eq:C_k}
C_k:=\begin{cases} k^{\times} & \text{ if }k\text{ is local},\\
\A_k^{\times}/k^{\times} & \text{ if }k\text{ is global}.\end{cases}
\end{align}
For a number field $F$ and a bounded complex of $F$-tori $T^{\bullet}=(\cdots\rightarrow T^i \rightarrow \cdots)$ (concentrated in non-positive degrees) and $i\in\Z$, we define
\begin{align*}
\mathbb{H}^i(\A_F/F,T^{\bullet})&:=\mathbb{H}^i(F,T^{\bullet}(C_{\overline{F}})) \\
& =\varinjlim_{E}\mathbb{H}^i(E/F,T^{\bullet}(C_E)),
\end{align*}
where $E$ runs through finite Galois extensions of $F$, 
$T^{\bullet}(C_{\overline{F}})$ denotes the complex of discrete $\Gamma_F$-modules $\cdots\rightarrow T^i(\A_{\overline{F}})/T^i(\overline{F})\rightarrow\cdots$ and $T^{\bullet}(C_E)$ is defined similarly.

%%%%%%%%%%%%%%%%%%%%
\begin{lem} \label{lem:identification_of_Kottwitz_A(H)}
(1) Let $H$ be a (connected) reductive group over a field, either global or local, and $T$ a maximal $k$-torus of $H$; set $\tilde{T}:=\rho_H^{-1}(T)$. If $k'$ is a finite Galois extension of $k$ splitting $T$, the two pairings are compatible:
\[ \xymatrix{ \pi_0(Z(\hat{H})^{\mathcal{G}}) \ar@{->>}[d] & \bigotimes & \pi_1(H)_{\mathcal{G},\tors}  \ar[r] & \Q/\Z \ar@{=}[d] \\
\widehat{\mathbb{H}}^{0}(k'/k,X^{\ast}(H_{\mathbf{ab}})) & \bigotimes & \widehat{\mathbb{H}}^{1}(k'/k,H_{\mathbf{ab}}(C_{k'})) \ar[r] \ar@{^(->}[u] & \Q/\Z }, \]
where the bottom pairing is the local/global Tate-Nakayama pairings for the complex $H_{\mathbf{ab}}=(\tilde{T}\rightarrow T)$ (cf. \cite[A.2.2]{KottwitzShelstad99}).%% 
\footnote{By our convention, $H_{\mathbf{ab}}=\tilde{T}\rightarrow T$ and $X^{\ast}(H_{\mathbf{ab}})=X^{\ast}(T)\rightarrow X^{\ast}(\tilde{T})$ are both placed in degree $[-1,0]$.}%% 

Moreover, the two vertical maps are bijections if the extension $k'/k$ has sufficiently large degree.
\end{lem}

\begin{proof} 
(1) First, we note that the bottom local/global Tate-Nakayama pairings are perfect pairings: its proof can be easily reduced to the classical Tate-Nakayama duality for tori (cf. proof of (A.2.2) of \cite{KottwitzShelstad99}). 
In view of the canonical maps (\ref{eq:H^0_H^-1}), the existence of the vertical maps and the compatibility of the pairings deduced from the canonical isomorphisms
\begin{align*}
\widehat{\mathbb{H}}^{0}(k'/k,Z(\hat{H}))^D &=\widehat{\mathbb{H}}^{-1}(k'/k,\hat{H}_{\mathbf{ab}})^D \\
&\stackrel{(a)}{=}\widehat{\mathbb{H}}^{0}(k'/k,X^{\ast}(T)\rightarrow X^{\ast}(\tilde{T}))^D \\
&\stackrel{(b)}{=}\widehat{\mathbb{H}}^{1}(k'/k,H_{\mathbf{ab}}(C_{k'})) \\
&\stackrel{(c)}{=}\widehat{\mathbb{H}}^{-1}(k'/k,\pi_1(H)), 
\end{align*}
where the equalities are as follows:

(a) For any two-term complex $D_{-1}\rightarrow D_0$ of diagonalizable $\C$-groups with $\Gal(k'/k)$-action, the exponential sequence $0\rightarrow X_{\ast}(D_i)\rightarrow \mathrm{Lie}(D_i) \rightarrow D_i\rightarrow 0$ gives an isomorphism, for all $i\in\Z$,
\begin{equation} \label{eq:connecting_isom_diagonal_gp}
\widehat{\mathbb{H}}^i(k'/k,D_{-1}\rightarrow D_0) = \widehat{\mathbb{H}}^{i+1}(k'/k,X_{\ast}(D_{-1})\rightarrow X_{\ast}(D_0)),
\end{equation} (First, establish this for a single term complex using the cohomology long exact sequence of the associated exponential sequence. Then, the general case follows from it by an appropriate application of the five lemma). 

(b) For a single term complex, this is the classical Tate-Nakayma duality. Then, one uses the five lemma in an obvious way.

(c) Apply the five lemma to the isomorphisms from the long cohomology exact sequence attached to the short exact sequence $0\rightarrow X_{\ast}(\tilde{T})\rightarrow X_{\ast}(T)\rightarrow \pi_1(H)\rightarrow 0$ to the long exact sequence 
 \[ \cdots \rightarrow \widehat{\mathbb{H}}^{i}(k'/k,\tilde{T}(C_{k'})) \rightarrow  \widehat{\mathbb{H}}^{i}(k'/k,T(C_{k'})) \rightarrow \widehat{\mathbb{H}}^{i}(k'/k,H_{\mathbf{ab}}(C_{k'})) \rightarrow  \widehat{\mathbb{H}}^{i+1}(k'/k,\tilde{T}(C_{k'})) \rightarrow \cdots, \]
provided by cup-product with the fundamental class in $H^2(k'/k,C_{k'})$, to see that the same cup-product gives an isomorphism
\[ \widehat{\mathbb{H}}^{i-1}(k'/k,\pi_1(H)) \isom \widehat{\mathbb{H}}^{i+1}(k'/k,H_{\mathbf{ab}}(C_{k'})). \]

Finally, the two maps (\ref{eq:H^0_H^-1}) become isomorphisms  if the degree $[k':k]$ is sufficiently large \cite[Prop.B.4]{Milne92}.
\end{proof}

Following Kisin \cite[(4.4)]{Kisin17}, we introduce some Galois cohomology notations. For a crossed module of algebraic $\Q$-groups $H'\rightarrow H$, we set
\begin{equation} \label{eq:Sha^{infty}}
\Sha^{\infty}(\Q,H'\rightarrow H)=\ker\left[H^1(\Q,H'\rightarrow H)\rightarrow H^1(\R,H'\rightarrow H)\right].
\end{equation}
When $H$ is a \emph{connected} reductive $\Q$-group, the natural map 
\[ \Sha^{\infty}(\Q,H):=\Sha^{\infty}(\Q,\{e\}\rightarrow H)\ \rightarrow\ \Sha^{\infty}(\Q,H^{\uc}\rightarrow H) \] 
is a bijection, thus $\Sha^{\infty}(\Q,H)$ is an abelian group in a natural way. For a general \emph{$\Q$-subgroup} $H$ of $G$, we introduce a pointed set
\[ \Sha^{\infty}_G(\Q,H):=\ker\left[ \Sha^{\infty}(\Q,H) \rightarrow \Sha^{\infty}(\Q,G)\right].\]
If $H$ is connected, by the above fact this is an abelian group.
We also define $\Sha^{\infty}_G(\Q,H_1)$ for an inner twist $H_1$ of a \emph{connected} reductive subgroup $H$ of $G$ as follows.
Let $I_0$ be a \emph{connected} reductive $\Q$-subgroup of $G$ and $I_1$ an inner twist of $I_0$.
Then, we define
\begin{equation} \label{eq:Sha^{infty}_G(I_1)}
\Sha^{\infty}_G(\Q,I_1):=\ker\left[ \Sha^{\infty}(\Q,I_1)\isom \Sha^{\infty}(\Q,I_0) \rightarrow \Sha^{\infty}(\Q,G)\right],
\end{equation}
where the isomorphism $\Sha^{\infty}(\Q,I_1)\isom \Sha^{\infty}(\Q,I_0)$ arises from the canonical isomorphisms (see the proof of Lemma \ref{lem:abelianization_exact_seq}, cf. \cite[Lem.4.4.3]{Kisin17}):
\begin{align} \label{eq:Kisin17_Lem.4.4.3}
\Sha^{\infty}(\Q,I_1)\isom \Sha^{\infty}(\Q,I_{1\mathbf{ab}}) \isom \Sha^{\infty}(\Q,I_{0\mathbf{ab}})\isom \Sha^{\infty}(\Q,I_0)
\end{align}
In general, the definition of $\Sha^{\infty}_G(\Q,I_1)$ will depend on the auxiliary group $I_0$ as well (but, depends only on its $G(\Qb)$-conjugacy class). But, in the case of our interest (\ref{eq:Sha^{infty}_G}), it will depend only on $I_1$ in a natural manner.

%%%%%%%%%%%%%%%%%%%%
%%%%%%%%%%%%%%%%%%%%

\section{Proof of Kottwitz formula: Cardinality of fixed point set of Frobenius-twisted Hecke correspondence} 
\label{sec:Pf_Kottwitz_formula}

In this section, we prove the main theorem, Theorem \ref{thm_intro:Kottwitz_formula}, i.e. the Kottwitz formula \cite[(3.1)]{Kottwitz90}, \cite[(19.6)]{Kottwitz92} on the expression of the Lefschetz number of a Frobenius-twisted Hecke correspondence acting on an automorphic lisse sheaf on a Shimura variety of Hodge-type with hyperspecial level as a certain sum of products of orbital and Twisted orbital integrals. 

Since in this section we will work exclusively with the mod-$\wp$ reductions of $\sS_{\mathbf{K}}$ and $\sS:=\varprojlim_{H^p}\sS_{H^p}$, we denote such $\kappa(\wp)$-schemes again by $\sS_{\mathbf{K}}$ and $\sS$.
Any element $g\in G(\A_f^p)$ gives rise to a Hecke correspondence (from $\sS_{\mathbf{K}}$ to itself): 
\begin{equation} \label{eq:Hecke_corr_f}
\sS_{\mathbf{K}}\stackrel{p_1'}{\longleftarrow} \sS_{\mathbf{K}_g} \stackrel{p_2}{\longrightarrow} \sS_{\mathbf{K}}.
\end{equation}
where $\mathbf{K}_g:=\mathbf{K}\cap g\mathbf{K} g^{-1}$, the right-hand map $p_2$ is the natural projection induced by the inclusion $\mathbf{K}_g\subset \mathbf{K}$ and the left-hand map $p_1'$ is the composite of the natural projection $\sS_{\mathbf{K}_g}\stackrel{p_1}{\rightarrow} \sS_{g\mathbf{K} g^{-1}}$ (induced by the inclusion $\mathbf{K}_g\subset g\mathbf{K} g^{-1}$) with the (right) action by $g$: $\sS_{g\mathbf{K} g^{-1}} \isom \sS_{\mathbf{K}}$; we denote by $f$ this Hecke correspondence. 
We are interested in the fixed point set of the composition $\Phi^m\circ f$ of the morphism $\Phi^m$ and the correspondence $f$, namely the fixed point set of the correspondence: 
\begin{equation} \label{eq:Hecke_corr_twisted_by_Frob}
\sS_{\mathbf{K}}\stackrel{p_1'}{\longleftarrow} \sS_{\mathbf{K}_g} \stackrel{p_2'=\Phi^m\circ p_2}{\longrightarrow} \sS_{\mathbf{K}}.
\end{equation}
By definition, a fixed point of this correspondence $\Phi^m\circ f$ is a point in $\sS_{\mathbf{K}_g}(\Fpb)$ whose images in $\sS_{\mathbf{K}}(\Fpb)$ under $p_1'$ and $p_2'=\Phi^m\circ p_2$ coincide. 

The main part of the Kottwitz conjecture is concerned with obtaining a group-theoretic description of the cardinality of the fixed point set of such a Frobenius-twisted Hecke correspondence on each isogeny class.
In \cite[$\S$1]{Kottwitz84b}, Kottwitz discusses how one can arrive at such description from a description of the set of mod-$p$ points as provided by the Langlands-Rapoport conjecture (but at that time his deduction was based on an incomplete early version of it suggested by Langlands \cite{Langlands76}, \cite{Langlands79}). This argument by Kottwitz was one of the motivations for the work \cite{LR87} and will also serve as a guide for our proof here.

Our proof relies heavily on geometric results obtained by Kisin \cite{Kisin17}. But we do not use his version of Langlands-Rapoport conjecture \cite[Thm.(0.3)]{Kisin17}, for that matter, in fact no Galois-gerb theory, except in the proof of Theorem \ref{thm:stable_isogeny_diagram} through Theorem \ref{thm:LR_special_admissible_morphism} which summarizes the necessary facts from \cite{LR87} (Theorem \ref{thm:stable_isogeny_diagram} is needed to prove Corollary \ref{cor:stable_conjugacy_from_nice_tori}, i.e. the ingredient (S) below). As mentioned in the introduction, Langlands and Rapoport \cite{LR87} provided group-theoretic arguments (formulated in terms of Galois-gerb theory) deriving Kottwitz conjecture from their conjecture, but Kisin's version of Langlands-Rapoport conjecture is too weak to invoke those arguments for, and we just emulate the (Galois-gerb theoretic) arguments in \cite{LR87}. We proceed in parallel to \cite{LR87}, establishing analogous, but formulated in more familiar terms, group-theoretic statements. Several ingredients are needed to substitute the Langlands-Rapoport's description of $\sS_{\mathbf{K}}(\Fpb)$ and to imitate their deduction arguments. 

The geometric ingredients (from Kisin's work \cite{Kisin17}) are as follows:
\begin{itemize}
\item[(I)] Definition/description of isogeny classes (in terms of affine Deligne-Lusztig varieties) (ibid., Prop.2.1.3) and their moduli interpretation (ibid., Prop. 1.4.15), and the resulting description of $\sS_{\mathbf{K}_p}(G,X)(\Fpb)$ as disjoint union of isogeny classes;
\item[(SCM)] Strong CM-lifting theorem (ibid., Cor. 2.2.5);
\item[(Ta)] Generalization of the Tate's theorem on endomorphisms of abelian varieties over finite fields (ibid., Cor. 2.3.2);
\item[(Tw)] Twisting method of isogeny classes or CM points (ibid., Prop. 4.4.8, 4.4.13);
\end{itemize}

Notable group-theoretic ingredients are as follows:
\begin{itemize}
\item[(TO)] Non-vanishing of $TO_{\delta}(\phi_p)$ implies Langlands-Rapoport condition \ref{itm:ast(epsilon)} (Theorem \ref{thm:LR-Satz5.21} (2));
\item[(S)] Stable Kottwitz triples attached to admissible K-pairs (Corollary \ref{cor:stable_conjugacy_from_nice_tori}).
\end{itemize}
The fact (TO) is a crucial ingredient in the proof of the effectivity criterion (\ref{eq:E}) explained in the introduction. The ingredient (S), which appears for the first time in this work, is invoked only when the derived group of $G$ is not necessarily simply connected, and is the only argument in this work which uses (a bit of basic) Galois gerb theory.

%%%%%%%%%%%%%%%%%%%%
\subsection{Kisin's results \cite{Kisin17}: geometric description of the $\F$-point set of Hodge type Shimura varieties}

Here we give a very brief review of some (geometric) results in \cite{Kisin17} insofar as relevant and necessary for our proof of the Kottwitz formula on Lefschetz numbers. We will use its notation freely.
Recall that we have fixed an embedding $\sigma_p:\Qpb\hookrightarrow \C$ of $E(G,X)$-algebras.

(A) There exists a $\Z_{(p)}$-lattice $V_{\Z_{(p)}}$ of $V$ and a set of tensors $\{s_{\alpha}\}_{\alpha}$ on it which defines the reductive closed $\Zp$-subgroup scheme $G_{\Zp}$ of $G_{\Qp}$ giving the hyperspecial subgroup $\mathbf{K}_p$ (i.e. $\mathbf{K}_p=G_{\Zp}(\Zp)$) \cite[Prop.1.3.2]{Kisin10}. 
Let $\pi:\mathcal{A}\rightarrow \sS_K$ be the universal abelian scheme over $\sS_K$ (for sufficiently small $\mathbf{K}^p$) and $\mathcal{V}=R^1\pi_{\ast}\Omega^{\bullet}$ the first relative de Rham chomology (algebraic vector bundle).
The tensors $\{s_{\alpha}\}_{\alpha}$ give rise to horizontal sections $\{s_{\alpha,\Betti}\}_{\alpha}$ on the local system $R^1\pi^{\mathrm{an}}_{\ast}(\Z_{(p)})$ over $Sh_{\mathbf{K}}(G,X)$ (with $\pi^{\mathrm{an}}$ denoting the analytification of $\pi$) and sections $\{s_{\alpha,\dR}\}_{\alpha}$ on $\mathcal{V}=R^1\pi_{\ast}\Omega^{\bullet}$ which correspond to each other (for the same $\alpha$) via the de Rham isomorphism over $\C$ \cite[1.3.6]{Kisin17}.

(B) Suppose $x\in \sS_{\mathbf{K}_p}(k)$ for a finite extension $k\subset \Fpb$ of $\kappa(\wp)$. Let $\mathcal{A}_x$ be the underlying abelian variety over $k$ and $\bar{x}$ the $\Fpb$-point induced by $x$.
Let $H^1_{\et}(\mathcal{A}_{\bar{x}}/\Ql)$ and $H_{\cris}^1(\mathcal{A}_x/K_0)$ be respectively the $l$-adic \'etale and the cristalline cohomology groups of $\mathcal{A}_{\bar{x}}$ (for $l\neq p$) and $\mathcal{A}_x$, where $K_0:=W(k)[1/p]\ (\subset\Qpb)$. 
We let $H^{\et}_1(\mathcal{A}_{\bar{x}},\Ql):=\Hom_{\Ql}(H^1_{\et}(\mathcal{A}_{\bar{x}}/\Ql),\Ql)$ and $H^{\cris}_1(\mathcal{A}_x/K_0):=\Hom_{K_0}(H^1_{\cris}(\mathcal{A}_{x}/K_0),K_0)$ denote their linear dual (homology) groups. The latter group is equipped with the Frobenius operator $\phi$: $\phi(f)(v):=p^{-1}\cdot{}^{\sigma}f(Vv)$ for $f\in H_1^{\cris}$ and $v\in H^1_{\cris}$, where $V$ is the Verschiebung on $H^1_{\cris}$. 
Recall that the relative Frobenius morphism  $\Fr_{\mathcal{A}_{x}/k}$ of $\mathcal{A}_{x}/k$ acts on these homology spaces ($\Ql$ and $K_0$-linearly) by the geometric $p^{[k:\Fp]}$-Frobenius $\Fr_k^{-1}$ in $\Gal(\Fpb/k)$ and by the inverse of $\phi^{[k:\Fp]}$, respectively.

Then, there exist tensors $\{s_{\alpha,l,x}\}_{\alpha}$ on $H^1_{\et}(\mathcal{A}_{\bar{x}}/\Ql)\ (l\neq p)$ and tensors $\{s_{\alpha,0,x}\}_{\alpha}$ on $H_{\cris}^1(\mathcal{A}_x/K_0)$ which are Frobenius invariant (for the geometric Frobenius acting on $H^1_{\et}(\mathcal{A}_{\bar{x}}/\Ql)$ and the absolute Frobenius automorphism acting on $H^1_{\cris}(\mathcal{A}_x/K_0)$); for the construction of $s_{\alpha,l,x}$ and $s_{\alpha,0,x}$, see \cite[(2.2)]{Kisin10} and \cite[Prop.1.3.9, 1.3.10]{Kisin17} respectively.
Also, there exist isomorphisms matching the respective tensors for each $\alpha$:
\begin{equation} \label{eq:isom_eta}
\begin{split}
\eta_l:(V,\{s_{\alpha}\})\otimes\Ql &\lisom (H^{\et}_1(\mathcal{A}_{\bar{x}},\Ql),\{s_{\alpha,l,x}\}) , \\
\eta_p:(V,\{s_{\alpha}\})\otimes K_0 &\lisom (H^{\cris}_1(\mathcal{A}_x/K_0),\{s_{\alpha,0,x}\}).
\end{split}
\end{equation}
Most of the time, we are just contented with an isomorphism defined over $\Qpnr$ (or even over $\mfk$):
\begin{equation} \label{eq:isom_eta_nr}
\eta_p^{\nr}:(V,\{s_{\alpha}\})\otimes \Qpnr \lisom (H^{\cris}_1(\mathcal{A}_x/\Qpnr),\{s_{\alpha,0,x}\}).
\end{equation}
For almost all $l\neq p$ (in particular, such that $G_{\Ql}$ is unramified), we may assume that the following conditions hold:
there exist a $\Z_{(l)}$-lattice $V_{\Z_{(l)}}$ of $V$ such that the tensors $\{s_{\alpha}\}_{\alpha}$ live on it and defines a reductive $\Z_l$-subgroup scheme $G_{\Z_l}$ of $G_{\Ql}$, and a similar statement holds true of the lattice $H^{\et}_1(\mathcal{A}_{\bar{x}},\Z_l)$ and the tensors $\{s_{\alpha,l,x}\}$. Further, for these $\Z_l$-structures, there exists an $\Z_l$-isomorphism extending $\eta_l$; we denote it again by $\eta_l$.

(C) We define $\gamma_l\in G(\Ql)$ ($l\neq p$), $\gamma_p\in G(K_0)$ by
\begin{equation}  \label{eq:gamma_v}
\gamma_v^{-1}:=\Int(\eta_v^{-1})(\Fr_{\mathcal{A}_{x}/k}).
\end{equation}
These elements $\gamma_l$, $\gamma_p$ are well-defined up to conjugacy in $G(\Ql)$ and $G(K_0)$, respectively.
We also define $\delta\in \mathrm{GL}(V_{K_0})$ by 
\[\delta(1_V\otimes\sigma):=\Int(\eta_p^{-1})(\phi),\]
so one has $\gamma_p=\Nm_n\delta:=\delta\sigma(\delta)\cdots\sigma^{n-1}(\delta)$ ($n=[K_0:\Qp]$).

Let $k'\subset\Fpb$ be a finite extension of $k$. For each finite place $l\neq p$ of $\Q$, let $I_{l/k'}$ be the centralizer of $\gamma_l^{[k':k]}\in G(\Ql)$ and define $I_{p/k'}$ to be the $\sigma$-centralizer $G_{\delta\sigma}$ of $\delta$ in $G(K_0')$, where $K_0':=W[k'][\frac{1}{p}]$ \cite[p.802]{Kottwitz82}: $G_{\delta\sigma}$ is a closed subgroup of $\Res_{K_0'/\Qp}(G)$ such that $G_{\delta\sigma}(\Qp)=\{ y\in G(K_0')\ |\ y(\delta\sigma)=(\delta\sigma)y\}$. One has $(G_{\delta\sigma})_{K_0'}\simeq Z_{G_{K_0'}}(\Nm_n\delta)$.

The increasing sequence of subgroups $\{I_{l/k'}\}_{k'\subset\Fpb}$ of $G_{\Ql}$ stabilizes to a subgroup $I_l$, which also equals the centralizer of $\gamma_l^n$ in $G_{\Ql}$ for (any) sufficiently large $n\in\N$. By similar reasoning (cf. \cite[(2.1.2)]{Kisin17}), there exists a $\Qp$-subgroup $I_p$ of $J_{\delta}$ which equals $I_{p/k'}$ for all sufficiently large $k'\subset\Fpb$.
Write $\Aut_{\Q}(\mathcal{A}_{\bar{x}})$ for the automorphism group of $\mathcal{A}_{\bar{x}}$ in the isogeny category, and let $I_{x}\subset \Aut_{\Q}(\mathcal{A}_{\bar{x}})$ denote the subgroup consisting of elements fixing all the tensors $\{s_{\alpha,l,x}\}_{\alpha}\ (l\neq p)$, $\{s_{\alpha,0,x}\}_{\alpha}$, regarded as an algebraic $\Q$-subgroup of $\Aut_{\Q}(\mathcal{A}_{\bar{x}})$.

%%%%%%%%%%%%%%%%%%%%
%%%%%%%%%%%%%%%%%%%%
\begin{thm} \label{thm:Kisin17_Cor.2.3.2;Tate_isom} \cite[Cor.2.1.7, Cor.2.3.2]{Kisin17}
For every finite place $v$ of $\Q$, the restriction of $\Int(\eta_v)^{-1}$ if $v\neq p$ or of $\Int(\eta_p^{\nr})^{-1}$ if $v=p$ to $(I_{x})_{\Qv}$ induces an isomorphism 
\[\Int(\eta_v)^{-1}(v\neq p),\ \Int(\eta_p^{\nr})^{-1}\ :\ (I_{x})_{\Qv}\isom I_v.\] 
In particular, for all $x\in \mathscr{I}$, the $\Q$-group $I_{x}$  has the same $\Qb$-rank as $G$.
\end{thm}

For each $x\in \sS_{\mathbf{K}_p}(\Fpb)$, if $x$ is defined over a finite field $k\subset \Fpb$, 
any choice of isomorphisms $\eta_l$ (\ref{eq:isom_eta}), $\eta_p^{\nr}$ (\ref{eq:isom_eta_nr}) gives an element $b\in G(\Qpnr)$ by $\Int(\eta_p^{\nr})(b\sigma)=\phi$, and
embeddings defined over $\Ql$ and $\Qp$
\[ i_{x,l}:(I_{x})_{\Ql} \hookrightarrow G_{\Ql},\quad i_{x,p}:(I_{x})_{\Qp} \hookrightarrow J_b.\]
A different choice of $\eta^p:=\prod'\eta_l$, $\eta_p^{\nr}$ gives rise to another triple $(b',{i_{x}^p}':=\prod'i_{x,l},i_{x,p}')$ such that $(b'\sigma,i_{x,p}')=g_p(b\sigma,i_{x,p})g_p^{-1}$ for some $g_p\in G(\mfk)$ and ${i_{x}^p}'=\Int(g^p)\circ i_{x}^p$ for some $g^p\in G(\A_f^p)$; we regard two such triples as being \emph{equivalent}. For almost all $l\neq p$, so that, among others, one can find a $\Z_l$-isomorphism $\eta_l$ for the $\Z_l$-structures explained above, $i_{x,l}$ extends to an embedding $(I_x)_{\Z_l}\hookrightarrow G_{\Z_l}$ of $\Z_l$-group schemes.

Then, there exists a $\langle\Phi\rangle\times Z_G(\Qp)\rtimes G(\A_f^p)$-equivariant injective map 
\[\iota_x:X(\{\mu_X^{-1}\},b)_{\mathbf{K}_p}\times G(\A_f^p)\rightarrow \sS_{\mathbf{K}_p}(\Fpb),\] 
where $\langle\Phi\rangle$ is the cyclic group generated by $\Phi$ \cite[Cor. 1.4.13, Prop. 2.1.3]{Kisin17}; the image of $\iota_x$ is called the \emph{isogeny class} containing $x$. 
According to \cite[Prop.1.4.15]{Kisin17}, two isogeny classes coincide if they have non-empty intersection, so we obtain a decomposition 
\begin{equation} \label{eq:isogeny_decomp:geometric}
 \bigsqcup_{x} I_{x}(\Q)\backslash [X(\{\mu_X^{-1}\},b)_{\mathbf{K}_p}\times G(\A_f^p)] \isom \sS_{\mathbf{K}_p}(\Fpb),
\end{equation}
where $x$ runs through a set of representatives of (i.e. points lying in) the isogeny classes and the quotient by $I_x(\Q)$ is done by $i_{x,p}$ and $i_x^p$. 

It is shown in \cite[(4.4.6)]{Kisin17} that to each isogeny class $\mathscr{I}$, one can attach a Kottwitz triple $(\gamma_0;\gamma,\delta)$: in Kisin's definition \cite[(4.3.1),(4.3.2)]{Kisin17}, this is well-defined up to ``multiplicative'' equivalence relation as well as up to the usual equivalence relation (among the triples of the same level) (also, $\gamma_0$ is by construction $i_T(\Fr_{\mathcal{A}_{x}/k})$ for the relative Frobenius endomorphism $\Fr_{\mathcal{A}_{x}/k}$ of $\mathcal{A}_x$ over a sufficiently big finite field $\F_q\subset\F$ and any embedding $i_T$ of a maximal $\Q$-torus of $I_x$ to be discussed in the next subsection).%%
\footnote{We may regard such Kottwitz triple as a Kottwitz triple attached to $(\mathscr{I},\Fr_{\mathcal{A}_{x}/k})$ as discussed in the introduction.}
We may choose a representative $(\gamma_0;\gamma,\delta)$ with sufficiently large level so that the centralizer $G_{\gamma_0}$ of $\gamma_0$ is connected (which is then independent of the choice of that representative in its multiplicative equivalence class). This $\Q$-group $I_0$ is an inner form of $I_x$ which satisfies the condition (iv$'$) of (\ref{subsubsec:pre-Kottwitz_triple}) with $I=I_x$ \cite[(2.3.5)]{Kisin17} and thus whose class in $H^1(\Q,I_x^{\ad})$ is uniquely determined by the isogeny class only (or its Kottwitz triple), and for any two such representatives with the same level, their rational components $\gamma_0$, $\gamma_0'$ are stably conjugate, from which it follows that their centralizers are inner forms of each other with an inner twist being given by conjugation by an element of $G(\Qb)$, hence (by \cite[4.4.3]{Kisin17}) that the abelian group 
\begin{equation} \label{eq:Sha^{infty}_G}
 \Sha^{\infty}_G(\Q,I_{x}):=\ker\left[\Sha^{\infty}(\Q,I_{x}) \isom \Sha^{\infty}(\Q,I_0) \rightarrow \Sha^{\infty}(\Q,G)\right],
 \end{equation}
where the first isomorphism is induced by any inner twisting $(I_{x})_{\Qb}\isom (I_0)_{\Qb}$ satisfying the condition (iv$'$) of (\ref{subsubsec:pre-Kottwitz_triple}) with $I=I_{x}$ and the second map by the inclusion $I_0\subset G$, is well defined, independent of the choice of a representative of the Kottwiz triple.

%%%%%%%%%%%%%%%%%%%%
%%%%%%%%%%%%%%%%%%%%
\subsection{Group-theoretic description of the $\Fpb$-point set of Hodge type Shimura varieties} \label{subsec:uncond_proof_K-formula}

Here, we emulate the Galois-gerb theoretic deduction arguments of Langlands and Rapoport \cite{LR87} with Kisin's geometric results (which are formulated in terms of abelian varieties with extra structure (of weak polarization and Hodge cycles) and their \'etale and crystalline cohomological realizations). To this end, we 
present purely group-theoretic, making no use of Galois gerb theory, statements which will serve as inputs for our arguments. The geometric ingredients from \cite{Kisin17} will appear only in the proofs of those statements (and nowhere else).
These statements are the four theorems in this subsection: Theorem \ref{thm:Kisin17_Cor.1.4.13,Prop.2.1.3,Cor.2.2.5} Theorem \ref{thm:Kisin17_Cor.2.3.2}, Theorem \ref{thm:Kisin17_Prop.4.4.8}, Theorem \ref{thm:equiv_K-triples}.
Then, from these theorems it requires just some standard arguments in algebraic group theory (more specifically, Galois cohomology theory of algebraic groups over local/global fields) to deduce the key facts (\ref{eq:Tate_thm_SV}), (\ref{eq:E}), and (\ref{eq:C}) discussed in the introduction. Such ``formal'' deductions are given in three Corollaries to the preceding theorems: Corollary \ref{cor:Tate_thm2}, Corollary \ref{cor:K-triple_of_twisted_K-pair}, Corollary \ref{cor:LR-Satz5.25}.
As was explained before, this will complete the proof of the formula in question in the constant coefficient case (this will be done in the last subsubsection).
For general automorphic sheaves, we apply an idea of Kottwitz in \cite{Kottwitz92} which underlies
his notion of \textit{virtual abelian variety over a finite field}, but here our notion is more direct and purely group theoretic, not involving any geometric object.

We remind the reader that we have fixed an embedding $\Qb\hookrightarrow \Qpb$ and that $\Qb$ is a subfield of $\C$. We also choose an embedding $\sigma_p:\Qpb\hookrightarrow \C$ compatible with these choices.

We recall (\ref{eq:b_{tau}}) that every $b\in G(\mfk)$ defines a cocycle $\tau\mapsto b_{\tau}:=\Nm_{i(\tau)}b$ in $Z^1(W_{\Qp},G)$ ($\tau|_{\mfk}=\sigma^{i(\tau)}$) and that $J_b$ is the $\Qp$-group (\ref{eq:J_b}) such that
$J_b(\Qp)$ is equal to the $W_{\Qp}$-invariants of $G(\bar{\mfk})$ for the twisted Galois action $x\mapsto b_{\tau}\tau(x)b_{\tau}^{-1}$  (\cite[3.3]{Kottwitz97}).

%%%%%%%%%%%%%%%%%%%%
%%%%%%%%%%%%%%%%%%%%
\begin{thm} \label{thm:Kisin17_Cor.1.4.13,Prop.2.1.3,Cor.2.2.5} 
As a set with action by $\langle\Phi\rangle\times Z_G(\Qp)\rtimes G(\A_f^p)$, $\sS_{\mathbf{K}_p}(\Fpb)$ is a disjoint union of subsets $S(\mathscr{I})$
%, called \emph{isogeny classes}:
\begin{equation} \label{eq:isogeny_decomp}
 \sS_{\mathbf{K}_p}(\Fpb)\isom \bigsqcup_{\mathscr{I}}S(\mathscr{I}),
\end{equation}
wihere

(1) For each isogeny class $S(\mathscr{I})$, there exist a connected reductive $\Q$-group $I_{\mathscr{I}}$, an element $b\in G(\mfk)$ with $\kappa_G([b])=-\mu^{\natural}$ (\ref{eqn:mu_natural}), and embeddings of group schemes (over $\Qv$ for finite $v$) \[ i_p:(I_{\mathscr{I}})_{\Qp} \hookrightarrow J_b,\quad i_l:(I_{\mathscr{I}})_{\Ql} \hookrightarrow G_{\Ql}\] such that for almost all $l\neq p$, $i_l$ extends to an embedding $i_l:(I_{\mathscr{I}})_{\Z_l}\hookrightarrow G_{\Z_l}$ between reductive $\Z_l$-group schemes $(I_{\mathscr{I}})_{\Z_l}$, $G_{\Z_l}$, and in terms of which, one has \[S(\mathscr{I}):=I_{\mathscr{I}}(\Q)\backslash [X(\{\mu_X^{-1}\},b)_{\mathbf{K}_p}\times G(\A_f^p)].\] Here, $I_{\mathscr{I}}$ acts on $X(\{\mu_X^{-1}\},b)_{\mathbf{K}_p}\times G(\A_f^p)$ diagonally via $i_p\times i^p$, where $i^p$ denotes the restricted product $\prod_{l\neq p}'i_l$, and $\Phi$ acts on $S(\mathscr{I})$ via its action on $X(\{\mu_X^{-1}\},b)_{\mathbf{K}_p}$ by $(b\sigma)^r$ ($r=[\kappa(\wp):\Fp]$) while $g\in G(\A_f^p)$ acts on $S(\mathscr{I})$ via its right translation of $G(\A_f^p)$.

The following statements hold:

(2) The $\Q$-group $I_{\mathscr{I}}$ has the same $\Qb$-rank as $G$, 
and there exists an embedding $Z(G)\subset I_{\mathscr{I}}$ such that $(I_{\mathscr{I}}/Z(G))_{\R}$ is a subgroup of a compact inner form of $G^{\ad}_{\R}$.

(3) For every maximal $\Q$-torus $T$ of $I_{\mathscr{I}}$, there exists a stable conjugacy class of $\Q$-embeddings
\[i_T:T\hookrightarrow G \]
(\ref{defn:stable_conj_MT,SD}) with the following properties:
\setlist{nolistsep}
\begin{itemize} [noitemsep]
\item[(i)] It contains a member $i_T$ for which there exist a $G(\mfk)\times G(\A_f^p)$-conjugate of the triple of (1), denoted again by $(b,i_p, i^p)$, with $b\in i_T(T)(\mfk)$, and $h\in X\cap \Hom(\dS,i_T(T)_{\R})$ such that the embeddings $i_p:(I_{\mathscr{I}})_{\Qp}\hookrightarrow J_b$, $i^p:(I_{\mathscr{I}})_{\A_f^p}\hookrightarrow G_{\A_f^p}$
are $T_{\Qp}\times T_{\A_f^p}$-equivariant:
\begin{equation} \label{eq:(i_p,i^p)_adapted_to_i_T}
i_p|_{T_{\Qp}}=i_T|_{T_{\Qp}},\quad {i^p}|_{T_{\A_f^p}}=i_T|_{T_{\A_f^p}},
\end{equation}
and that 
the $\sigma$-conjugacy class $[b]$ in $B(i_T(T))$ equals $-\underline{\mu_h}\in X_{\ast}(i_T(T))_{\Gamma_p}$.%%{sign convention}
\footnote{The minus sign in $-\underline{\mu_h}$ is due to our sign convention resulting from working with homology instead of cohomology.}
\item[(ii)] Two $\Q$-embeddings $T\hookrightarrow I_{\mathscr{I}}$ which are stably conjugate under $I_{\mathscr{I}}(\Qb)$ give rise to the same stable conjugacy class of $\Q$-embeddings $i_T$ (Definition \ref{defn:stable_conj_MT,SD}).
\end{itemize}

(4) For any special Shimura subdatum $(T,h)$ of $(G,X)$, if $\mathscr{I}=\mathscr{I}_{T,h}$ is the isogeny class of the reduction of the special point $[h,1]\in Sh_{\mathbf{K}_p}(G,X)(\Qb)$, there exists an $I_{\mathscr{I}}(\Q)$-conjugacy class of $\Q$-embeddings 
\[j_{T,h}:T\hookrightarrow I_{\mathscr{I}}\] 
with the following properties: 
\begin{itemize} [noitemsep]
\item[(iii)] For any choice of $j_{T,h}$ in its conjugacy class, the associated (by the claim of (3)) stable conjugacy class of $\Q$-embeddings $i_T:T\hookrightarrow G$ contains the inclusion $i:T\subset G$. The pair $(i,h)$ qualifies as a special Shimura subdatum in (i) for $\mathscr{I}_{T,h}$ and $j_{T,h}$.

\item[(iv)] For any isogeny class $\mathscr{I}$ and each maximal $\Q$-torus $T\subset I_{\mathscr{I}}$, the isogeny class $\mathscr{I}_{i_T(T),h}$ attached to the special Shimura subdatum $(i_T(T),h\in X\cap\Hom(\dS,X_{\ast}(i_T(T))_{\R}))$ in (i) equals $\mathscr{I}$, and the $I_{\mathscr{I}}(\Q)$-conjugacy class of embeddings $T\stackrel{i_T}{\rightarrow}i_T(T)\stackrel{j_{i_T(T),h}}{\hookrightarrow} I_{\mathscr{I}}$ contains the inclusion $T\subset I_{\mathscr{I}}$. 
\end{itemize}
\end{thm}

Some explanations are in order.
%%%%%%%%%%%%%%%%%%%%
\begin{rem} \label{rem:isogeny_decomposition}
(1) The decomposition (\ref{eq:isogeny_decomp}) will be obtained from (\ref{eq:isogeny_decomp:geometric}) via the following choice:
For each isogeny class $\mathscr{I}$, we fix an $\Fpb$-point $x=x_{\mathscr{I}}$ lying in it and define $I_{\mathscr{I}}$ to be $I_{x}$. By choosing isomorphisms $\eta_l\ (l\neq p)$ (\ref{eq:isom_eta}), $\eta_p^{\nr}$ (\ref{eq:isom_eta_nr}) (defined over $\mfk$) such that for almost all $l\neq p$, $\eta_l$ extends over $\Z_l$, we will then obtain an element $b\in G(\mfk)$ by that $b\sigma$ is the absolute Frobenius element acting on $H^{\cris}_1(\mathcal{A}_x/\mfk)$ via $\eta_p^{\nr}$, and embeddings $i_l:=i_{x,l}:(I_{x})_{\Ql}\hookrightarrow G_{\Ql}\ (l\neq p)$, $i_p:=i_{x,p}:(I_{x})_{\Qp}\hookrightarrow J_b$ of group schemes over $\Ql$ and $\Qp$. 
In the rest of this subsection, in the proofs of \emph{theorems} we will remember and use this geometric origin of the decomposition (\ref{eq:isogeny_decomp}).%%
\footnote{Some further properties of the objects appearing in this theorem will be added later as three more theorems (Theorem \ref{thm:Kisin17_Cor.2.3.2}, \ref{thm:Kisin17_Prop.4.4.8}, \ref{thm:equiv_K-triples}). But these properties will not be deduced from the properties stated in this theorem. Rather they will be additional properties which are group-theoretic statements but are established by geometric means (just like the properties stated in this theorem). 
From a logical viewpoint, it might have been more desirable to combine all these theorems into one big theorem which presents all the objects and their complete properties at once. But for easy readability, we decided to break it into four theorems.}

(2) Regarding Theorem \ref{thm:Kisin17_Cor.1.4.13,Prop.2.1.3,Cor.2.2.5} (3), note that for a triple $(b,i_p,i^p)$ as in (1) an equivalent triple 
\[(b':=g_pb\sigma(g_p^{-1}),\ \Int(g_p)\circ i_p:(I_{\mathscr{I}})_{\Qp}\hookrightarrow J_b\isom J_{b'},\ \Int(g^p)\circ i^p),\] 
($(g_p,g^p)\in G(\mfk)\times G(\A_f^p)$) again satisfies (1): left multiplication by $g_p\times g^p$ induces a $\langle\Phi\rangle\times Z_G(\Qp)\rtimes G(\A_f^p)$-equivariant bijection
\[I_{\mathscr{I}}(\Q)\backslash X(\{\mu_X^{-1}\},b)_{\mathbf{K}_p}\times G(\A_f^p) \isom I_{\mathscr{I}}(\Q)\backslash X(\{\mu_X\},b')_{\mathbf{K}_p}\times G(\A_f^p).\]
Furthermore, as will be evident from the proof, each equivalent triple enjoys all the other properties; thus we will not distinguish them.

(3) In Theorem \ref{thm:Kisin17_Cor.1.4.13,Prop.2.1.3,Cor.2.2.5} (3) (i), we note that when $b\in  i_T(T)(\mfk)$, for any $\Qp$-algebra $R$, $i_T(T)(R) (\subset G(R)\subset G(\mfk\otimes R))$ lies in $J_{b}(R)=G(\mfk\otimes R)^{\langle\sigma\rangle}$ (\ref{eq:J_b}).
% which is regarded as a subset of $G(R)$ via the map $\Xi$ (\ref{eq:Xi2}) (which is induced by the canonical $\mfk$-algebra homomorphism $\mfk\otimes R\rightarrow R:l\otimes x\mapsto lx$); so the proper formulation for (\ref{eq:(i_p,i^p)_adapted_to_i_T}) will be $\Xi\circ i_p|_{T_{\Qp}}=i_T|_{T_{\Qp}}$. 
%(but, the restriction of $\Xi$ to $J_b(\Qp)$ equals the inclusion $G(\mfk)$).

We will say that a triple $(b,i_p,i^p)$ in (1) is \textit{adapted to} a given embedding $i_T:T\hookrightarrow G$ if it satisfies the conditions in Theorem \ref{thm:Kisin17_Cor.1.4.13,Prop.2.1.3,Cor.2.2.5} (3) (i).
For fixed $i_T$, any two triples $(b,i_p,i^p)$ adapted to $i_T$ (in fact, satisfying just the first condition of (i)) differ by conjugation by an element of $i_T(T)(\mfk)\times i_T(T)(\A_f^p)$ (since $\mathrm{rk}(T)=\mathrm{rk}(G)$ and $T(\Qp)$ is Zariski-dense in $T_{\mfk}$). In particular, the $\sigma$-conjugacy class of $b$ in $B( i_T(T))$ is uniquely determined by $\mathscr{I}$ once $i_T$ is fixed. 
%Also, when $(b,i_p,i^p)$ adapted to $(i_T, h)$, the two maps $i_T:T(R)\isom i_T(T)(R)\subset G(R)\subset G(\mfk\otimes R)$, $ i_p: T(R)\hookrightarrow J_b(R)\subset G(\mfk\otimes R)$ are equal. )

(4) Let us call a transfer of maximal torus $\Int(g):T\hookrightarrow G\ (g\in G(\Qb))$ \emph{$\R$-rational} if $\Int(g)|_{T_{\R}}=\Int(g_{\infty})|_{T_{\R}}$ for some $g_{\infty}\in G(\R)$.
Then, if $i_T$ and $(i_T(T),h)$ are the $\Q$-embedding and the special Shimura subdatum in (i), for any transfer $\Int(g):T\hookrightarrow G$ of the maximal torus $T$ which is $\R$-rational, the new $\Q$-embedding $i_T':=\Int(g)\circ i_T$ and $\Int(g)\circ h$ also fulfill the condition (i). For the $\Q$-embedding $i_T$ in (i), we allow any member in its $\R$-rational, stable conjugacy class (not just in its stable conjugacy class) as such embedding  (see the remark after Definition \ref{defn:stable_conj_MT,SD}). 
\end{rem}

Before moving on,
we derive the first (primitive) description of the fixed point set of the Frobenius-twsted Hecke correspondence $\Phi^m\circ f$ (\ref{eq:Hecke_corr_twisted_by_Frob}) acting on an isogeny class $S(\mathscr{I})_{\mathbf{K}}:=S(\mathscr{I})/\mathbf{K}^p\subset \sS_{K}(\Fpb)$: 
\[S(\mathscr{I})_{\mathbf{K}}\stackrel{p_1'}{\longleftarrow} S(\mathscr{I})_{\mathbf{K}_g} \stackrel{p_2'=\Phi^m\circ p_2}{\longrightarrow} S(\mathscr{I})_{\mathbf{K}}.\]
From now on, we assume $\mathbf{K}^p$ to be small enough so that the following conditions hold:
\begin{align} \label{item:Langlands-conditions}
(a) &\text{ If }\epsilon\in I_{\mathscr{I}}(\Q)\text{ and }\epsilon x= xg\text{ for some }x\in (X^p(\phi)/\mathbf{K}^p_g)\times X_p(\phi),\text{ then }\epsilon\in Z(\Q)_K:=Z(\Q)\cap K. \\
(b) &\ I_{\mathscr{I}}^{\der}\cap K\cap Z(G)(\Q)=\{1\}. \nonumber
\end{align}
This is possible by \cite[p.1171-1172]{Langlands79} (cf. \cite[1.3.7, 1.3.8]{Kottwitz84b}, \cite[Lem.5.5]{Milne92}).
Then, an elementary argument (see \cite[$\S$1.4]{Kottwitz84b}, \cite[Lem.5.3]{Milne92}) shows that the fixed point set decomposes into disjoint subsets (cf. \cite[1.4.3, 1.4.4]{Kottwitz84b}):
\begin{align} \label{eq:fixed_pt_set_of_Heck-corresp2}
S(\mathscr{I})_{\mathbf{K}}^{\Phi^m\circ f=\mathrm{Id}} &= \bigsqcup_{\epsilon} I_{\mathscr{I},\epsilon}(\Q)\backslash X(\mathscr{I},\epsilon)_{\mathbf{K}_g} \\
&:= \bigsqcup_{\epsilon} I_{\mathscr{I},\epsilon}(\Q)\backslash \bigl[ X_p(\mathscr{I},\epsilon)\times X^p(\mathscr{I},\epsilon,g)/\mathbf{K}^p_g \bigr] ,\nonumber
\end{align}
where the index $\epsilon$ runs through a set of representatives in $I_{\mathscr{I}}(\Q)$ for the conjugacy classes of $I_{\mathscr{I}}(\Q)/Z(\Q)_K$, $I_{\mathscr{I},\epsilon}$ is the centralizer of $\epsilon$ in $I_{\mathscr{I}}$ (regarded as an algebraic $\Q$-subgroup of $I_{\mathscr{I}}$), and 
\begin{align*}
X_p(\mathscr{I},\epsilon) &:=\{\ x_p\in X(\{\mu_X^{-1}\},b)_{\mathbf{K}_p}  \ \ |\ \  i_p(\epsilon) x_p=(b\sigma)^n x_p\ \}, \\
X^p(\mathscr{I},\epsilon,g) &:=\{\ x^p\in G(\A_f^p) \ |\ \  i^p(\epsilon) x^pg=x^p\text{ mod }\mathbf{K}^p\ \}.
\end{align*}

\begin{proof} (of Theorem \ref{thm:Kisin17_Cor.1.4.13,Prop.2.1.3,Cor.2.2.5})
(1) This is how we construct the decomposition.

(2) This is \cite[2.1.7]{Kisin17}

(3) It is shown in the proof of \cite[Thm. 2.2.3]{Kisin17} that for every maximal $\Q$-torus $T\subset I_{x}\subset \Aut_{\Q}(\mathcal{A}_{\bar{x}})$, and for any choice of a cocharacter $\mu_T\in X_{\ast}(T)$ satisfying the conditions of \cite[Lem.2.2.2]{Kisin17} (in particular, $\mu_T$ must lie in the conjugacy class $c(G,X)$),
there exist a point $y$ in the isogeny class of $x$, endowed with a quasi-isogeny $\mathcal{A}_{x}\rightarrow \mathcal{A}_{y}$ preserving the extra structures (defined over a finite extension of $k$), which lifts to a $K$-valued point $\tilde{y}$ of $\sS_{\mathbf{K}_p}$ for a finite extension $K\subset\Qpb$ of $K_0=\mathrm{Frac}(W(k))$ ($\mu_T$ must be defined over $K$). This lifting is determined from $\mu_T$ by the condition that the Hodge filtration on $H^1_{\cris}(\mathcal{A}_{y}/K)\cong H^1_{\dR}(\mathcal{A}_{\tilde{y}}/K)$ defined by $\mathcal{A}_{\tilde{y}}$ is the filtration induced by $\mu_T^{-1}$. This implies that the action of $T$ on $\mathcal{A}_{y}$ (in the isogeny category) lifts to $\mathcal{A}_{\tilde{y}}$, and $\sigma_p(\tilde{y})$ is a special (=CM) point on $Sh_{\mathbf{K}_p}(G,X)(\Qb)$ (as the ranks of $I_x$ and $G$ are the same).
%; the existence of $\eta_{\Betti}$ is due to the moduli interpretation of $\Sh_{\mathbf{K}}(G,X)(\C)$. 
That is, when we fix such $\mu_T$, via a choice of an isomorphism of $\Q$-vector spaces endowed with a set of tensors
\begin{equation} \label{eq:Betti-isom}
\eta_{\Betti}:(H^{\Betti}_1(\mathcal{A}_{\sigma_p(\tilde{y})},\Q),\{s_{\alpha,\Betti,\sigma_p(\tilde{y})}\}) \lisom (V,\{s_{\alpha}\})
\end{equation} 
which is well-defined up to action of $G(\Q)\subset\GL(V)$, we obtain an embedding and a special Shimura subdatum
\begin{equation} \label{eq:CM-lifting_via_T}
i_T:T\hookrightarrow G,\quad h\in X\cap \Hom(\dS,i_T(T)_{\R})
\end{equation}
such that $\sigma_p(\mu_T)=\mu_h$ under the canonical isomorphism $H^1_{\dR}(\mathcal{A}_{\tilde{y}}/K)\otimes_{\sigma_p}\C\cong H^1_{\Betti}(\mathcal{A}_{\sigma_p(\tilde{y})},\C)$ ($\mu_T$ is defined over $\Qpb$ and thus over $\C$ via $\sigma_p$; in fact, it is defined over $\Qb\subset\C$ as $\mu_T\in X_{\ast}(T)$).

For such embedding $i_T$ and $h$, we claim that there exist isomorphisms $\eta_l$ (\ref{eq:isom_eta}), $\eta_p^{\nr}$ (\ref{eq:isom_eta_nr}) which are \emph{$T$-equivariant} with respect to $i_T:T\hookrightarrow G$ and the action of $T$ on $\mathcal{A}_{y}$. 
By construction of $i_T$ via the choice of $\eta_B$ (\ref{eq:Betti-isom}), it suffices to find such $T$-equivariant isomorphisms with $(H^{\Betti}_1(\mathcal{A}_{\sigma_p(\tilde{y})},\Q),\{s_{\alpha,\Betti,\sigma_p(\tilde{y})}\})$ replacing $(V,\{s_{\alpha}\})$ (for the lifted action of $T$ on $\mathcal{A}_{\tilde{y}}$).
At $l\neq p$, this is clear since there exist \emph{canonical} isomorphisms of $\Ql$-vector spaces matching the respective tensors:
\begin{equation} \label{eq:isom_epsilon_l}
\begin{split}
\epsilon_l:(H^{\et}_1(\mathcal{A}_{\bar{y}},\Ql),\{s_{\alpha,l,y}\}) &\lisom (H^{\et}_1((\mathcal{A}_{\tilde{y}})_{\Qpb},\Ql),\{s_{\alpha,l,\tilde{y}}\})\\
& \lisom (H^{\Betti}_1(\mathcal{A}_{\sigma_p(\tilde{y})},\Q),\{s_{\alpha,\Betti,\sigma_p(\tilde{y})}\})\otimes\Ql.
\end{split}
\end{equation}
which are thus $T$-equivariant (by functoriality): in fact, the tensors $s_{\alpha,l,\tilde{y}}$ are constructed by these isomorphisms, cf. \cite[(2.2)]{Kisin10}.
At $p$, we also have canonical isomorphisms of $\C$-vector spaces matching the respective tensors
\begin{equation} \label{eq:isom_epsilon_C}
\begin{split}
(H^{\cris}_1(\mathcal{A}_{y}/K),\{s_{\alpha,0,y}\})\otimes_{\sigma_p}\C & \lisom (H^{\dR}_1(\mathcal{A}_{\tilde{y}}/K),\{s_{\alpha,\dR,\tilde{y}}\})\otimes_{\sigma_p}\C\\ 
&\lisom (H^{\Betti}_1(\mathcal{A}_{\sigma_p(\tilde{y})},\Q),\{s_{\alpha,\Betti,\sigma_p(\tilde{y})}\})\otimes\C 
\end{split}
\end{equation}
which are $T$-equivariant.
For the fact that the first isomorphism matches the respective tensors, see the proof of \cite[Prop.2.3.5]{Kisin10} (cf. \cite[Prop.1.3.9]{Kisin17}).
Then, we consider the functor which assigns to a $K_0$-algebra $R$ the set of $R$-linear isomorphisms $H^{\cris}_1(\mathcal{A}_{y}/K_0)\otimes R \isom H^{\Betti}_1(\mathcal{A}_{\sigma_p(\tilde{y})},\Q)\otimes R$
which take $s_{\alpha,0,y}$ to $s_{\alpha,\Betti,\sigma_p(\tilde{y})}$ for every $\alpha$ and are $T$-equivariant.
Since $i_T(T)$ is a maximal torus of $G$ by (2) and $T(\Qp)$ is Zariski-dense in $T_{K_0}$ (as $T$ is unirational), this is a pseudo-torsor under $T_{K_0}$, which must then be a torsor as it is non-empty (\ref{eq:isom_epsilon_C}). So by Steinberg's theorem $H^1(\Qpnr,T)=\{1\}$, one can find a $\Qpnr$-isomorphism 
\begin{equation} \label{eq:isom_epsilon_ur_p} 
\epsilon_p:(H^{\cris}_1(\mathcal{A}_{y}/K_0),\{s_{\alpha,0,y}\})\otimes_{K_0}\Qpnr \lisom (H^{\Betti}_1(\mathcal{A}_{\sigma_p(\tilde{y})},\Q),\{s_{\alpha,\Betti,\sigma_p(\tilde{y})}\})\otimes\Qpnr. 
\end{equation}
matching the tensors and intertwining the two $T$-actions.

For such $T$-equivariant $\eta_l$, $\eta_p^{\nr}$, the resulting embeddings $i_{x,p}:I_{x}(\Qp)\simeq I_y(\Qp)\hookrightarrow J_b(\Qp)$, $i^p_{x}:I_{x}(\A_f^p)\simeq I_{y}(\A_f^p)\hookrightarrow G(\A_f^p)$ clearly satisfy (\ref{eq:(i_p,i^p)_adapted_to_i_T}), where $b\in G(\Qpnr)$ is given by that $b\sigma$ is the absolute Frobenius automorphism acting on $H^{\cris}_1(\mathcal{A}_x/\Qpnr)$ via $\eta_p^{\nr}$.
Note that we have $b\in i_T(T)(\Qpnr)$ because $b$ commutes with $i_T(T)(\Qp)$, thus $b\in Z_{G(\Qpnr)}(i_T(T)(\Qp))=i_T(T)(\Qpnr)$ (the equality holds again since $T$ is unirational so that $T(\Qp)$ is Zariski dense in $T$).
Since $\mathcal{A}_{x}$ is the reduction of the CM point $[h,1]\in Sh_{K_T}(i_T(T),\{h\})\in Sh_{\mathbf{K}}(G,X)$ (for $K_{T}:=i_T(T)(\A_f)\cap K$), the claim on $[b]\in B(i_T(T))$ of (i) follows from \cite[Lem.3.1.1]{Lee18b}:
%%{sign convention}

\begin{lem} \label{lem:Lee18_Lem.3.1.1} Let $(T,h)$ be a special Shimura subdatum.
The $T_{\Qp}$-isocrystal of the reduction of $[h,1]\in \sS_{\mathbf{K}_p}(\Qb)$ is the image of $-\underline{\mu_h}\in X_{\ast}(T_{\Qp})_{\Gamma_p}$ under the isomorphism $\kappa_{T_{\Qp}}:B(T_{\Qp})\isom X_{\ast}(T_{\Qp})_{\Gamma_p}$.
\end{lem}

\begin{proof}
Although a weaker claim is asserted in the statement, this is what is proved in \cite[Lem.3.1.1]{Lee18b}.
\end{proof}

Next, we prove that

$(\star)$ for any maximal $\Q$-torus $T\subset I_x$,  the resulting stable conjugacy class of $\Q$-embeddings $T\hookrightarrow G$ is independent of all the choices made in its construction, i.e. of the cocharacter $\mu_T$, the point $y$, the quasi-isogeny $\mathcal{A}_{x}\rightarrow \mathcal{A}_{y}$, and the isomorphism (\ref{eq:Betti-isom}).
%\footnote{These choices are not independent of each other: rather, $y$, $\tilde{y}$, $\mathcal{A}_{x}\rightarrow \mathcal{A}_{y}$ are either determined or their choices are restricted by choice of $\mu_T$, but we will not be concerned with this matter here.} 

In fact, we will show the following stronger statement which will establish (ii): 

Suppose given three $\Qb$-quasi-isogenies $\theta_i:\mathcal{A}_x\otimes\Qb\rightarrow \mathcal{A}_{y_i}\otimes\Qb\ (i=0,1,2)$ matching the weak polarizations, \'etale and crystalline tensors, and whose induced $\Qb$-isomorphism $\theta_{i\ast}:(I_x)_{\Qb}\isom (I_{y_i})_{\Qb}$ of algebraic groups restrict to a $\Q$-isomorphism $T_i\isom T_{y_i}$ of $\Q$-maximal tori, where $y_0=x$, $T_1=T_0=:T$ and $T_2=\theta_0(T)$, and assume that for $i=1,2$, $y_i$ together with $T_{y_i}$-action (in the isogeny category) lifts to a point $\tilde{y}_i$ in characteristic zero, giving a $\Q$-embedding $T_{y_i}\hookrightarrow G$ via an isomorphism (\ref{eq:Betti-isom}) (which we fix one for each $i=1,2$). Then, the resulting embeddings $i_1,i_2\circ \theta_{0\ast}:T\isom T_{y_i}\hookrightarrow G$ are stably conjugate. 

For this, we consider the functor, defined on the category of $\Q$-algebras $R$, of $R$-linear isomorphisms between $\Q$-vector spaces
\begin{equation} \label{eq:pseudo-torsor1}
 H^{\Betti}_1(\mathcal{A}_{\sigma_p(\tilde{y}_1)},\Q)\otimes R \ \isom \ H^{\Betti}_1(\mathcal{A}_{\sigma_p(\tilde{y}_2)},\Q)\otimes R
 \end{equation} 
which preserve the Betti tensors $\{s_{\alpha,\tilde{y}_1}\}$, $\{s_{\alpha,\tilde{y}_2}\}$ and the $T$-actions $\theta_{1\ast}|_{T}$, $\theta_2\circ\theta_{0\ast}|_{T}$.
This functor is a pseudo-torsor under $T=Z_G(T)$, so will be a $T$-torsor if it is non-empty. The fact that $T$ is a maximal torus of $G$ also implies that any $\Qb$-valued point of this pseudo-torsor would be given by $\Int(g)$ $(g\in G(\Qb))$ such that $i_2\circ \theta_{0\ast}=\Int(g)\circ i_1$. But, the $\Qb$-isogeny $\theta_2\circ \theta_0\circ\theta_1^{-1}$ induces isomorphisms
\begin{align*}
 H^{\et}_1(\mathcal{A}_{\bar{y}_1},\Ql)\otimes\Qlb \ & \isom  H^{\et}_1(\mathcal{A}_{x},\Ql)\otimes\Qlb \isom \ H^{\et}_1(\mathcal{A}_{\bar{y}_2},\Ql)\otimes\Qlb,\\ 
 H^{\cris}_1(\mathcal{A}_{y_1}/\mfk)\otimes\mfkb \ & \isom  H^{\cris}_1(\mathcal{A}_{x}/\mfk)\otimes\mfkb \isom \ H^{\cris}_1(\mathcal{A}_{y_2}/\mfk)\otimes\mfkb
\end{align*}
which preserve the tensors $\{s_{\alpha,l,y_i}\}$, $\{s_{\alpha,0,y_i}\}$  $(i=1,2)$ and the $T$-actions, thus provide $\Qlb$ or $\mfkb$-points of (\ref{eq:pseudo-torsor1}) via the canonical isomorphisms (\ref{eq:isom_epsilon_l}), (\ref{eq:isom_epsilon_ur_p}), respectively.
We have established statement (ii), thus (3).

For (4), there exist a $K$-valued point $\tilde{y}$ of $\sS_{\mathbf{K}_p}$ defined over a finite extension $K\subset \Qpb$ of $\Qp$ such that $[h,1]=\sigma_p(\tilde{y})$ and an identification $H^{\Betti}_1(\mathcal{A}_{\sigma_p(\tilde{y})},\Q)=V$ matching the weak polarizations and the respective tensors $\{s_{\alpha,\Betti,\sigma_p(\tilde{y})}\}$, $\{s_{\alpha}\}$; let $y$ be its reduction. As $h$ factors through $T$, the composite map $T\subset G\subset \GL(V)=\Aut(H^{\Betti}_1(\mathcal{A}_{\sigma_p(\tilde{y})},\Q))$ induces a homomorphism $T\hookrightarrow \Aut_{\Q}(\mathcal{A}_{\tilde{y}})$, and the embedding $j_{y}:T\hookrightarrow \Aut_{\Q}(\mathcal{A}_{\tilde{y}})\hookrightarrow \Aut_{\Q}(\mathcal{A}_{y})$ factors through the subgroup $I_{y}$.
If $x$ is the prechosen point of $\mathscr{I}=\mathscr{I}_{T,h}$, we define $j_{T,h}$ to be the composite $T\hookrightarrow I_{y}\isom I_{x}$ for any isomorphism $I_{y}\isom I_{x}$ induced by a choice of a quasi-isogeny $\mathcal{A}_{y} \rightarrow \mathcal{A}_{x}$ preserving the (\'etale and crystalline) tensors and the weak polarizations, and hence this embedding is well-defined up to $I_x(\Q)$-conjugacy. Now, in view of the property $(\star)$, statements (iii) and (iv) are immediate from this construction: property $(\star)$ allows us to take any quasi-isogeny $\mathcal{A}_{y} \rightarrow \mathcal{A}_{x}$ (preserving the extra structures) to fix $j_{T,h}$, and in the construction of $i_T$ to take any $\mu_T\in X_{\ast}(T)$ giving a point in $\mathscr{I}$ which is endowed with a $T$-action and lifts (with $T$-action) to a CM-point. Thus we can take $\mu_T$ such that it gives $y$ and $\tilde{y}$ and use the accompanying quasi-isogeny $\mathcal{A}_{y} \rightarrow \mathcal{A}_{x}$, in which case $i_T$ and the special Shimura subdatum in (\ref{eq:CM-lifting_via_T}) become just the given embedding $T\subset G$ and the original datum $(T,h)$, when one uses the identification $H^{\Betti}_1(\mathcal{A}_{\sigma_p(\tilde{y})},\Q)=V$ for the isomorphism (\ref{eq:Betti-isom}). 
Statement (v) is proved by a similar method as (iii). 
\end{proof}

%%%%%%%%%%%%%%%%%%%%
\begin{defn} \label{defn:admissible_pair2}
A \emph{K-pair} is a pair $(\mathscr{I},\epsilon)$ consisting of an isogeny class $\mathscr{I}\subset\sS_{\mathbf{K}}(\Fpb)$ and an element $\epsilon$ of $I_{\mathscr{I}}(\Q)$. 

A K-pair $(\mathscr{I},\epsilon)$ is said to be \emph{admissible} of level $n=m[\kappa(\wp):\Fp]$ ($m\in\N$) if 
for a triple $(b\in G(\mfk),i_p,i^p)$ in Theorem \ref{thm:Kisin17_Cor.1.4.13,Prop.2.1.3,Cor.2.2.5},
and there exists $x_p\in G(\mfk)/\tilde{\mathbf{K}}_p$ such that 
\[i_p(\epsilon)x_p=(b\sigma)^nx_p.\]
According to \cite[Lemma 1.4.9]{Kottwitz84b}, this condition is the same as the existence of $x\in G(\mfk)$ such that \[ i_p(\epsilon)x \sigma^n=(b\sigma)^nx \] 
(as usual, an equality in $G(\mfk)\rtimes\langle\sigma\rangle$), namely $x^{-1}i_p(\epsilon)x =\Nm_n (x^{-1}b \sigma(x))$. Two K-pairs $(\mathscr{I},\epsilon)$, $(\mathscr{I}',\epsilon')$ are said to be \emph{equivalent} if $\mathscr{I}=\mathscr{I}'$ and $\epsilon'=\Int(g)(\epsilon)$ for some $g\in I_{\mathscr{I}}(\Q)$.

A K-pair $(\mathscr{I},\epsilon)$ is said to be \emph{$\mathbf{K}_p$-admissible} if one can find $x_p$ above further satisfying that $x_p\in X(\{\mu_X^{-1}\},b)_{\mathbf{K}_p}$ (\ref{eq:X_p(mu_X,b)}).
\end{defn}
Clearly, the admissibility condition in this definition does not depend on the choice of a representative $(b,i_p)$ in its equivalence class. 
%Note that any rational element $\epsilon$ of $I_{\mathscr{I}}$ is semisimple).

Next, with any admissible K-pair $(\mathscr{I},\epsilon\in I_{\mathscr{I}}(\Q))$, say, of level $n=mr$, we associate a Kottwitz triple by imitating a similar recipe in \cite{Kottwitz84b}, \cite{LR87} (for what is called an admissible pair there). Choose a representative of the triple $(b,i_p,i^p)$ in Theorem \ref{thm:Kisin17_Cor.1.4.13,Prop.2.1.3,Cor.2.2.5} (1).
First, for any $c\in G(\mfk)$ such that 
\begin{equation} \label{eq:(epsilon,b,c)->delta2}
c(i_p(\epsilon)^{-1}(b\sigma)^n)c^{-1}=\sigma^n,
\end{equation}
one readily checks that $\delta:=cb\sigma(c^{-1})\in G(L_n)$ and $\Nm_n\delta=ci_p(\epsilon)c^{-1}$ (see \cite[$\S$1.4]{Kottwitz84b}): the $\sigma$-conjugacy class of $\delta$ in $G(L_n)$ depends only on the K-pair $(\mathscr{I},\epsilon)$ (i.e. is independent of the choice of either $c$ or a representative $(b,i_p)$ in its equivalence class).
We put $\gamma=(\gamma_l)_{l\neq p}:=i^p(\epsilon)$: its $G(\A_f^p)$-conjugacy class is also uniquely determined by $(\mathscr{I},\epsilon)$. 
To define $\gamma_0\in G(\Q)$, we choose a maximal $\Q$-torus $T\subset I_{\mathscr{I}}$ with $\epsilon\in T(\Q)$, and fix an embedding $i_T:T\hookrightarrow G$ in Theorem \ref{thm:Kisin17_Cor.1.4.13,Prop.2.1.3,Cor.2.2.5} (3). Then, we obtain a triple of elements in $G(\Q)\times G(L_n)\times G(\A_f^p)$
\begin{equation} \label{eq:K-triple_attached_to_admissible_K-pair}
(\gamma_0;\gamma,\delta):=(i_T(\epsilon);i^p(\epsilon),cb\sigma(c^{-1})),
\end{equation} 
where $n$ is the level of $(\mathscr{I},\epsilon)$.

%%%%%%%%%%%%%%%%%%%%
\begin{lem} \label{lem:K-triple_attached_to_admissible_K-pair}
The triple (\ref{eq:K-triple_attached_to_admissible_K-pair}) is a Kottwitz triple (Definition \ref{defn:Kottwitz_triple}). When the triple $(b,i_p,i^p)$ is adapted to $i_T:T\hookrightarrow G$ (i.e. the condition (i) of Theorem \ref{thm:Kisin17_Cor.1.4.13,Prop.2.1.3,Cor.2.2.5} (3) holds), it satisfies the conditions of Definition \ref{defn:stable_Kottwitz_triple} except for (i$'$) and has trivial Kottwitz invariant. In general, the stable conjugacy class of $\gamma_0$ depends possibly on $T$, but not on the choice of $i_T$. 
\end{lem}

We call any such triple $(\gamma_0;\gamma,\delta)$ the \emph{Kottwitz triple attached to} the admissible K-pair $(\mathscr{I},\epsilon)$.
This should be distinguished from the notion of ``\emph{stable} Kottwitz triple attached to $(\mathscr{I},\epsilon)$'', to be defined later in Definition \ref{defn:stable_K-triple_attached_to_admissible_K-pair}.

\begin{proof}
That the triple (\ref{eq:K-triple_attached_to_admissible_K-pair}) is a Kottwitz triple is obvious: $\gamma_l$ is $G(\Qlb)$-conjugate to $\gamma_0$ since there is a $G(\Ql)$-conjugate of $\gamma_l=i_l(\epsilon)$ which is so (i.e. when $(b,i_p,i^p)$ is adapted to $i_T$). So, the $G(\Qb)$-conjugacy class of $\gamma_0$ does not depend on the choices of either $T$ or $i_T$ since the same is true of the $G(\Qlb)$-conjugacy class of $\gamma_l$. The condition $\ast(\delta)$ follows from Theorem \ref{thm:Kisin17_Cor.1.4.13,Prop.2.1.3,Cor.2.2.5} (i). 
Next, assume that $(b,i_p,i^p)$ is adapted to $i_T$. It is then obvious that $\gamma_0$ and $\gamma(=\gamma_0)$ are stably conjugate under $G(\bar{\A}_f^p)$, and from $\Nm_n\delta=ci_p(\epsilon)c^{-1}$, we see that $c^{-1}\delta\sigma(c)=b\in i_T(T)(\mfk)\subset I_0(\mfk)$.
To see that for suitable elements $(g_v)_v\in G(\bar{\A}_f^p)\times G(\mfk)$, the Kottwitz invariant $\alpha(\gamma_0;(\gamma_l)_{l\neq p},\delta;(g_v)_v)$ vanishes, we take $g_v=1$ for $v=l\neq p$ and $g_p=c$ above. Then, we have $\alpha_l=1$ for $l\neq p$, and $\alpha_v$'s for $v=p,\infty$ are $\pm1$ times the restrictions of the character of $\hat{T}$ defined by $\mu_h$ with different signs, from which the vanishing of $\alpha(\gamma_0;(\gamma_l)_{l\neq p},\delta;(g_v)_v)$ is obvious;
we remark that as observed before \cite[p.172]{Kottwitz90}, this further establish condition (iv$'$) of Definition \ref{defn:stable_Kottwitz_triple} (this fact will be also established by direct means in Corollary \ref{cor:Tate_thm2}, according to which for the group $I$ in (iv$'$), one can take $I_{\mathscr{I},\epsilon}^{\mathrm{o}}$, the neutral component of the centralizer of $\epsilon$ in $I_{\mathscr{I}}$).
Finally, the stable conjugacy class of $\gamma_0=i_T(\epsilon)$ obtained from an embedding $i_T:T\hookrightarrow G$ in Theorem \ref{thm:Kisin17_Cor.1.4.13,Prop.2.1.3,Cor.2.2.5} (3) does not depend on $i_T$, since two such embeddings $i_T$ are stably conjugate to each other (in the sense of Definition \ref{defn:stable_conj_MT,SD}) by \emph{loc. cit.} 
\end{proof}

%%%%%%%%%%%%%%%%%%%%
\begin{rem}
The arguments sketched in the introduction (\ref{subsec:sketch_of_proof}) work equally well for geometrically effective K-pairs $(\mathscr{I},\epsilon)$ (i.e. pairs with $S(\mathscr{I},\epsilon)_{\mathbf{K}}\neq \emptyset$) and their associated ``geometrically effective'' Kottwitz triples, instead of admissible K-pairs and effective Kottwitz triples. And it even appears more natural to do so for the purpose of proving (\ref{eq:E}) and (\ref{eq:C}). But the main utility of admissible K-pairs is to produce Kottwitz triple and the ``geometric effectivity'' is too rigid for that purpose.
\end{rem}

%%%%%%%%%%%%%%%%%%%%
\begin{lem}  \label{lem:canonical_decomp_of_epsilon} 
For any $\Q$-group $T$ of multiplicative type such that $\Q\operatorname{-}\mathrm{rk}(T)=\R\operatorname{-}\mathrm{rk}(T)$, there exist a positive integer $s$ and elements $\pi_0$, $t\in T(\Q)$ such that
\begin{align*} 
%\label{eq:virtual_admissible_pair} 
(a)\quad & \epsilon^s=\pi_0 t,\\ 
(b)\quad &\pi_0\in K_l\text{ for all }l\neq p, \nonumber \\ 
(c)\quad &t\in K_p, \nonumber 
\end{align*}
where for each finite place $v$ of $\Q$, $K_v$ denotes the maximal compact subgroup of $T^{\mathrm{o}}(\Q_v)$.
The pair $(\pi_0,t)$ is uniquely determined by $\epsilon$, up to taking simultaneous powers. In particular, the construction of $(\pi_0,t;s)$ is functorial in $(T,\epsilon)$: $f:(T,\epsilon) \rightarrow (T',\epsilon')$ is a morphism of pairs as above, $f$ matches the corresponding elements $(\pi_0,t)$, $(\pi_0',t')$ (for the same $s$)
\end{lem}

Clearly, if one wants, one may further assume that $\pi_0,t\in T^{\mathrm{o}}(\Q)$.

\begin{proof}
(cf. \cite[Lem.10.12]{Kottwitz92})
The uniqueness up to taking simultaneous powers of a pair of elements $(\pi_0,t)$ satisfying the properties (a) - (c) is an easy consequence of the property that $T(\Q)$ is discrete in $T(\A_f)$ (implied by $\Q\operatorname{-}\mathrm{rk}(T)=\R\operatorname{-}\mathrm{rk}(T)$). For existence, it follows from the same property that the canonical map $\varphi:T(\Q)\rightarrow X:=\oplus_{v\neq\infty} (T(\Qv)/K_v)$ has finite kernel and cokernel, so if we consider the two elements $(a_v)_v,\ (b_v)_v\in X$ where $a_v=1$, $b_v=\epsilon\mod K_v$ for $v\neq p$, and $a_v=\epsilon\mod K_p$, $b_p=1$, then some (common) power $(a_v)^r$, $(b_v)^r$ are the images of elements $a$, $b$ of $T(\Q)$. As $\epsilon^r$ and $ab$ have the same images under $\varphi$, some powers of them are equal: $\epsilon^s=\pi_0t$ where $\pi_0:=a^{s/r}, t:=b^{s/r}$.
\end{proof}

%%%%%%%%%%%%%%%%%%%%
\begin{prop} \label{prop:canonical_decomp_of_epsilon2} 
Assume that $\Q\operatorname{-}\mathrm{rk}(Z(G))=\R\operatorname{-}\mathrm{rk}(Z(G))$.
For any admissible K-pair $(\mathscr{I},\epsilon)$, there exist $s\in\N$ and $\pi_0, t\in T_{\epsilon}(\Q)$,
where $T_{\epsilon}$ is the $\Q$-subgroup (of multiplicative type) of $I_{\mathscr{I}}$ generated by $\epsilon\in I_{\mathscr{I}}(\Q)$, satisfying the properties of Lemma \ref{lem:canonical_decomp_of_epsilon}: we have 
\[(a)\ \epsilon^s=\pi_0 t\ ;\qquad (b)\ \pi_0\in K_l \text{ for all }l\neq p\ ;\qquad (c)\ t\in K_p,\] 
for the maximal compact subgroup $K_v$ of $T_{\epsilon}^{\mathrm{o}}(\Q_v)$ (for each finite place $v$). Also, the K-pair $(\mathscr{I},\pi_0^k)$ is admissible for any $k\gg1$.
The pair $(\pi_0,t)$ is uniquely determined by $\epsilon$, up to taking simultaneous powers. 
\end{prop}

\begin{proof}
Any $\Q$-torus of $I_{\mathscr{I}}$ has the same ranks over $\Q$ and $\R$, by Theorem \ref{thm:Kisin17_Cor.1.4.13,Prop.2.1.3,Cor.2.2.5} (2) and our assumption on $Z(G)$. Under that condition,
the construction of $\pi_0, t\in T_{\epsilon}(\Q)$ satisfying the properties (a) - (c) of Lemma \ref{lem:canonical_decomp_of_epsilon} uses only the $\Q$-group structure of the $\Q$-group $T_{\epsilon}$ (of multiplicative type) and $\epsilon$. 
To see the second statement, we choose a maximal $\Q$-torus $T$ of $I_{\mathscr{I}}$ containing $\epsilon$, and the data $i_T:T\hookrightarrow G$, $h\in X\cap \Hom(\dS,i_T(T)_{\R})$, $(b\in i_T(T)(\mfk),i_p,i^p)$ in Theorem \ref{thm:Kisin17_Cor.1.4.13,Prop.2.1.3,Cor.2.2.5}; thus, $T_{\epsilon}\subset T$.
In this situation, we establish the statement, assuming that $t\in T(\Qp)_1=T(\Qp)\cap\ker(w_{T})$, which is allowed, as $T(\Qp)_1$ is a subgroup with finite index of the maximal compact subgroup of $T(\Qp)$.
Since $(\mathscr{I},\epsilon)$ is admissible, there exists $c\in G(\mfk)$ such that $c(\epsilon^{-1}\cdot(b\sigma)^n)c^{-1}=\sigma^n$, where $n$ is the level of $(\phi,\epsilon)$. Since $t$ lies in the parahoric subgroup $T(\Qp)_1$ of $T(\Qp)$, we can find $t_p\in T(\mfk)$ such that $t=t_p^{-1}\sigma^{ns}(t_p)$ \cite[Prop.3]{Greenberg63}. So, we have $(ct_p)(\pi_0^{-1}\cdot(b\sigma)^{ns})(ct_p)^{-1}=\sigma^n$, as was asserted. 
\end{proof}

For a group $A$, we consider the following equivalence relation $\sim$ among the set of pairs $(a,n)$ with $a\in A$, $n\in\N$: $(a,n)\sim (a',n')$ if there exists $N\in\N$ such that $a^{n'N}=(a')^{nN}$. We define a \emph{germ of an element} of $A$ to be an equivalence class of pairs $(a,n)$ for this equivalence relation. For an algebraic group $G$ over a field $k$ (say, of characteristic zero) and a germ $\pi$ of an element of $G(k)$, if $(\pi_n,n)$ is a representative of $\pi$, the Zariski closures in $G$ of the subgroups generated by $\pi_n^e\ (e\in\N)$ form a decreasing sequence (with $e$'s being ordered multiplicatively) of subgroups, thus stabilizes to a $k$-subgroup. It is then easy to see that this $k$-group depends only on the given germ $\pi$, not on the choice of representative $(\pi_n,n)$; we call this subgroup of $G$ the \emph{subgroup generated by} the germ $\pi$. Also we call the centralizer in $G$ of this subgroup the \emph{centralizer} (in $G$) of the germ $\pi$ and denote it by $G_{\pi}$. Note that the subgroup generated by a germ of an element and thus its centralizer are always \emph{connected}.

%%%%%%%%%%%%%%%%%%%%
%%%%%%%%%%%%%%%%%%%%
\begin{thm} \label{thm:Kisin17_Cor.2.3.2} 
Assume that $\Q\operatorname{-}\mathrm{rk}(Z(G))=\R\operatorname{-}\mathrm{rk}(Z(G))$. Let $\mathscr{I}\subset\sS_{\mathbf{K}_p}(\Fpb)$ be an isogeny class.

There exists a unique germ $\pi_{\mathscr{I}}$ of an element in $Z(I_{\mathscr{I}})(\Q)$, called the \emph{germ of Frobenius endomorphism} of $\mathscr{I}$, with the following properties:

(a) There exist a representative $i_p:(I_{\mathscr{I}})_{\Qp}\hookrightarrow J_b$ of the $G(\mfk)$-conjugacy class of Theorem \ref{thm:Kisin17_Cor.1.4.13,Prop.2.1.3,Cor.2.2.5} and a representative $(\pi_N,N)$ of $\pi=\pi_{\mathscr{I}}$ such that one has $b\in G(L_N)$ and $i_p(\pi_N)=\Nm_Nb$;

(b) The embeddings $i_l:(I_{\mathscr{I}})_{\Ql} \hookrightarrow G_{\Ql}\ (l\neq p)$, $i_p:(I_{\mathscr{I}})_{\Qp}\hookrightarrow J_b$ in Theorem \ref{thm:Kisin17_Cor.1.4.13,Prop.2.1.3,Cor.2.2.5} induce isomorphisms of group schemes
\[ i_l:(I_{\mathscr{I}})_{\Ql} \isom I_{i_l(\pi)},\quad i_p:(I_{\mathscr{I}})_{\Qp}\isom I_{i_p(\pi)},\] 
where $I_{i_l(\pi)}$ and $I_{i_p(\pi)}$ denote respectively the centralizer in $G_{\Ql}$ of the germ $i_l(\pi)$ and the centralizer in $J_b$ of the germ $i_p(\pi)$. Also for almost all $l\neq p$, the $\Z_l$-embedding $i_l:(I_{\mathscr{I}})_{\Z_l} \hookrightarrow G_{\Z_l}$ of Theorem \ref{thm:Kisin17_Cor.1.4.13,Prop.2.1.3,Cor.2.2.5} induces an isomorphism  $(I_{\mathscr{I}})_{\Z_l} \isom (G_{\Z_l})_{i_l(\pi)}$.

(c) For any representative $(\pi_n,n)$ of $\pi_{\mathscr{I}}$, the K-pair $(\mathscr{I},\pi_n^k)$ is admissible for all $k\gg1$ and $\pi_n$ lies in a compact subgroup of $I_{\mathscr{I}}(\Ql)$ for every finite place $l\neq p$.
For any admissible K-pair $(\mathscr{I},\epsilon)$ of level $n$ and every triple $(\pi_0,t\in T_{\epsilon}(\Q);s\in\N)$ attached to $\epsilon$ as in Proposition \ref{prop:canonical_decomp_of_epsilon2}, $(\pi_0,ns)$ represents the germ $\pi_{\mathscr{I}}$.
\end{thm}

\begin{proof}
This is basically \cite[Cor.2.3.2]{Kisin17}.
There is a natural candidate for $\pi_n\ (n\gg1)$, i.e. the $p^n$-th relative Frobenius endomorphism in $Z(I_{\mathscr{I}})(\Q)$ of the isogeny class $\mathscr{I}$. The uniqueness (as a germ of element in $Z(I_{\mathscr{I}})(\Q)$) will follow from (c).

Property (b) is Theorem \ref{thm:Kisin17_Cor.2.3.2;Tate_isom}, i.e. Kisin's generalization of the Tate's theorem on endomorphisms of abelian varieties over finite fields.
We choose a point $x$ in $\mathscr{I}$ defined over a finite field $\F_q$ and via some $K_0:=W(\F_q)[1/p]$-linear isomorphism $V\otimes K_0 \isom H^{\cris}_1(\mathcal{A}_x/K_0)$ matching the tensors $s_{\alpha}$, $s_{\alpha,0,x}$ on both sides, we identify the Frobenius automorphism on $H^{\cris}_1(\mathcal{A}_x/K_0)$ with $\delta\sigma$ for $\delta\in G(K_0)$. From this choice, we obtain a datum $(b=\delta,i_p:(I_{\mathscr{I}})_{\Qp}\hookrightarrow J_b)$ in Theorem \ref{thm:Kisin17_Cor.1.4.13,Prop.2.1.3,Cor.2.2.5}; note that in this case, we have 
\begin{equation} \label{eq:(delta,i_p)}
i_p(\pi_N)=\Nm_N\delta
\end{equation}
for any $N\in N$ with $\F_q\subset \F_{p^N}$; this establishes (a).
Now, property (b) for $l\neq p$ is immediate from Theorem \ref{thm:Kisin17_Cor.2.3.2;Tate_isom}. In the $p$-adic case, \textit{loc. cit.} says that for any $N\gg1$, $i_p$ identifies $(I_{\mathscr{I}})_{\Qp}$ with the $\sigma$-centralizer $G_{\delta\sigma}(\subset \Res_{L_N/\Qp}(G))$ inside $J_{\delta}$. But, since $\Nm_N\delta\rtimes\sigma^N=(\delta\rtimes\sigma)^N$, for any $\Qp$-algebra $R$, an element $g\in G(\mfk\otimes_{\Qp} R)$ commutes simultaneously with $\delta\sigma$ and $\sigma^N$ if and only if it does so with $\delta\sigma$ and $\Nm_N\delta=i_p(\pi_N)$, thus one has
\begin{align*}
G_{\delta\sigma}(R) &=\{ g\in G(L_N\otimes_{\Qp} R)\ |\ \delta\sigma(g)=g\delta \} \\
&=\{g\in G(\mfk\otimes_{\Qp} R)\ |\ \sigma^Ng=g,\ \delta\sigma(g)=g\delta \} \\
&=J_{\delta,i_p(\pi_N)}(R),
\end{align*}
as was asserted.
The second statement on the extension of $i_l$ over $\Z_l$ is clear, since $\pi_N$ is a relative Frobenius endomorphism so that for almost all $l\neq p$, $i_l(\pi_N)\in G_{\Z_l}(\Z_l)$ for all $N\gg1$.

Next, for property (c), the admissibility of $(\mathscr{I},\pi_N)$ (for $N\gg1$) follows from (a), since for the choice of $(b,i_p)$ there, one has $i_p(\pi_N)^{-1}(b\sigma)^N=\sigma^N$. Also, being a relative Frobenius endomorphism,  for every $l\neq p$, the eigenvalues of $\pi_N$ are all $l$-adic units, which implies that the subgroup of $Z(I_{\mathscr{I}})(\Ql)$ generated by $\pi_N$ is bounded (thus, its closure is compact). The last claim of (c) will be proved in Proposition \ref{prop:phi(delta)=gamma_0_up_to_center} below.
\end{proof}

%%%%%%%%%%%%%%%%%%%%
\begin{prop} \label{prop:phi(delta)=gamma_0_up_to_center}
Let $(\mathscr{I},\epsilon)$ be an admissible K-pair, say of level $n$. Suppose that the image of $\epsilon$ in $(I_{\mathscr{I}}/Z(G))(\Q)$ lies in a compact subgroup of $(I_{\mathscr{I}}/Z(G))(\A_f^p)$.
 
(1) For any sufficiently large $k\in\N$ divisible by $n$, if $(\pi_k,k)$ represents the germ $\pi_{\mathscr{I}}$ of Frobenius endomorphism of $\mathscr{I}$, the element $\pi_k\cdot \epsilon^{-\frac{k}{n}}$ of $I_{\mathscr{I}}(\Q)$ lies in the subgroup $Z(G)(\Q)$. 

(2) Assume one of the following two conditions: 

(a)  the weight homomorphism $w_X=(\mu_h\cdot\iota(\mu_h))^{-1}\ (h\in X)$ is rational, and $\epsilon$ is a Weil $q=p^n$-element of weight $-w=-w_X$, in the sense that for the $\Q$-subgroup $S\subset I_{\mathscr{I}}$ (of multiplicative type) generated by $\epsilon$ and any character $\chi$ of $S$, $\chi(\epsilon)\in\Qb$ is a Weil $q=p^n$-number of weight $-\langle \chi,w_X\rangle\in\Z$ in the usual sense; or,

(b)  the anisotropic kernel of $Z(G)$ remains anisotropic over $\R$ and $\epsilon\in I_{\mathscr{I}}(\A_f^p)$ itself lies in a compact subgroup of $I_{\mathscr{I}}(\A_f^p)$.

Then one has $\epsilon^{\frac{k}{n}}=\pi_k$ for any sufficiently large $k\in\N$ divisible by $n$.
\end{prop}

Obviously, the assumption is equivalent to that for a triple $(b,i_p,i^p)$ as in Theorem \ref{thm:Kisin17_Cor.1.4.13,Prop.2.1.3,Cor.2.2.5} (1), the image of $i^p(\epsilon)$ in $G^{\ad}(\A_f^p)$ lies in a compact subgroup. Similarly, the condition in (b) that $\epsilon\in I_{\mathscr{I}}(\A_f^p)$ itself lies in a compact subgroup of $I_{\mathscr{I}}(\A_f^p)$ is the same as that for any such $i^p$, $i^p(\epsilon)$ itself lies in a compact subgroup of $G(\A_f^p)$.

\begin{proof} 
For the proof, we first need a fact.

%%%%%%%%%%%%%%%%%%%%
\begin{lem} \label{lem:equality_of_two_Newton_maps}
(1) Let $\mathscr{I}$ be an isogeny class and $T\subset I_{\mathscr{I}}$ a maximal $\Q$-torus with $i_T:T\hookrightarrow G$ and $(b,i_p,i^p)$ an associated $\Q$-embedding and triple as in Theorem \ref{thm:Kisin17_Cor.1.4.13,Prop.2.1.3,Cor.2.2.5} (3) and satisfying the condition (i).
Then, for all $k\gg1$ the two Newton homomorphisms $\nu_{\pi_k}, \nu_b\in X_{\ast}(T)_{\Q}^{\Gamma_p}$ are related (via $i_T$) by: $\nu_{\pi_k}=k\nu_b$.

(2) Let $\gamma_0\in G(\Qp)$ be a semisimple element and suppose that there exists $\delta\in G(L_n)$ such that $\Nm_n\delta=c'\gamma_0 c'^{-1}$ for some $c'\in G(\mfk)$; let $I_0:=G_{\gamma_0}^{\mathrm{o}}$. If $b':=c'^{-1}\delta\sigma(c')\in G_{\gamma_0}(\mfk)$ belongs to $I_0$, the two Newton homomorphisms $\nu_{\gamma_0}$, $\nu_{b'}\in \Hom_{\mfk}(\mathbb{D},I_0)$ are related by: $\nu_{\gamma_0}=n\nu_{b'}$. 
In particular, in this case $[b']_{I_0}\in B(I_0)$ is basic \cite[(5.1)]{Kottwitz85}. If $b'$ lies in a maximal $\Qp$-torus $T_1$ of $G_{\gamma_0}$, the same statement holds true with $I_0$ replaced by $T_1$.

(3) For every admissible K-pair $(\mathscr{I},\epsilon)$ and for any maximal $\Q$-torus $T$ of $I_{\mathscr{I}}$ with $\epsilon\in T(\Q)$, we have equality of quasi-cocharacters of $T$: $\frac{1}{k}\nu_{\pi_k}=\frac{1}{n}\nu_{\epsilon}\ (k\gg1)$.
\end{lem}

\begin{proof}
(1) This follows from (\ref{eq:(delta,i_p)}) (cf. \cite[4.3]{Kottwitz85}).

(2) This is proved in Lemma 5.15 of \cite{LR87} (under the assumption $I_0=G_{\gamma_0}$). We briefly sketch its arguments. 

First, we observe that for any $c'\in G(\mfk)$ as in the statement and $n'\in \N$ divisible by $n$, $c_{n'}:=c'^{-1}\sigma^{n'}(c')$ belongs to $G_{\gamma_0}(\mfk)$ and lies in any small neighborhood of $1$ in $I_0(\mfk)$ as $n'$ becomes large (in fact, even becomes $1$ if $c'\in G(\Qpnr)$). 

Secondly, if we choose a maximal $\Qp$-torus $T_1$ of $G_{\Qp}$ containing $\gamma_0$, then the Newton quasi-cocharacter $\nu_{\gamma_0}$ of $\gamma_0\in T_1(\Qp)$ satisfies \cite[4.4]{Kottwitz85} that for every $\Qp$-rational $\lambda\in X_{\ast}(T_1)$, one has
\[|\lambda(\gamma_0)|_p=p^{-\langle\lambda,\nu_{\gamma_0}\rangle}.\]
It follows from this equation that $\nu_1:=\frac{1}{n}\nu_{\gamma_0}\in X_{\ast}(T_1)_{\Q}^{\Gamma_p}$ maps into the (connected) center of $G_{\gamma_0}$ and that $p^{-n\nu_1}\gamma_0\in T_1(\Qp)$ lies in the maximal compact subgroup of $T_1(\Qp)$. Especially, $(p^{-n\nu_1}\gamma_0)^k$ also lies in any small neighborhood of $1$ in $T_1$ if $k$ becomes large. Therefore, according to \cite[Prop.3]{Greenberg63}, for sufficiently large $k\in\N$, there exists $d\in I_0(\mfk)$ such that with $n'=nk$, 
\[\gamma_0^k c_{n'}=p^{n'\nu_1}d^{-1}\sigma^{n'}(d).\]

Finally, from $\Nm_n\delta=c'\gamma_0 c'^{-1}$, one easily checks that
\begin{equation} \label{eq:geom_conj_at_p}
\Nm_{n'}b'=\gamma_0^k c_{n'}(=p^{n'\nu_1}d^{-1}\sigma^{n'}(d)).
\end{equation}
It follows from this equality and the definition \cite[4.3]{Kottwitz85} that when $b'\in I_0(\mfk)$, $\nu_1\in\Hom_{\mfk}(\mathbb{D},I_0)$ is the Newton homomorphism of $b'$. 
If $b'\in T_1(\mfk)$, then $c_{n'}\in T_1(\mfk)$ and one can find $d$ in $T_1(\mfk)$, so the same conclusion holds true of $T_1$ instead of $I_0$.

(3) Let $i_T:T\hookrightarrow G$ and $(b,i_p,i^p)$ be an associated $\Q$-embedding and triple as in Theorem \ref{thm:Kisin17_Cor.1.4.13,Prop.2.1.3,Cor.2.2.5} (3) and satisfying the condition (i). 
Then, for $\gamma_0:=i_T(\epsilon)$, if we choose $c\in G(\mfk)$ such that $c(\gamma_0^{-1}(b\sigma)^n)c^{-1}=\sigma^n$, we have $\delta:=cb\sigma(c)^{-1}\in G(L_n)$ and $\Nm_n\delta=c\gamma_0c^{-1}$, so by choosing $c'$ to be $c$ in (2), we have $b'=b\in i_T(T)(\mfk)\subset G_{\gamma_0}^{\mathrm{o}}(\mfk)$.
Hence, the claim follows from (1) and (2).
\end{proof}

%%%%%%%%%%%%%%%%%%%%
\textsc{Proof of Proposition \ref{prop:phi(delta)=gamma_0_up_to_center} continued.} 
By assumption, $\epsilon\in T(\Q)$ lies in the maximal compact subgroup of $T/Z(G)(\A_f^p)$.

(1) As $\pi_{kd}=\pi_k^{d}$, it is enough to show that for some $k\in\N$ (divisible by $n$), the image of $\pi_k\cdot \epsilon^{-\frac{k}{n}}$ in $T/Z(G)(\Q)$ is a torsion element. For that, we use the fact that for any linear algebraic group $G$ over a number field $F$, $G(F)$ is discrete in $G(\A_F)$, so for any compact subgroup $K\subset G(\A_F)$, $G(F)\cap K$ will be finite, particularly, a torsion group. 
We will check that for every place $v$ of $\Q$, the image of $\pi_k\cdot \epsilon^{-\frac{k}{n}}$ in $T/Z(G)(\Qv)$ lies in the maximal compact (open) subgroup of $T/Z(G)(\Qv)$. Recall that for an $F$-torus $T$ and any finite place $v$ of $F$, the maximal compact subgroup $H$ of $T(F_v)$ equals 
\[\bigcap_{\chi\in X^{\ast}(T),\ F_v-\text{rational}}\ker(\mathrm{val}_v\circ\chi),\]
where $\mathrm{val}_v$ is the (normalized) valuation on $F_v$. 
For every finite place $l\neq p$, the image of $\epsilon$ in $T/Z(G)(\Ql)$ is a unit (i.e. lies in a compact subgroup) by assumption, and so is $\pi_k$ by definition of $\pi_k$ (in fact, $\pi_k$ is itself a unit in $T(\Ql)$ for every $l\neq p$). As $T/Z(G)$ is anisotropic over $\R$, the claim is trivial for the archimedean place. Hence, it suffices to show that for every $\Qp$-rational character $\chi$ of $T/Z(G)$, $|\chi(\pi_k)|_p=|\chi(\epsilon^{\frac{k}{n}})|_p$. In fact, we will show this for $\Qp$-rational characters $\chi$ of $T$.
Indeed, for all sufficiently large $k\in\N$ divisible by $n$ and for any $\Qp$-rational character $\chi$ of $T$, one has 
\[ |\chi(\pi_k)|_p=p^{-\langle\chi,\nu_{\pi_k}\rangle}=p^{-\frac{k}{n}\langle\chi,\nu_{\epsilon}\rangle}=|\chi(\epsilon)|_p^{\frac{k}{n}}. \]
due to Lemma \ref{lem:equality_of_two_Newton_maps} (3) (for the second equality). 

(2) In case (a), the additional assumption tells us that $|\chi(\pi_k)|_{\infty}=|\chi(\epsilon^{\frac{k}{n}})|_{\infty}$ for every $\Q_{\infty}$-rational character $\chi$ of $S$, and also implies that $\epsilon\in I_{\mathscr{I}}(\A_f^p)$ itself lies in a compact open subgroup of $I_{\mathscr{I}}(\A_f^p)$. In case (b), it is well-known that the stated condition implies that for any maximal $\Q$-torus $T_0$ of $G$ which is elliptic over $\R$, $T_0(\Q)$ is discrete in $T_0(\A_f)$; this is the condition which Kisin called \emph{the Serre condition for $T_0$}, \cite[(3.7.3)]{Kisin17}. 

In both cases, we can again repeat the argument of (1).
\end{proof}

%%%%%%%%%%%%%%%%%%%%
\begin{lem} \label{lem:Zariski_group_closure}
Let $G$ be a reductive group $G$ over a field $k$ and $\epsilon\in G(k)$ a semisimple element; let $S$ be the $k$-subgroup of $G$ generated by $\epsilon$.

(1) If $\epsilon':=g\epsilon g^{-1}\in G(k)$ for $g\in G(\bar{k})$, the $k$-subgroup of $G$ generated by $\epsilon'$ equals $\Int(g)(S)$.

(2) For any field extension $k'/k$, the $k'$-subgroup of $G_{k'}$ generated by $\epsilon$ equals $S_{k'}$,
and the elements $\{\epsilon^n\}_{n\in\N}$ are Zariski dense in $S_{k'}$.
\end{lem}

\begin{proof} 
(1) As $G_{\epsilon}=Z_G(S)$, the map $\Int(g):S\rightarrow \Int(g)(S)$ is a $k$-isomorphism of $k$-groups, and thus $\Int(g)(S)$ is the $k$-subgroup of $G$ generated by $\epsilon'$.
(2) This follows from the following easy fact:
if one takes a $k$-torus $T\subset G$ containing $\epsilon$, the $k$-subgroup $S$ (of multiplicative type) of $T$ generated by $\epsilon$ is the kernel of the surjective homomorphism $T\rightarrow T'$ of $k$-tori, where $T'$ is defined by $X^{\ast}(T')=\{\chi\in X^{\ast}(T)\ |\ \chi(\epsilon)=1\}$.
\end{proof}

%%%%%%%%%%%%%%%%%%%%
\begin{cor} \label{cor:Tate_thm2}
Assume that $\Q\operatorname{-}\mathrm{rk}(Z(G))=\R\operatorname{-}\mathrm{rk}(Z(G))$.
Let $(\mathscr{I},\epsilon)$ be an admissible K-pair and $(\gamma_0;\gamma,\delta)$ the associated Kottwitz triple defined by (\ref{eq:K-triple_attached_to_admissible_K-pair}): $\gamma_l=i_l(\epsilon)$ and $\delta=cb\sigma(c)^{-1}$ for a triple $(b,i_p,i^p)$ in Theorem \ref{thm:Kisin17_Cor.1.4.13,Prop.2.1.3,Cor.2.2.5} and an element $c\in G(\mfk)$ satisfying (\ref{eq:(epsilon,b,c)->delta2}).

(1) The embedding $i^p=\prod{}'_{l\neq p,\infty}i_l:(I_{\mathscr{I}})_{\A_f^p} \hookrightarrow G_{\A_f^p}$ induces an isomorphism of $\A_f^p$-group scheme
\begin{equation} \label{eq:isom_i^p}
i^p\ :\ (I_{\mathscr{I},\epsilon})_{\A_f^p} \isom G_{\gamma},
\end{equation} 
where $G_{\gamma}$ is the centralizer of $\gamma$ in $G_{\A_f^p}$. 

(2) The embedding $\Int(c)\circ i_p:(I_{\mathscr{I}})_{\Qp}\hookrightarrow J_b\isom J_{\delta}$ induces an isomorphism of $\Qp$-groups
\begin{equation} \label{eq:isom_Int(c)circi_p}
\Int(c)\circ i_p\ :\ (I_{\mathscr{I},\epsilon})_{\Qp}\isom G_{\delta\sigma},
\end{equation}
where $G_{\delta\sigma}(\subset \Res_{L_n/\Qp}(G))$ is the $\sigma$-centralizer of $\delta\in G(L_n)$ (\ref{eq:G_{xtheta}}). 
The functorial composite map $\Xi\circ i_p(=\Xi\circ \Int(c^{-1})\circ \Int(c)\circ i_p) :I_{\mathscr{I}}(R)\hookrightarrow J_b(R)\subset G(\mfk\otimes R)\isom G(R)$ for $\mfk$-algebras $R$ induces a $\mfk$-isomorphism
\begin{equation} \label{eq:Xicirci_p}
\Xi \circ i_p\ :\ (I_{\mathscr{I},\epsilon}^{\mathrm{o}})_{\mfk}\isom (I_0)_{\mfk},
\end{equation}
where $\Xi:G(\mfk\otimes R)\isom G(R)$ is the map (\ref{eq:Xi2}) induced by the $\mfk$-algebra homomorphism $\mfk\otimes R\rightarrow R:l\otimes x\mapsto lx$. 
Some $I_0(\mfk)$-conjugate is an inner twisting defined over $\Qpnr$ with class $[b^{\ad}]\in H^1(\Qpnr/\Qp,I_0)\subset B(I_0)$ (Proposition \ref{prop:psi_p} (3)). 

(3) Suppose that $(b,i_p,i^p)$ is adapted to a pair $(i_T:T\hookrightarrow G, h\in X\cap \Hom(\dS,i_T(T)_{\R}))$ as in Theorem \ref{thm:Kisin17_Cor.1.4.13,Prop.2.1.3,Cor.2.2.5} (3) ($T\subset I_{\mathscr{I},\epsilon}^{\mathrm{o}}$ being some maximal $\Q$-torus) and $\gamma_0=i_T(\epsilon)$.
Then, there exists an inner twisting 
\[\varphi\ :\ (I_{\mathscr{I},\epsilon}^{\mathrm{o}})_{\Qb}\isom (I_0)_{\Qb}\] 
which restricts to $i_T$ and 
such that for each finite place $v$ of $\Q$, $\varphi_{\Qvb}$ is \emph{conjugate} to $i_v$ in (\ref{eq:isom_i^p}) if $v\neq p,\infty$, to $\Xi\circ i_p$ in (\ref{eq:Xicirci_p}) if $v=p$.
\end{cor}

\begin{proof}
(1) By Theorem \ref{thm:Kisin17_Cor.2.3.2} and Proposition \ref{prop:canonical_decomp_of_epsilon2}, $i_l$ induces isomorphisms (for any $k\gg1$)
\[(I_{\mathscr{I},\epsilon})_{\Ql} \isom Z_{G_{\Ql}}(i_l(\pi_{nsk}),i_l(\epsilon)) \isom Z_{G_{\Ql}}(i_l(\pi_0^k),i_l(\epsilon)) \isom Z_{G_{\Ql}}(i_l(\epsilon)),\] 
where $(\pi_{nsk},nsk)$ is a representative of the germ $\pi_{\mathscr{I}}$ and $(s\in\N,\pi_0\in T_{\epsilon}(\Q))$ are as in Proposition \ref{prop:canonical_decomp_of_epsilon2} (the groups in the middle are the simultaneous centralizers of the elements inside the round brackets); $i_l(\pi_0^k)$ lies in the subgroup of $G_{\Ql}$ generated by $i_l(\epsilon)$, Lemma \ref{lem:Zariski_group_closure}.

(2) Again, as in the case $l\neq p$, it follows from Theorem \ref{thm:Kisin17_Cor.2.3.2} and Proposition \ref{prop:canonical_decomp_of_epsilon2} that $i_p$ induces $\Qp$-isomorphisms (for any $k\gg1$)
\begin{equation} \label{eq:i_{p,epsilon}}
(I_{\mathscr{I},\epsilon})_{\Qp} \isom Z_{J_b}(i_p(\pi_{nsk}),i_p(\epsilon)) \isom Z_{J_b}(i_p(\pi_0^k),i_p(\epsilon)) \isom Z_{J_b}(i_p(\epsilon)).
\end{equation}
($i_p(\pi_0^k)$ lies in the subgroup of $J_b$ generated by $i_p(\epsilon)$, Lemma \ref{lem:Zariski_group_closure}.)
Then, since $c(i_p(\epsilon)^{-1}(b\sigma)^n)c^{-1}=\sigma^n$ (\ref{eq:(epsilon,b,c)->delta2}), 
$\Int(c)$ induces a $\Qp$-isomorphism 
\begin{equation} \label{eq:Int(c)}
 \Int(c): Z_{J_b}(i_p(\epsilon)) \isom G_{\delta\sigma}
\end{equation}
(for any $\Qp$-algebra $R$, an element $g\in G(\mfk\otimes_{\Qp} R)$ commutes simultaneously with $b\sigma$ and $i_p(\epsilon)$ if and only if $\Int(c)(g)$ does so with $\delta\sigma=c(b\sigma)c^{-1}$ and $\sigma^n$, by the relation (\ref{eq:(epsilon,b,c)->delta2})). This proves the first claim.

The second and the last claims follows from the first one and the facts (Proposition \ref{prop:psi_p}) that $\Xi\circ \Int(c^{-1}):G(\mfk\otimes R)\isom  G(\mfk\otimes R)\isom G(R)$ induces an isomorphism $(G_{\delta\theta}^{\mathrm{o}})_{\mfk} \isom (J_b^{I_0})_{\mfk} \isom (I_0)_{\mfk}$ and that for any $h\in I_0(\mfk)$ as in Lemma \ref{lem:Xi_{b'}_inner-twisting} (1) (i.e. such that $\Nm_n(b')=\nu(p)$ holds for $b':=h^{-1}b\sigma(h)$ and $\nu:=\nu_{I_0}(b)$), $\Int(h^{-1})\circ \Xi\circ i_p$ is a  conjugate with the required property.

(3) The situation is that we are given two connected reductive $\Q$-groups $I_{\mathscr{I},\epsilon}^{\mathrm{o}}$, $I_0$ which share the same maximal $\Q$-torus $T\isom i_T(T)$. We will show that the two canonical homomorphisms $\Gal(\Qb/\Q)\rightarrow \Aut(G,T)/N_G(T)$ attached to these groups $G$ endowed with a maximal $\Q$-torus $T$ are the same. Since this homomorphism gives the canonical Galois action on the based root datum $(X^{\ast}(T),\Delta,X_{\ast}(T),\Delta^{\vee})$, where $\Delta\subset X^{\ast}(T)$ is a set of simple roots, this will imply the existence of an inner twisting $\varphi:(I_{\mathscr{I},\epsilon}^{\mathrm{o}})_{\Qb}\isom (G_{\gamma_0}^{\mathrm{o}})_{\Qb}$ which restricts to $i_T$ (\cite[2.4, 2.9]{Springer79}).
By Chebotarev density theorem, it is enough to check this locally for places in a set of Dirichlet density $1$. For a place $v$, the equality of the two restrictions $\Gal(\Qvb/\Qv)\rightarrow \Aut(G,T)/N_G(T)$ in turn will follow from the existence of similar (i.e. \emph{$T$-equivariant}) inner twistings \[ \varphi_v: (I_{\mathscr{I},\epsilon}^{\mathrm{o}})_{\Qvb}\isom (G_{\gamma_0}^{\mathrm{o}})_{\Qvb}.\] 
For finite $v$, the inner twistings provided by $i_v(v\neq p)$ (1) and $\Xi\circ i_p$ (2) are $i_T$-equivariant as the triple $(b,i_p,i^p)$ is adapted to $i_T$ (Theorem \ref{thm:Kisin17_Cor.1.4.13,Prop.2.1.3,Cor.2.2.5} (3) (i)): to see that $\Xi \circ i_p$ restricts to $i_T$, recall that for $\mfk$-algebras $R$, $i_T(T)(R)\subset G(R)$ was regarded as a subgroup of $J_b(R)=G(\mfk\otimes R)^{\langle \sigma\rangle}$ via the natural map $G(R)\subset G(\mfk\otimes R)$, so the restriction to $i_T(T)(R)$ of $\Xi: J_b(R)\subset G(\mfk\otimes R)\isom G(R)$ is the inclusion. Furthermore, the $T$-equivariances also imply that the two local inner twistings $\varphi_v$, $i_v$ ($\Xi\circ i_p$ when $v=p$) are conjugate since any automorphism of a connected reductive group fixing (pointwise) a maximal torus is inner. 
\end{proof}

Next, we propose a recipe of associating a \emph{stable} Kottwitz triple $(\gamma_0,\gamma,\delta)$, not just a Kottwitz triple, with any admissible K-pair $(\mathscr{I},\epsilon)$. The (only) problem with the recipe (\ref{eq:K-triple_attached_to_admissible_K-pair}) for that purpose is that we do not know yet whether the stable conjugacy class of $\gamma_0$ does not depend on the choice of the maximal $\Q$-torus $T\subset I_{\mathscr{I},\epsilon}^{\mathrm{o}}$ (Lemma \ref{lem:K-triple_attached_to_admissible_K-pair}).%%
\footnote{The Galois-gerb description of Langlands-Rapoport conjecture \cite{LR87} suggests that it should not be so, but at the moment we were not establish this.}
Our solution to this problem is to use some restricted class of maximal tori $T$ (which we call ``\emph{nice}''): see Definition \ref{defn:stable_K-triple_attached_to_admissible_K-pair}.

%%%%%%%%%%%%%%%%%%%%
\begin{defn} \label{defn:germ_of_Frobenius_elt}
Let $(T,h)$ be a special Shimura subdatum. Choose a finite CM, Galois extension $L\subset \Qb$ splitting $T$ and $x\in L$ such that the principal ideal $(\varpi)$ is some power $\mathfrak{P}^s$ of $\mathfrak{P}$, the prime ideal of $\mathcal{O}_L$ induced by the chosen embedding $\Qb\hookrightarrow \Qpb$.

The \emph{germ of Frobenius element} attached to $(T,h)$ is the germ of an element in $T(\Q)$ with a representative $(\Nm(\mu_h)(\varpi^{-1}),s[L_{\mathfrak{P},0}:\Qp])$, where $\Nm(\mu_h)$ is the ``reciprocity map''
\[ \Nm(\mu_h):\Res_{L/\Q}(\Gm) \stackrel{\Res_{L/\Q}(\mu_T)}{\ra} \Res_{L/\Q}T_L \stackrel{\Nm_{L/\Q}}{\ra} T, \]
and $L_{\mathfrak{P},0}$ is the maximal unramified subfield of $L_{\mathfrak{P}}$.
\end{defn}

We know (cf. \cite[4.3.9]{Kisin17}) that for any special Shimura subdatum $(T,h)$, if $\mathscr{I}$ is the isogeny class of the reduction of any associated CM point, the germ of Frobenius element attached to $(T,h)$ equals $i_T(\pi_{\mathscr{I}})$, where $\pi_{\mathscr{I}}$ is the germ of Frobenius endomorphism of $\mathscr{I}$ (Theorem \ref{thm:Kisin17_Cor.2.3.2}) and $i_T:T\hookrightarrow G$ is the inclusion, regarded as a $\Q$-embedding attached to $j_{T,h}:T\subset I_{\mathscr{I}}$ as in Theorem \ref{thm:Kisin17_Cor.1.4.13,Prop.2.1.3,Cor.2.2.5} (4).

We use the following fact which is deduced from a basic theory of Galois gerb presented in \cite{LR87} (cf. \cite{Lee18a}). 
As was mentioned before, this theorem is the only part in the current work which appeal to Galois gerb theory at all.

%%%%%%%%%%%%%%%%%%%%
%%%%%%%%%%%%%%%%%%%%
\begin{thm} \label{thm:LR_special_admissible_morphism}
For every special Shimura subdatum $(T,h)$ and a finite CM, Galois extension $L\subset \Qb$ splitting $T$, there exists a $T(\Qb)$-valued cochain $a^L=(a^L_{\tau})_{\tau}$ on $\Gamma=\Gal(\Qb/\Q)$ with the following properties:

(i) Its coboundary $(\partial a^L)_{\rho,\tau}:=a^L_{\rho}\rho(a^L_{\tau})(a^L_{\rho\tau})^{-1}\ (\rho,\tau\in\Gamma)$ lies in the subgroup generated by the Frobenius element $\pi=\Nm(\mu_h)(\varpi^{-1})\ (s\gg1)$ attached to $(T,h,L)$ (Definition \ref{defn:germ_of_Frobenius_elt}).

(ii) If $\mathscr{I}$ is the isogeny class of the reduction of an(y) associated CM point and if $G_{\pi} \subset G$ is the (connected) centralizer of $\pi$, the inner-twist $I_a$ of $G_{\pi}$ by $\bar{a}\in Z^1(\Q,G_{\pi}^{\ad})$ is $\Q$-isomorphic to $I_{\mathscr{I}}$.

(iii) For every prime $l\neq p$, there exist $t_l\in T(\Qlb)$ such that $a^L(l)_{\tau}=t_l^{-1}\tau(t_l)$.
There exists $t_p\in T(\Qpb)$ such that $a^L(p)^{\nr}_{\tau}:=t_pa^L_{\tau}\tau(t_p)^{-1}$ is unramified in the sense that $a^L(p)^{\nr}_{\tau}\in T(\Qpnr)$ for all $\tau\in\Gamma_p$ and $a^L(p)^{\nr}_{\tau}=1$ on $\Gal(\Qpb/\Qpnr)$, and that $w_T(a^L(p)^{\nr}_{\sigma})=-\underline{\mu_h}\in X_{\ast}(T)_{\Gal(\Qpb/\Qpnr)}$. 

(iv) For every $\R$-elliptic maximal $\Q$-torus $T_0$ of $G$ such that for some $g_0\in G(\Qb)$, $(a^L_0)_{\tau}:=g_0a^L_{\tau}\tau(g_0)^{-1}\in T_0(\Qb)$ for all $\tau\in\Gamma$, there exists $g_1\in G(\Qb)$ such that $\Int(g_1)$ induces a $\Q$-isomorphism of maximal tori over $\Q$, $T_0\isom T_1:=g_1T_0g_1^{-1}$ and the $T_1$-valued cochain $(a^L_1)_{\tau}:=g_1(a^L_0)_{\tau}\tau(g_1)^{-1}$ is the cochain attached to some special Shimura subdatum $(T_1,h_1)$.
\end{thm}

Note that by (i), the image $\bar{a}$ of the cochain $a$ in $G_{\pi}^{\ad}$ becomes a cocycle.

\begin{proof}
We review some rudiments of Galois gerb theory developed in Langlands-Rapoport \cite{LR87} (cf. \cite{Lee18a}), only a bare minimum necessary to deduce these results. For an algebraic $H$ over $\Q$, we denote by $\mathscr{G}_H$ the extension group $H(\Qb)\rtimes\Gal(\Qb/\Q)$ with the action of $\Gamma$ on $H(\Qb)$ given by the $\Q$-structure of $H$ (this is called the \textit{neutral Galois gerb} attached to $H$).

The basic ingredient of \textit{loc. cit.} is the so-called pseduo-motivic Galois gerb: 
for each CM field $L$, finite and Galois over $\Q$, there exist a maximal $\Q$-torus $P^L$ which splits over $L$ and a topological group $\mathfrak{P}^L$, an extension of $\Gamma=\Gal(\Qb/\Q)$ by $P^L(\Qb/\Q)$:
\[ 1\ra P^L(\Qb/\Q) \ra \mathfrak{P}^L \ra \Gal(\Qb/\Q) \ra 1\]
which are characterized by certain properties \cite[$\S$3]{LR87} (cf. \cite[3.2.7]{Lee18a}). Here the action of $\Gamma$ on $P^L(\Qb/\Q)$ is via the $\Q$-structure of $P^L$, and $\mathfrak{P}^L$ is defined by a continuous cocycle on $\Gamma$ with values in $P^L(\Qb)$.%%
%\footnote{In fact, $\mathfrak{P}^L$ is obtained via a pull-back and push-out from another extension $1\ra P^L(L/\Q) \ra \mathfrak{P}^L_L \ra \Gal(L/\Q)\ra 1$, \cite[$\&$2]{LR87}} 

Next, we construct the objects in this theorem. For each datum $(T,\mu)$ of a $\Q$-torus split over $L$ and a cocharacter $\mu\in X_{\ast}(T)$, there exists a (continuous) map of extensions compatible with the projections to $\Gamma$
\[ \psi_{T,\mu}:\mathfrak{P}^L \ra \mathscr{G}_T \]
whose restriction to $P^L(\Qb)$ is the map $P^L(\Qb)\ra T(\Qb)$ induced by a morphism of $\Q$-tori $P^L\ra T$ and such that $\psi_{T,\mu}(P^L(\Q))$ is generated by the germ of Frobenius element attached to $(T,h)$ (Definition \ref{defn:germ_of_Frobenius_elt}) (see \cite{LR87}, Satz. 2.3 and the discussion on p.143-144, cf. \cite[(3.1)]{Kisin17}, \cite[3.3]{Lee18a}).

We choose a continuous section $q_{\rho}$ to the projection $\mathfrak{P}^L\ra \Gamma$, and define $a^L_{\rho}\in T(\Qb)$ by 
\[ \psi_{T,\mu_h}(q_{\rho})=a^L_{\rho}\rtimes\rho. \] 
Then, (i) follows since $(\partial a^L)_{\rho,\tau}=\psi_{T,\mu_h}(q_{\rho}q_{\tau}q_{\rho\tau}^{-1})\in \psi_{T,\mu}(P^L(\Q))$.
(iii) follows from the construction of $\psi_{T,\mu_h}$ (\cite[Satz.2.3]{LR87} and \cite[3.3.9]{Lee18a}), and (ii) is a consequence of (iii), in view of Theorem \ref{thm:Kisin17_Cor.2.3.2} and Hasse principle for adjoint groups over number fields.
(iv) follows from the part in \cite{LR87} beginning from Lemma 5.11 until the end of the proof of Satz 5.3: a slight improvement of this argument (and in English) is given in \cite[Prop.4.1.5]{Lee18a}. The discussion in \textit{loc. cit.} concerns certain maps $\mathfrak{P}\ra \mathscr{G}_G$ called ``admissible morphisms'', of which the above $\psi_{T,\mu_h}$ composed with the natural map $\mathscr{G}_T\ra  \mathscr{G}_G$ is a special example (called \textit{special admissible morphism}).
\end{proof}

For a $p$-adic or real field $F$ and a connected reductive group $G$ over $F$, a maximal $F$-torus $T$ of $G$ is said to be \emph{fundamental} if its $F$-rank is minimal \cite[$\S$10]{Kottwitz86}, \cite[5.3.1]{Borovoi98}; in the $p$-adic case, fundamental tori are elliptic (since every semisimple group over a $p$-adic field always contains anisotropic maximal tori). We also note that for any admissible K-pair $(\mathscr{I},\epsilon)$, fundamental tori in $(I_{\mathscr{I},\epsilon}^{\mathrm{o}})_{\R}$ are also elliptic by Theorem \ref{thm:Kisin17_Cor.1.4.13,Prop.2.1.3,Cor.2.2.5} (2).

%%%%%%%%%%%%%%%%%%%%
\begin{defn} \label{defn:nice_tori}
For an admissible K-pair $(\mathscr{I},\epsilon)$, a maximal $\Q$-torus $T$ of $I_{\mathscr{I},\epsilon}^{\mathrm{o}}$ is said to be ``\emph{nice}'' if $T_{\Qv}$ is fundamental in $(I_{\mathscr{I}})_{\Qv}$ at $v=p$ and $\infty$, plus at every place $v$ of $\Q$ where $(I_{\mathscr{I},\epsilon}^{\mathrm{o}})_{\Qv}$ is not quasi-split. 
%either $(I_{\mathscr{I},\epsilon}^{\mathrm{o}})_{\Qv}$ or $(I_{\mathscr{I}})_{\Qv}$ is not quasi-split. 
\end{defn}

%%%%%%%%%%%%%%%%%%%%
\begin{thm} \label{thm:stable_isogeny_diagram}
Let $(\mathscr{I},\epsilon)$ be an admissible K-pair and $T$, $T'$ two nice maximal $\Q$-tori of $I_{\mathscr{I}}$. Let $(i_T:T\hookrightarrow G, h_0\in \Hom(\dS,i_T(T)_{\R})\cap X)$ be a $\Q$-embedding and special Shimura subdatum as in Theorem \ref{thm:Kisin17_Cor.1.4.13,Prop.2.1.3,Cor.2.2.5} (3). Choose a representative $(\pi_N,N)$ of the germ $\pi_{\mathscr{I}}$ of Frobenius endomorphism of $\mathscr{I}$ with $N\in\N$ big enough so that $\pi_N$ lies in the subgroup generated by $\epsilon$, the centralizer $G_{\pi_0}$ of $\pi_0:=i_T(\pi_N)\in i_T(T)(\Q)$ is connected and the new K-pair $(\mathscr{I},\pi_N)$ is also admissible (Theorem \ref{thm:Kisin17_Cor.2.3.2}). 

Then, there exists a commutative diagram
\begin{equation} \label{eq:stable_isogeny_diagram}
\xymatrix{ I_{\mathscr{I}} \ar[r]^{\sim}_{\iota} \ar[d]_{\varphi} & I_{a} & I_{a,\gamma_0}^{\mathrm{o}}  \ar@{_(->}[l] \ar[r]^{\sim}_{\iota'} & I_{a',\gamma_0'}^{\mathrm{o}} \\
G_{\pi_0} \ar[r]^{\mathrm{Id}}  & G_{\pi_0} \ar[u]_{\psi}  & G_{\gamma_0}^{\mathrm{o}} \ar@{_(->}[l] \ar[r]^{\Int(g)} \ar[u]_{\psi}   & G_{\gamma_0'}^{\mathrm{o}} \ar[u]_{\psi'} }
\end{equation}
with the following properties:

(i) All the maps are $\Qb$-isomorphisms except for the middle horizontal ones (which are inclusions). The two horizontal isomorphisms $\iota$, $\iota'$ are $\Q$-rational and the vertical maps are inner twistings;

(ii) The inner twisting $\varphi$ is an inner twisting attached, as in Corollary \ref{cor:Tate_thm2} (3),  to the admissible pair $(\mathscr{I},\pi_N)$ and the $\Q$-embedding $i_T$ (so, $\varphi|_T=i_T$);

(iii) The two vertical maps denoted by the same $\psi$ are the inner twistings defined by a cochain $a_{\tau}=g_0(a_0)_{\tau}\tau(g_0)^{-1}$ ($g_0\in G_{\pi_0}(\Qb)$), where $(a_0)_{\tau}$ is the cochain attached, as in Theorem \ref{thm:LR_special_admissible_morphism}, to the special Shimura subdatum $(i_T(T),h_0)$;
% (so, $\psi|_{T_1}$ is $\Q$-rational);

(iv) The map $\Int(g)\circ \psi^{-1}\ (g\in G(\Qb))$ induces a $\Q$-isomorphism of maximal $\Q$-tori $\iota(T')\isom T'_1$, and $\gamma_0':=\Int(g)(\gamma_0)$.
The vertical map $\psi'$ is the inner twisting defined by the cochain $a'_{\tau}=ga_{\tau}\tau(g)^{-1}$. This cochain is attached, as in Theorem \ref{thm:LR_special_admissible_morphism}, to some special Shimura subdatum $(T'_1,h'_1)$.
% (so, $\psi_1|_{T'_1}$ is $\Q$-rational);
\end{thm}

Some remarks are in order.
Here, we fixed a finite CM, Galois extension $L\subset \Qb$ big enough to split $T$ and $T'$, and omit $L$ in the notations for the cochains $a_0$, $a'$ attached to the special Shimura subdata $(i_T(T),h_0)$, $(T'_1,h'_1)$. We regard $I_a$ as $(G_{\pi_0})_{\Qb}$ with the $\Q$-structure given by $x\mapsto a_{\tau}\tau(x)a_{\tau}^{-1}$, so the map on $\Qb$-points underlying the inner twisting $\psi:(G_{\pi_0})_{\Qb}\isom (I_a)_{\Qb}$ is the identity. The same remark applies to $\psi'$ too. Also, $I_{a,\gamma_0}^{\mathrm{o}}$ denotes the neutral component of the centralizer $I_{a,\gamma_0}$ of $\psi(\gamma_0)$ in $I_{a}$. 
We note that as $\pi_N$ lies in the subgroup generated by $\epsilon$, the image of of the cochain $a$ in $(G_{\gamma_0}^{\mathrm{o}})^{\ad}$ becomes also a cocycle.

\begin{proof}
\textbf{Step 1.} Construction of the leftmost square diagram. This consists of three subsquares:
\begin{equation} \label{eq:stable_isogeny_diagram1}
\xymatrix{ I_{\mathscr{I}} \ar[r]^{\sim}_{\iota_0} \ar[d]_{\varphi} & I_{a_0} \ar[r]^{\sim}_{\Int(g_0)} & I_{a}  \\
G_{\pi_0} \ar[r]^{\Int(g_0^{-1})}  & G_{\pi_0} \ar[u]_{\psi_0} \ar[r]^{\Int(g_0)} & G_{\pi_0} \ar[u]_{\psi} }
\end{equation}

Put $T_0:=i_T(T)$. Let $(T_0,h_0)$ be a special Shimura subdatum provided by Theorem \ref{thm:Kisin17_Cor.1.4.13,Prop.2.1.3,Cor.2.2.5} (3) and denote by $\psi_0:G_{\pi_0}\ra I_{a_0}$ the inner twisting defined by the cochain $a_0$ attached to $(T_0,h_0)$ as in Theorem \ref{thm:LR_special_admissible_morphism}.
By Corollary \ref{cor:Tate_thm2}, Theorem \ref{thm:LR_special_admissible_morphism} and Hasse principle for adjoint groups, the two inner twistings $\varphi^{-1}:G_{\pi_0}\isom I_{\mathscr{I}}$, $\psi_0:G_{\pi_0}\ra I_{a_0}$ are defined by the same class in $H^1(\Q,G_{\pi_0}^{\ad})$ (in fact, they belong to $H^1(\Q,T_0/Z(G_{\pi_0}))$). Hence, there exist a $\Q$-isomorphism $\iota_0: I_{\mathscr{I}}\isom I_{a_0}$ and $g_0\in G_{\pi_0}(\Qb)$ such that $\iota_0=\psi_0\circ\Int(g_0^{-1})\circ\varphi$. If we define a new cochain $a_{\tau}:=g_0(a_0)_{\tau}\tau(g_0)^{-1}\in G_{\pi_0}(\Qb)$ on $\Gamma$ and $\psi:G_{\pi_0}\isom I_{a}$ is the inner twisting by $\bar{a}\in Z^1(\Q,G_{\pi_0}^{\ad})$, the map $\Int(g_0):(I_{a_0})_{\Qb}\isom (I_{a})_{\Qb}$ induces a $\Q$-isomorphism $I_{a_0}\isom I_{a}$ which also equals $\psi\circ\Int(g_0)\circ\psi_0^{-1}$. This completes the construction of the first left two subsquares, with $\iota:=\Int(g_0)\circ\iota_0$, . 

We note that as $\Int(g_0)\circ \iota_0$  and $\varphi|_T$ are both $\Q$-rational, so is $\psi|_{T_0}$. This is the same condition as that 
\[ a_{\tau}\in T_0(\Qb) \] for all $\tau\in \Gamma$. Indeed, the first condition holds if and only if for every $x\in T_0(\Qb)$ and $\tau\in\Gamma$, $\tau(\psi(x))=\psi(\tau(x))$, i.e. $\psi^{-1}\circ{}^{\tau}\psi(\tau(x))=\tau(x).$ As $\psi^{-1}\circ{}^{\tau}\psi=\Int(a_{\tau})$ and $Z_G(T_0)=T_0$, the claim follows. 

Note that $G_{\gamma_0}\subset G_{\pi_0}$ since $\pi_N$ lies in the subgroup generated by $\gamma_0$.

\textbf{Step 2.} Construction of the rightmost square diagram. This consists of two subsquares:
\begin{equation} \label{eq:stable_isogeny_diagram2}
\xymatrix{ I_{a,\gamma_0}^{\mathrm{o}} \ar[r]^{\sim}_{\Int(g_1)} & I_{a'_0,\gamma_0}^{\mathrm{o}} \ar[r]^{\sim}_{\Int(g_2)} & I_{a',\gamma_0'}^{\mathrm{o}} \\
G_{\gamma_0}^{\mathrm{o}} \ar[r]^{\Int(g_1)} \ar[u]_{\psi}  & G_{\gamma_0}^{\mathrm{o}} \ar[r]^{\Int(g_2)} \ar[u]_{\psi_0'} & G_{\gamma_0'}^{\mathrm{o}} \ar[u]_{\psi'} }
\end{equation}
Here, $\gamma_0=\iota(\epsilon)\in T_0(\Q)$ which also equals $\iota(\epsilon)\in \iota(T')(\Q)$ via $\psi:T_0\isom \iota(T)$.

\textbf{Step 2.a.} We claim that there exists $g_1\in G_{\gamma_0}^{\mathrm{o}}(\Qb)$ such that $\Int(g_1)\circ\psi^{-1}:I_{a,\gamma_0}^{\mathrm{o}} \isom G_{\gamma_0}^{\mathrm{o}}$ becomes $\Q$-rational when restricted to $\iota(T')$. 
Since $\psi$ is an inner twisting and $\iota(T')_{\Qv}$ is elliptic in $(I_{a,\gamma_0}^{\mathrm{o}})_{\Qv}$ at least at one place $v$ (i.e. at $v=\infty$), according to \cite[Lem.5.6]{LR87} (for which one does not need $G^{\ast}$ to be quasi-split; see also the discussion in $\S$9 of \cite{Kottwitz84a}, particularly 9.4.1, 9.5), it suffices to check that $\iota(T')$ transfers to $G_{\gamma_0}^{\mathrm{o}}$ (with respect to the inner twisting $\psi^{-1}$) locally everywhere. Note that for any finite place $l\neq p$, as $\psi_{\Qlb}$ is conjugate to a $\Ql$-isomorphism, we have $(I_{a,\gamma_0}^{\mathrm{o}})_{\Ql}\approx (G_{\gamma_0}^{\mathrm{o}})_{\Ql}$. Let $S$ be the set of places $v$ where $(I_{a,\gamma_0}^{\mathrm{o}})_{\Ql}\approx (G_{\gamma_0}^{\mathrm{o}})_{\Ql}$ is not quasi-split, plus $v=p$ and $\infty$; by our assumption that $T$ is nice, for every $v\in S$, $T_{\Qv}$ is fundamental in $(I_{\mathscr{I},\epsilon}^{\mathrm{o}})_{\Qv}\approx (I_{a,\gamma_0}^{\mathrm{o}})_{\Qv}$. This implies the existence of local transfers. In more detail, when $v\in S$, we use the fact \cite[$\S$10]{Kottwitz86} (or \cite[Lem.5.8, 5.9]{LR87}) that fundamental maximal tori of a connected reductive group over a local field transfer to \emph{all} inner forms. When $v\notin S$, we resort to the fact \cite[p.340]{PR94} that any maximal torus in a connected reductive group transfers to the quasi-split inner form.

If we put $(a'_0)_{\tau}:=g_1a_{\tau}\tau(g_1)^{-1}$, its image in $(G_{\gamma_0}^{\mathrm{o}})^{\ad}$ becomes a cocycle and defines an inner twisting $\psi_0':G_{\gamma_0}^{\mathrm{o}} \isom I_{a'_0,\gamma_0}^{\mathrm{o}}$. The $\Qb$-morphism $\Int(g_1): (I_{a,\gamma_0}^{\mathrm{o}})_{\Qb} \isom (I_{a'_0,\gamma_0}^{\mathrm{o}})_{\Qb}$ induces a $\Q$-isomorphism $\iota'_0:I_{a,\gamma_0}^{\mathrm{o}}\isom I_{a'_0,\gamma_0}^{\mathrm{o}}$, and satisfies $\iota'_0\circ \psi=\psi_0'\circ \Int(g_1)$, giving the left square in the above diagram. 

\textbf{Step 2.b.} Note that as $\Int(g_1)\circ \iota$  and $\Int(g_1)\circ\psi^{-1}|_{\iota(T')}$ are both $\Q$-rational, so is $\psi'_0|_{T'_0}$ for $T'_0:=\Int(g_1)\circ\psi^{-1}(\iota(T'))\subset G_{\gamma_0}^{\mathrm{o}}$. As has been observed above, this is the same condition as that the cochain $(a'_0)_{\tau}=g_1a_{\tau}\tau(g_1)^{-1}=g_1g_0(a_0)_{\tau}\tau(g_1g_0)^{-1}$ is $T'_0$-valued. Therefore, by Theorem \ref{thm:LR_special_admissible_morphism} (iv), there exist $g_2\in G(\Qb)$ and a special Shimura subdatum $(T'_1,h'_1)$ such that $\Int(g_2)$ induces a transfer of maximal torus $T'_0\isom T'_1=g_2T'_0(g_2)^{-1}$ and the cochain $a'_{\tau}:=g_2(a'_0)_{\tau}\tau(g_2)^{-1}$ is the one attached to $(T'_1,h'_1)$ as in Theorem \ref{thm:LR_special_admissible_morphism}; then, we obtain the required commutative diagram for the inner twisting $\psi':G_{\gamma_0'}^{\mathrm{o}} \isom I_{a',\gamma_0}^{\mathrm{o}}$ defined by $\bar{a}'\in Z^1(\Q,T'_1/Z(G))$ and $\iota':=\Int(g')$ for $g':=g_2g_1$.
\end{proof}

Following Kottwitz \cite[2.1,2.2]{Kottwitz97}, we define a \emph{$\sigma$-$\Qpnr$-space} to be a finite dimensional vector space $U$ over $\Qpnr$ and a $\sigma$-linear bijection $\Phi:U\ra U$, and a \emph{$W_{\Qp}$-$\Qpb$-space} to be a finite dimensional $\Qpb$-vector space $U$ equipped with a discrete, semilinear action of the Weil group $W_{\Qp}$ (\textit{semilinear} means that $\tau(\alpha u)=\tau(\alpha)\tau(u)$ for all $\tau\in W_{\Qp}$, $\alpha\in \Qpb$, $u\in U$).
There is an obvious (tensor) functor from $\sigma$-$\Qpnr$-spaces to $W_{\Qp}$-$\Qpb$-spaces $U\mapsto \Qpb\otimes_{\Qpnr}U$, the action of $W_{\Qp}$ on $\Qpb\otimes_{\Qpnr}U$ being given by 
\[ \tau(\alpha\otimes u)=\tau(\alpha)\otimes \Phi^ju \] 
for all $\alpha\in\Qpb$, $u\in U$ and $\tau\in W_{\Qp}$ mapping to $\sigma^j\in\langle\sigma\rangle$.
There is a (tensor) functor quasi-inverse to this functor: $V\mapsto V^I$ (invariants of inertia).
The natural extensions of $\sigma$-$\Qpnr$ spaces to $L=\mathrm{Frac}(W(\Fpb))$ are isocrystals over $\Fpb$ (which were also called $\sigma$-$L$-space by Kottwitz, loc. cit.) which in turn extend to $W_{\Qp}$-$\bar{L}$-space by the same formula as above. The passage from $\sigma$-$L$-spaces to $W_{\Qp}$-$\bar{L}$-spaces induces an equivalence of (tensor) categories. 
For a linear algebraic group $G$ over $\Qp$, and any of theses spaces, there is a notion of ``that space with $G$-structure''. The isomorphism classes of $\sigma$-$L$ spaces or $W_{\Qp}$-$\bar{L}$-spaces with $G$-structure are both classified by $B(G)=H^1(W_{\Qp},G)$.

%%%%%%%%%%%%%%%%%%%%:
\begin{cor} \label{cor:stable_conjugacy_from_nice_tori}
For any two nice maximal $\Q$-tori $T$, $T'$ of $I$ with embeddings $i_T:T\hookrightarrow G$, $i_{T'}:T'\hookrightarrow G$ as in Theorem \ref{thm:Kisin17_Cor.1.4.13,Prop.2.1.3,Cor.2.2.5} (3), the elements $i_T(\epsilon)$ and $i_{T'}(\epsilon)$ are stably conjugate.
\end{cor}

\begin{proof}
The $\Q$-embedding $i_T$ is provided by a CM lifting $\tilde{x}$ of a point $x$ in the isogeny class $\mathscr{I}$. Suppose that $x$ is defined over a finite field $\F_q$ and
let $H_1^{\cris}:=H_1^{\cris}(\mathcal{A}_x/W(\F_q))\otimes\Qpnr$ denote its crystalline homology group, regarded as a $\sigma$-$\Qpnr$-space; fixing an isomorphism $I_{\mathscr{I}}\isom I_x$, we regard $T'$ as a subgroup of $\mathrm{Aut}(\mathcal{A}_x)$. Recall the special Shimura subdatum $(T'_1,h'_1)$ in Theorem \ref{thm:stable_isogeny_diagram} (iv); we have a $\Q$-isomorphism $\iota'_{T'}:=\psi'^{-1}\circ\iota'\circ\iota: T' \isom T'_1$ (\ref{eq:stable_isogeny_diagram}).
It suffices to show that 

\emph{There exists an isomorphism $\eta_p':H_1^{\cris}\otimes\Qpb \isom (V_{\Qpnr},\delta'\sigma)_{\Qpb}$ of $W_{\Qp}$-$\Qpb$-spaces for some $\delta'\in T'_1(\Qpnr)$ with $[\delta']= -\underline{\mu}_{h'_1}\in B(T'_1)\simeq X_{\ast}(T'_1)_{\Gamma_p}$ which preserves the extra structures (quasi-polarizations and the relevant tensors) and is compatible with the actions of $\iota'_{T'}:T'\isom T'_1$.} 

Indeed, then $\mu_{T'}:=(\iota'_{T'})^{-1}(\mu_{h'_1})\in X_{\ast}(T')$ would satisfy the Kisin's conditions \cite[Lem.2.2.2]{Kisin17} with respect to $\eta_p'$, which shows that the reduction of any special point defined by $(T'_1,h'_1)$ is $\mathscr{I}$ and the $\Q$-embedding $\iota'_{T'}$ qualifies as a $\Q$-embedding $T'\hookrightarrow G$ in Theorem \ref{thm:Kisin17_Cor.1.4.13,Prop.2.1.3,Cor.2.2.5} (3). The same theorem also says that $\iota'_{T'}(\epsilon)$ is stably conjugate to $i_{T'}(\epsilon)$. But, $\iota'_{T'}(\epsilon)=\Int(g)\circ\varphi(\epsilon)=\Int(g)(i_T(\epsilon))$ by (\ref{eq:stable_isogeny_diagram}) ($\epsilon\in T(\Q)$) and $\Int(g)(i_T(\epsilon))$ is stably conjugate to $i_T(\epsilon)$, as $\Int(g):T'\isom T'_1$ is a transfer of maximal tori, which proves the corollary.

Let $(T_0:=i_T(T),h_0)$ be the special Shimura subdatum giving the CM lifting $\tilde{x}$. 
When $a_0$ is the $T_0$-valued cochain attached to this datum as in Theorem \ref{thm:LR_special_admissible_morphism}, it was shown in the course of construction of (\ref{eq:stable_isogeny_diagram1}) that $a_{\tau}=g_0(a_0)_{\tau}\tau(g_0)^{-1}$ for some $g_0\in G_{\pi_0}(\Qb)$.

For each cochain $c=a_0,a,a'$, let $c(p)$ be its restriction to $\Gamma_p=\Gal(\Qpb/\Qp)$.
Choose $t_{p0}\in T_0(\Qpb)$ and $t'_p\in T'_1(\Qpb)$ such that $(a_0(p)^{\nr})_{\tau}:=t_{p0}a_{0}(p)_{\tau}\tau(t_{p0})^{-1}$ and $(a'(p)^{\nr})_{\tau}:=t'_pa'_{\tau}\tau(t'_p)^{-1}$ are unramified (Theorem \ref{thm:stable_isogeny_diagram} (iii)). Note that 
\[ (a(p)^{\nr})_{\tau}:=t_pa_{\tau}\tau(t_p)^{-1}=(a_0(p)^{\nr})_{\tau} \] 
for $t_p:=t_{p0}g_0^{-1}\in G_{\pi_0}(\Qpb)$.
These unramified cochains give $\sigma$-$\Qpnr$-spaces $(V_{\Qpnr},\delta\sigma)$, $(V_{\Qpnr},\delta'\sigma)$, where
\[ \delta:=(a(p)^{\nr})_{\sigma} \in T_0(\Qpnr),\quad \delta':=(a'(p)^{\nr})_{\sigma} \in T'_1(\Qpnr) \]  
(so, $[\delta']= -\underline{\mu}_{h'_1}\in B(T'_1)\simeq X_{\ast}(T'_1)_{\Gamma_p}$, by the same theorem). 
The $\sigma$-conjugacy classes $[\delta]$, $[\delta']$ in $B(G)=H^1(W_{\Qp},G)$ depend only on the cochains $a(p)$, $a'(p)$ (not on the choice of $t_p$, $t_p'$ giving unramified cochains $a(p)^{\nr}$, $a'(p)^{\nr}$).

From Theorem \ref{thm:stable_isogeny_diagram}, we obtain a commutative diagram consisting of isomorphisms of $\Qpb$-vector spaces which are endowed with either cochains on $\Gamma_p$ (upstairs) or structures of $W_{\Qp}$-$\Qpb$-spaces (downstairs):
\begin{equation} \label{eq:stable_isom_diagram}
\xymatrix{  &  (V_{\Qpb},a(p)) \ar[r]^{g}_{(1)} \ar[d]_{t_p}^{(2)} &  (V_{\Qpb},a'(p)) \ar[d]_{t'_p}^{(3)} \\ 
H_1^{\cris}\otimes_{\Qpnr}\Qpb  \ar[r]^{\eta_p} \ar@(dl,dr)[]_{T'}  &  (V_{\Qpnr},\delta\sigma)_{\Qpb} \ar[r]^{f} \ar@(dl,dr)[]_{\iota(T')}   &  (V_{\Qpnr},\delta'\sigma)_{\Qpb} \ar@(dl,dr)[]_{T'_1\simeq \psi'(T'_1)}  }
\end{equation}
First, we explain the square diagram. The isomorphism (1) (of $\Qpb$-vector spaces) is the left multiplications by $g\in G(\Qpb)$. It carries the cochain $a(p)_{\tau}$ to $a'(p)_{\tau}:=ga(p)_{\tau}\tau(g)^{-1}$. In the vertical isomorphisms (2), (3) which are multiplications by $t_p,t'_p\in G(\Qpb)$, the cochains upstairs and the $\sigma$-$\Qpnr$-structures downstairs are compatible in the following sense: for example, $t_p$ carries the cochain $a(p)_{\tau}$ to the unramified cochain $a(p)^{\nr}_{\tau}=t_pa(p)_{\tau}\tau(t_p)^{-1}$, which amounts to the $\sigma$-$\Qpnr$-space $(V_{\Qpnr},\delta\sigma)$, or the associated $W_{\Qp}$-$\Qpb$-space $(V_{\Qpnr},\delta\sigma)_{\Qpb}$. These vertical isomorphisms induce $\Qp$-embeddings %%
%\footnote{In fact, these are isomorphisms, which would follow from some finer information that in Theorem \ref{thm:LR_special_admissible_morphism} (i), the coboundary and the Frobenius germ generate the same subgroup.}
%% 
\[ (I_{a})_{\Qp}= (G_{\pi_0})_{a(p)} \isom (G_{\pi_0})_{a(p)^{\nr}} \hra  G_{\delta\sigma},\quad (I_{a'})_{\Qp}\isom ((G_{\pi_0'})_{\Qp})_{a'(p)^{\nr}}  \hra G_{\delta'\sigma}. \] 
where for $c=a(p), a(p)^{\nr}$, $((G_{\pi_0})_{\Qp})_{c}$ is the inner twist of $(G_{\pi_0})_{\Qp}$ by the cocycle $\overline{c}\in Z^1(\Qp, G_{\pi}^{\ad})$.
The group $(G_{\delta\sigma})_{\Qpb}$ acts on $V_{\Qpb}$ via $\Res_{L_n/\Qp}(G_{L_n})_{\Qpb}\cong\prod_{L_n\hra \Qpb}G_{\Qpb}\ra G_{\Qpb}$ (projection onto the identity factor), and the induced action of $(G_{\pi_0})_{a(p)^{\nr}}(\Qpb)$ is via the identification $(G_{\pi_0})_{a(p)^{\nr}}(\Qpb)=G_{\pi_0}(\Qpb)\subset \GL(V_{\Qpb})$. Hence, we see that $(I_a)_{\Qp}$, $(I_{a'})_{\Qp}$ act respectively on the $\sigma$-$\Qpnr$-spaces $(V_{\Qpnr},\delta\sigma)$, $(V_{\Qpnr},\delta'\sigma)$ via 
\[ \Int(t_p)\circ\psi^{-1},\quad \Int(t'_p)\circ\psi'^{-1};\]
this gives the actions in the bottom of the square diagram.
Then, the commutativity of the diagrams (\ref{eq:stable_isogeny_diagram}) implies that multiplication by $f:=t'_pgt_p^{-1}$ is compatible with these actions via the $\Qp$-isomorphism $\iota':\iota(T'_{\Qp})\isom \psi'((T'_1)_{\Qp})$, i.e. as $t'_p\in T'_1(\Qpb)$, we have the equality:
\begin{equation} \label{eq:itoa(T')=T'_1}
\Int(ft_p)\circ\psi^{-1}|_{\iota(T')}=\psi'^{-1}\circ \iota'|_{\iota(T')}\ :\ \iota(T'_{\Qp})\isom (T'_1)_{\Qp}.
\end{equation}
Here, it is most important that the $(T'_1)_{\Qpb}$-action on $V_{\Qpb}$ via $\psi'(T'_1)_{\Qb}\hra (I_{a'})_{\Qb}\subset G_{\Qb}$ arises from the $\Q$-structure $T'_1\subset G\subset \GL(V)$, (this is not true for the $\iota(T')_{\Qb}$-action on $V_{\Qb}$ via $\iota(T')\subset I_{a}$) and that $[\delta']= -\underline{\mu}_{h'_1} \in B(T'_1)=X_{\ast}(T'_1)_{\Gamma_p}$ by Theorem \ref{thm:LR_special_admissible_morphism}.

Next, we explain the isomorphism $\eta_p$. By Theorem \ref{thm:LR_special_admissible_morphism} (iii), there exists $t_0\in T_0(\mfk)$ such that 
\[ \Nm_{L_{\mathfrak{P}}/L_{\mathfrak{P},0}}(\mu_{h_{0}}(\varpi^{-1}))=t_0a_{0}(p)^{\nr}_{\sigma}\sigma(t_0)^{-1}=:\delta_0, \] 
since $\kappa_{T_0}(\Nm_{L_{\mathfrak{P}}/L_{\mathfrak{P},0}}(\mu_{h_{0}}(\varpi^{-1})))=-\underline{\mu}_{h_{0}}=\kappa_{T_0}(a_{0}(p)^{\nr}_{\sigma})$; so, $\Nm_n(\delta_0)=\Nm(\mu_{h_{0}})(\varpi^{-1})$, where $n=s[L_{\mathfrak{P},0}:\Qp]$ ($s\gg1$) (the right side is a representative of the germ of Frobenius endomorphism).
By construction (Corollary \ref{cor:Tate_thm2} (3)), the inner twisting $\varphi:(I_{\mathscr{I}})_{\Qb}\isom (G_{\pi_0})_{\Qb}$ in (\ref{eq:stable_isogeny_diagram}) restricts to $i_T$, and over $\Qpb$, $\varphi$ is conjugate to the composite $(I_{\mathscr{I}})_{\Qpb}\isom (G_{\delta_0\sigma})_{\Qpb}\isom (G_{\pi_0})_{\Qpb}$. The former isomorphism $(I_{\mathscr{I}})_{\Qp}\isom G_{\delta_0\sigma}$ is obtained by an isomorphism $H_1^{\cris} \isom (V_{\Qpnr},\delta_0\sigma)$ of $\sigma$-$\Qpnr$-spaces (with extra structures), and the latter isomorphism $(G_{\delta_0\sigma})_{\Qpb}\isom (G_{\pi_0})_{\Qpb}$ is conjugate over $\Qpb$ to $p_{\delta_0}:(G_{\delta_0\sigma})_{L_n}\isom (G_{L_n})_{\Nm\delta_0}$ which is obtained by restriction from the projection onto the identity factor $\Res_{L_n/\Qp}(G_{L_n})_{L_n}\cong\prod_{L_n\hra L_n}G_{L_n}\ra G_{L_n}$ (Proposition \ref{prop:psi_p}).
Since $\Int(g)$ and $\psi$ in (\ref{eq:stable_isogeny_diagram}) are both induced by isomorphisms of vector spaces $V_{\Qb}\isom V_{\Qb}$, it follows that $\Int(t_p)\circ\psi^{-1}\circ \iota_{\Qpb}:(I_{\mathscr{I}})_{\Qpb}\isom (I_{a})_{\Qpb} \hra (G_{\delta\sigma})_{\Qpb}$ is also induced by an isomorphism $W_{\Qpb}$-$\Qpb$-spaces
\[\eta_p: H_1^{\cris}\otimes\Qpb \isom (V_{\Qpnr},\delta\sigma)_{\Qpb} \]
(which preserves the extra structures of quasi-polarizations and the tensors): namely, $\eta_p$ is compatible with the actions of $\Int(t_p)\circ\psi^{-1}\circ\iota:T'\isom \iota(T')\hra G_{\delta\sigma}\hra \GL(V_{\Qpb})$. 

Now, $\eta_p':=f\circ \eta_p$ is an isomorphism  $W_{\Qp}$-$\Qpb$-spaces, compatible with the actions of $ T'$, $T'_1$ via the isomorphism $\Int(f)\circ\Int(t_p)\circ\psi^{-1}\circ\iota=\psi'^{-1}\circ\iota'\circ\iota=\iota'_{T'}$ (\ref{eq:itoa(T')=T'_1}), which we looked for.
\end{proof}

This corollary and Lemma \ref{lem:K-triple_attached_to_admissible_K-pair} allow for the following definition:

%%%%%%%%%%%%%%%%%%%%
\begin{defn} \label{defn:stable_K-triple_attached_to_admissible_K-pair}
For any admissible K-pair $(\mathscr{I},\epsilon)$, we call the stable Kottwitz triple $(\gamma_0;\gamma,\delta)$ obtained by the recipe (\ref{eq:K-triple_attached_to_admissible_K-pair}) for any \emph{nice} maximal $\Q$-torus $T\subset I_{\mathscr{I}}$ with $\epsilon\in T(\Q)$ the \emph{stable Kottwitz triple attached to} the admissible K-pair $(\mathscr{I},\epsilon)$.
\end{defn}

%%%%%%%%%%%%%%%%%%%%
To prove the effectivity criterion of stable Kottwitz triple (\ref{eq:E}) and establish a formula for the number of admissible pairs giving a fixed effective stable Kottwitz triple, we follow the idea in \cite{LR87} of twisting an ``admissible pair'' (a Galois gerb theoretic version of admissible K-pair: it consists of an ``admissible morphism'' $\phi$ and an element $\epsilon\in \underline{\mathrm{Aut}}(\phi)(\Q)$(a subgroup of $G(\Qb)$) by certain cohomology classes in $H^1(\Q,\underline{\mathrm{Aut}}(\phi,\epsilon))$).
For that we adapt a similar technique of Kisin of twisting an isogeny class $\mathscr{I}$ by cohomology classes in $\Sha^{\infty}_G(\Q,I_{\mathscr{I}})$, and twist arbitrary admissible K-pair $(\mathscr{I},\epsilon\in I_{\mathscr{I}}(\Q))$ by cohomology classes in 
\[ \im[\Sha^{\infty}_G(\Q,I_{\mathscr{I},\epsilon}^{\mathrm{o}})\rightarrow H^1(\Q,I_{\mathscr{I},\epsilon})],\]
where the group $\Sha^{\infty}_G(\Q,I_{\mathscr{I},\epsilon}^{\mathrm{o}})$ is defined by
\begin{equation} \label{eq:Sha_G^{infty}(Q,I_{sI,epsilon}^{mathrm{o}})}
 \Sha^{\infty}_G(\Q,I_{\mathscr{I},\epsilon}^{\mathrm{o}}):=\ker\left[\Sha^{\infty}(\Q,I_{\mathscr{I},\epsilon}^{\mathrm{o}})\rightarrow \Sha^{\infty}(\Q,I_{\mathscr{I}}) \stackrel{\ast}{\rightarrow} \Sha^{\infty}(\Q,G)\right].
\end{equation}
(the map $\ast$ is defined by any inner twist $I$ of $I_{\mathscr{I}}$ such that $I\hookrightarrow G$ \cite[(4.4.3)]{Kisin17}.)

Suppose given an isogeny class $\mathscr{I}$ and $T\subset I_{\mathscr{I}}$ a maximal $\Q$-torus. Let $i_T:T\hookrightarrow G$ and $h\in X\cap \Hom(\dS,i_T(T)_{\R})$ be as in Theorem \ref{thm:Kisin17_Cor.1.4.13,Prop.2.1.3,Cor.2.2.5} (3).
Let $\tilde{\beta}\in \Sha^{\infty}(\Q,T)$ and assume that the image of $i_T(\tilde{\beta})\in \Sha^{\infty}(\Q,i_T(T))$ in $H^1(\Q,G)$ is trivial. So, there exists $\tilde{\omega}\in G(\Qb)$ such that the cochain $\tau\mapsto \tilde{\omega}^{-1}\tau(\tilde{\omega})$ on $\Gal(\Qb/\Q)$ belongs to $Z^1(\Q,i_T(T))$ and as such one has 
\begin{equation} \label{eq:tilde{omega}}
[\tilde{\omega}^{-1}\tau(\tilde{\omega})]=i_T(\tilde{\beta}) 
\end{equation}
in $H^1(\Q,i_T(T))$ (thus, $\tilde{\omega}$ is defined up to right translate under $i_T(T)(\Qb)$). Equivalently, $\Int(\tilde{\omega}):G_{\Qb}\isom G_{\Qb}$ becomes $\Q$-rational when restricted to the maximal torus $i_T(T)$ (i.e. is a transfer of maximal torus $i_T(T)\hookrightarrow G$)
and its base-change $\Int(\tilde{\omega})_{\R}:i_T(T)_{\R}\hookrightarrow G_{\R}$ is induced from conjugation by an element in $G(\R)$ (i.e. there exists $g\in G(\R)$ such that $\Int(\tilde{\omega})(t)=gtg^{-1}$ for all $t\in i_T(T)(\R)$).
Therefore, for such $\tilde{\omega}$, the pair 
\begin{equation} \label{eq:stable-conj._of_special_SD}
(i_T(T)^{\tilde{\beta}}, h^{\tilde{\beta}}):=(\Int(\tilde{\omega})(i_T(T)),\Int(\tilde{\omega})(h))
\end{equation}
is another special Shimura subdatum of $(G,X)$, which is easily seen to depend only on 
$(i_T(T),h)$ and the cohomology class $\tilde{\beta}\in \Sha^{\infty}(\Q,T)$, but not on the choice of $\tilde{\omega}$.

%%%%%%%%%%%%%%%%%%%%
\begin{thm} \label{thm:Kisin17_Prop.4.4.8}
(1) Let $\mathscr{I}$ be an isogeny class and $\beta_1,\beta_2\in \Sha^{\infty}_G(\Q,I_{\mathscr{I}})$. 
For each $i=1,2$, choose a maximal $\Q$-torus $T_i\subset I_{\mathscr{I}}$ such that $\beta_i\in H^1(\Q,I_{\mathscr{I}})$ is the image of some $\tilde{\beta}_i\in H^1(\Q,T_i)$ (which exists \cite[Thm.5.11]{Borovoi98} and must lie in $\Sha^{\infty}(\Q,T_i)$ \cite[4.4.5]{Kisin17}), and fix $i_{T_i}:T_i\hookrightarrow G$, $h_i\in X\cap \Hom(\dS,i_{T_i}(T_i)_{\R})$, and $\tilde{\omega}_i\in G(\Qb)$ as above. 

Then, the isogeny classes $\mathscr{I}_i\ (i=1,2)$ of the reductions of $[\Int(\tilde{\omega}_i)(h_i),1]$ are the same if and only if $\beta_1=\beta_2$.
Namely, for any $\beta\in \Sha^{\infty}_G(\Q,I_{\mathscr{I}})$, the corresponding isogeny class depends only on $\beta$ (not on the choices of auxiliary data $T$, $\tilde{\beta}$, $i_{T}$, $h$, $\tilde{\omega}$), and when we denote it by $\mathscr{I}^{\beta}$ (the \emph{twist} of $\mathscr{I}$ by $\beta$), the assignment $\beta\mapsto \mathscr{I}^{\beta}$ defines an inclusion of $\Sha^{\infty}_G(\Q,I_{\mathscr{I}})$ into the set of isogeny classes in $\sS_{\mathbf{K}_p}(G,X)(\Fpb)$.

(2) Let $\mathscr{I}$ be an isogeny class and $\beta\in \Sha^{\infty}_G(\Q,I_{\mathscr{I}})$. Then, there exist a $\Qb$-isomorphism 
\begin{equation*} \label{eq:inner-twist:xi}
 \xi:(I_{\mathscr{I}})_{\Qb}\isom (I_{\mathscr{I}^{\beta}})_{\Qb}
\end{equation*}
and a representative $(a_{\tau})_{\tau}\in Z^1(\Q,I_{\mathscr{I}})$ of $\beta$ such that
\setlist{nolistsep} \begin{itemize} [noitemsep] 
\item[(i)] the image of $a_{\tau}$ in $Z^1(\Q,I_{\mathscr{I}}^{\ad})\hookrightarrow Z^1(\Q,\mathrm{Aut}(I_{\mathscr{I}}))$ equals the cocycle $(\xi^{-1}\circ{}^{\tau}\xi)_{\tau}$;
\item[(ii)] if $(b,i_p,i^p)$ is a triple attached to $\mathscr{I}$ as in Theorem \ref{thm:Kisin17_Cor.1.4.13,Prop.2.1.3,Cor.2.2.5}, for any $h_l\in G(\Qlb)\ (l\neq p)$ and $h_p\in J_b(\Qpb)$ such that $i_v(a_{\tau})=h_v^{-1}\tau(h_v)$ for all $\tau\in\Gamma_v$, 
the datum $(b',i_p',{i^p}')$ defined by that $b'=b$ and 
\begin{equation} \label{eq:i_v'} i_v'=\Int(h_v)\circ i_v \circ \xi^{-1} \end{equation} 
is a triple attached likewise to $\mathscr{I}^{\beta}$.
\end{itemize}

Moreover, for any maximal $\Q$-torus $T$ of $I_{\mathscr{I}}$ admitting $\tilde{\beta}\in H^1(\Q,T)$ mapping to $\beta$,
there exist such an isomorphism
\begin{equation} \label{eq:inner-twist:xi2} \xi:(I_{\mathscr{I}})_{\Qb}\isom (I_{\mathscr{I}^{\beta}})_{\Qb} \end{equation}
and a representative $(a_{\tau})_{\tau}\in Z^1(\Q,I_{\mathscr{I}})$ of $\beta$, 
with the additional properties that 
\begin{itemize} [noitemsep]
\item[(iii)] $\xi$ induces a $\Q$-embedding $T\hookrightarrow I_{\mathscr{I}^{\beta}}$;
\item[(iv)] if $(b,i_p,i^p)$ is a triple attached to $\mathscr{I}$ and adapted to a pair $(i_T:T\hookrightarrow G, h\in X\cap \Hom(\dS,i_T(T)_{\R}))$ and $\tilde{\omega}\in G(\Qb)$ is any element satisfying (\ref{eq:tilde{omega}}), there exists a $G(\A_f^p)\times G(\mfk)$-conjugate $(b'\sigma,i_p',{i^p}')$ of a triple in (ii) which is adapted to the pair
\[ i_{T^{\tilde{\beta}}}:T^{\tilde{\beta}} \stackrel{\xi^{-1}}{\longrightarrow} T \stackrel{\Int(\tilde{\omega})\circ i_T}{\hookrightarrow} G,\quad  h^{\tilde{\beta}}\in X\cap\Hom(\dS,i_{T^{\tilde{\beta}}}(T^{\tilde{\beta}})_{\R}),\]
where $T^{\tilde{\beta}}:=\xi(T)$ and $h^{\tilde{\beta}}$ is (\ref{eq:stable-conj._of_special_SD}).
%In particular, $\mathscr{I}^{\beta}$ is the reduction of the special point $[h^{\tilde{\beta}},1]$.
The cocycle $a_{\tau}$ is the image of $\tilde{a}_{\tau}:=i_T^{-1}(\tilde{\omega}^{-1}\tau(\tilde{\omega}))\in Z^1(\Q,T)$.
\end{itemize}

(3) Let $(\mathscr{I},\epsilon)$ be an admissible K-pair and $\beta\in \im[\Sha^{\infty}_G(\Q,I_{\mathscr{I},\epsilon}^{\mathrm{o}})\rightarrow H^1(\Q,I_{\mathscr{I},\epsilon})]$. 
Then, there exists a K-pair $(\mathscr{I}^{\beta},\epsilon^{\beta})$ with the following properties: 
\begin{itemize} [noitemsep]
\item[(v)] as an isogeny class, $\mathscr{I}^{\beta}$ is the twist of $\mathscr{I}$ by the image of $\beta$ in $\Sha^{\infty}_G(\Q,I_{\mathscr{I}})$ defined in (1); 
\item[(vi)] for any maximal $\Q$-torus $T$ of $I_{\mathscr{I},\epsilon}^{\mathrm{o}}$ such that $\beta$ is the image of some $\tilde{\beta}\in H^1(\Q,T)$, 
one has
\begin{equation} \label{eq:epsilon^{beta}=xi(epsilon)}
\epsilon^{\beta}=\xi(\epsilon)
\end{equation}
up to $I_{\mathscr{I}^{\beta}}(\Q)$-conjugacy, where $\xi:(I_{\mathscr{I}})_{\Qb}\isom (I_{\mathscr{I}^{\beta}})_{\Qb}$ is (\ref{eq:inner-twist:xi2}). Also, its $I_{\mathscr{I}^{\beta}}(\Q)$-conjugacy class contains 
\begin{equation} \label{eq:epsilon^{P}}
\Int(\tilde{\omega})(i_T(\epsilon)),
\end{equation}
where $i_T:T\hookrightarrow G$ and $\tilde{\omega}$ are as in (iv); here, $\Int(\tilde{\omega})(i_T(\epsilon))\in i_T(T)^{\tilde{\beta}}(\Q)$ is regarded as an element of $I_{\mathscr{I}^{\beta}}(\Q)$ via any embedding $j_{i_T(T)^{\tilde{\beta}},h^{\tilde{\beta}}}:i_T(T)^{\tilde{\beta}}\hookrightarrow I_{\mathscr{I}^{\beta}}$ in Theorem  \ref{thm:Kisin17_Cor.1.4.13,Prop.2.1.3,Cor.2.2.5} (4).
\item[(vii)] the assignment $\beta\mapsto (\mathscr{I}^{\beta},\epsilon^{\beta})$ gives a well-defined inclusion of $\im[\Sha^{\infty}_G(\Q,I_{\mathscr{I},\epsilon}^{\mathrm{o}})\rightarrow H^1(\Q,I_{\mathscr{I},\epsilon})]$ into the set of equivalence classes of admissible K-pairs.
\end{itemize}
\end{thm}

In the rest of this article, for an element $\beta$ of $H^1(\Q,I_{\mathscr{I},\epsilon}^{\mathrm{o}})$, its images in $H^1(\Q,I_{\mathscr{I},\epsilon})$ or $H^1(\Q,I_{\mathscr{I}})$ will be also denoted by $\beta$ by abuse of notation (illustrations: $(\mathscr{I}^{\beta},\epsilon^{\beta})$, $\mathscr{I}^{\beta}$).

%%%%%%%%%%%%%%%%%%%%
\begin{rem} \label{rem:issue_with_twisting_method}
(1) For the proof of (\ref{eq:E}) and (\ref{eq:C}), it is critical to know the twisted group $I_{\mathscr{I}^{\beta}}$ ``on the nose'', i.e. not just by its inner class (the image of $\beta$ in $H^1(\Q,I_{\mathscr{I}}^{\ad})$), but as a group endowed with an isomorphism $\xi:(I_{\mathscr{I}})_{\Qb}\isom (I_{\mathscr{I}^{\beta}})_{\Qb}$ and an explicit relation, as in (ii), between the associated data $(b,i_p,i^p)$, $(b',i_p',{i^p}')$. 
It is noteworthy that $\xi$ (\ref{eq:xi_a_{tau}}) is constructed via an isomorphism (\ref{eq:Qb-isogeny_theta}) of abelian varieties, i.e. of motives (which can be treated as ``vector spaces'' in this situation).

(2) In the proof of statement (ii), it will be shown that for each finite place $v$, there always exists an element $h_v$ satisfying the cocycle equation $i_v(a_{\tau})=h_v^{-1}\tau(h_v)$. When $(b,i_p,i^p)$ and $a_{\tau}$ are fixed, such an element is uniquely determined by this equation, up to left-translate under $G(\Ql)$ or $J_b(\Qp)$. In particular, the triples $(b',i_p',{i^p}')$ (\ref{eq:i_v'}) obtained thus are indeed all equivalent.

(3) In the statement (iv), the $G(\Q)$-conjugacy class of $i_{T^{\tilde{\beta}}}$ does not depend on the choice of $\tilde{\omega}$ and $i_T(T)^{\tilde{\beta}}=i_{T^{\tilde{\beta}}}(T^{\tilde{\beta}})$. 
\end{rem}

\begin{proof}
(1) This is \cite[Prop.4.4.8]{Kisin17}.

(2) and (3): We prove (2) and (3) together, recalling ingredients from the proof of (1) along the way necessary. The isomorphism $\xi:(I_{\mathscr{I}})_{\Qb}\isom (I_{\mathscr{I}^{\beta}})_{\Qb}$ will be constructed as $\Int(\theta_{\tilde{\omega}}^{-1})$ for a $\Qb$-isogeny $\theta_{\tilde{\omega}}:\mathcal{A}_{x}^{\mathcal{P}}\otimes_{\Q}\Qb\isom \mathcal{A}_{x}\otimes_{\Q}\Qb\ (x\in \mathscr{I})$, where $\mathcal{A}_{x}^{\mathcal{P}}$ is the abelian variety corresponding to some point in the twisted isogeny class $\mathscr{I}^{\beta}$ which will be obtained from $\mathcal{A}_x$ through ``twisting'' by $\beta$. To construct these objects, we need a moduli interpretation of $\mathscr{I}^{\beta}$ (which was also used in the proof of (1)). 

\textbf{Construction of $\xi$ and proof of (i)}
For any $\beta\in \Sha^{\infty}_G(\Q,I_{\mathscr{I}})\subset H^1(\Q,I_{\mathscr{I}})$, one can find a maximal $\Q$-torus $T$ of $I_{\mathscr{I}}$ such that $\beta$ is the image of some $\tilde{\beta}\in H^1(\Q,T)$ \cite[Thm.5.11]{Borovoi98}, which then must lie in $\Sha^{\infty}(\Q,T)$ \cite[4.4.5]{Kisin17}. We fix $i_T:T\hookrightarrow G$, $\tilde{\omega}$ as in (1).
If $x\in \mathscr{I}$ denotes the reduction of $[h,1]$, the reduction of $[\Int(\tilde{\omega})(h),1]$ is $\mathbf{i}_{\omega}(x)$ in the notation of (4.4.7) of \cite{Kisin17}. Its underlying abelian variety $\mathcal{A}_{\mathbf{i}_{\omega}(x)}$ is isomorphic to the twist $\mathcal{A}_{x}^{\mathcal{P}}$ of $\mathcal{A}_{x}$ by the $T$-torsor $\mathcal{P}$ corresponding to $\tilde{\beta}\in H^1(\Q,T)$ (ibid., (4.1.6)), and there exists a natural $\Qb$-isogeny 
\begin{equation} \label{eq:Qb-isogeny_theta}
\theta_{\tilde{\omega}}:\mathcal{A}_{x}^{\mathcal{P}}\otimes_{\Q}\Qb\isom \mathcal{A}_{x}\otimes_{\Q}\Qb
\end{equation}
preserving the weak polarizations and the set of (crystalline and \'etale) cycles (the cycles on $ \mathcal{A}_{x}$ are $s_{\alpha,0,x}$, $s_{\alpha,l,x}\ (l\neq p)$ in the notation of the proof of Theorem \ref{thm:Kisin17_Cor.1.4.13,Prop.2.1.3,Cor.2.2.5}, but they can be also regarded as cycles on $\mathcal{A}_{x}^{\mathcal{P}}$ in a natural way (ibid., Lemma 4.1.5, 4.1.7)). 
This $\Qb$-isogeny $\theta_{\tilde{\omega}}$ is unique up to $T(\Qb)$-conjugacy, once we fix an identification $I_{\mathscr{I}}(\Q)= I_x(\Q)$. Indeed, by construction (ibid., 4.1.6, Lemma 4.1.2), it is obtained from some universal isomorphism (in the isogeny category)
\begin{equation} \label{eq:univ_isogeny_theta}
\mathcal{A}_x^{\mathcal{P}}\otimes_{\Q}\mathcal{O}_{\mathcal{P}}\isom \mathcal{A}_{x}\otimes_{\Q}\mathcal{O}_{\mathcal{P}}
\end{equation}
(preserving the same extra structures) by specialization via any point of $\mathcal{P}(\Qb)=\Hom(\mathcal{O}_{\mathcal{P}},\Qb)$, and any two points $p_1,p_2\in\Hom(\mathcal{O}_{\mathcal{P}},\Qb)$ differ by a point $t\in T(\Qb)$ acting on $\mathcal{O}_{\mathcal{P}}$ (via $(tf)(x)=f(xt)$ for $f\in \mathcal{O}_{\mathcal{P}}$, $x\in \mathcal{P}$); it is then easy to see that the resulting $\Qb$-isogenies $\mathcal{A}_{x}^{\mathcal{P}}\otimes_{\Q}\Qb\isom \mathcal{A}_{x}\otimes_{\Q}\Qb$ differ by composition with $t\in I_{x}(\Qb)$. Fixing one $\theta_{\tilde{\omega}}$ (\ref{eq:Qb-isogeny_theta}), we put 
\begin{equation} \label{eq:xi_a_{tau}} 
\xi:=\Int(\theta_{\tilde{\omega}}^{-1}):(I_{\mathscr{I}})_{\Qb}\isom (I_{\mathscr{I}^{\beta}})_{\Qb},\quad (a_{\tau})_{\tau}:=(\theta_{\tilde{\omega}}\circ{}^{\tau}\theta_{\tilde{\omega}}^{-1})_{\tau}\in Z^1(\Q,I_{\mathscr{I}}).
\end{equation}
This completes the construction of $\xi$, $(a_{\tau})_{\tau}$, and the proof of (i).

\textbf{Construction of K-pair $(\mathscr{I}^{\beta},\epsilon^{\beta})$ and proof of (3)}
Now, when $\beta\in \im[\Sha^{\infty}_G(\Q,I_{\mathscr{I},\epsilon}^{\mathrm{o}})\rightarrow H^1(\Q,I_{\mathscr{I},\epsilon})]$, as before
we can find a maximal $\Q$-torus $T$ of $I_{\mathscr{I},\epsilon}^{\mathrm{o}}$ such that $\beta$ is the image of some $\tilde{\beta}\in \Sha^{\infty}(\Q,T)$.
Since $\epsilon\in T(\Q)$, it is immediate from the construction of $\mathcal{A}_{x}^{\mathcal{P}}$ (ibid., (4.1.6)) that $\epsilon$ gives a self-isogeny
\[ \epsilon^{\mathcal{P}} \in \mathrm{Aut}_{\Q}(\mathcal{A}_x^{\mathcal{P}}). \]  
Moreover, one readily checks (using ibid., Lemma 4.1.5, 4.1.7) that the universal isomorphism (\ref{eq:univ_isogeny_theta}) takes $\epsilon^{\mathcal{P}}$ to $\epsilon$, so for any $\Qb$-isogeny as in $\theta_{\tilde{\omega}}$ (\ref{eq:Qb-isogeny_theta}) one has 
\begin{equation*} 
\epsilon^{\mathcal{P}} =\Int(\theta_{\tilde{\omega}}^{-1})(\epsilon)
\end{equation*}
which thus lies in $I_{\mathbf{i}_{\omega}(x)}(\Q)$. As an element of $\mathrm{Aut}_{\Q}(\mathcal{A}_x^{\mathcal{P}})$, $\epsilon^{\mathcal{P}}$ still depends on $\mathcal{P}$, i.e. on $\tilde{\beta}$.

To show that the $I_{\mathscr{I}^{\beta}}(\Q)$-conjugacy class of $\epsilon^{\mathcal{P}}$ depends only on $\beta$, not on $\tilde{\beta}$ (so, justifying the notation $\epsilon^{\beta}$) and the statement (vii), 
let us put $I_{\mathbf{i}_{\omega}(x),\epsilon}:=Z_{I_{\mathbf{i}_{\omega}(x)}}(\epsilon)$ and consider the $\Q$-scheme $\mathcal{Q}_{\epsilon}$ of isomorphisms 
\[\mathcal{A}_x^{\mathcal{P}}\lisom \mathcal{A}_{x}, \]
preserving the weak polarizations plus the crystalline and the etale tensors $\{s_{\alpha,0,x}\}$, $\{s_{\alpha,l,x}\}\ (l\neq p)$ and which further takes $\epsilon^{\mathcal{P}}$ to $\epsilon$. 
This is a $I_{\mathscr{I},\epsilon}$-torsor which admits a $T$-equivariant map $\mathcal{P}\rightarrow\mathcal{Q}_{\epsilon}$, hence is isomorphic to the $I_{\mathscr{I},\epsilon}$-torsor associated with $\beta\in H^1(\Q,I_{\mathscr{I},\epsilon})$ (cf. proof of ibid., Prop. 4.4.8). 
This implies that the assignment $\beta\mapsto (\mathscr{I}^{\mathcal{P}},\epsilon^{\mathcal{P}})$ gives a well-defined inclusion of $\im[\Sha^{\infty}_G(\Q,I_{\mathscr{I},\epsilon}^{\mathrm{o}})\rightarrow H^1(\Q,I_{\mathscr{I},\epsilon})]$ into the set of equivalence classes of K-pairs, where $\mathcal{P}$ is the $T$-torsor corresponding to any choice of $\tilde{\beta}\in H^1(\Q,T)$ mapping to $\beta$ (for a maximal $\Q$-torus $T$ of $I_{x,\epsilon}$).
In particular, the $I_{\mathbf{i}_{\omega}(x)}(\Q)$-conjugacy class of $\epsilon^{\mathcal{P}}$ depends only on $\beta$; we let $\epsilon^{\beta}$ denote any representative of the associated $I_{\mathscr{I}^{\beta}}(\Q)$-conjugacy class (via some identification $I_{\mathscr{I}^{\beta}}=I_{\mathbf{i}_{\omega}(x)}$) and write $(\mathscr{I}^{\beta},\epsilon^{\beta})$ for $(\mathscr{I}^{\mathcal{P}},\epsilon^{\mathcal{P}})$. This completes the construction of K-pair $(\mathscr{I}^{\beta},\epsilon^{\beta})$ and we have proved (3), except for the property (\ref{eq:epsilon^{P}}).

\textbf{Proof of (iii), (\ref{eq:epsilon^{P}})} We prove the property (\ref{eq:epsilon^{P}}), i.e. the equality (up to $I_{\mathscr{I}^{\beta}}(\Q)$-conjugacy):
 \[ j^{\tilde{\beta}}(i_T(\epsilon)^{\tilde{\beta}}) = \epsilon^{\beta}, \] 
where $i_T(\epsilon)^{\tilde{\beta}}:=\Int(\tilde{\omega})(i_T(\epsilon))$ and $j^{\tilde{\beta}}:=j_{i_T(T)^{\tilde{\beta}},h^{\tilde{\beta}}}:i_T(T)^{\tilde{\beta}}\hookrightarrow I_{\mathscr{I}^{\beta}}$ is any $\Q$-embedding  (in its $I_{\mathscr{I}^{\beta}}(\Q)$-conjugacy class) specified in Theorem \ref{thm:Kisin17_Cor.1.4.13,Prop.2.1.3,Cor.2.2.5} (4). 

Recall that we have fixed identifications $I_{\mathscr{I}}=I_x$, $I_{\mathscr{I}^{\beta}}=I_{\mathbf{i}_{\omega}(x)}$.
Via $\sigma_p(\tilde{x})=[h,1]$, the special point $[h,1]$ gives a point $\tilde{x}$ in $\sS_{\mathbf{K}_p}(K)$ for a finite extension $K\subset\Qpb$ of $\Qp$. 
Let us fix $\tilde{\omega}\in G(\Qb)$ such that $[\tilde{\omega}^{-1}\tau(\tilde{\omega}))]=i_T(\tilde{\beta})$ in $H^1(\Q,i_T(T))$. Then, the special point $[h^{\tilde{\beta}},1]$ ($h^{\tilde{\beta}}=\Int(\tilde{\omega})(h)$) corresponds to the twist $(\mathcal{A}_{\sigma_p(\tilde{x})}^{\mathcal{P}},\lambda^{\mathcal{P}},\eta_{\mathbf{K}'^{\omega}}^{\omega})$ by the $T$-torsor $\mathcal{P}$ of the triple $(\mathcal{A}_{\sigma_p(\tilde{x})},\lambda,\eta_{\mathbf{K}'})$ associated with $\sigma_p(\tilde{x})$, as constructed in ibid., (4.2.2) (ibid., Prop. 4.2.6). Also, there exists a $\Qb$-isogeny (canonical up to $T(\Qb)$-conjugacy, being determined by choice of a point in $\mathcal{P}(\Qb)$)
\begin{equation} \label{eq:Qb-isogeny_theta_in_char0}
\theta_{\tilde{\omega}}:\mathcal{A}_{\sigma_p(\tilde{x})}^{\mathcal{P}}\otimes_{\Q}\Qb\isom \mathcal{A}_{\sigma_p(\tilde{x})}\otimes_{\Q}\Qb
\end{equation}
which preserves the Betti tensors and is compatible with the weak polarizations; this $\Qb$-isogeny provides, via mod-$p$ reduction, the $\Qb$-isogeny (\ref{eq:Qb-isogeny_theta}). If we choose an isomorphism $\eta_{\Betti}:H^{\Betti}_1(\mathcal{A}_{\sigma_p(\tilde{x})},\Q)\isom V$ as in (\ref{eq:Betti-isom}), the Hodge structure on $H^{\Betti}_1(\mathcal{A}_{\sigma_p(\tilde{x})}^{\mathcal{P}},\Q)$ is given by $h^{\tilde{\beta}}=\Int(\tilde{\omega})(h)$ via some $\Q$-isomorphism 
$H^{\Betti}_1(\mathcal{A}_{\sigma_p(\tilde{x})}^{\mathcal{P}},\Q) \isom V$ which can be constructed from $\eta_{\Betti}$ and $\mathcal{P}$ in a certain natural manner.

The two isogenies $\theta_{\tilde{\omega}}$ (\ref{eq:Qb-isogeny_theta}), (\ref{eq:Qb-isogeny_theta_in_char0}) respectively produce inner twistings
\[ \xi=\Int(\theta_{\tilde{\omega}}^{-1})\ :\ \Aut_{\Q}(\mathcal{A}_{\sigma_p(\tilde{x})})_{\Qb} \isom \Aut_{\Q}(\mathcal{A}_{\sigma_p(\tilde{x})}^{\mathcal{P}})_{\Qb},\quad (I_{x})_{\Qb}\isom (I_{\mathbf{i}_{\omega}(x)})_{\Qb}, \]
and cocycles $(\theta_{\tilde{\omega}}\circ {}^{\tau}\theta_{\tilde{\omega}}^{-1})_{\tau}$ in either $Z^1(\Q, \Aut_{\Q}(\mathcal{A}_{\sigma_p(\tilde{x})}))$ or $Z^1(\Q,I_x)$.
According to \cite[Lem.4.1.2]{Kisin17}, for each $\theta_{\tilde{\omega}}$, we have equalities of cochains in $Z^1(\Q,i_T(T))=Z^1(\Q,T)$:
\begin{equation} \label{eq:Kisin17_Lem.4.1.2}
(\theta_{\tilde{\omega}}\circ {}^{\tau}\theta_{\tilde{\omega}}^{-1})_{\tau}=(\tilde{\omega}^{-1}\tau(\tilde{\omega}))_{\tau}.
\end{equation}
For this, we note that any $\tilde{\omega}\in G(\Qb)$ satisfying $[\tilde{\omega}^{-1}\tau(\tilde{\omega})]=i_T(\tilde{\beta})$ (\ref{eq:tilde{omega}}) gives the same image $\omega$ under the natural quotient map $G\rightarrow G/i_T(T)$ and the $i_T(T)$-torsor $\mathcal{P}$ corresponding to $i_T(\tilde{\beta})$ is the fiber over $\omega$ of this map: in this way, $\tilde{\omega}$ can be regarded as a point of $\mathcal{P}$, cf. \cite[(4.3.2)]{Kisin17}.
In particular, $\theta_{\tilde{\omega}}\circ {}^{\tau}\theta_{\tilde{\omega}}^{-1}\in T$, and thus $\Int(\theta_{\tilde{\omega}}^{-1}):(I_{x})_{\Qb}\isom (I_{\mathbf{i}_{\omega}(x)})_{\Qb}$ induces a $\Q$-isomorphism $T\isom T^{\tilde{\beta}}$ for some maximal $\Q$-torus $T^{\tilde{\beta}}\subset  I_{\mathbf{i}_{\omega}(x)}$; this proves the statement (iii). Also, if we define a $\Q$-embedding 
\[ i_{T^{\tilde{\beta}}}:T^{\tilde{\beta}}\isom i_T(T)^{\tilde{\beta}} \subset G \] 
by the commutativity of the diagram
\begin{equation} \label{eq:new_embedding_T-beta}
\xymatrix{ i_T(T) \ar[r]^{\Int(\tilde{\omega})} & i_T(T)^{\tilde{\beta}} \\ T \ar[u]^{i_T} \ar[r]_{\Int(\theta_{\tilde{\omega}}^{-1})} & T^{\tilde{\beta}} \ar@{->}[u]_{i_{T^{\tilde{\beta}}}}, } 
\end{equation}
its stable conjugacy class is the one specified in Theorem \ref{thm:Kisin17_Cor.1.4.13,Prop.2.1.3,Cor.2.2.5} (3), because we have $h^{\tilde{\beta}}=\Int(\tilde{\omega})(h)$ which implies that the lifting $\mathcal{A}_{\sigma_p(\tilde{x})}^{\mathcal{P}}$ of $\mathcal{A}_{\mathbf{i}_{\omega}(x)}$ is determined by the cocharacter $\mu_{h^{\tilde{\beta}}}=\Int(\theta_{\tilde{\omega}}^{-1})(\mu_h)\in X_{\ast}(T^{\tilde{\beta}})$ by the procedure explained in the proof of Theorem \ref{thm:Kisin17_Cor.1.4.13,Prop.2.1.3,Cor.2.2.5} (3) and the embedding $i_{T^{\tilde{\beta}}}:T^{\tilde{\beta}}\hookrightarrow G $ is the induced lifting of
the action $T^{\tilde{\beta}}\subset I_{\mathbf{i}_{\omega}(x)}$ to $\mathcal{A}_{\sigma_p(\tilde{x})}^{\mathcal{P}}$.

Now, as $\epsilon^{\beta}=\Int(\theta_{\tilde{\omega}}^{-1})(\epsilon)$ up to $I_{\mathscr{I}^{\beta}}(\Q)$-conjugacy (\ref{eq:epsilon^{P}}), we have
\[i_T(\epsilon)^{\tilde{\beta}}=i_{T^{\tilde{\beta}}}(\epsilon^{\beta}).\]
Then, since the $I_{\mathscr{I}^{\beta}}(\Q)$-conjugacy class of $j_{i_{T}(T)^{\tilde{\beta}},h^{\tilde{\beta}}}\circ i_{T^{\tilde{\beta}}}:T^{\tilde{\beta}}\hookrightarrow I_{\mathbf{i}_{\omega}(x)}$ contains the inclusion $T^{\tilde{\beta}} \subset I_{\mathbf{i}_{\omega}(x)}$ (Theorem \ref{thm:Kisin17_Cor.1.4.13,Prop.2.1.3,Cor.2.2.5} (4)), we have $j^{\tilde{\beta}}(i_T(\epsilon)^{\tilde{\beta}}) = \epsilon^{\beta}$ up to $I_{\mathscr{I}^{\beta}}(\Q)$-conjugacy; this proves the statement (\ref{eq:epsilon^{P}}) and thus completes the proof of (vi).

\textbf{Proof of (ii), (iv)} The $\Qb$-isogeny $\theta_{\tilde{\omega}}:\mathcal{A}_{x}^{\mathcal{P}}\otimes_{\Q}\Qb\isom \mathcal{A}_{x}\otimes_{\Q}\Qb$ (\ref{eq:Qb-isogeny_theta}) induces isomorphisms
\begin{align*}
(\theta_{\tilde{\omega}})_{\ast}:H_1^{\et}(\mathcal{A}_{x}^{\mathcal{P}},\Ql)\otimes_{\Ql}\Qlb &\isom H_1^{\et}(\mathcal{A}_{x},\Ql)\otimes_{\Ql}\Qlb, \\
(\theta_{\tilde{\omega}})_{\ast}:H_1^{\cris}(\mathcal{A}_{x}^{\mathcal{P}}/\mfk)\otimes_{\Qp}\mfkb &\isom H_1^{\cris}(\mathcal{A}_{x}/\mfk)\otimes_{\Qp}\mfkb
\end{align*}
preserving the extra structures (of weak polarizations and cycles). 

As $\beta=[\theta_{\tilde{\omega}}\circ {}^{\tau}\theta_{\tilde{\omega}}^{-1}]\in \Sha^{\infty}_G(\Q,I_{x})$ (\ref{eq:Sha^{infty}_G}), its image in $H^1(\Q,G_{\mathbf{ab}})$ under the map
\[\Sha^{\infty}(\Q,I_{x}) \isom \Sha^{\infty}(\Q,I') \rightarrow \Sha^{\infty}(\Q,G)=\Sha^{\infty}(\Q,G_{\mathbf{ab}})\] 
is trivial by assumption, where $I'=G_{\gamma_k}$ is the (connected) centralizer of a semisimple element $\gamma_k$ in $G(\Q)$ (coming from a relative Frobenius endomorphism of $\mathcal{A}_x$ over a sufficiently large finite field) and is an inner form of $I_x$. 
Thus, the localization of $\beta(l)$ at $l$ lies in
\[ \beta(l)\in \ker[H^1(\Ql,I_{\mathscr{I}}) \stackrel{i_{l\ast}}{\rightarrow} H^1(\Ql,G)], \] 
(by $H^1(\Ql,G^{\uc})=\{1\}$ and Lemma \ref{lem:abelianization_exact_seq}) 
and so there exists $h_l\in G(\Qlb)$ such that $i_l(\theta_{\tilde{\omega}}\circ {}^{\tau}\theta_{\tilde{\omega}}^{-1})=h_l^{-1}\cdot \tau(h_l)$.
Recall that $i_l=\Int(\eta_l^{-1}):(I_{x})_{\Ql}\hookrightarrow G_{\Ql}$ is given by a $\Ql$-isomorphism 
 $\eta_l:V\otimes\Ql \isom H^{\et}_1(\mathcal{A}_{\bar{x}},\Ql)$ as in (\ref{eq:isom_eta}) matching the tensors $\{s_{\alpha}\}$, $\{s_{\alpha,l,x}\}$.
Then, for any choice of $h_l$ above, the isomorphism 
\[ \eta_l':=(\theta_{\tilde{\omega}})_{\ast}^{-1} \circ \eta_l \circ h_l^{-1}\ :\ V\otimes\Qlb \ \isom \ H^{\et}_1(\mathcal{A}_{x}^{\mathcal{P}},\Ql)\otimes\Qlb \] 
(mathcing the tensors $\{s_{\alpha}\}$, $\{s_{\alpha,l,x}\}$) is defined over $\Ql$  (i.e. is the base-change of a $\Ql$-isomorphism $V\otimes\Ql \isom H^{\et}_1(\mathcal{A}_{x}^{\mathcal{P}},\Ql)$, again denoted by $\eta_l'$), and thus
\[ i_l':=\Int({\eta_l'}^{-1})=\Int(h_l)\circ i_l\circ \Int(\theta_{\tilde{\omega}}): (I_{\mathbf{i}_{\omega}(x)})_{\Ql} \hookrightarrow G_{\Ql} \] qualifies as a $\Ql$-embedding attached to the isogeny class $\mathscr{I}^{\beta}$ in Theorem \ref{thm:Kisin17_Cor.1.4.13,Prop.2.1.3,Cor.2.2.5}.

At $p$, the given embedding $i_p=\Int(\eta_p^{-1}):(I_{x})_{\Qp}\hookrightarrow J_b$ is given by a $\mfk$-isomorphism $\eta_p: V_{\mfk} \isom H^{\cris}_1(\mathcal{A}_x/\mfk)$ as in (\ref{eq:isom_eta_nr}) matching the tensors $\{s_{\alpha}\}$, $\{s_{\alpha,0,x}\}$. Since we have $H^1(\Qp,\tilde{J}_b)=\{1\}$ by Lemma \ref{lem:LR87_Lem.5.18} below, where $\tilde{J}_b:=\rho^{-1}(J_b)$ for the canonical morphism $\rho:G^{\uc}\rightarrow G$, it follows from Lemma \ref{lem:abelianization_exact_seq} that there exists $h_p\in J_b(\Qpb)$ such that $i_p(\theta_{\tilde{\omega}}\circ {}^{\tau}\theta_{\tilde{\omega}}^{-1})=h_p^{-1}\cdot \tau(h_p)$ for all $\tau\in\Gamma_p$. Then, the isomorphism 
\[ \eta_p':=(\theta_{\tilde{\omega}})_{\ast}^{-1} \circ \eta_p \circ h_p^{-1}\ :\ V_{\mfk}\otimes_{\Qp}\mfkb \ \isom \ H^{\cris}_1(\mathcal{A}_{x}^{\mathcal{P}}/\mfk)\otimes_{\Qp}\mfkb \] 
(mathcing the tensors $\{s_{\alpha}\}$, $\{s_{\alpha,0,x}\}$) is defined over $\mfk$ (i.e. is the base change of a $\mfk$-isomorphism $V_{\mfk}\isom H^{\cris}_1(\mathcal{A}_{x}^{\mathcal{P}}/\mfk)$, again denoted by $\eta_p'$), since it is invariant under the Weil group $W_{\Qp}\subset \Gamma_p$, which acts on $V_{\mfk}\otimes_{\Qp}\mfkb$, $H^{\cris}_1(\mathcal{A}_{x}^{\mathcal{P}}/\mfk)\otimes_{\Qp}\mfkb$ and $J_b(\mfkb)= G(\mfk\otimes\mfkb)^{\langle\sigma\rangle}$ all via its actions on the tensor factors $\mfkb$ (cf. (\ref{eq:Xi2})).
Further, as $h_p\in J_b$, we have the relation 
\[ {\eta_p'}^{-1}\circ\phi_{\mathcal{A}_{x}^{\mathcal{P}}}\circ\eta_p'=h_p\circ i_p(\phi_{\mathcal{A}_{x}})\circ h_p^{-1}=h_p (b\sigma) h_p^{-1}
=b\sigma.\] Therefore, the map
\begin{equation} \label{eq:i_p'}
i_p':=\Int({\eta_p'}^{-1})=\Int(h_p)\circ i_p\circ \Int(\theta_{\tilde{\omega}}): (I_{\mathbf{i}_{\omega}(x)})_{\Qp} \hookrightarrow J_{b} 
\end{equation}
qualifies as a $\Qp$-embedding attached to the isogeny class $\mathscr{I}^{\beta}$ in Theorem \ref{thm:Kisin17_Cor.1.4.13,Prop.2.1.3,Cor.2.2.5} with $b'=b$; this proves (ii).

To prove statement (iv), we use this the fact (\ref{eq:Kisin17_Lem.4.1.2}) that $i_T(\theta_{\tilde{\omega}}\circ {}^{\tau}\theta_{\tilde{\omega}}^{-1})=\tilde{\omega}^{-1}\tau(\tilde{\omega})$ for all $\tau\in\Gamma$. When $(b,i_p,i^p)$ adapted to $(i_T, h)$, at any finite place $v\neq p$, we may take $h_v=\tilde{\omega}$. Then, from (\ref{eq:new_embedding_T-beta}), one sees that
\[ i_v'|_{T^{\tilde{\beta}}(\Qv)}= i_{T^{\tilde{\beta}}}. \] 
At $p$, recall that the injection $\Xi:J_b(\mfkb)\subset G(\mfk\otimes\mfkb)\rightarrow G(\mfkb)$ (\ref{eq:Xi2}) induced by the canonical $\mfk$-algebra homomorphism $\mfk\otimes R\rightarrow R:x\otimes y\mapsto xy$ in the case $R=\mfkb$ is $W_{\Qp}$-equivariant for the standard action on $J_b(\mfkb)$ and the twisted action on $G(\mfkb)$: $z\mapsto b_{\tau}\tau(z) b_{\tau}^{-1}$. 
When $(b,i_p,i^p)$ is adapted to $(i_T, h)$, we have $\bar{h}_p\cdot b_{\tau}\tau(\bar{h}_p)b_{\tau}^{-1}=\tilde{\omega}^{-1}\tau(\tilde{\omega})$ for $\bar{h}_p:=\Xi(h_p)$ and all $\tau\in W_{\Qp}$,  which implies that $g_p:=\bar{h}_p\cdot  \tilde{\omega}^{-1} \in G(\mfk)$ as $b_{\tau}=1$ for $\tau\in \Gal(\mfkb/\mfk)$. For $i_p'$ (\ref{eq:i_p'}), it follows from (\ref{eq:new_embedding_T-beta}) that
 \[ \Xi\circ i_p'|_{T^{\tilde{\beta}}}=\Int(g_p)\circ i_{T^{\tilde{\beta}}}, \] 
which shows that
\[ \Int(g_p^{-1})\circ i_p'= \Int(g_p^{-1}h_p)\circ i_p\circ \xi^{-1}: (I_{\mathbf{i}_{\omega}(x)})_{\Qp} \hookrightarrow J_b\isom J_{b'} \] is adapted to $i_{T^{\tilde{\beta}}}$, where $b':=g_p^{-1}b\sigma(g_p)$. 
This proves (iv).
\end{proof}

%%%%%%%%%%%%%%%%%%%%
\begin{lem} \label{lem:LR87_Lem.5.18}
Suppose given $b\in G(\mfk)$ whose Newton quasi-cocharacter $\nu_b:=\nu_G(b)\in\Hom_{\mfk}(\mathbb{D},G)$ is $\Qp$-rational.
For $\tilde{J}_b:=\rho^{-1}(J_b)$ for the morphism $\rho:G^{\uc}\rightarrow G$, we have
\[H^1(\Qp,\tilde{J}_{b})=\{1\}.\]
\end{lem}

\begin{proof}
This is proved in the discussion proceeding Lemma 5.18 of \cite{LR87} when $G^{\der}=G^{\uc}$, whose proof carries over here. 
Since $\tilde{J}_{b}$ is an inner twist of the connected reductive group $\tilde{J}:=Z_{G^{\uc}}(\nu_{b})$ over the $p$-adic local field $\Qp$ \cite[3.3]{Kottwitz97}, the claim is equivalent to vanishing of $H^1(\Qp,\tilde{J})$. Let $A$ be the connected center of $J\cap G^{\der}_{\Qp}$ (a split $\Qp$-torus) and $T\subset G^{\der}_{\Qp}$ the centralizer of a maximal $\Qp$-split torus of $G^{\der}_{\Qp}$ containing $A$. The simple roots $\alpha_1,\cdots,\alpha_s$ of $T$ in $G$ can be divided into two sets, the roots $\alpha_1,\cdots,\alpha_r$ of $T$ in $J$ and the others $\alpha_{r+1},\cdots,\alpha_s$. Let $\{\omega_1,\cdots,\omega_s\}$ be the corresponding fundamental weights which form a basis of $X^{\ast}(\tilde{T})$ for $\tilde{T}=\rho^{-1}(T)$, and $R\subset \tilde{T}$ be the kernel of $\{\omega_1,\cdots,\omega_r\}$. Note that $R$ is an induced torus since the basis $\{\omega_{r+1},\cdots,\omega_s\}$ of $X^{\ast}(R)$ is permuted under the canonical action of $\Gal(\Qpb/\Qp)$ which equals the naive action as $G^{\der}_{\Qp}$ is quasi-split. Now, the claim follows from the fact that $\tilde{J}$ is the semi-direct product of a simply connected semi-simple group over $\Qp$ (i.e. $J^{\uc}$) and $R$.
\end{proof}

%%%%%%%%%%%%%%%%%%%%
\begin{lem} \label{lem:b(gamma_0;delta;g_p)_{tau}}
If $(\gamma_0,\delta,g_p)\in G(\Qp)\times G(L_n)\times G(\mfkb)$ satisfies that $g_p \gamma_0 g_p^{-1}=\Nm_n\delta$, then
\begin{equation} \label{eq:b(gamma_0;delta;g_p)_{tau}}
\tau\mapsto b(\gamma_0;\delta;g_p)_{\tau}:=g_p^{-1}\cdot \Nm_{i(\tau)}\delta \cdot \tau(g_p)
\end{equation}
is a cocycle in $Z^1(W_{\Qp},G_{\gamma_0}(\mfkb))$,
where $i:W_{\Qp}\ra \langle\sigma\rangle=\Z$ is the canonical surjection and $\Nm_{i(\tau)}\delta$ is defined as in (\ref{eq:b_{tau}}). If $g_p\in G(\mfk)$, this cocycle is the same as the one (\ref{eq:b_{tau}}) defined by $b=b(\gamma_0;\delta;g_p):=g_p^{-1}\delta\sigma(g_p)$.
\end{lem}

\begin{proof}
It is immediate that the cochain (\ref{eq:b(gamma_0;delta;g_p)_{tau}}) is a cocycle (valued in $G(\mfkb)$).
That it is valued in $G_{\gamma_0}(\mfkb)$ follows from the following equality obtained by applying $\tau$ to both sides of $g_p \gamma_0 g_p^{-1}=\Nm_n\delta$:
\begin{align*} \tau(g_p)\gamma_0\tau(g_p)^{-1} = & (\Nm_{i(\tau)}\delta)^{-1}\Nm_{i(\tau)+n}\delta =(\Nm_{i(\tau)}\delta)^{-1} \cdot \Nm_n\delta \cdot \sigma^n(\Nm_{i(\tau)}\delta) \\ =& (\Nm_{i(\tau)}\delta)^{-1} g_p \gamma_0 g_p^{-1} \Nm_{i(\tau)}\delta.
\end{align*}
The last claim is clear and left to readers.
\end{proof}

%%%%%%%%%%%%%%%%%%%%
\begin{lem} \label{lem:existence_of_nice_tori}
Let $(\mathscr{I},\epsilon)$ an admissible K-pair and $I=I_{\mathscr{I},\epsilon}^{\mathrm{o}}$. Then, for any finite set $S$ of places $v$ of $\Q$ containing those where $I_{\Qv}$ is not quasi-split and the infinite place $\infty$, and for any finite subset $\Xi\subset H^1(\Q,I)$, there exists a maximal $\Q$-torus of $I_{\mathscr{I},\epsilon}^{\mathrm{o}}$ such that
$\Xi$ lies in the image of $H^1(\Q,T)\ra H^1(\Q,I)$ and $T_{\Qv}$ is fundamental in $I_{\Qv}$ for every $v\in S$.
\end{lem}

\begin{proof}
%. We claim that such $T$ works for the purpose at hand. 
In fact, this is what the proof of \cite{Borovoi98}, Theorem 5.11 establishes. In more detail, we take the same finite set $\Sigma$ of places of $\Q$ in that proof ($\Sigma$ is chosen by certain local properties of $I$ only). Then in the remaining part of \textit{loc. cit.}, it is shown that for any maximal $\Q$-torus $T'$ of $I^{\uc}$ presented in ibid., Lemma 5.6.5, $T:=T'\cdot Z(I)^{\mathrm{o}}$ works for the purpose at hand; in this part of the proof the only property required of $T$ is that it is fundamental at all archimedean places. Hence, it suffices to check that we can take a maximal $\Q$-torus $T'$ of $G^{\uc}$ as in that lemma to further satisfy our condition that $T_{\Qv}$ is fundamental in $I_{\Qv}$ for every $v\in S$. But, the only condition used about $T'$ to prove the properties in the statement of that lemma is that $T'_{\Qv}$ is fundamental in $I^{\uc}_{\Qv}$ for every $v\in \Sigma$ ($\Sigma$ is enlarged include all archimedean places). Therefore, any maximal $\Q$-torus of $I$ which is fundamental in $I_{\Qv}$ for every $v\in S\cup\Sigma$ serves our purpose, and it is well-known that such $\Q$-torus exists (for each $v\in S\cup\Sigma$, choose a fundamental maximal torus $T_v$ of $(I_{\mathscr{I},\epsilon})_{\Qv}$ which exists for $v\neq\infty$ by \cite[Thm.6.21]{PR94} and for $v=\infty$ by Theorem \ref{thm:Kisin17_Cor.1.4.13,Prop.2.1.3,Cor.2.2.5}. Then, there exists a maximal $\Q$-torus $T$ of $I_{\mathscr{I},\epsilon}$ such that for each $v\in S\cup\Sigma$, $T_{\Qv}$ is conjugate to $T_v$ by \cite[$\S$7.1,Cor.3]{PR94}).
\end{proof}

Recall (Subsection \ref{subsec:Kottwitz_invariant}) that a stable Kottwitz triple $(\gamma_0;\gamma,\delta)$ together with a choice of elements $(g_v)_v\in G(\bar{\A}_f^p)\times G(\mfk)$ satisfying (\ref{eq:stable_g_l}) and Definition \ref{defn:stable_Kottwitz_triple} (iii$'$) (with $g_p=c$ there) gives rise to an adelic class \[(\alpha_v(\gamma_0;\gamma,\delta;g_v))_v\in  \bigoplus_v \pi_1(I_0)_{\Gamma_v} \] 
(whose finite part $(\alpha_v)_{v\neq\infty}$ measures the difference between $\gamma_0$ and $(\gamma,\Nm_n\delta)$, $n$ being the level of the triple) 
and that the triple $(\gamma_0;\gamma,\delta)$ having trivial Kottwitz invariant means that one can find such a tuple $(g_v)_v$ such that the invariant $\alpha(\gamma_0;\gamma,\delta;(g_v)_v)$ (\ref{eq:Kottwitz_invariant}) defined by this adelic class vanishes.

%%%%%%%%%%%%%%%%%%%%
\begin{cor} \label{cor:K-triple_of_twisted_K-pair}
Let $(\mathscr{I},\epsilon)$ be an admissible K-pair and $(\gamma_0;\gamma,\delta)$ the associated stable Kottwitz triple (\ref{defn:stable_K-triple_attached_to_admissible_K-pair}) with trivial Kottwitz invariant and $(g_{v})_v\in G(\bar{\A}_f^p)\times G(\mfk)$ a tuple of elements satisfying (\ref{eq:stable_g_l}) and Definition \ref{defn:stable_Kottwitz_triple} (iii$'$) giving trivial Kottwitz invariant; set $I_0:=G_{\gamma_0}^{\mathrm{o}}$.

Then, for any $\beta\in\Sha^{\infty}_G(\Q,I_{\mathscr{I},\epsilon}^{\mathrm{o}})$, the stable Kottwitz triple $(\gamma_0;\gamma',\delta')$ attached to $(\mathscr{I}^{\beta},\epsilon^{\beta})$ has $\gamma_0$ as the rational component, and for this rational component $\gamma_0$, there exists a tuple of elements $(g_{v}')_v\in G(\bar{\A}_f^p)\times G(\mfk)$ satisfying (\ref{eq:stable_g_l}), Definition \ref{defn:stable_Kottwitz_triple} (iii$'$) whose associated Kottwitz invariant is trivial, such that at every finite place $v$ of $\Q$, we have an equality in $\ker[H^1(\Qv,I_0)\rightarrow H^1(\Qv,G)]$:
\[ \alpha_v(\gamma_0;\gamma'_v,\delta';g_{v}') -\alpha_v(\gamma_0;\gamma_{v},\delta;g_{v})=\beta(v). \]
Here, $\alpha_v(\gamma_0;\gamma_{v},\delta;g_{v})$ is $\alpha_l(\gamma_0;\gamma_{l};g_{l})$ (\ref{eq:alpha_l}) for $v=l\neq p$ and $\kappa_{I_0}([b(\gamma_0;\delta;g_p)]_{I_0})$ (\ref{eq:alpha_p}) for $v=p$.
\end{cor}

\begin{rem} \label{rem:kappa(b)-kappa(b')_is_torsion}
(1) As $I_{\mathscr{I},\epsilon}^{\mathrm{o}}$ is an inner form of $I_0$ (Corollary \ref{cor:Tate_thm2}), there exist canonical isomorphisms
$H^1(\Qv,I_{\mathscr{I},\epsilon}^{\mathrm{o}}) \cong \pi_1(I_{\mathscr{I},\epsilon}^{\mathrm{o}})_{\Gamma_v,\mathrm{tor}} \cong \pi_1(I_0)_{\Gamma_v,\mathrm{tor}} \cong H^1(\Qv,I_0)$  (\ref{eq:H^1(Qv,G)=H^1(Qv,G_1)}).

(2) We note that if $b:=b(\gamma_0;\delta;g_p)$ and $b':=b(\gamma_0;\delta';g_p')$, although $\kappa_{I_0}([b]_{I_0})$ and $\kappa_{I_0}([b']_{I_0})$ are each elements of $\pi_1(I_0)_{\Gamma_p}$, their difference lies in $\pi_1(I_0)_{\Gamma_p,\mathrm{tor}}$, in fact in
 \[ \ker[\pi_1(I_0)_{\Gamma_p,\mathrm{tor}} \rightarrow \pi_1(G)_{\Gamma_p,\mathrm{tor}}]=\ker[H^1(\Qp,I_0)\rightarrow H^1(\Qp,G)]. \] 
Indeed, the equality of the Newton quasi-cocharacters $\nu_{I_0}(b)=\frac{1}{n}\nu_{I_0}(\gamma_0)=\nu_{I_0}(b')$ (Proposition \ref{prop:psi_p}) implies \cite[3.5]{Kottwitz97} that $[b]_{I_0}$ and $[b']_{I_0}$ are in the image of $j_b^{I_0}:H^1(\Qp,J_b^{I_0})\rightarrow B(I_0)_{basic}$. This image maps onto $\pi_1(I_0)_{\Gamma_v,\mathrm{tor}}+\kappa_{I_0}([b])$ under $\kappa_{I_0}$  (Lemma \ref{lem:Xi_{b'}_inner-twisting} (2)). Also, we have $\kappa_{G}([b]_{G})=\kappa_{G}([b']_{G})$ by condition $\ast(\delta)$.
\end{rem}

\begin{proof}
Let $(\mathscr{I}^{\beta},\epsilon^{\beta})$ be the admissible K-pair obtained from $(\mathscr{I},\epsilon)$ by twisting by $\beta\in H^1(\Q,I_{\mathscr{I},\epsilon})$ (Theorem \ref{thm:Kisin17_Prop.4.4.8} (3)) and $S$ the finite set of places $v$ of $\Q$ where either of $(I_{\mathscr{I}^{\beta},\epsilon^{\beta}}^{\mathrm{o}})_{\Qv}$, $(I_{\mathscr{I}^{\beta}})_{\Qv}$, $(I_{\mathscr{I},\epsilon}^{\mathrm{o}})_{\Qv}$ or $(I_{\mathscr{I}})_{\Qv}$ is not quasi-split.
For $\Xi=\{\beta\}$ and $S$, fix a maximal $\Q$-torus $T\subset I_{\mathscr{I},\epsilon}^{\mathrm{o}}$ as in Lemma \ref{lem:existence_of_nice_tori}, which thus admits a class $\tilde{\beta}\in H^1(\Q,T)$ mapping to $\beta$.
Then, by construction (\ref{eq:(epsilon,b,c)->delta2}), Definition \ref{defn:stable_K-triple_attached_to_admissible_K-pair}, we may assume that $\gamma_0=i_T(\epsilon)$ for an embedding $i_T:T\hookrightarrow G$ as in Theorem \ref{thm:Kisin17_Cor.1.4.13,Prop.2.1.3,Cor.2.2.5} (3).
Take an inner twisting $\xi:(I_{\mathscr{I}})_{\Qb}\isom (I_{\mathscr{I}^{\beta}})_{\Qb}$ and a cocycle representative $\tilde{a}_{\tau}=\tilde{\omega}^{-1}\tau(\tilde{\omega})\in Z^1(\Q,T)$ of $\tilde{\beta}$ as in (\ref{eq:inner-twist:xi2}), so that the embedding $i_{T^{\tilde{\beta}}}:T^{\tilde{\beta}} \stackrel{\xi^{-1}}{\rightarrow} T\stackrel{i_T}{\hookrightarrow} G \stackrel{\Int(\tilde{\omega})}{\rightarrow} G$ is one attached to the maximal torus $T^{\tilde{\beta}}=\xi(T)\subset I_{\mathscr{I}^{\beta}}$ as in Theorem \ref{thm:Kisin17_Prop.4.4.8} (iv), (cf. (\ref{eq:new_embedding_T-beta})).
Since $I_{\mathscr{I}^{\beta},\epsilon^{\beta}}^{\mathrm{o}}$ is an inner twist of $I_{\mathscr{I},\epsilon}^{\mathrm{o}}$, $T^{\tilde{\beta}}$ is also nice (in the sense of Definition \ref{defn:nice_tori}) for $(\mathscr{I}^{\beta},\epsilon^{\beta})$ for every $v\in S$). Hence, taking $\epsilon^{\beta}:=\xi(\epsilon)$ (\ref{eq:epsilon^{beta}=xi(epsilon)}), we obtain a rational component $i_{T^{\tilde{\beta}}}(\epsilon^{\beta})=\Int(\tilde{\omega})(\gamma_0)$ of the Kottwitz triple $(\gamma_0',\gamma',\delta')$ attached to $(\mathscr{I}^{\beta},\epsilon^{\beta})$, and
since $\Int(\tilde{\omega})$ is a transfer of maximal tori, we may assume that $\gamma_0'=\gamma_0$. 
Below, we will denote by $a_{\tau}$ the image of $\tilde{a}_{\tau}$ in $Z^1(\Q,H)$ for any of $H=I_{\mathscr{I},\epsilon}^{\mathrm{o}}$, $I_{\mathscr{I},\epsilon}$ and $I_{\mathscr{I}}$.

At $l\neq p$, by construction (\ref{eq:(epsilon,b,c)->delta2}), we have $\gamma_{l}=i_l(\epsilon)$ for some embedding $i_l:(I_{\mathscr{I}})_{\Ql}\hookrightarrow G_{\Ql}$ as in Theorem \ref{thm:Kisin17_Cor.1.4.13,Prop.2.1.3,Cor.2.2.5}. The element $g_{l}\in G(\Qlb)$ satisfies that $g_{l}\gamma_0g_{l}^{-1}=\gamma_{l}$ and that $g_{l}^{-1} \tau(g_l)\in I_0(\Qlb)$ for every $\tau\in\Gamma_l$; then, $\alpha_l(\gamma_0;\gamma_{l};g_{l})=[g_{l}^{-1}\tau(g_l)]\in H^1(\Ql,I_0)$.
We choose $h_l\in G(\Qlb)$ trivializing the cocycle $i_l(a(l))$:
\[ i_l(a(l))=h_l^{-1}\tau(h_l)\ (\tau\in\Gamma_l),\]
so that $i_l':=\Int(h_l)\circ i_l\circ \xi^{-1}$ is an embedding $(I_{\mathscr{I}^{\beta}})_{\Ql}\hookrightarrow G_{\Ql}$ attached to $\mathscr{I}^{\beta}$ (Theorem \ref{thm:Kisin17_Prop.4.4.8} (2)): it was observed in the proof of this theorem that $i_l(\beta(l))\in H^1(\Ql,G)$ is trivial. Then, for $\gamma_{l}':=i_l'(\epsilon^{\beta})=\Int(h_{l})(\gamma_{l})$ (which gives the $l$-component of a Kottwitz triple attached to $(\mathscr{I}^{\beta},\epsilon^{\beta})$), we have the relation 
\[\gamma'_l=\Int(g_{l}')(\gamma_0)\ \in G(\Ql)\] 
where $g_{l}':=h_lg_{l}$. 
Since one has 
\begin{equation*} \label{eq:cocycles_at_l}
\quad g_{l}'^{-1}\tau(g_l')=j_l(a(l)) \cdot g_{l}^{-1}\tau(g_l),
\end{equation*}
for $j_l:=\Int(g_{l}^{-1})\circ i_l:(I_{\mathscr{I},\epsilon})_{\Qlb}\isom (I_0)_{\Qlb}$ which is an inner twist 
such that $j_l\circ{}^{\tau}j_l^{-1}=\Int(g_{l}^{-1}\tau(g_l))$ (Corollary \ref{cor:Tate_thm2}),
the cohomology classes $\alpha_l(\gamma_0;\gamma_{l};g_{l})=[g_{l}^{-1}\tau(g_l)]$, $\alpha_l(\gamma_0;\gamma_{l}';g_{l}')=[g_{l}'^{-1}\tau(g_l')]$ in $H^1(\Ql,I_0)$ satisfy the relation:
\begin{equation} \label{eq:difference_of_alpha_l's}
\alpha_l(\gamma_0;\gamma'_l;g_{l}')=\alpha_l(\gamma_0;\gamma_{l};g_{l})+\beta(l)
\end{equation}
under the canonical isomorphism $H^1(\Ql,I_{\mathscr{I},\epsilon}^{\mathrm{o}})=H^1(\Ql,I_0)$ (\ref{eq:H^1(Qv,G)=H^1(Qv,G_1)}) (apply \cite{Borovoi98}, 3.15 and Lemma 3.15.1 to $(I_{\mathscr{I},\epsilon})_{\mathbf{ab}}={}_{\psi}(I_0)_{\mathbf{ab}}$, $\psi:=g_{l}^{-1}\tau(g_l)$, $\psi':=a(l)$; cf. Lemma \ref{lem:beta_v=beta_v'+zeta} below.)

At $v=p$, the chosen element $g_p\in G(\mfk)$ satisfies that $g_p\gamma_0g_p^{-1}=\Nm_n\delta$ and 
\begin{equation} \label{eq:b=b(gamma_0;delta;g_p)}
b:=b(\gamma_0;\delta;g_p):=g_p^{-1}\delta\sigma(g_p)\in I_0(\mfk),
\end{equation}
which gives a cocycle $b_{\tau}=b(\gamma_0;\delta;g_p)_{\tau}\in Z^1(W_{\Qp},I_0)$ (\ref{eq:b_{tau}}).
On the other hand, by construction (\ref{eq:(epsilon,b,c)->delta2}), $\delta\in G(L_n)$ is defined as 
\[ \delta=c b_1\sigma(c^{-1}) \] 
from some embedding $i_p:(I_{\mathscr{I}})_{\Qp}\rightarrow J_{b_1}\ (b_1\in G(\mfk))$ attached to $\mathscr{I}$ in Theorem \ref{thm:Kisin17_Cor.1.4.13,Prop.2.1.3,Cor.2.2.5} and some $c\in G(\mfk)$ such that $c(i_p(\epsilon)^{-1}(b_1\sigma)^n)c^{-1}=\sigma^n$, i.e. $c i_p(\epsilon) c^{-1}=\Nm_n(\delta)$. 
We may assume that $(b_1,i_p)$ is adapted to $i_T:T\hookrightarrow G$ so that $\gamma_0=i_p(\epsilon)$. 
Then, by conjugation (of $(b_1\sigma,i_p)$), we may further assume that $b_1=b$ and $i_p$ is an embedding $(I_{\mathscr{I}})_{\Qp}\rightarrow J_{b}$ and also that $\gamma_0=i_p(\epsilon)$ again; we obtain the same $\delta$ by taking $c=g_p$.
Indeed, as $b_1=c^{-1}\delta\sigma(c)=d^{-1}b\sigma(d)$ for $d:=g_p^{-1}c$, we conjugate $(b_1\sigma,i_p)$ by $\Int(d)$, which gives the relation $g_p i_p(\epsilon) g_p^{-1}=\Nm_n(\delta)$ for the new $i_p$ (as $cd^{-1}=g_p$), so still $\gamma_0=i_p(\epsilon)$. 

Now, we choose $h_p\in J_b(\Qpb)$ satisfying $i_p(a_{\tau})=h_p^{-1}\tau(h_p)$ so that $i_p':=\Int(h_p)\circ i_p \circ \xi^{-1}$ (\ref{eq:i_v'}) is an embedding $(I_{\mathscr{I}^{\beta}})_{\Qp}\hookrightarrow J_b$ attached to $\mathscr{I}^{\beta}$ (Theorem \ref{thm:Kisin17_Prop.4.4.8} (2)): note that we are using the same $b$. 
By admissibility of the K-pair $(\mathscr{I}^{\beta},\epsilon^{\beta})$, there exists $c'\in G(\mfk)$  such that $c'(i_p'(\epsilon^{\beta})^{-1}(b\sigma)^n)c'^{-1}=\sigma^n$, and we define $\delta':=c' b\sigma(c'^{-1})$ so that $c' i_p'(\epsilon^{\beta}) c'^{-1}=\Nm_n(\delta')$.
Then, since
\[ i_p'(\epsilon^{\beta})=i_p'\circ\xi(\epsilon)= \bar{h}_p \gamma_0 \bar{h}_p^{-1},\]
where $\bar{h}_p=\Xi(h_p)$ for the injection $\Xi:J_b(\mfkb)=G(\mfk\otimes\bar{\mfk})^{\langle\sigma\rangle} \hookrightarrow G(\bar{\mfk})$ (\ref{eq:Xi2}) induced by the $\mfk$-algebra homomorphism $\mfk\otimes\bar{\mfk}\rightarrow \bar{\mfk}:l\otimes x\mapsto lx$,
we obtain
\[ x\gamma_0x^{-1}=\Nm_n(\delta') \] for $x:=c'\bar{h}_p$. 
The associated cocycle $b(\gamma_0;\delta';x)_{\tau}$ (\ref{eq:b(gamma_0;delta;g_p)_{tau}}) in $Z^1(W_{\Qp},G_{\gamma_0}(\mfkb))$ becomes
\begin{align*} \label{eq:a_{tau}b_{tau}=b_{tau}'}
b(\gamma_0;\delta';x)_{\tau}&:=x^{-1}\Nm_{i(\tau)}\delta'\tau(x) \\
&= \bar{h}_p^{-1}\cdot b_{\tau}\cdot \tau(\bar{h}_p) \nonumber \\
&= (\bar{h}_p^{-1}\cdot b_{\tau}\tau(\bar{h}_p) b_{\tau}^{-1}) \cdot b_{\tau} \nonumber \\
&= \Xi(i_p(a_{\tau}))\cdot b_{\tau}, \nonumber
\end{align*}
which thus lies in $Z^1(W_{\Qp},I_0(\mfkb))$ as $a_{\tau}\in Z^1(\Qp,I_{\mathscr{I},\epsilon}^{\mathrm{o}})$ ($\Xi$ is $W_{\Qp}$-equivariant for the standard action on $J_b(\mfkb)$ and the twisted action on $G(\mfkb)$). Since $H^1(\Gal(\mfkb/\mfk),I_0)=\{1\}$, we can find $g\in I_0(\mfkb)$ such that $g_p':=xg\in G(\mfk)$, so the cocycle $b(\gamma_0;\delta';g_p')_{\tau}$ (\ref{eq:b(gamma_0;delta;g_p)_{tau}}), which is cohomologous to $b(\gamma_0;\delta';x)_{\tau}$, is valued in $I_0(\mfk)$ and equal to $b'_{\tau}$ (\ref{eq:b_{tau}}) for $b':=g_p'^{-1}\delta'\sigma(g_p')$.
So, we have $[b(\gamma_0;\delta';g_p')]_{I_0}= j_b^{I_0}([i_p(a_{\tau})]_{I_0})$
for the injection 
\[ j_b^{I_0}: H^1(\Qp,J_b^{I_0})\hookrightarrow B(I_0)=H^1(W_{\Qp},I_0) : [x_{\tau}] \mapsto [\Xi(x_{\tau})\cdot b_{\tau}], \]
thus the equality
\[ \kappa_{I_0}([b(\gamma_0;\delta';g_p')])=\kappa\circ j_b^{I_0}([i_p(a_{\tau})])=[i_p(a_{\tau})]+\kappa_{I_0}([b(\gamma_0;\delta;g_p)_{\tau}]) \]
by Lemma \ref{lem:Xi_{b'}_inner-twisting} (2), which proves the claim at $p$.

Now, it is immediate from the formulae proved that if the tuple $(g_v)_v$ satisfies the condition that the associated Kottwitz invariant $\alpha(\gamma_0;\gamma,\delta;(g_v)_v)$ (\ref{eq:Kottwitz_invariant}) is trivial, so does the tuple $(g_v')_v$ that we have constructed.
\end{proof}

%%%%%%%%%%%%%%%%%%%%
%%%%%%%%%%%%%%%%%%%%
\begin{thm} \label{thm:equiv_K-triples}
Let $(\mathscr{I},\epsilon)$ and $(\mathscr{I}',\epsilon')$ be two admissible K-pairs whose associated stable Kottwitz triples $(\gamma_0;\gamma,\delta)$, $(\gamma_0';\gamma',\delta')$ are (stably) equivalent.

Then, for any nice maximal $\Q$-torus $T$ of $I_{\mathscr{I},\epsilon}^{\mathrm{o}}$ and any special Shimura subdatum $(i_T(T),h\in X\cap\Hom(\dS,i_T(T)_{\R}))$ as in Theorem \ref{thm:Kisin17_Cor.1.4.13,Prop.2.1.3,Cor.2.2.5} (3), there exist an embedding $\rho:T\hookrightarrow I_{\mathscr{I}',\epsilon'}^{\mathrm{o}}$ and $g\in G(\Qb)$ such that 
the pair
\[(i_T',h'):=\Int(g)(i_T,h)\] is a pair attached likewise to $\rho:T\hookrightarrow I_{\mathscr{I}',\epsilon'}^{\mathrm{o}}$.

So, the cohomology class $\tilde{\beta}=[i_T^{-1}(g^{-1}\tau(g))]\in H^1(\Q,T)$ lies in $\Sha^{\infty}_G(\Q,i_T(T))$ and
the K-pair $(\mathscr{I}',\epsilon')$ is the twist $(\mathscr{I}^{\beta},\epsilon^{\beta})$ of $(\mathscr{I},\epsilon)$ by the image $\beta$ of $\tilde{\beta}$ in $H^1(\Q,I_{\mathscr{I},\epsilon})$ (Theorem \ref{thm:Kisin17_Prop.4.4.8} (3)).
\end{thm}

\begin{proof}
First, we observe that the equivalence of the Kottwiz triples attached to $(\mathscr{I},\epsilon)$, $(\mathscr{I}',\epsilon')$ implies that the $\Q$-groups $I_{\mathscr{I},\epsilon}^{\mathrm{o}}$, $I_{\mathscr{I}',\epsilon'}^{\mathrm{o}}$ are isomorphic locally at every place (in fact, one can prove that they are even globally isomorphic, but we will not  need this stronger fact here). This follows from Corollary \ref{cor:Tate_thm2}. In particular, we note that $T'$ is also nice for $(\mathscr{I}',\epsilon')$.

The $\Q$-torus $i_T(T)$ in the special Shimura subdatum $(i_T(T),h)$ is the image of a $\Q$-embedding $i_{T}:T\hookrightarrow G$. Choose $(b\in i_T(T)(\mfk),\ i_p:(I_{\mathscr{I}})_{\Qp}\hra J_b,\ i^p:(I_{\mathscr{I}})_{\A_f^p}\hra G_{\A_f^p})$ which is adapted to $i_T$ (Theorem \ref{thm:Kisin17_Cor.1.4.13,Prop.2.1.3,Cor.2.2.5} (3)). We may assume that the triple $(\gamma_0;\gamma,\delta)$ is given by (\ref{eq:K-triple_attached_to_admissible_K-pair}) from theses data and some $c$ satisfying (\ref{eq:(epsilon,b,c)->delta2}): $\gamma_l=i_l(\epsilon)=i_T(\epsilon)=\gamma_0\in G(\Ql)$, $i_p(\epsilon)=i_T(\epsilon)=\gamma_0\in G(\Qp)$. 
Recall (see the proof of Theorem \ref{thm:Kisin17_Cor.1.4.13,Prop.2.1.3,Cor.2.2.5}) that $i_T$ is given by a lifting of a point $x\in \mathscr{I}$ as being endowed with the action of $T\subset I_x=I_{\mathscr{I}}$ to a special point $\tilde{x}$ of $\sS_{\mathbf{K}_p}$. Such lifting is determined by any cocharacter $\mu_T$ of $T$ satisfying certain conditions \cite[Lem.2.2.2]{Kisin17} and is characterized by that the special Shimura datum $h\in X\cap \Hom(\dS,i_{T}(T)_{\R})$ defining the CM point $\tilde{x}$ has the Hodge cocharacter $\mu_{h}=i_T\circ\sigma_p(\mu_T)$, where $\sigma_p(\mu_T)$ is the base-change to $\C$ via $\sigma_p:\Qpb\hookrightarrow \C$ of $\mu_T\in \Hom_{\Qpb}(\Gm,T)$. 

For $(\mathscr{I}',\epsilon')$, we may also assume that $(\gamma_0';\gamma',\delta')$ is given by (\ref{eq:K-triple_attached_to_admissible_K-pair}) from similar data: some (nice) maximal $\Q$-torus $T'$ of $I_{\mathscr{I}'}$ containing $\epsilon'$, a $\Q$-embedding $i_{T'}:T'\hookrightarrow G$ which is produced by a CM lifting $\tilde{x}'$ of a point $x'$ in $\mathscr{I}'$, and a tuple  $(b',\{i_v'\}_v)$ as in Theorem \ref{thm:Kisin17_Cor.1.4.13,Prop.2.1.3,Cor.2.2.5} (3) that is adapted to $i_{T'}$. Again, we fix an identification $I_{\mathscr{I}'}=I_{x'}$. 

\textbf{Step 1.}  
Let $\varphi:(I_{x,\epsilon}^{\mathrm{o}})_{\Qb}\isom (I_0)_{\Qb}$, $\varphi':(I_{x',\epsilon'}^{\mathrm{o}})_{\Qb}\isom (I_0')_{\Qb}$ be the inner twistings in Corollary \ref{cor:Tate_thm2} (3), constructed by the choices of $(i_T,b,i_v)$, $(i_{T'},b',i_v')$, where $I_0:=G_{\gamma_0}^{\mathrm{o}}$, $I_0'=G_{\gamma_0'}^{\mathrm{o}}$; in particular, $\varphi|_T=i_T$.
Then, we claim that there exists $g_1\in G(\Qb)$ such that $\Int(g_1)$ induces an inner twisting $(I_0)_{\Qb}\isom (I_0')_{\Qb}$ and the composite inner twisting 
\begin{equation} \label{eq:rho_vectorsp_isom} 
\rho:=\varphi'^{-1}\circ\Int(g_1)\circ\varphi\ :\ (I_{x,\epsilon}^{\mathrm{o}})_{\Qb}\isom (I_0)_{\Qb}\isom (I_0')_{\Qb}\isom (I_{x',\epsilon'}^{\mathrm{o}})_{\Qb} \end{equation} 
becomes $\Q$-rational when restricted to $T$.

By the assumption that $\gamma_0$ and $\gamma_0'$ are stably conjugate, there exists $g_0\in G(\Qb)$ such that $\Int(g_0)$ induces an inner twisting $(I_0)_{\Qb}\isom (I_0')_{\Qb}$, namely that $\gamma_0'=\Int(g_0)(\gamma_0)$ and $\zeta:=(\tau\mapsto g_0^{-1}\tau(g_0))$ is a $I_0$-valued cocycle on $\Gamma_{\Q}$. 
We will prove that some conjugate (under $I_{x',\epsilon'}^{\mathrm{o}}(\Qb)$) of $\rho_0:=\varphi'^{-1}\circ\Int(g_0)\circ\varphi$ becomes $\Q$-rational upon restriction to $T$. In this case that $\rho$ (\ref{eq:rho_vectorsp_isom}) is of the form $\rho=\Int(h)\circ \rho_0$ for some $h\in I_{x',\epsilon'}^{\mathrm{o}}(\Qb)$, one says that the maximal torus $T$ of the connected reductive group $I_{x,\epsilon}^{\mathrm{o}}$ \emph{transfers to} the inner form $I_{x',\epsilon'}^{\mathrm{o}}$ for (the conjugacy class of) the inner twisting $\rho_0$, and the resulting $\Q$-embedding $\rho|_T:T\hookrightarrow I_{x',\epsilon'}^{\mathrm{o}}$ of $T$ is called an \emph{admissible embedding} relative to (the conjugacy class of) the inner twisting $\rho_0$, cf. \cite[$\S$9]{Kottwitz84a}.

Since $T_{\Qv}$ is elliptic in $(I_{x,\epsilon}^{\mathrm{o}})_{\Qv}$ at least at one place $v$ (i.e. at $v=\infty$), according to \cite[Lem.5.6]{LR87} (for which one does not need $G^{\ast}$ to be quasi-split; see also the discussion in $\S$9 of \cite{Kottwitz84a}, particularly 9.4.1, 9.5), it suffices to check that $T$ transfers to $I_{x',\epsilon'}^{\mathrm{o}}$ for $\rho$ locally everywhere. Let $S$ be the set of places $v$ of $\Q$ where $(I_{x',\epsilon'}^{\mathrm{o}})_{\Qv}\approx (I_{x,\epsilon}^{\mathrm{o}})_{\Qv}$ is not quasi-split.
When $v\in S$, this follows from the assumption of $T$ being nice and the fact \cite[$\S$10]{Kottwitz86} (or \cite[Lem.5.8, 5.9]{LR87}) that fundamental maximal tori of a connected reductive group over a local field transfer to \emph{all} inner forms. When $v\notin S$, this follows from the fact \cite[p.340]{PR94} that any maximal torus in a connected reductive group transfers to the quasi-split inner form.

\textbf{Step 2.} 
We first verify that via the $\Q$-embedding $\rho:T\hookrightarrow I_{x',\epsilon'}^{\mathrm{o}}$, the \emph{same} cocharacter $\mu_T\in X_{\ast}(T)$ still satisfies the two conditions of \cite[Lem.2.2.2]{Kisin17}. 

The two conditions for $\mu_T\in X_{\ast}(T)$ when $T$ is regarded as a subgroup of $I_{x,\epsilon}^{\mathrm{o}}$ concern the map 
\[ j_p: (I_{x,\epsilon}^{\mathrm{o}})_{\Qpb} \isom (G_{\delta\sigma}^{\mathrm{o}})_{\Qpb} \isom (G_{\Qpb})_{\Nm\delta}^{\mathrm{o}}. \]
where the first map $(I_{x,\epsilon}^{\mathrm{o}})_{\Qpb} \isom (G_{\delta\sigma}^{\mathrm{o}})_{\Qpb}$ is the base-change to $\Qpb$ of the $\Qp$-isomorphism $\Int(c)\circ\iota_p$ (\ref{eq:isom_Int(c)circi_p}) and the second map is the base-change to $\Qpb$ of $p_{\delta}:(G_{\delta\sigma})_{L_n}\isom (G_{L_n})_{\Nm\delta}$, which is induced by the projection onto the identity factor: $\Res_{L_n/\Qp}(G_{L_n})\otimes L_n\cong\prod_{L_n\hra L_n}G_{L_n}\ra G_{L_n}$ (Proposition \ref{prop:psi_p}).
On the other hand, by Corollary \ref{cor:Tate_thm2} (3), $\varphi_{\Qpb}$ is conjugate to the composite of these two maps $\Int(c)\circ\iota_p$, $p_{\delta}$, followed by the isomorphism $\Int(c^{-1}):(G_{\Qpb})_{\Nm\delta}^{\mathrm{o}} \isom (I_0)_{\Qpb}$, where $\Nm\delta=c\gamma_0c^{-1}$.
Hence, $\rho_{\Qpb}:(I_{x,\epsilon}^{\mathrm{o}})_{\Qpb}\isom (I_{x',\epsilon'}^{\mathrm{o}})_{\Qpb} $ (\ref{eq:rho_vectorsp_isom}) equals $j_p'^{-1}\circ\Int(g_p)\circ j_p$, where $j_p':  (I_{x',\epsilon'}^{\mathrm{o}})_{\Qpb} \isom (G_{\Qpb})_{\Nm\delta'}^{\mathrm{o}}$ is defined similarly to $j_p$ and $g_p\in G(\Qpb)$ is some element satisfying $\Nm\delta'=g_p\Nm\delta g_p^{-1}$ (so $\Int(g_p):(G_{\Qpb})_{\Nm\delta}^{\mathrm{o}} \isom (G_{\Qpb})_{\Nm\delta'}^{\mathrm{o}}$).
Thus, the conditions to be checked for $\mu_T$ when $T$ is regarded as a subgroup of $I_{x',\epsilon'}^{\mathrm{o}}$ via $\rho$ concern the composite $\Int(g_p)\circ j_p$.
In view of this, the first condition of \cite[Lem.2.2.2]{Kisin17} is obvious, and the second condition follows from Lemma \ref{lem:equality_of_two_Newton_maps} which implies that 
$\Int(g_p)(\nu_{b})=\nu_{b'}$.

Therefore, $\rho\circ \mu_T$ determines a special point $\tilde{x}''$ lifting a point $x''$ in $\mathscr{I}'$, to which the $T$-action on $x''$ (via $\rho$ and a suitable isogeny $x'\ra x''$ which we fix) also lifts. Then, the resulting embedding 
\[ i_{T}':T\hookrightarrow G \]
(obtained via a choice of an isomorphism (\ref{eq:Betti-isom})) is, by construction, an embedding in the stable conjugacy in Theorem \ref{thm:Kisin17_Cor.1.4.13,Prop.2.1.3,Cor.2.2.5} (3) (i) for the maximal torus $\rho:T\subset I_{\mathscr{I}'}$. Let $b''\in i_T'(T)(\mfk)$ be an element provided in that theorem. The pair $(i_T',b'')$ gives a new triple $(\gamma_0'';\gamma'',\delta'')$ by the recipe (\ref{eq:K-triple_attached_to_admissible_K-pair}). By Corollary \ref{cor:Tate_thm2} (3) applied to these data  $(i_T',\gamma_0'':=i_T'(\epsilon'),b''\in I_0''(\mfk))$, we obtain a new inner twisting \[\varphi'':(I_{\mathscr{I}',\epsilon'}^{\mathrm{o}})_{\Qb} \isom (I_0'')_{\Qb}\] that is $T$-equivariant with respect to $\rho:T\hookrightarrow I_{\mathscr{I}',\epsilon'}^{\mathrm{o}}$ and $i_T':T\hookrightarrow I_0'':=G_{\gamma_0''}^{\mathrm{o}}$. 

\textbf{Step 3.} 
It remains to see that there exists $g\in G(\Qb)$ such that \[i_T'=\Int(g)\circ i_T;\] then, it will follow that $h''=\Int(g)(h)$ for the morphism $h''\in X\cap\Hom(\dS,i_T'(T)_{\R})$ giving the CM lifting $\tilde{x}''$, since $h''$ and $h$ are defined by $\mu_{h''}=i_T'\circ \mu_T$ and $\mu_h=i_T\circ \mu_T$. 

For the proof, we consider the functor, defined on the category of $\Q$-algebras $R$, of $R$-linear isomorphisms
\begin{equation} \label{eq:pseudo-torsor2}
H^{\Betti}_1(\mathcal{A}_{\sigma_p(\tilde{x})},\Q) \otimes R \lisom H^{\Betti}_1(\mathcal{A}_{\sigma(\tilde{x}'')},\Q) \otimes R
\end{equation} 
preserving the tensors $\{s_{\alpha,\tilde{x}}\}$, $\{s_{\alpha,\tilde{x}''}\}$ and the $T$-actions $i_T$, $i_T'$ (of course, these identifications of $T$-actions with $i_T$, $i_T'$ are made via the choice of suitable $\Q$-isomorphisms $H^{\Betti}_1(\mathcal{A}_{\sigma_p(\tilde{x})},\Q)\isom V$, $H^{\Betti}_1(\mathcal{A}_{\sigma_p(\tilde{x}'')},\Q)\isom V$, and we will omit to mention them in the rest of the proof).
Since this functor is a pseudo-torsor under $T$ ($\mathrm{rk}_{\Qb}T =\mathrm{rk}_{\Qb}G$), it suffices to show that it is non-empty; then, it will have a $\Qb$-point which would give the required $g\in G(\Q)$ (This argument already appeared in the proof of Theorem \ref{thm:Kisin17_Cor.1.4.13,Prop.2.1.3,Cor.2.2.5} (3) (ii)). 

For $l\neq p$, $\varphi_{\Qlb}:(I_{\mathscr{I},\epsilon}^{\mathrm{o}})_{\Qlb}\isom (I_0)_{\Qlb}$ is \emph{conjugate} to a $\Ql$-isomorphism $i_l:(I_{\mathscr{I},\epsilon}^{\mathrm{o}})_{\Ql}\isom G_{\gamma_l}^{\mathrm{o}}=(I_0)_{\Ql}$ (Corollary \ref{cor:Tate_thm2} (3)) which itself arises from the canonical $\Ql$-isomorphism
\[ \epsilon_l : H^{\et}_1(\mathcal{A}_{\bar{x}},\Ql)\lisom H^{\Betti}_1(\mathcal{A}_{\sigma_p(\tilde{x})},\Q)_{\Ql} \]
matching the tensors $\{s_{\alpha,l,x}\}$, $\{s_{\alpha,\Betti,\sigma_p(\tilde{x})}\}$ and the $T_{\Ql}$-actions $T\subset I_x$, $i_T$ (i.e. $i_l=\Int(\epsilon_l)$ via the choice of suitable isomorphisms $H^{\et}_1(\mathcal{A}_{\bar{x}},\Ql)\isom V_{\Ql}$, $H^{\Betti}_1(\mathcal{A}_{\sigma_p(\tilde{x})},\Q)\isom V$).
The same is true of $\varphi_{\Qlb}':(I_{\mathscr{I}',\epsilon'}^{\mathrm{o}})_{\Qlb}\isom (I_0')_{\Qlb}$ and $i_l':(I_{\mathscr{I}',\epsilon'}^{\mathrm{o}})_{\Ql}\isom G_{\gamma_l'}^{\mathrm{o}}=(I_0')_{\Ql}$, except that the corresponding $\Ql$-isomorphism 
\[ \epsilon_l': H^{\et}_1(\mathcal{A}_{\bar{x}'},\Ql)\lisom H^{\Betti}_1(\mathcal{A}_{\sigma_p(\tilde{x}')},\Q)_{\Ql} \] 
are \emph{not} $T$-equivariant with respect to $\rho:T\hookrightarrow I_{x',\epsilon'}^{\mathrm{o}}$ and $i_T'$: instead, the $\Ql$-isomorphism
\[ \epsilon_l'': H^{\et}_1(\mathcal{A}_{\bar{x}''},\Ql)\lisom H^{\Betti}_1(\mathcal{A}_{\sigma_p(\tilde{x}'')},\Q)_{\Ql} \]
has such property.
Choose $f_l,f_l'\in I_0(\Qlb)$ such that $\varphi_{\Qlb}=\Int(f_l)\circ i_l$ and $\varphi'_{\Qlb}=\Int(f_l')\circ i_l'$.

Now, the base-change $\rho_{\Qlb}=i_l'^{-1}\circ\Int(f_l'^{-1}g_1f_l)\circ i_l$ of $\rho=\varphi'^{-1}\circ\Int(g_1)\circ\varphi$ (\ref{eq:rho_vectorsp_isom}) arises from a composite $\Qlb$-isomorphism of $\Ql$-vector spaces $\epsilon_l'^{-1}\circ g_1\circ \epsilon_l$:
\[ H^{\et}_1(\mathcal{A}_{\bar{x}},\Ql)_{\Qlb} \isom H^{\Betti}_1(\mathcal{A}_{\sigma_p(\tilde{x})},\Q)_{\Qlb} \isom H^{\Betti}_1(\mathcal{A}_{\sigma_p(\tilde{x}')},\Q)_{\Qlb} \isom H^{\et}_1(\mathcal{A}_{\bar{x}'},\Ql)_{\Qlb}=H^{\et}_1(\mathcal{A}_{\bar{x}''},\Ql)_{\Qlb}, \] 
which matches the tensors and is $T$-equivariant with respect to $T\subset I_x$, $\rho:T\hookrightarrow I_{x'}=I_{x''}$, where the middle isomorphism is induced by $g_1\in \Aut(V_{\Qlb})$  (via $H^{\Betti}_1(\mathcal{A}_{\sigma_p(\tilde{x})},\Q)\isom V \stackrel{\sim}{\leftarrow} H^{\Betti}_1(\mathcal{A}_{\sigma_p(\tilde{x}')},\Q)$); note that the second and the third isomorphisms may not be $T$-equivariant individually, but their composite is so.
Then, composing this with $\epsilon_l''$, we obtain an $\Qlb$-isomorphism 
\[ H^{\Betti}_1(\mathcal{A}_{\sigma_p(\tilde{x})},\Q)_{\Qlb} \isom H^{\Betti}_1(\mathcal{A}_{\sigma_p(\tilde{x}'')},\Q)_{\Qlb} \]
which match the tensors and the $T$-actions, as required.

Therefore, we obtain a cocycle $(\tilde{\beta}_{\tau})_{\tau}:=i_T^{-1}(g^{-1}\tau(g))\in Z^1(\Q,T)$.
As $h'=\Int(g)(h)$ and $h$ are $G(\R)$-conjugate, the image of $i_T(\tilde{\beta})$ in $H^1(\R,K_{\infty})$ is trivial, where $K_{\infty}$ is the centralizer of $h$ (compact-modulo-center inner form of $G_{\R}$) which implies \cite[Lem.4.4.5]{Kisin17} that $\tilde{\beta}$ lies in $H^1(\R,T)$, and so in $\Sha^{\infty}_G(\Q,i_T(T))$ (via $i_T$).

Finally, we have $(\mathscr{I}',\epsilon')=(\mathscr{I}^{\beta},\epsilon^{\beta})$ by Theorem \ref{thm:Kisin17_Prop.4.4.8} (3).
\end{proof}

%%%%%%%%%%%%%%%%%%%%
\begin{rem} \label{rem:Kisin17's_error_Prop.4.4.13}
This gives a correct proof of \cite{Kisin17}, Proposition 4.4.13. In the proof of \textit{loc. cit.}, there is no guarantee that the diagram in the bottom on p.80 commutes up to conjugation:
\[ \xymatrix{ I_0\otimes_{\Q}\Qb \ar[rr]^{\sim}_{\varphi} \ar[d]^{\Int(g)} \ar@{->}[rrd]_{\psi} & &  I\otimes_{\Q}\Qb \ar[d]^{\sim}_{\theta} \\ I_0'\otimes_{\Q}\Qb \ar[rr]^{\sim}_{\varphi'} & & I'\otimes_{\Q}\Qb} \]
We use the notation of \textit{loc. cit.}, except that $\psi,\theta$ are introduced below by us. 
Here, $\varphi:I_0\otimes_{\Q}\Qb \ra I\otimes_{\Q}\Qb$ is an inner-twisting whose local inner classes in $H^1(\Qv,I_0^{\ad})$ are determined by the Kottwitz triple $(\gamma_0;\gamma,\delta)$ attached to $\mathscr{I}$ (\textit{ibid.} Corollary 2.3.5) (Kottwitz triples in Kisins' work have the same definition as ours, except that their equivalence classes are coarser than ours \cite[4.3.1]{Kisin17}). It has the additional property that the restriction of $\varphi^{-1}$ to any maximal torus $T$ of $I$ (Kisin fixes one in the beginning) is $\Q$-rational (in fact, $\gamma_0$ is determined via this $\varphi^{-1}:T\hra G$ as $\varphi^{-1}(\pi_{\mathscr{I}})$ for some $\F_q$-Frobenius $\pi_{\mathscr{I}}$ attached to $\mathscr{I}$). The inner-twisting $\varphi':I_0'\otimes_{\Q}\Qb \ra I'\otimes_{\Q}\Qb$ $\varphi'$ is a similar inner-twisting defined by the Kottwitz triple $(\gamma_0';\gamma',\delta')$ attached to $\mathscr{I}'$; but, in the proof it will be constructed after a $\Q$-isomorphism $\theta$ is constructed first (see below).

Before going into the explanation of the issue, we make some preliminary observations: from the assumption of the equivalence $(\gamma_0;\gamma,\delta)\sim(\gamma_0';\gamma',\delta')$, it follows that $I_{\Qv}\sim I'_{\Qv}$ for every place $v$ of $\Q$ and that for any $g\in G(\Qb)$ satisfying $\gamma_0'=g\gamma_0g^{-1}$, the classes in $H^1(\Qv,I_0^{\ad})$ of the inner twistings $\varphi_{\Qvb},\varphi'_{\Qvb}\circ\Int(g)$ are the same for all $v$ (cf. Lemma \ref{lem:beta_v=beta_v'+zeta} below. So, $\varphi,\varphi'':=\varphi'\circ\Int(g)$ have the same (global) inner classes by Hasse principle for adjoint groups, and there exists $f\in I_0^{\ad}(\Qb)$ such that $(\varphi\circ \Int(f))^{-1}\circ{}^{\tau}(\varphi\circ\Int(f))=\varphi''^{-1}\circ{}^{\tau}\varphi''$ for all $\tau\in\Gal(\Qb/\Q)$, namely we know that there exists a $\Q$-isomorphism $\theta':=\varphi''\circ\Int(f^{-1})\circ\varphi^{-1}:I\isom I'$ making the diagram above (with $\theta'$ replacing $\theta$) commute up to conjugation.

For Kisin's proof, it is critical that with the diagram as above, there exists a maximal torus $T$ of $I_{\mathscr{I}}$ such that $\varphi^{-1}|_T$ and $\varphi'^{-1}\circ \theta |_T$ are \emph{both} $\Q$-rational. The $\Q$-rationality for $\varphi^{-1}|_T$ can be realized (for arbitrary $T$) simply by the construction of $\varphi$ (\textit{ibid.} Corollary 2.3.5). But the same property for $\varphi'^{-1}\circ \theta |_T$ involves both $\varphi'$ and $\theta$ and the constructions of these maps to attain that property require the construction of each other first: $\theta$ is to be determined (up to $I^{\ad}(\Q)$) by $\varphi'$ (to provide the commutative diagram), while $\varphi'$ is constructed so as to satisfy, among others, the condition that the restriction of $\varphi'^{-1}$ to the pre-chosen maximal torus $\theta(T)$ of $I'$ is $\Q$-rational. With some thinking (and as argued below), it does not seem possible to achieve these two goals simultaneously by this method (of constructing inner twistings \emph{by cohomology classes}).

We continue to explain the rest of the proof of \cite[Prop. 4.4.13]{Kisin17} after the construction of $\varphi$.
For the purpose just stated, from the equivalence of the Kottwitz triples $(\gamma_0;\gamma,\delta)\sim (\gamma_0;\gamma',\delta')$ ($\gamma_0'$ is not involved yet), Kisin first deduces an inner-twisting $\psi:I_0\otimes_{\Q}\Qb \ra I'\otimes_{\Q}\Qb$ whose local inner classes in $H^1(\Qv,I_0^{\ad})$ are the same everywhere as those of $\varphi$ (\textit{ibid.} Corollary 2.3.5). Then, by the same argument as above, this fixes (up to $I'^{\ad}(\Q)$) a $\Q$-isomorphism $\theta:I\isom I'$ (of inner twistings $\varphi$, $\psi$) such that $\theta\circ\varphi$ and $\psi$ are conjugate. Next, \emph{using this $\theta$}, one regards $T$ as a $\Q$-torus of $I'$ and applying the same corollary again (relying on \textit{ibid.} Cor. 2.2.5) he obtains an inner twisting $\varphi':I_0'\otimes_{\Q}\Qb \ra I'\otimes_{\Q}\Qb$ \emph{whose restriction to $\theta(T)$ is $\Q$-rational} (again, this latter condition gives $\gamma_0':=\varphi'^{-1}(\pi_{\mathscr{I}'})\in G(\Qb)$). 
%For the later part of Kisin's proof, the order is crucial that one first constructs a $\Q$-isomorphism $\theta:I\isom I'$ and then, regarding $T$ as a $\Q$-torus of $I'$ via this $\theta$, construct an inner twisting $\varphi'$ whose restriction to $\theta(T)$ is $\Q$-rational. 
Then, from the discussion above, we know that the inner-twistings $\varphi'\circ\Int(g)$, $\psi$ are the same, and thus
there exists a $\Q$-automorphism $\Xi:I'\ra I'$ such that $\Xi\circ \psi$ and $\varphi'\circ\Int(g)$ commute up to conjugation: but, only when $\Xi$ is itself an inner automorphism, $\psi,\varphi'\circ\Int(g)$ commute up to conjugation.

Kisin's proposition can be proved by taking, in our proof above, $\epsilon$, $\epsilon'$ to be suitable elements of the Frebenius germs of $\mathscr{I}$, $\mathscr{I}'$ in Theorem \ref{thm:Kisin17_Cor.2.3.2} (such that their centralizers become connected, i.e. equal to $I_{\mathscr{I}}$, $I_{\mathscr{I}'}$). Note that in our proof of this theorem, we do not assert the existence of a commutative (up to conjugation) diagram of inner twistings as Kisin does and mostly importantly that we construct an inner-twisting $\rho:I_{\Qb}\isom I'_{\Qb}$ (\ref{eq:rho_vectorsp_isom}) (not a $\Q$-isomorphism $\theta$ as above) which becomes $\Q$-rational only upon restriction to a given maximal torus and is conjugate over $\Qvb$ to an inner twisting which \emph{arises from an isomorphism of vector spaces}.%%
\footnote{This principle that inner twistings of reductive groups should arise from isomorphisms of underlying vector spaces lies at the heart of our method adapting the Galois gerb theory in \cite{LR87}, where inner twistings always come from twisting some subgroups of the general linear group with cochains (which become cocycles only modulo centers), thus arise from isomorphisms of vector spaces. Such twistings with cochains are also demonstrated here, in the proofs of Theorem \ref{thm:stable_isogeny_diagram} and Corollary \ref{cor:stable_conjugacy_from_nice_tori} (the former is the only place in this work using Galois gerb theoretic results in \cite{LR87}). }
\end{rem}

%%%%%%%%%%%%%%%%%%%%
\begin{lem} \label{lem:beta_v=beta_v'+zeta}
For a finite place $v$ of $\Q$, let $\gamma_0\in G(\Qv)_{\mathrm{ss}}$ and $g\in G(\Qvb)$ be such that $\zeta:=(g^{-1}\tau(g))_{\tau}\in Z^1(\Qv,I_0)$; so, $\gamma_0':=\Int(g)(\gamma_0)\in G(\Qv)$ and $\Int(g):(I_0)_{\Qvb}\isom (I_0')_{\Qvb}$ is an inner twisting, where $I_0':=G_{\gamma_0'}^{\mathrm{o}}$.

If $\varphi:I_{\Qvb}\isom (I_0)_{\Qvb}$ is an inner twisting (for a $\Qv$-group $I$) with the corresponding cohomology class $\beta_v\in H^1(\Qv,I_0^{\ad})$, then 
the cohomology class $\beta_v'\in H^1(\Qv,I_0'^{\ad})$ of $\Int(g)\circ\varphi$ is given by
\[ \beta_v'=\beta_v- [\zeta^{\ad}] \]
under the canonical isomorphism $H^1(\Qv,I_0^{\ad}) \cong \pi_1(I_0^{\ad})_{\Gamma_v,\mathrm{tor}}\cong H^1(\Qv,I_0'^{\ad})$, where $[\zeta^{\ad}]$ is the class of the image $\zeta^{\ad}\in Z^1(\Qv,I_0^{\ad})$  of $\zeta$.
\end{lem}

\begin{proof}
Via the map $\Int(g^{-1}):I_0'^{\ad}(\Qvb)\isom I_0^{\ad}(\Qvb)$, the Galois action on $I_0'^{\ad}(\Qvb)$ is identified with the twisted Galois action on $I_0^{\ad}(\Qvb)$ defined by $x\mapsto \Int(\zeta_{\tau})(\tau(x))$: $\Int(g^{-1})(\tau x)=\zeta_{\tau}\tau(g^{-1}xg)\zeta_{\tau}^{-1}$, for $x\in I_0'^{\ad}(\Qv)$. Namely, in the notation of \cite[5.3, Examples 2)]{Serre02}, $I_0'^{\ad}={}_{\zeta}I_0^{\ad}$ via $\Int(g^{-1})$. So, by ibid. Prop.35bis, there exists a bijection
\[ t_{\zeta} \ :\ Z^1(\Qv,I_0'^{\ad}) \isom Z^1(\Qv,I_0^{\ad})\quad : \quad x_{\tau} \mapsto \Int(g^{-1})(x_{\tau})\cdot \zeta_{\tau}^{\ad}=(g^{\ad})^{-1}x_{\tau}\tau(g^{\ad}) \]
which induces a bijection 
\[ t_{\zeta}: H^1(\Qv,I_0'^{\ad}) \isom H^1(\Qv,I_0^{\ad}). \] 

On the other hand, if $x'_{\tau}\in Z^1(\Qv,I_0'^{\ad})$ represents the inner twisting $\varphi':=\Int(g)\circ\varphi$ (i.e. $\varphi'\circ{}^{\tau}\varphi'^{-1}=\Int(x'_{\tau})$ for every $\tau\in\Gamma_v$), one readily finds that the inner twisting $\varphi$ is represented by $t_{\zeta}(x'_{\tau})$.
Thus, the claim follows since under the canonical isomorphism $H^1(\Qv,G)\cong H^1(\Qv,G_{\mathbf{ab}})\cong \pi_1(G)_{\Gamma_p,\mathrm{tor}}$ for connected reductive $\Qv$-groups $G$, (for nonarchimedean $v$), this map $t_{\zeta}$ is the translation by the cohomology class $ [\zeta^{\ad}]$  (\cite[Prop.3.16]{Borovoi98}). 
\end{proof}

%%%%%%%%%%%%%%%%%%%%
\begin{lem} \label{lem:|widetilde{Sha}_G(Q,I_{phi,epsilon})^+|}
For any admissible pair $(\mathscr{I},\epsilon)$ with an associated Kottwitz triple $(\gamma_0;\gamma,\delta)$, the set 
\[ \widetilde{\Sha}_G(\Q,I_{\mathscr{I},\epsilon})^+:=\ker[\Sha^{\infty}_G(\Q,I_{\mathscr{I},\epsilon}^{\mathrm{o}})\rightarrow H^1(\A,I_{\mathscr{I},\epsilon})]\] 
(cf. (\ref{eq:Sha_G^{infty}(Q,I_{sI,epsilon}^{mathrm{o}})})) is canonically identified with $\ker[\Sha^{\infty}_G(\Q,I_0)\rightarrow \bigl(H^1(\A_f^p,G_{\gamma})\oplus H^1(\Qp,G_{\delta\sigma})\bigr)]$ ($I_0:=G_{\gamma_0}^{\mathrm{o}}$ as usual). In particular, it is a finite set which depends only on the associated equivalence class of Kottwitz triple. Its cardinality $\tilde{i}(\gamma_0;\gamma,\delta)$ is given by
\[  |\widetilde{\Sha}_G(\Q,I_{\mathscr{I},\epsilon})^+|=|\ker[ \ker^1(\Q,I_0)\rightarrow \ker^1(\Q,G)]| \cdot |\mathfrak{D}(\gamma_0;\gamma,\delta)|,\]
with
\[ \mathfrak{D}(\gamma_0;\gamma,\delta) :=\im[\Sha^{\infty}_G(\Q,I_0)\rightarrow H^1(\A_f,I_0)] \cap 
\ker [H^1(\A_f,G_{\gamma,\delta}^{\mathrm{o}})\rightarrow H^1(\A_f,G_{\gamma,\delta})], \]
where $H^1(\A_f,G_{\gamma,\delta}^{\mathrm{o}}):=\oplus_l H^1(\Qv,G_{\gamma_l}^{\mathrm{o}})\oplus H^1(\Qp,G_{\delta\sigma}^{\mathrm{o}})$ and $H^1(\A_f,G_{\gamma,\delta})$ is defined similarly.
\end{lem}

Recall that for any connected reductive $\Q$-group $I_0$ and its inner twist $I_1$, one has canonical isomorphisms $\Sha^{\infty}_G(\Q,I_0)=\Sha^{\infty}_G(\Q,I_{0\mathbf{ab}})$, $H^1(\A_f,I_0)=H^1(\A_f,I_{0\mathbf{ab}})$, and that the abelianization complexes $I_{0\mathbf{ab}}$ and $I_{1\mathbf{ab}}$ are canonically isomorphic in the derived category of commutative algebraic $\Q$-group schemes, both being quasi-isomorphic to $Z(I_0^{\uc})\rightarrow Z(I_0)$, cf. Appendix \ref{sec:abelianization_complex}.
Since $I_{\mathscr{I},\epsilon}^{\mathrm{o}}$ is an inner form of $I_0$ (Corollary \ref{cor:Tate_thm2}), the identification in the first claim of the lemma is made via the corresponding canonical isomorphism $\Sha^{\infty}_G(\Q,I_{\mathscr{I},\epsilon}^{\mathrm{o}})\cong \Sha^{\infty}_G(\Q,I_0)$ (independent of the choice of the inner twisting).
Similarly, the localization map 
\[\Sha^{\infty}_G(\Q,I_0)\rightarrow \bigl(H^1(\A_f^p,G_{\gamma})\oplus H^1(\Qp,G_{\delta\sigma})\bigr)\] 
and the intersection in the definition of $\mathfrak{D}(\gamma_0;\gamma,\delta)$ make sense because there exist inner twists $(I_0)_{\Qvb}\isom (G_{\gamma_v}^{\mathrm{o}})_{\Qvb}\ (v\neq p)$, $(I_0)_{\Qpb}\isom(G_{\delta\sigma}^{\mathrm{o}})_{\Qpb}$ (Corollary \ref{cor:Tate_thm2}), thereby canonical identifications
\[H^1(\A_f,I_0) \cong H^1(\A_f,(I_0)_{\mathbf{ab}}) \cong H^1(\A_f,(G_{\gamma,\delta}^{\mathrm{o}})_{\mathbf{ab}}) \cong H^1(\A_f,G_{\gamma,\delta}^{\mathrm{o}}).\]

%%%%%%%%%%%%%%%%%%%%
\begin{rem} \label{rem:|widetilde{Sha}_G(Q,I_{phi,epsilon})^+|}
When $I_{\mathscr{I},\epsilon}$ is connected (e.g., if $G^{\der}=G^{\uc}$), one has $|\mathfrak{D}(\gamma_0;\gamma,\delta)|=1$ and the constant $|\widetilde{\Sha}_G(\Q,I_{\mathscr{I},\epsilon})^+|$ depends only on the stable conjugacy class of $\gamma_0$, rather than on the whole triple $(\gamma_0;\gamma,\delta)$; but, we do not know whether this is true in general, since for a non-connected group $H$ over a local field $\Qv$, $H^1(\Qv,H)$ may not be invariant under inner twists of $H$. 

The constant $\ker^1(\Q,I_0)$ also appears in \cite[Lem.17.2]{Kottwitz92} with a similar interpretation (see the discussion on p.441-442 of ibid. for the appearance of the constant $|\ker[ \ker^1(\Q,I_0)\rightarrow \ker^1(\Q,G)]|$). \end{rem}

\begin{proof}
This set is equal to $\ker[\Sha^{\infty}_G(\Q,I_{\mathscr{I},\epsilon}^{\mathrm{o}})\rightarrow H^1(\A_f,I_{\mathscr{I},\epsilon})]$ (no $\infty$-component in the target). There exists a commutative diagram with isomorphic vertical maps:
\[ \xymatrix{\Sha^{\infty}_G(\Q,I_{\mathscr{I},\epsilon}^{\mathrm{o}}) \ar[r]^j \ar[d]_i &H^1(\A_f,I_{\mathscr{I},\epsilon}^{\mathrm{o}}) \ar[r]^{k}  \ar[d]^f &  H^1(\A_f,I_{\mathscr{I},\epsilon}) \ar[d]^f \\ 
\Sha^{\infty}_G(\Q,I_0) \ar[r]^{j_0} & H^1(\A_f,G_{\gamma,\delta}^{\mathrm{o}}) \ar[r]^{k_0} &H^1(\A_f,G_{\gamma,\delta}), }\]
where the last two vertical maps are given by any isomorphism $f:(I_{\mathscr{I},\epsilon})_{\A_f}\isom G_{\gamma,\delta}$, e.g. as in Corollary \ref{cor:Tate_thm2}. The first vertical map $i$ and $j_0:\Sha^{\infty}_G(\Q,I_0)\rightarrow H^1(\A_f,I_0)\cong H^1(\A_f,G_{\gamma,\delta}^{\mathrm{o}})$ are both induced by (some) inner twists $(I_{\mathscr{I},\epsilon}^{\mathrm{o}})_{\Qb}\isom (I_0)_{\Qb}$, $(I_0)_{\bar{\A}_f} \isom (G_{\gamma,\delta}^{\mathrm{o}})_{\bar{\A}_f}$, thus in fact are \emph{canonical}, independent of the choice of the inner twists. This shows that the left square of the diagram commutes. 
Note that although the last vertical map is not canonical, depending on the choice of an isomorphism $f:(I_{\mathscr{I},\epsilon})_{\A_f}\isom G_{\gamma,\delta}$ (which itself is not determined by the Kottwitz triple alone), the middle vertical map induced by it is canonical, and the commutativity of the right square implies that $\widetilde{\Sha}_G(\Q,I_{\mathscr{I},\epsilon})^+=\ker(k\circ j)$ (although the map $k\circ j$ itself may not) is canonically identified with $\ker(k_0\circ j_0)$; in particular, for two admissible pairs $(\mathscr{I},\epsilon)$, $(\mathscr{I}',\epsilon')$ with the same associated equivalence classes of Kottwitz triples, the corresponding sets $\widetilde{\Sha}_G(\Q,I_{\mathscr{I},\epsilon})^+$, $\widetilde{\Sha}_G(\Q,I_{\mathscr{I}',\epsilon'})^+$ are also canonically isomorphic. This proves the first claim.

Next, we see that our set $\widetilde{\Sha}_G(\Q,I_{\mathscr{I},\epsilon})^+$ is the disjoint union of $j^{-1}(\beta)$'s with $\beta$ running through
\begin{align} \label{eq:D(gamma_0;gamma,delta)}
&\ \ker[H^1(\A_f,I_{\mathscr{I},\epsilon}^{\mathrm{o}})\rightarrow H^1(\A_f,I_{\mathscr{I},\epsilon})] \cap \im(j) \\ \cong &\ \ker [H^1(\A_f,G_{\gamma,\delta}^{\mathrm{o}})\rightarrow H^1(\A_f,G_{\gamma,\delta})]  \cap \im[\Sha^{\infty}_G(\Q,I_0)\rightarrow H^1(\A_f,I_0)] \nonumber \\
=: &\ \mathfrak{D}(\gamma_0;\gamma,\delta) \nonumber. 
\end{align}
(again, by the same reasoning, this identification is canonical.) For $\alpha\in H^1(\Q,I_{\mathscr{I},\epsilon}^{\mathrm{o}})$, let $H^1(\Q,I_{\mathscr{I},\epsilon}^{\mathrm{o}})_{\alpha}$ denote the subset of $H^1(\Q,I_{\mathscr{I},\epsilon}^{\mathrm{o}})$ consisting of elements having the same image in $H^1(\A,I_{\mathscr{I},\epsilon}^{\mathrm{o}})$ as $\alpha$. 
Then, we claim that for each $\beta=j(\alpha)$ with $\alpha\in \Sha_G^{\infty}(\Q,I_{\mathscr{I},\epsilon}^{\mathrm{o}})$, there are bijections
\begin{equation} \label{eq:j^{-1}(beta)}
 j^{-1}(\beta)\ \isom\ \ker[H^1(\Q,I_{\mathscr{I},\epsilon}^{\mathrm{o}})_{\alpha}\stackrel{\mathbf{ab}^1}{\longrightarrow} H^1(\Q,G_{\mathbf{ab}})]\ \isom\ \ker[\ker^1(\Q,I_0)\rightarrow \ker^1(\Q,G)], 
\end{equation}
where $\mathbf{ab}^1$ is the composite map
\[\mathbf{ab}^1:H^1(\Q,I_{\mathscr{I},\epsilon}^{\mathrm{o}})\rightarrow H^1(\Q,(I_{\mathscr{I},\epsilon}^{\mathrm{o}})_{\mathbf{ab}})=H^1(\Q,(I_0)_{\mathbf{ab}})\rightarrow H^1(\Q,G_{\mathbf{ab}})\] 
(or its restriction $\Sha^{\infty}(\Q,I_{\mathscr{I},\epsilon}^{\mathrm{o}})=\Sha^{\infty}(\Q,(I_{\mathscr{I},\epsilon}^{\mathrm{o}})_{\mathbf{ab}})\rightarrow \Sha^{\infty}(\Q,G_{\mathbf{ab}})=\Sha^{\infty}(\Q,G)$). 
The first bijection is clear. For the second bijection, we first recall the fact (\cite[Prop.35bis]{Serre02}) that for any cocycle $a\in Z^1(\Q,I_{\mathscr{I},\epsilon}^{\mathrm{o}})$ representing $\alpha$, if $I_a^{\mathrm{o}}$ is the inner twist of $I_{\mathscr{I},\epsilon}^{\mathrm{o}}$ via $a$, there exists a commutative diagram 
\[ \xymatrix@R-1pc{ H^1(\Q,I_{\mathscr{I},\epsilon}^{\mathrm{o}}) \ar[r]^{\sim} \ar[d] & H^1(\Q,I^{\mathrm{o}}_a) \ar[d] \\ H^1(\A,I_{\mathscr{I},\epsilon}^{\mathrm{o}}) \ar[r]^{\sim} & H^1(\A,I^{\mathrm{o}}_a) } \] where horizontal bijections are given by $[a']\mapsto [a'a^{-1}]$ (so send $[a]$ to the neutral element, in each row), which thus induces a bijection
\[H^1(\Q,I_{\mathscr{I},\epsilon}^{\mathrm{o}})_{\alpha} \isom \ker[H^1(\Q,I_a^{\mathrm{o}})\rightarrow H^1(\A,I_a^{\mathrm{o}})]. \] 
As $\alpha$ maps to zero in $H^1(\Q,G_{\mathbf{ab}})$, this bijection commutes with the respective natural maps $\mathbf{ab}^1$ into $H^1(\Q,G_{\mathbf{ab}})$ (see \cite[3.15]{Borovoi98}).
Hence, the second bijection in (\ref{eq:j^{-1}(beta)}) follows from this bijection and the well-known fact that for any \emph{connected} reductive group $H$ over $\Q$, the finite abelian group $\ker^1(\Q,H)$ is unchanged under inner twists.
Note that the set $\ker[H^1(\A_f,I_{\mathscr{I},\epsilon}^{\mathrm{o}})\rightarrow H^1(\A_f,I_{\mathscr{I},\epsilon})]$, being the image of the finite set $\pi_0(I_{\mathscr{I},\epsilon})(\A_f)$ in $H^1(\A_f,I_{\mathscr{I},\epsilon}^{\mathrm{o}})$, is also finite, which implies the same property for $\mathfrak{D}(\gamma_0;\gamma,\delta)$ and $\widetilde{\Sha}_G(\Q,I_{\mathscr{I},\epsilon})^+$. The formula for $|\widetilde{\Sha}_G(\Q,I_{\mathscr{I},\epsilon})^+|$ is immediate from this discussion.
\end{proof}

%%%%%%%%%%%%%%%%%%%%
\begin{cor} \label{cor:LR-Satz5.25} 
Keep the assumptions of Theorem \ref{thm:LR-Satz5.21}.
Let $(\gamma_0;\gamma=(\gamma_l)_{l\neq p},\delta)$ be a stable Kottwitz triple with trivial Kottwitz invariant.
Suppose that one of the following two conditions holds: (a) $Z(G)$ has same ranks over $\Q$ and $\R$, or (b) the weight homomorphism $w_X$ is rational and $\gamma_0$ is a Weil $p^n$-element of weight $-w_X$, where $n$ is the level of $(\gamma_0;(\gamma_l)_{l\neq p},\delta)$ (cf. Proposition \ref{prop:phi(delta)=gamma_0_up_to_center}).

If $\mathrm{O}_{\gamma}(f^p)\cdot \mathrm{TO}_{\delta}(\phi_p)\neq 0$, the stable Kottwitz triple $(\gamma_0;\gamma,\delta)$ is \emph{effective}, that is, there exists an admissible pair $(\mathscr{I},\epsilon)$ giving rise to it, in which case the set of such admissible pairs has the same cardinality as the set
\[ \Sha_G(\Q,I_{\mathscr{I},\epsilon})^+ :=\mathrm{im}[\widetilde{\Sha}_G(\Q,I_{\mathscr{I},\epsilon})^+\rightarrow H^1(\Q,I_{\mathscr{I},\epsilon})]. \]
In fact, the twisting $\beta\in \Sha_G(\Q,I_{\mathscr{I},\epsilon})^+ \mapsto (\mathscr{I}^{\beta},\epsilon^{\beta})$ of Theorem \ref{thm:Kisin17_Prop.4.4.8} gives a bijection.
\end{cor}

%%%%%%%%%%%%%%%%%%%%
\begin{rem} \label{rem:LR-Satz5.25} 
From this theorem, we see that the number $|\Sha_G(\Q,I_{\mathscr{I},\epsilon})^+|$ depends only on the effective stable Kottwitz triple $(\gamma_0;\gamma,\delta)$, not on the choice of an admissible pair giving rise to it. When $G_{\gamma_0}$ is connected so that $I_{\mathscr{I},\epsilon}$ is also connected (as they are inner forms of each other), the set $\Sha_G(\Q,I_{\mathscr{I},\epsilon})^+$ is identical to the group denoted by $\Sha_G(\Q,I_{\mathscr{I},\epsilon})$ in \cite[4.4.9]{Kisin17}. So, in this case, one knows that this is invariant under inner twist of $I_{\mathscr{I},\epsilon}$, thereby equal to the finite abelian group $\ker[ \ker^1(\Q,I_0)\rightarrow \ker^1(\Q,G)]$ which depends even only on $\gamma_0$ (this last set is the one used in \cite[Lem.5.24]{LR87} for the same purpose when $G^{\der}=G^{\uc}$).
But, in general cases the author does not know whether any of these properties holds for our set $\Sha_G(\Q,I_{\mathscr{I},\epsilon})^+$ (except for finiteness which follows from the same property of $\widetilde{\Sha}_G(\Q,I_{\mathscr{I},\epsilon})^+$).
\end{rem}

\begin{proof} 
One easily checks (cf. \cite[$\S$1.4, $\S$1.5]{Kottwitz84b}) that non-vanishing of $\mathrm{TO}_{\delta}(\phi_p)$ is equivalent to non-emptiness of the set $Y_p(\delta)$ (\ref{eq:Y_p(delta)}).
Hence, by Theorem \ref{thm:LR-Satz5.21}, there exist a special Shimura subdatum $(T,h)$ and an element $\epsilon_0\in T(\Q)$ such that $\epsilon_0$ is stably conjugate to $\gamma_0$ and $(T,h,\epsilon_0)$ is admissible in the sense of Definition \ref{defn:Frobenius_pair}. Let $(\mathscr{I}_1,\epsilon_1):=(\mathscr{I}_{T,h},j_{T,h}(\epsilon_0))$ be the corresponding K-pair. This K-pair is admissible (in the sense of Definition \ref{defn:admissible_pair2}) by Theorem \ref{thm:Kisin17_Cor.1.4.13,Prop.2.1.3,Cor.2.2.5} (i) and (iii): $\kappa_T([\Nm_{K/K_0}(\mu_h(\pi))^{-1}])=-\underline{\mu_h}$ (\ref{eq:Kottwitz97-(7.3.1)}).
Then, to prove the effectivity of the given stable Kottwitz triple $(\gamma_0;\gamma,\delta)$, it suffices to find a class $\beta\in \Sha^{\infty}(\Q,I_{\mathscr{I}_1,\epsilon_1}^{\mathrm{o}})$ such that the twist  $(\mathscr{I}^{\beta},\epsilon^{\beta})$ of $(\mathscr{I}_1,\epsilon_1)$ by (the image in $H^1(\Q,I_{\mathscr{I}_1,\epsilon_1})$ of) $\beta$ produces the triple. 
In fact, we will prove a stronger statement.

For the stable Kottwitz triple $(\gamma_0;\gamma,\delta)$ (with trivial Kottwitz invariant), we fix a tuple of elements $(g_v)_v\in G(\bar{\A}_f^p)\times G(\mfk)$ satisfying (\ref{eq:stable_g_l}) and Definition \ref{defn:stable_Kottwitz_triple} (iii$'$) such that the associated Kottwitz invariant is trivial. This choice gives an invariant $\alpha_v(\gamma_0;\gamma,\delta;g_v)\in X^{\ast}(Z(\hat{I}_0)^{\Gamma_v})$ for each place $v$ of $\Q$ (notation: $\alpha_l(\gamma_0;\gamma,\delta;g_l)=\alpha_l(\gamma_0;\gamma_l;g_l)$ (\ref{eq:alpha_l}) for finite $l\neq p$, $\alpha_p(\gamma_0;\gamma,\delta;g_p)=\alpha_p(\gamma_0;\delta;g_p)$  (\ref{eq:alpha_p})). 
Let $(\gamma_0;\gamma_1,\delta_1)$ be the stable Kottwitz triple attached to the admissible pair $(\mathscr{I}_1,\epsilon_1)$ (with the same rational component $\gamma_0$), and also choose a similar tuple of elements $(g_{1v})_v\in G(\bar{\A}_f^p)\times G(\mfk)$, thereby obtaining invariants $\alpha_v(\gamma_0;\gamma_1,\delta;g_{1v})$, whose associated Kottwitz invariant is trivial. Given these choices, we will find $\beta\in \Sha^{\infty}(\Q,I_{\mathscr{I}_1,\epsilon_1}^{\mathrm{o}})$ such that for each place $v$, its localization $\beta(v)$ equals
\[ (\alpha(v):=\alpha_v(\gamma_0;\gamma,\delta;g_v)-\alpha_v(\gamma_0;\gamma_1,\delta_1;g_{1v}))_v\in  \bigoplus_{v} \mathbb{H}^1(\Qv,I_{0\mathbf{ab}}): \]
$\alpha(v)$ belongs to the subgroup $\mathbb{H}^1(\Qv,I_{0\mathbf{ab}})\subset \pi_1(I_0)_{\Gamma_v,\tors}$, since these groups are equal for nonarchimedean $v$ \cite[Thm.1.2]{Kottwitz86} and $[\alpha(p)]=\kappa_{I_0}(b(\gamma_0;\delta;g_p))-\kappa_{I_0}(b(\gamma_0;\delta_1;g_{1p}))\in \pi_1(I_0)_{\Gamma_v,\tors}$ (cf. Remark \ref{rem:kappa(b)-kappa(b')_is_torsion}), while $\alpha(\infty)=0$.
Then, by Corollary \ref{cor:K-triple_of_twisted_K-pair}, for the new admissible K-pair $(\mathscr{I}^{\beta},\epsilon^{\beta})$, there will exist an auxiliary tuple $(g_v')_v\in G(\bar{\A}_f^p)\times G(\mfk)$ such that $\alpha_v(\gamma_0;\gamma_{v}',\delta';g_{v}')=\alpha_v(\gamma_0;\gamma_v,\delta;g_{v})$. Of course, this implies that $(\gamma_0;\gamma',\delta')$ is equivalent to $(\gamma_0;\gamma,\delta)$. 

We set  $I_1:=I_{\mathscr{I},\epsilon}^{\mathrm{o}}$. With the choice of a maximal $\Q$-torus of $T$ of $I_{\mathscr{I}}$ with $\epsilon\in T(\Q)$ and an embedding $T\hookrightarrow G$, the inner-twist $\varphi:(I_1)_{\Qb}\isom (I_0)_{\Qb}$ in Corollary \ref{cor:Tate_thm2} satisfies the condition of Lemma \ref{lem:abelianization_exact_seq}; let $\tilde{I}_0=\rho^{-1}(I_0)$ and $\tilde{I}_1$ be as in that lemma.
Then, in view of Lemma \ref{lem:abelianization_exact_seq}, we need to find a global class $[\widetilde{a}]\in H^1(\Q,\tilde{I}_{1})$ whose localizations $([\tilde{a}(v)])$ go over to $(\alpha(v))$ under the map $\oplus f_v:=\oplus \mu_v\circ \mathbf{ab}_v$ in the commutative diagram:
\begin{equation}  \label{eq:H^1_{ab}_for_(I'->I)}
 \xymatrix{ \bigoplus_v H^1(\Qv,I_{1}) \ar[r] & \bigoplus_v \mathbb{H}^1(\Qv,I_{1\mathbf{ab}}) & \mathfrak{K}(I_{1}/\Q)^D  \\ 
 \bigoplus_v H^1(\Qv,\tilde{I}_{1}) \ar[u] \ar[r]^(.45){\mathbf{ab}} \ar[ru]^{\oplus f_v} & \bigoplus_v \mathbb{H}^1(\Qv,\tilde{I}_{1\mathbf{ab}}) \ar[u]_{\mu=\oplus_v\mu_v} \ar[r]^(.5){\lambda} & \mathbb{H}^1(\A/\Q,\tilde{I}_{1\mathbf{ab}}) \ar@{->>}[u]_{\nu}  \\
 \bigoplus_v \mathbb{H}^0(\Qv,\tilde{I}_{1}\rightarrow I_{1}) \ar[r]^(.5){\overline{\mathbf{ab}}} \ar[u] & \bigoplus_v \mathbb{H}^0(\Qv,G_{\mathbf{ab}})  \ar[u]_{\xi} \ar@{->>}[r] & \im(\lambda\circ\xi) \ar@{^(->}[u] }
 \end{equation} 
Here, the first two columns (consisting of adelic cohomology groups) are each a part of the cohomology long exact sequence for the crossed module $\tilde{I}_1\rightarrow I_1$ \cite[(3.4.3.1)]{Borovoi98} and for the distinguished triangle (\ref{eq:DT_of_CX_of_tori}) attached to the morphism $\tilde{I}_{1\mathbf{ab}}\rightarrow I_{1\mathbf{ab}}$; so, they are exact. The two horizontal maps for $H^1$ between these two columns are the abelianization maps for $I_1$, $\tilde{I}_1$, and the left lower commutative diagram (containing the map $\overline{\mathbf{ab}}$) is induced from an obvious commutative diagram of crossed modules (replace the complexes of tori $\tilde{I}_{1\mathbf{ab}}[1]$, $G_{\mathbf{ab}}$ by the quasi-isomorphic complexes $(\tilde{I}^{\uc}_1\rightarrow \tilde{I}_1)[1]$, $\tilde{I}_1\rightarrow I_1$); in particular, the map $\overline{\mathbf{ab}}$ is an isomorphism. The map $\lambda$ comes from the long exact sequence for $\tilde{I}_{1\mathbf{ab}}$:
\[ \cdots \longrightarrow \mathbb{H}^1(\Q,\tilde{I}_{1\mathbf{ab}}) \longrightarrow \mathbb{H}^1(\A,\tilde{I}_{1\mathbf{ab}}) \stackrel{\lambda}{\longrightarrow} \mathbb{H}^1(\A/\Q,\tilde{I}_{1\mathbf{ab}})\longrightarrow \mathbb{H}^{2}(\Q,\tilde{I}_{1\mathbf{ab}}) \longrightarrow \cdots, \] 
and $\nu$ is the dual of the inclusion $\mathfrak{K}(I_1/F) \hookrightarrow \pi_0(Z(\hat{\tilde{I}}_1)^{\Gamma})$ under the identification $\mathbb{H}^1(\A/\Q,\tilde{I}_{1\mathbf{ab}})=\pi_0(Z(\hat{\tilde{I}}_1)^{\Gamma})^D$ (Lemma \ref{lem:identification_of_Kottwitz_A(H)}). 
Since $\im[H^1(\Q,\tilde{I}_{1})\rightarrow \bigoplus_v H^1(\Q_v,\tilde{I}_{1})]=\ker(\lambda\circ\mathbf{ab})$ \cite[Thm.5.16]{Borovoi98} (cf. \cite[Prop.2.6]{Kottwitz85}), hence it suffices to find an adelic cohomology class 
\[ [\tilde{a}(v)]\in \oplus_vH^1(\Qv,\tilde{I}_{1}) \] 
such that $\alpha(v)=f_v([\tilde{a}(v)])$ for each $v$ and $\lambda\circ\mathbf{ab}([\tilde{a}(v)])=0$.
We find such class $([\tilde{a}(v)])_v$ as follows. We first note that $(\alpha(v))$ lifts to an element $\tilde{\alpha}=(\tilde{\alpha}(v))_v\in \bigoplus_v \mathbb{H}^1(\Qv,\tilde{I}_{1\mathbf{ab}})$. In view of the cohomology exact sequence attached to the distinguished triangle (\ref{eq:DT_of_CX_of_tori}), this amounts to vanishing of the image of $(\alpha(v))_v$ in $\bigoplus_v \mathbb{H}^1(\Qv,G_{\mathbf{ab}}) (\subset \bigoplus_v X^{\ast}(Z(\hat{G})^{\Gamma_v})_{\tors})$, which then follows from the fact (\ref{eq:restriction_of_alpha_to_Z(hatG)}) that for any Kottwitz triple $(\epsilon;\gamma,\delta)$, the image of the invariant $\alpha_v(\epsilon;\gamma,\delta;g_v)$ in $X^{\ast}(Z(\hat{G})^{\Gamma_v})$ is independent of the triple.
Now, by assumption (of vanishing of Kottwitz invariants), we have $\nu\circ\lambda(\tilde{\alpha})=0$ in $\mathfrak{K}(I_1/\Q)^D$. 
Hence, by Lemma \ref{lem:proof_of_Kottwitz86_Thm.6.6} below (which asserts exactness of the third column of (\ref{eq:H^1_{ab}_for_(I'->I)})), we may further assume that $(\tilde{\alpha}(v))_v\in \ker(\lambda)$.
Since the (local) abelianization map $\mathbf{ab}_v$ is surjective for every place $v$ (even bijective for finite $v$) \cite[Thm.5.4]{Borovoi98}, we can find $([\tilde{a}(v)])_{v\neq\infty} \in \bigoplus_{v\neq\infty} H^1(\Qv,\tilde{I}_{1})$ mapping to $(\tilde{\alpha}(v))_{v\neq\infty}$. For $v=\infty$, we use the condition that $\alpha(\infty)=0$, by which and the exactness of the vertical sequence for $\tilde{I}_{1\mathbf{ab}}\rightarrow I_{1\mathbf{ab}}$, we find an element $\bar{\alpha}(\infty)\in \mathbb{H}^0(\R,\tilde{I}_{1\mathbf{ab}}\rightarrow I_{1\mathbf{ab}})$ mapping to $\tilde{\alpha}(\infty)$. Then, using that $\overline{\mathbf{ab}}$ is an isomorphism, we define $\tilde{a}(\infty)$ to be the the image of $\bar{\alpha}(\infty)$ under the map $\mathbb{H}^0(\R,\tilde{I}_{1}\rightarrow I_{1})\rightarrow H^1(\R,\tilde{I}_1)$. The cohomology class $([\tilde{a}(v)])_v$ thus found is the one that we are looking for. 

Next, we observe that for any admissible pair $(\mathscr{I},\epsilon)$ and each 
\[ \beta\in \Sha_G(\Q,I_{\mathscr{I},\epsilon})^+=\mathrm{im}[\widetilde{\Sha}_G(\Q,I_{\mathscr{I},\epsilon})^+\rightarrow H^1(\Q,I_{\mathscr{I},\epsilon})],\] 
the stable Kottwitz triples $(\gamma_0;\gamma,\delta)$, $(\gamma_0';\gamma',\delta')$ attached to $(\mathscr{I},\epsilon)$, $(\mathscr{I}^{\beta},\epsilon^{\beta})$ are the same; recall (Lemma \ref{lem:|widetilde{Sha}_G(Q,I_{phi,epsilon})^+|}) that $\widetilde{\Sha}_G(\Q,I_{\mathscr{I},\epsilon})^+=\ker[\Sha^{\infty}_G(\Q,I_{\mathscr{I},\epsilon}^{\mathrm{o}})\rightarrow H^1(\A,I_{\mathscr{I},\epsilon})]$. This follows from Corollary \ref{cor:K-triple_of_twisted_K-pair}. Indeed, one may lift $\beta$ to $\Sha^{\infty}_G(\Q,I_{\mathscr{I},\epsilon}^{\mathrm{o}})$ and assume $\gamma_0'=\gamma_0$. Then, with the notation of that corollary, the ``difference'' of the equivalence class of $\gamma_l' (l\neq p)$ (resp. of $\delta'$) and that of $\gamma_l$ (resp. of $\delta$) is measured by the image in $H^1(\Qv,G_{\gamma_0})$ of the cohomology class $\alpha_l(\gamma_0;\gamma'_l;g_{l}') -\alpha_l(\gamma_0;\gamma_{l};g_{l})$ (resp. $\kappa_{I_0}([b(\gamma_0;\delta';g_p')]_{I_0})-\kappa_{I_0}([b(\gamma_0;\delta;g_p)]_{I_0})$) in $H^1(\Qv,I_0)$ (see Subsection \ref{subsec:Kottwitz_invariant}).

It remains to show that the image of the inclusion (Theorem \ref{thm:Kisin17_Prop.4.4.8} (vii))
\[ \Sha_G(\Q,I_{\mathscr{I},\epsilon})^+\hookrightarrow \{\text{ admissible K-pairs }\}/\sim \ :\  \beta\mapsto (\mathscr{I}^{\beta},\epsilon^{\beta})\] 
exhausts all the admissible K-pairs whose associated stable Kottwitz triples are equivalent to that of $(\mathscr{I},\epsilon)$, namely that if $(\mathscr{I},\epsilon)$ and $(\mathscr{I}',\epsilon')$ are two admissible K-pairs having equivalent associated stable Kottwitz triples, one has $(\mathscr{I}',\epsilon')=(\mathscr{I}^{\beta},\epsilon^{\beta})$ for some $\beta\in \Sha_G(\Q,I_{\mathscr{I},\epsilon})^+$. By Theorem \ref{thm:equiv_K-triples}, there exist a nice maximal $\Q$-torus $T\subset I_{\mathscr{I},\epsilon}^{\mathrm{o}}$, $\Q$-embeddings $i_T:T\hookrightarrow G$, $ i_T'=\Int(g)\circ i_T:T\hookrightarrow G$ such that the cohomology class $\tilde{\beta}=[i_T^{-1}(g^{-1}\tau(g))]\in \Sha^{\infty}_G(\Q,T)$ belongs to $\Sha_G^{\infty}(\Q,i_T(T))$ (since for some $h\in X$, $h'=\Int(g)(h)\in X$), and that the admissible K-pair $(\mathscr{I}',\epsilon')$ is the twist $(\mathscr{I}^{\beta},\epsilon^{\beta})$ of $(\mathscr{I},\epsilon)$ by the image $\beta$ of $\tilde{\beta}$ in $\im[\Sha^{\infty}_G(\Q,I_{\mathscr{I},\epsilon}^{\mathrm{o}})\rightarrow H^1(\Q,I_{\mathscr{I},\epsilon})]$.
Given this, the fact that $\beta$ vanishes in $\prod_{v\neq\infty} H^1(\Q_v,I_{\mathscr{I},\epsilon})$ again follows by the same reasoning as above (deduced from  Corollary \ref{cor:K-triple_of_twisted_K-pair}): the image of $\beta$ in $\prod_{v\neq\infty}H^1(\Q_v,I_{\mathscr{I},\epsilon})$ measures the ``difference'' of the equivalence class of $(\gamma_0;\gamma',\delta')$ and that of $(\gamma_0;\gamma,\delta)$.
This completes the proof.
\end{proof}

%%%%%%%%%%%%%%%%%%%%
\begin{lem} \label{lem:proof_of_Kottwitz86_Thm.6.6}
Let $G$ be a connected reductive group over a number field $F$ and $I$ an $F$-Levi subgroup of $G$.
Set $\tilde{I}:=\rho^{-1}(I)$ for the canonical homomorphism $\rho:G^{\uc}\rightarrow G$, $\Gamma:=\Gal(\overline{F}/F)$, and $\Gamma_v:=\Gal(\bar{F}_v/F_v)$ for any place $v$ of $F$.
The kernel of the natural map 
\begin{equation}
\nu:\mathbb{H}^1(\A_F/F,\tilde{I}_{\mathbf{ab}})=\pi_0(Z(\hat{\tilde{I}})^{\Gamma})^D\rightarrow \mathfrak{K}(I/F)^D
\end{equation}
equals the image of $\lambda:\mathbb{H}^1(\A_F,\tilde{I}_{\mathbf{ab}}) \rightarrow \mathbb{H}^1(\A_F/F,\tilde{I}_{\mathbf{ab}})$ of the kernel of
\[\mu:\mathbb{H}^1(\A_F,\tilde{I}_{\mathbf{ab}})\rightarrow  \mathbb{H}^1(\A_F,I_{\mathbf{ab}}). \]
In other words, we have an equality
\[\mathrm{coker}[\mathbb{H}^0(\A_F,G_{\mathbf{ab}})\rightarrow \mathbb{H}^1(\A_F/F,\tilde{I}_{\mathbf{ab}})]=\mathfrak{K}(I/F)^D.\]
\end{lem}

%%%%%%%%%%%%%%%%%%%%
\begin{rem} \label{eq:K_Labesse_Kottwitz}
(1) It is easy to see that for a reductive $H$ over $F$ and for each place $v$ of $F$, the map $\mathbb{H}^1(F_v,H)\rightarrow \pi_0(Z(\hat{H})^{\Gamma_v})^D$ that was constructed by Kottwitz in \cite[Thm.1.2]{Kottwitz86} factors through the abelianization map $\mathbb{H}^1(F_v,H)\rightarrow \mathbb{H}^1(F_v,H_{\mathbf{ab}})$ and the induced map 
\[\mathbb{H}^1(F_v,H_{\mathbf{ab}})\rightarrow \pi_0(Z(\hat{H})^{\Gamma_v})^D.\]
is a monomorphism (and an isomorphism for non-archimedean $v$), equal to the map in Lemma \ref{lem:identification_of_Kottwitz_A(H)} (the point is that both monomorphisms are induced by Tate-Nakayama duality); this also equals the monomorphism constructed by Borovoi \cite[Prop.4.1, 4.2]{Borovoi98}. 
If $H^{\der}=H^{\uc}$, this map for archimedean $v$ is also surjective, because then $H_{\mathbf{ab}}$ is quasi-isomorphic to an \emph{$F$-torus} whose dual torus is $Z(\hat{H})$ (i.e. $H^{\ab}=H/H^{\der}$), so that Tate-Nakayama duality for tori applies.

(2) This lemma is mentioned, without proof, in \emph{Remarque} on p.43 of \cite{Labesse99} ($\mathbb{H}^1(\A_F/F,\tilde{I}_{\mathbf{ab}})$ is the same as $\mathbb{H}^0_{\mathbf{ab}}(\A_F/F,I\backslash G):=\mathbb{H}^0(\A_F/F,I_{\mathbf{ab}}\rightarrow G_{\mathbf{ab}})$ defined there, cf. \cite[Prop.1.8.1]{Labesse99} and (\ref{eq:DT_of_CX_of_tori})). Earlier, this was also observed in the proof of Theorem 6.6 of \cite{Kottwitz86} when $G^{\der}=G^{\uc}$, in which case, according to (1), the maps $\lambda$, $\mu$ become the natural maps: 
\[\lambda:\oplus_v \pi_0(Z(\hat{\tilde{I}})^{\Gamma_v})^D \rightarrow \pi_0(Z(\hat{\tilde{I}})^{\Gamma})^D,\quad \mu: \oplus_v \pi_0(Z(\hat{\tilde{I}})^{\Gamma_v})^D \rightarrow \oplus_v \pi_0(Z(\hat{I})^{\Gamma_v})^D\]
which are defined only in terms of $Z(\hat{\tilde{I}})$, $Z(\hat{I})$, and which are the definitions for these maps used in \cite[Thm.6.6]{Kottwitz86}.
\end{rem}

Since we could not find a proof (for the general case) in literatures, we present a proof. In our proof, we work systematically with the abelianized cohomology groups $H^1(k,H_{\mathbf{ab}}(C_k))$ (for $H=I,\tilde{I},G$) instead of $\pi_0(Z(\hat{H})^{\Gal(\bar{k}/k)})^D$. Also, we use Galois (hyper)cohomology of complexes of $\Q$-tori of length $2$, especially their Poitou-Tate-Nakayama local/global dualities, as expounded in \cite{Borovoi98}, \cite{Demarche11}, \cite[Appendix A]{KottwitzShelstad99}, \cite[Ch.1]{Labesse99}. 

\begin{proof} (of Lemma \ref{lem:proof_of_Kottwitz86_Thm.6.6})
We have two diagrams (the left diagram defines $\mathfrak{K}(I/F)$ by being cartesian and the right one is its dual):
\begin{equation}  \label{eq:dual_of_connecting_map_Z(hat{G})}
\xymatrix{ \pi_0(Z(\hat{\tilde{I}})^{\Gamma}) \ar[r]^(.44){\partial} & H^1(F,Z(\hat{G})) & \pi_0(Z(\hat{\tilde{I}})^{\Gamma})^D \ar@{->>}[d]_{\nu} & H^1(F,Z(\hat{G}))^D \ar[l]_{\partial^D} \ar@{->>}[d]^{i^D} \\ \mathfrak{K}(I/F) \ar@{^(->}[u] \ar[r] & \ker^1(F,Z(\hat{G})) \ar@{^(->}[u]_{i} & \mathfrak{K}(I/F)^D & \ker^1(F,Z(\hat{G}))^D \ar[l] },
\end{equation}
from which we see the equality:
\begin{align*}
\ker(\nu)&=\mathrm{im}(\partial^D)(\ker(i^D)).
\end{align*}
We need to describe the maps $\partial^D$, $i^D$ in terms of the complexes $G_{\mathbf{ab}}$, $\tilde{I}_{\mathbf{ab}}$. For that, we recall some facts.

For a connected reductive group $G$ over a field $k$ and a Levi $k$-subgroup $I$ of $G$, the map (\ref{eq:boundary_map_for_center_of_dual})
\[\partial:\pi_0(Z(\hat{\tilde{I}})^{\Gamma})\rightarrow H^1(F,Z(\hat{G}))\] 
is identified in a natural manner with the connecting homomorphism 
\[\widehat{\mathbb{H}}^{-1}(k'/k,\hat{\tilde{I}}_{\mathbf{ab}}) \rightarrow \widehat{\mathbb{H}}^0(k'/k,\hat{G}_{\mathbf{ab}}),\]
in the long exact sequence of Tate-cohomology arising from the distinguished triangle (\ref{eq:DT_of_dualCX_of_tori}), where $k'$ is a finite Galois extension of $k$ splitting $T$.
Indeed, we recall \cite[Lem.2.2, Cor.2.3]{Kottwitz84b} that the given map
is the connecting homomorphism $\mathrm{Ext}^1_{\Gamma}(X^{\ast}(Z(\hat{\tilde{I}})),\Z)\rightarrow \mathrm{Ext}^2_{\Gamma}(X^{\ast}(Z(\hat{G})),\Z)$ in the long exact sequence of the cohomology $\mathrm{Ext}^{\bullet}_{\Gamma}(-,\Z)$ for the exact sequence of $\Gamma$-modules: $0\rightarrow X^{\ast}(Z(\hat{\tilde{I}}))\rightarrow X^{\ast}(Z(\hat{I}))\rightarrow X^{\ast}(Z(\hat{G}))\rightarrow0$. As was shown in \textit{loc. cit.} (and using that $\pi_0(D^{\Gamma})=\widehat{\mathbb{H}}^{0}(\Gamma,D)$), one has a functorial isomorphism $\mathrm{Ext}^n_{\Gamma}(X^{\ast}(D),\Z)=\widehat{H}^{n-1}(\Gamma,D)\ (n\geq1)$ for any diagonalizable $\C$-group $D$ with $\Gamma$-action. 
Since $Z(\hat{\tilde{I}})[1]=\hat{\tilde{I}}_{\mathbf{ab}}$, $Z(\hat{G})[1]=\hat{G}_{\mathbf{ab}}$ (\ref{eq:center_of_complex_dual}), the claim follows.

For a connected reductive group $H$ over $F$ (number field), we $j^i_{H}$ denote the injection (arising from a suitable cohomology long exact sequence for $H_{\mathbf{ab}}$):
\[j^i_{H}:\mathbb{H}^i(\A_F,H_{\mathbf{ab}})/\mathbb{H}^i(F,H_{\mathbf{ab}}) \hookrightarrow \mathbb{H}^i(\A_F/F,H_{\mathbf{ab}}).\]

Next, there exist natural isomorphisms \cite[11.2.2, 4.2.2]{Kottwitz84b}, \cite[D.2.C]{KottwitzShelstad99}, \cite[Thm.5.13]{Borovoi98}
\begin{align} \label{eq:identification_of_TS-gps}
\ker^1(F,Z(\hat{G}))^D=\ker^1(F,G) =\ker^1(F,G_{\mathbf{ab}}),
\end{align}
where for $i\geq0$, we define (cf. \cite[(C.1)]{KottwitzShelstad99}, \cite[1.4]{Labesse99})
\begin{align*}
\ker^i(F,G_{\mathbf{ab}})&:=\ker(H^i(F,G_{\mathbf{ab}})\rightarrow H^i(\A_F,G_{\mathbf{ab}})) \\
&=\ker(H^i(F,G_{\mathbf{ab}})\rightarrow \oplus_vH^i(F_v,G_{\mathbf{ab}})) \ \text{ if }i\geq1
\end{align*}
(When $i\geq1$, there exists a canonical isomorphism $H^i(\A_F,G_{\mathbf{ab}})= \oplus_vH^i(F_v,G_{\mathbf{ab}})$, \cite[Lem.4.5]{Borovoi98}, \cite[Lem.C.1.B]{KottwitzShelstad99}, \cite[Prop.1.4.1]{Labesse99}).

%%%%%%%%%%%%%%%%%%%%
\begin{lem} \label{lem:commutativity_of_duality_diagram} 
There exists a commutative diagram induced by the global and local Tate-Nakayama dualities:
\[ \xymatrix{ \ker^1(F,G_{\mathbf{ab}}) \ar[r] & \ker^2(F,\tilde{I}_{\mathbf{ab}}) \\
H^1(F,Z(\hat{G}))^D \ar[r]^(0.55){\partial^D} \ar@{->>}[u]^{i^D} & \pi_0(Z(\hat{\tilde{I}})^{\Gamma})^D  \ar@{->>}[u] \\ 
\mathbb{H}^0(\A_F,G_{\mathbf{ab}})/\mathbb{H}^0(F,G_{\mathbf{ab}}) \ar[r]^{\bar{\xi}} \ar@{^(->}[u]^{j^0_G} & \mathbb{H}^1(\A_F,\tilde{I}_{\mathbf{ab}})/\mathbb{H}^1(F,\tilde{I}_{\mathbf{ab}})  \ar@{^(->}[u]_{j^1_{\tilde{I}}} } \] 
where $\bar{\xi}$ is induced by the connecting homomorphism $\xi:\mathbb{H}^0(\A_F,G_{\mathbf{ab}})\rightarrow \mathbb{H}^1(\A_F,\tilde{I}_{\mathbf{ab}})$ in the long exact sequence of cohomology attached to the distinguished triangle (\ref{eq:DT_of_CX_of_tori}).

The two vertical sequences are short exact sequences. 
\end{lem}

\begin{proof} In view of the discussion above and the isomorphism (\ref{eq:connecting_isom_diagonal_gp}), $\partial^D$ equals the dual of the obvious connecting homomorphism
\[\widehat{\mathbb{H}}^1(k'/k,X^{\ast}(G_{\mathbf{ab}}))^D \longrightarrow \widehat{\mathbb{H}}^{0}(k'/k,X^{\ast}(\tilde{I}_{\mathbf{ab}}))^D, \]
Hence, the existence and exactness of each vertical sequence is easily deduced by reading the relevant parts of the Poitou-Tate exact sequence for two-term complexes of tori \cite[Thm.6.1]{Demarche11} (take $\hat{C}=X^{\ast}(H_{\mathbf{ab}})$ for $H=\tilde{I}, G$ in \textit{loc. cit.}): to identify the upper diagram, use the identification (\ref{eq:identification_of_TS-gps}) and the fact \cite[Lem.C.3.B, C.3.C]{KottwitzShelstad99}) that
\[\ker^1(F,G_{\mathbf{ab}})=[\ker^0(F,\hat{G}_{\mathbf{ab}})_{\mathrm{red}}]^D=[\ker^1(F,X_{\ast}(\hat{G}_{\mathbf{ab}}))]^D.\]
The commutativity of the daigram follows from the compatibility of the global and local Tate-Nakayama dualities.
\end{proof}

Now, as $\mathrm{im}(\xi)=\ker(\mu)$, by Lemma \ref{lem:commutativity_of_duality_diagram} we see that
\begin{align*}
\ker(\nu)&=\mathrm{im}(\partial^D\circ j^0_G)=\mathrm{im}(\nu\circ\xi)\\
&=\mathrm{im}(\nu)(\ker\mu) 
\end{align*}
as was asserted. 
\end{proof}

%%%%%%%%%%%%%%%%%%%%
\begin{rem} \label{rem:Zhou}
In fact, the geometric ingredients used to establish the four theorems in this subsection, thus their three corollaries as well, are also available for more general parrahoric levels and the Kisin-Pappas integral model \cite{KisinPappas18}
under the assumption considered by Zhou \cite{Zhou18} that $G_{\Qp}$ is tamely ramified, $p\nmid \pi_1(G^{\der})|$ (which he called ($\ast$)) and $G_{\Qp}$ is residually split. More explicitly, 
\begin{itemize}
\item[(I)] \cite[Prop.2.1.3]{Kisin17} and \cite[Prop.1.4.15]{Kisin17} generalize to \cite[Prop.6.4]{Zhou18} and \cite[Cor.6.3]{Zhou18}, respectively;
\item[(SCM)] \cite[Cor.2.2.5]{Kisin17} is generalized in the argument of proving \cite[Thm.9.4]{Zhou18};
\item[(Ta)] \cite[Cor.2.3.2]{Kisin17} generalizes to \cite[Cor.9.5]{Zhou18}.
\end{itemize}
The twisting method is not worked out in \cite{Zhou18} (though some basic one is done in \cite[4.4]{KisinPappas18}), but it is not difficult to see that the twisting arguments from \cite{Kisin17} that we needed in this subsection require, as geometric inputs, only the moduli interpretation of isogeny classes, namely \cite[Prop.1.4.15]{Kisin17} which is generalized to \cite[Cor.6.3]{Zhou18}. Then, given these four ingredients and Corollary \ref{cor:stable_conjugacy_from_nice_tori}, we deduced (by standard arguments) the effectivity criterion of Kottwitz triple from Theorem \ref{thm:LR-Satz5.21} which is also generalized in \cite{Lee18a} for \emph{special maximal parahoric} level subgroups under the assumption that $G_{\Qp}$ is tamely ramified, quasi-split.
Consequently, the results in this section hold under the conditions satisfying both \cite{Zhou18} and \cite{Lee18a}.

On the other hand, to have a complete description of the Hasse-Weil zeta function of a Shimura variety (and for some other applications, such as construction of Galois representations, cf. \cite{ScholzeShin13}), one also needs formulae at bad reductions for the semi-simple Lefschetz number similar to (\ref{eq_intro:Kottwitz_formula}) but with the test function $\phi_p$ there being replaced by the trace function of Frobenius automorphism acting on a relevant nearby cycle sheaf (cf. \cite[$\S$10]{Rapoport05}, \cite[Conj.4.30, 4.31]{Haines14}). 
Then, we believe that the above results would be key ingredients in obtaining such formulae under relevant assumptions. 
%But, we can already utilize those partial results in \cite{Lee18a} to extend the scope of some previous works, for example, we can relax the ramification condition imposed on PEL datum in the main result of \cite{Scholze13}. 
\end{rem}

%%%%%%%%%%%%%%%%%%%%
%%%%%%%%%%%%%%%%%%%%
\subsection{Proof of Kottwitz formula}

Now we give an expression for the number of fixed points of the correspondence $\Phi^m\circ f$ (\ref{eq:Hecke_corr_twisted_by_Frob}), and, more generally, a weighted sum over the same fixed point set with weight being given by the trace of that correspondence acting on the stalk of an automorphic lisse sheaf. We briefly recall the setup. For more details, see \cite[$\S$6, 16]{Kottwitz92}.
We fix a rational prime $l\neq p$. Let $\xi$ be a finite-dimensional representation of $G$ on a vector space $W$ over a number field $L$, and let $\lambda$ be a place of $L$ lying above $l$. We consider the $\lambda$-adic representation $W_{\lambda}:=W\otimes_L{L_{\lambda}}$ of $G(\A_f^p)$ induced from the natural one of $G(\Ql)$ via the projection $G(\A_f^p)\rightarrow G(\Ql)$. Since $Z(\Q)$ is discrete in $Z(\A_f)$ by our assumption ($Z(G)$ has same ranks over $\Q$ and $\R$), it gives rise to lisse sheaves $\sF_{\mathbf{K}^p}$ on the spaces $\sS_{\mathbf{K}^p}$ for varying $\mathbf{K}^p$'s: the projective limit $\sS=\varprojlim_{H^p}\sS_{H^p}$ is a Galois covering of $\sS_{\mathbf{K}^p}$ with Galois group $\mathbf{K}^p$ if $\mathbf{K}^p$ is small enough such that $K\cap Z(G)(\Q)=\{1\}$ \cite[2.1.9-12]{Deligne79}, and we have $\sF_{\mathbf{K}^p}=\sS\times W_{\lambda}/\mathbf{K}^p$, where $k=(k_v)\in \mathbf{K}^p$ acts on $W_{\lambda}$ via $\xi_{L_{\lambda}}(k_l^{-1})$. It is clear that there are canonical isomorphisms 
\[\cdot g:\sF_{g\mathbf{K}^pg^{-1}}\isom (\cdot g)^{\ast}\sF_{\mathbf{K}^p},\quad p_2^{\ast}\sF_{\mathbf{K}^p} =\sF_{\mathbf{K}^p_g},\quad \Phi:\Phi^{\ast}\sF_{\mathbf{K}^p}\isom \sF_{\mathbf{K}^p},\]
where for any $g\in G(\A_f^p)$, $\cdot g$ denotes the right action $\sS_{g\mathbf{K}^pg^{-1}}\isom\sS_{\mathbf{K}^p}$.
Thus, the associated correspondence $\Phi^m\circ f$ (\ref{eq:Hecke_corr_twisted_by_Frob}) extends in a natural manner to the sheaf $\sF_{\mathbf{K}^p}$. In particular, for any fixed point $x'\in \sS_{\mathbf{K}^p_g}(\F)$ of $\Phi^m\circ f$ and $x:=p_1'(x')=p_2'(x')\in \sS_{\mathbf{K}^p}(\F)$, the correspondence $\Phi^m\circ f$ gives an automorphism of the stalk $\sF_x$ (we write $\sF$, $\sF'$ for $\sF_{\mathbf{K}^p}$, $\sF_{\mathbf{K}^p_g}$): 
\begin{equation} \label{eq:Frob-Hecke_corr_at_stalk}
\sF_x=(\Phi^{m})^{\ast}(\sF)_{p_2(x')}\  \stackrel{\Phi^m}{\longrightarrow}\ \sF_{p_2(x')} =(p_2^{\ast}\sF)_{x'}=\sF'_{x'}\ \stackrel{p_1^{\ast}(\cdot g)}{\longrightarrow}\ (p_1'^{\ast}\sF)_{x'}=\sF_x.
\end{equation}
We are interested in computing the sum:
\begin{equation} \label{eq:fixed-pt_set_of_Frob-Hecke_corr}
T(m,f):=\sum_{x'\in\mathrm{Fix}} \mathrm{tr}(\Phi^m\circ f;\sF_x),
\end{equation}
where $\mathrm{Fix}$ denotes the set of fixed points of $\Phi^m\circ f$.

Set $\mathbf{K}_q:=G(L_n)\cap \tilde{\mathbf{K}}_p$ ($q=p^n$).
Define $\phi_p$ to be the characteristic function of the subset (union of double cosets) of $\mathbf{K}_q\backslash G(L_n)/\mathbf{K}_q$ corresponding to the set:
\begin{equation} \label{eq:Adm_K_q(mu)}
\mathrm{Adm}_{\mathbf{K}_q}(\{\mu\}):=\{w\in\mathbf{K}_q\backslash G(L_n)/\mathbf{K}_q\ |\ w\leq \tilde{\mathbf{K}}_p\mu(p^{-1})\tilde{\mathbf{K}}_p\}.
\end{equation}
Here, in the inequality, $w$ also denotes its image in $\tilde{\mathbf{K}}_p\backslash G(\mfk)/\tilde{\mathbf{K}}_p$.
Let $dy_p$ (resp. $dy^p$) denote the Haar measure on $G(L_n)$ (resp. on $G(\A_f^p)$) giving measure $1$ on $\mathbf{K}_q$ (resp. on $\mathbf{K}^p$).

%%%%%%%%%%%%%%%%%%%%
\begin{lem} \label{lem:fixed-pt_subset_of_Frob-Hecke_corr}
Assume that $Z(G)$ has same ranks over $\Q$ and $\R$.
We take $\mathbf{K}^p$ small enough so that conditions of (\ref{item:Langlands-conditions}) hold and $Z(G)(\Q)\cap K=\{1\}$.
Let $(\gamma_0;\gamma,\delta)$ be a Kottwitz triple attached to an admissible K-pair $(\mathscr{I},\epsilon)$.
Then, the cardinality of the set $I_{\mathscr{I},\epsilon}(\Q)\backslash X(\mathscr{I},\epsilon)_{\mathbf{K}_g}$ (\ref{eq:fixed_pt_set_of_Heck-corresp2}) is given by
\[|I_{\mathscr{I},\epsilon}(\Q)\backslash X(\mathscr{I},\epsilon)_{\mathbf{K}_g}|= \frac{\mathrm{vol}(I_{\mathscr{I},\epsilon}^{\mathrm{o}}(\Q)\backslash I_{\mathscr{I},\epsilon}^{\mathrm{o}}(\A_f))}{[I_{\mathscr{I},\epsilon}(\Q):I_{\mathscr{I},\epsilon}^{\mathrm{o}}(\Q)]} \cdot \mathrm{O}_{\gamma}(f^p)\cdot \mathrm{TO}_{\delta}(\phi_p),\]
where $\mathrm{TO}_{\delta}(\phi_p)$ (twisted orbital integral) and $\mathrm{O}_{\gamma}(f^p)$ (orbital integral) are defined by: 
\begin{align} \label{eq:(twisted-)orbital_integral}
\mathrm{TO}_{\delta}(\phi_p)&=\int_{G_{\delta\sigma}^{\mathrm{o}}(\Qp)\backslash G(L_n)}\phi_p(y_p^{-1}\delta\sigma(y_p)) d\bar{y}_p,\\
\mathrm{O}_{\gamma}(f^p)&=\int_{G_{\gamma}^{\mathrm{o}}(\A_f^p)\backslash G(\A_f^p)}f^p((y^p)^{-1}\gamma y^p) d\bar{y}^p, \nonumber
\end{align}
Here, $f^p$ is the characteristic function of $\mathbf{K}^pg^{-1}\mathbf{K}^p$ in $G(\A_f^p)$, and the Haar measure $di_p$ on $G_{\delta\sigma}^{\mathrm{o}}(\Qp)$ (resp. $di^p$ on $G_{\gamma}^{\mathrm{o}}(\A_f^p)$) is obtained from a Haar measure on $I_{\mathscr{I},\epsilon}^{\mathrm{o}}(\Qp)$ (resp. on $I_{\mathscr{I},\epsilon}^{\mathrm{o}}(\A_f^p)$) that gives rational measures to compact open subgroups of $I_{\mathscr{I},\epsilon}^{\mathrm{o}}(\Qp)$ (resp. of $I_{\mathscr{I},\epsilon}^{\mathrm{o}}(\A_f^p)$) via the isomorphism $\Int(c)\circ i_p$ (\ref{eq:isom_Int(c)circi_p}) (resp. via $i^p$ (\ref{eq:isom_i^p})). We write $d\bar{y}_p$ $d\bar{y}^p$ for the quotient of $dy_p$ by $di_p$ and that of $dy^p$ by $di^p$, respectively.

\end{lem}

\begin{proof}
The argument given in \cite[$\S$1.4, $\S$1.5]{Kottwitz84b} (cf. \cite[p.432]{Kottwitz92}) works without change, since the necessary isomorphisms (1.4.7), (1.4.8) in \cite{Kottwitz84b} follow from Corollary \ref{cor:Tate_thm2}. 
Then, the cardinality in question becomes the triple product
\[ \mathrm{vol}(I_{\mathscr{I},\epsilon}(\Q)\backslash I_{\mathscr{I},\epsilon}(\A_f))\cdot \int_{G_{\gamma}(\A_f^p)\backslash G(\A_f^p)}f^p((y^p)^{-1}\gamma y^p) d\bar{y}^p\cdot \int_{G_{\delta\sigma}(\Qp)\backslash G(L_n)}\phi_p(y_p^{-1}\delta\sigma(y_p)) d\bar{y}_p \]
(the quotient measures $d\bar{y}^p$, $d\bar{y}_p$ being defined similarly).
The statement follows from this.
\end{proof}

Let $A_G$ denote the maximal $\Q$-split torus in the center of $G$. 
Recall that we assumed $G$ to be the smallest algebraic (connected reductive) $\Q$-group such that every $h\in X$ factors through $G_{\R}$ (i.e. $G$ is the Mumford-Tate group of a generic $h\in X$) so that in particular, $(G,X)$ satisfies the condition that $Z(G)$ has same ranks over $\Q$ and $\R$ (thus, $Z_G(\R)/A_G(\R)$ is also compact).

%%%%%%%%%%%%%%%%%%%%
%%%%%%%%%%%%%%%%%%%%
\begin{thm} \label{thm:Kottwitz_formula:Kisin}
Let $(G,X)$ be a Shimura datum of Hodge type such that $G_{\Qp}$ is unramified.
Fix a hyperspecial subgroup $\mathbf{K}_p$ and take $\mathbf{K}^p$ to be sufficiently small such that conditions (a), (b) of (\ref{item:Langlands-conditions}) hold and $K\cap Z(G)(\Q)=\{1\}$.

(1) We have the following expression for (\ref{eq:fixed-pt_set_of_Frob-Hecke_corr}):
\[T(m,f)=\sum_{(\gamma_0;\gamma,\delta)} c(\gamma_0;\gamma,\delta)\cdot \mathrm{O}_{\gamma}(f^p)\cdot \mathrm{TO}_{\delta}(\phi_p)\cdot \mathrm{tr}\xi(\gamma_0),\]
with
\[ c(\gamma_0;\gamma,\delta):=\tilde{i}(\gamma_0;\gamma,\delta)\cdot |\pi_0(G_{\gamma_0})(\Q)|^{-1} \cdot \tau(I_0)\cdot \mathrm{vol}(A_G(\R)^{\mathrm{o}}\backslash I_0(\infty)(\R))^{-1} \] 
where $I_0:=G_{\gamma_0}^{\mathrm{o}}$, $\tilde{i}(\gamma_0;\gamma,\delta)=|\widetilde{\Sha}_G(\Q,I_{\phi,\epsilon})^+|$ (Lemma \ref{lem:|widetilde{Sha}_G(Q,I_{phi,epsilon})^+|}), $\tau(I_0)$ is the Tamagawa number of $I_0$, and $I_0(\infty)$ is the (unique) inner form of $(I_0)_{\R}$ having compact adjoint group. Also, the sum is over a set of representatives $(\gamma_0;\gamma,\delta)$ of all equivalence classes of \emph{stable} Kottwitz triples of level $n=m[\kappa(\wp):\Fp]$ having trivial Kottwitz invariant.

(2) Then, for any $f^p$ in the Hecke algebra $\mathcal{H}(G(\A_f^p)/\!\!/ \mathbf{K}^p)$, there exists $m(f^p)\in\N$, depending on $f^p$, such that for each $m\geq m(f^p)$, we have
\begin{align} \label{eq:Lef-number1}
\sum_{i}(-1)^i\mathrm{tr}( & \Phi^m\times f^p | H^i_c(Sh_{\mathbf{K}}(G,X)_{\Qb},\sF_{\mathbf{K}})) \\
& = \sum_{(\gamma_0;\gamma,\delta)} c(\gamma_0;\gamma,\delta)\cdot \mathrm{O}_{\gamma}(f^p)\cdot \mathrm{TO}_{\delta}(\phi_p) \cdot \mathrm{tr}\xi(\gamma_0), \nonumber
\end{align}
where the sum is over a set of representatives $(\gamma_0;\gamma,\delta)$ of \emph{all} (stable) equivalence classes of \emph{stable} Kottwitz triples of level $n$ having trivial Kottwitz invariant. If $G^{\ad}$ is anisotropic or $f^p$ is the identity, we can take $m(f^p)$ to be $1$ (irrespective of $f^p$).
\end{thm}

In the definition of $c(\gamma_0;\gamma,\delta)$, the volume of (the quotient of) $\R$-groups is defined with respect to the unique Haar measure $di_{\infty}$ on $I_{\mathscr{I},\epsilon}^{\mathrm{o}}(\R)$ (or equivalently on $I_0(\infty)(\R)$ by transfer of measure) such that the product of $di^p,di_p,di_{\infty}$ is the canonical measure on $I_{\mathscr{I},\epsilon}^{\mathrm{o}}(\A)$ that is used to define the Tamagawa number $\tau(I_{\mathscr{I},\epsilon}^{\mathrm{o}})$  (cf. \cite[$\S$1.2]{Labesse01}). Then, as $Z_G(\R)/A_G(\R)$ is compact, one has 
\begin{equation} \label{eq:Tamagawa_number}
\mathrm{vol}(I_{\mathscr{I},\epsilon}^{\mathrm{o}}(\Q)\backslash I_{\mathscr{I},\epsilon}^{\mathrm{o}}(\A_f))=\tau(I)\cdot \mathrm{vol}(A_{G}(\R)^{\mathrm{o}}\backslash I_0(\R))^{-1}.
\end{equation}

\begin{proof}
We may assume that $f^p$ is the characteristic function of $\mathbf{K}^pg^{-1}\mathbf{K}^p$ of some $g\in G(\A_f^p)$.
For an isogeny class $\mathscr{I}$ and a comact open subgroup $\mathbf{K}^p$ of $G(\A_f^p)$, the fixed point subset
\[\mathrm{Fix}_{\mathscr{I}}:=S(\mathscr{I})_{\mathbf{K}}^{\Phi^m\circ f=\mathrm{Id}}=\{x'\in S_{\mathbf{K}^p_g}(\mathscr{I})\ |\ p_1'(x')=p_2'(x')\}\]
is a disjoint union $\sqcup_{\epsilon\in I_{\mathscr{I}}(\Q)} I_{\mathscr{I},\epsilon}(\Q)\backslash X(\mathscr{I},\epsilon)_{\mathbf{K}_g}$ (\ref{eq:fixed_pt_set_of_Heck-corresp2}). 
We claim that for any $x'\in \mathrm{Fix}_{\mathscr{I}}$ and $x:=p_1'(x')=p_2'(x')\in S_{\mathbf{K}^p}(\mathscr{I})$, there is an equality: 
\[\mathrm{tr}(\Phi^m\circ f;\sF_x)=\mathrm{tr}\xi(\gamma_0),\]
where $\gamma_0$ is the rational component of any Kottwitz triple attached to the admissible K-pair $(\mathscr{I},\epsilon)$ such that $x'\in I_{\mathscr{I},\epsilon}(\Q)\backslash X(\mathscr{I},\epsilon)_{\mathbf{K}_g}$ (i.e. $\gamma_0$ is stably conjugate to $i_T(\epsilon)$ for a $\Q$-embedding $i_T:T\hookrightarrow G$ of a maximal torus $T$ of $I_{\mathscr{I},\epsilon}$). In particular, this trace depends only on (the stable conjugacy class of) $\gamma_0$, so on the (equivalence class of) Kottwitz triples attached to the pair $(\mathscr{I},\epsilon)$.
Indeed (cf. \cite[$\S$16]{Kottwitz92}), choose $x_p\in X_p(\mathscr{I},\epsilon)$, $x^p\in X^p(\mathscr{I},\epsilon;g)$ such that $x'=[x_p,x^p]\in I_{\mathscr{I},\epsilon}(\Q)\backslash X(\mathscr{I},\epsilon)_{\mathbf{K}_g}$. It also gives a point $\tilde{x}=[x_p,x^p]$ of 
\[S(\mathscr{I}):=I_{\mathscr{I}}(\Q)\backslash X_p(\mathscr{I})\times X^p(\mathscr{I})=\varprojlim_{H^p}(\sS_{H^p}(\F)\cap S_{H^p}(\mathscr{I}))\]  
lying above $x'$, and $i^p(\epsilon) x^pgk=x^p$ for some $k\in \mathbf{K}^p$. So, one has
\[\Phi^m(\tilde{x})=[\Phi^mx_p,x^p]=[x_p,i^p(\epsilon) ^{-1} x^p]=\tilde{x}gk.\] 
If we use this point $\Phi^m(\tilde{x})$ of $S(\mathscr{I})$ to identify the stalk $\sF_x$ with $W_{\lambda}$,
we have $\beta(w)=\xi(k^{-1}g^{-1})w$: the automorphism (\ref{eq:Frob-Hecke_corr_at_stalk}) becomes
\[ [\Phi^m(\tilde{x}),w] \mapsto [\tilde{x},w] \mapsto [\tilde{x}g,\xi(g^{-1})w]= [\tilde{x}gk,\xi(k^{-1}g^{-1})w], \] 
hence $\mathrm{tr}(\Phi^m\circ f;\sF_x)=\mathrm{tr}\xi(\gamma_0)$ as $k^{-1}g^{-1}=(x^p)^{-1}i^p(\epsilon) x^p\in G(\A_f^p)$ is conjugate to $\gamma_0$ under $G(\bar{\A}_f^p)$.

Now, we have the following successive equalities:
\begin{align} \label{eq:T(m,f)1}
T(m,f) &=\sum_{\mathscr{I}} \sum_{x'\in\mathrm{Fix}_{\mathscr{I}}} \mathrm{tr}(\Phi^m\circ f;\sF_{p_1'(x')}) \\
& =\sum_{(\mathscr{I},\epsilon)} |I_{\mathscr{I},\epsilon}(\Q)\backslash X(\mathscr{I},\epsilon)_{\mathbf{K}_g}| \cdot \mathrm{tr}\xi(\gamma_0)   \nonumber \\ & =\sum_{(\mathscr{I},\epsilon)} \frac{\mathrm{vol}(I_{\mathscr{I},\epsilon}^{\mathrm{o}}(\Q)\backslash I_{\mathscr{I},\epsilon}^{\mathrm{o}}(\A_f))}{[I_{\mathscr{I},\epsilon}(\Q):I_{\mathscr{I},\epsilon}^{\mathrm{o}}(\Q)]} \cdot \mathrm{O}_{\gamma}(f^p)\cdot \mathrm{TO}_{\delta}(\phi_p) \cdot \mathrm{tr}\xi(\gamma_0) \nonumber
\end{align}
Here, in the first line, $\mathscr{I}$ runs through the set of isogeny classses, so the first equality results from (\ref{eq:isogeny_decomp}), (\ref{eq:fixed-pt_set_of_Frob-Hecke_corr}).
In the second line, $(\mathscr{I},\epsilon)$ runs through a set of representatives for the equivalence classes of admissible K-pairs. We have just seen that if $\epsilon\in I_{\mathscr{I}}(\Q)$ is such that $x'\in\mathrm{Fix}_{\mathscr{I}}$ belongs to the subset $I_{\mathscr{I},\epsilon}(\Q)\backslash X(\mathscr{I},\epsilon)_{\mathbf{K}_g}$ in the decomposition (\ref{eq:fixed_pt_set_of_Heck-corresp2}), we have $\mathrm{tr}(\Phi^m\circ f;\sF_{p_1'(x')})=\mathrm{tr}\xi(\gamma_0)$ for any Kottwitz triple $(\gamma_0;\gamma,\delta)$ attached to $(\mathscr{I},\epsilon)$, which gives the second equality. The third equality is Lemma \ref{lem:fixed-pt_subset_of_Frob-Hecke_corr}.

Next, we rewrite the last expression of (\ref{eq:T(m,f)1}) using (equivalence classes of) \emph{effective} stable Kottwitz triples as a new summation index.
For each equivalence class of effective stable Kottwitz triple $(\gamma_0;\gamma,\delta)$, we fix an admissible K-pair $(\mathscr{I}_1,\epsilon_1)$ giving rise to it. Then, the set of admissible pairs producing the same (stable) equivalence class of stable Kottwitz triple is in bijection with 
\[\Sha_G(\Q,I_{\mathscr{I}_1,\epsilon_1})^+=\im \left[ \widetilde{\Sha}_G(\Q,I_{\mathscr{I}_1,\epsilon_1})^+ \rightarrow H^1(\Q,I_{\mathscr{I}_1,\epsilon_1}) \right]\] (Corollary \ref{cor:LR-Satz5.25}, Lemma \ref{lem:|widetilde{Sha}_G(Q,I_{phi,epsilon})^+|}). More explicitly, for each 
\[ \beta\in \widetilde{\Sha}_G(\Q,I_{\mathscr{I}_1,\epsilon_1})^+=\ker\left[\Sha^{\infty}_G(\Q,I_{\mathscr{I}_1,\epsilon_1}^{\mathrm{o}})\rightarrow H^1(\A,I_{\mathscr{I}_1,\epsilon_1})\right],\] 
the admissible pair corresponding to the image of $\beta$ in $H^1(\Q,I_{\mathscr{I}_1,\epsilon_1})$ is the twist $(\mathscr{I},\epsilon):=(\mathscr{I}_1^{\beta},\epsilon_1^{\beta})$ (Theorem \ref{thm:Kisin17_Prop.4.4.8} (3))%%
\footnote{Again, this is an abuse of notation which will keep occurring in this proof: $\beta$ should be its image in $H^1(\Q,I_{\mathscr{I}_1,\epsilon_1})$.}
and the associated groups $I_{\mathscr{I},\epsilon}^{\mathrm{o}}$, $I_{\mathscr{I},\epsilon}$ are the (simultaneous) inner twists of $I_{\mathscr{I}_1,\epsilon_1}^{\mathrm{o}}$ and $I_{\mathscr{I}_1,\epsilon_1}$ via $\beta$ (Theorem \ref{thm:Kisin17_Prop.4.4.8} (2)); for simplicity, let us write $I_{\beta}^{\mathrm{o}}$ and $I_{\beta}$ for these twists.
Then, the last line of (\ref{eq:T(m,f)1}) becomes the first line of the following identity:
\begin{align} \label{eq:T(m,f)-2}
T(m,f) & =\sum_{(\gamma_0;\gamma,\delta)}\sum_{\beta\in \widetilde{\Sha}_G(\Q,I_{\mathscr{I}_1,\epsilon_1})^+} \frac{1}{|\ker[H^1(\Q,I_{\beta}^{\mathrm{o}})\rightarrow H^1(\Q,I_{\beta})]|} \cdot \frac{\mathrm{vol}(I_{\beta}^{\mathrm{o}}(\Q)\backslash I_{\beta}^{\mathrm{o}}(\A_f))}{[I_{\beta}(\Q):I_{\beta}^{\mathrm{o}}(\Q)]} \\  
& \qquad \qquad \qquad \qquad \qquad \qquad \qquad \qquad \qquad \qquad \cdot \mathrm{O}_{\gamma}(f^p)\cdot \mathrm{TO}_{\delta}(\phi_p) \cdot \mathrm{tr}\xi(\gamma_0) \nonumber \\
&=\sum_{(\gamma_0;\gamma,\delta)} c_1(\gamma_0)\cdot |\pi_0(G_{\gamma_0})(\Q)|^{-1} \cdot \tilde{i}(\gamma_0;\gamma,\delta) \cdot \mathrm{O}_{\gamma}(f^p)\cdot \mathrm{TO}_{\delta}(\phi_p) \cdot  \mathrm{tr}\xi(\gamma_0), \nonumber
\end{align}
where $c_1(\gamma_0):=\tau(I_0)\cdot \mathrm{vol}(A_{G}(\R)^{\mathrm{o}}\backslash I_0(\infty)(\R))^{-1}$.
Here, in the first line, the fist sum is over a set of representatives of the equivalence classes of \emph{effective} stable Kottwitz triples $(\gamma_0;\gamma,\delta)$ and in the second sum $\beta$ runs through the set $\widetilde{\Sha}_G(\Q,I_{\mathscr{I},\epsilon})^+$. Two elements $\beta,\beta'$ of $\widetilde{\Sha}_G(\Q,I_{\mathscr{I}_1,\epsilon_1})^+$ give equivalent admissible pairs $(\mathscr{I}_1^{\beta},\epsilon_1^{\beta})$, $(\mathscr{I}_1^{\beta'},\epsilon_1^{\beta'})$ if and only if $\beta$, $\beta'$ map to a same element in $H^1(\Q,I_{\mathscr{I}_1,\epsilon_1})$ (Theorem \ref{thm:Kisin17_Prop.4.4.8}). 
So, to establish the first equality, it suffices to prove that for each $\beta\in \widetilde{\Sha}_G(\Q,I_{\mathscr{I}_1,\epsilon_1})^+$, the set $S_{\beta}$ of those elements in $ \widetilde{\Sha}_G(\Q,I_{\mathscr{I}_1,\epsilon_1})^+$ that have the same image in $H^1(\Q,I_{\mathscr{I}_1,\epsilon_1})$ as $\beta$ is in bijection with $\ker[H^1(\Q,I_{\beta}^{\mathrm{o}})\rightarrow H^1(\Q,I_{\beta})]$; this will then also show that the latter set has the same size for all the elements in $S_{\beta}$, a set which really depends only on the image of $\beta$ in $H^1(\Q,I_{\mathscr{I}_1,\epsilon_1})$.
First, it is known \cite[Prop.35bis]{Serre02} that the set $\ker[H^1(\Q,I_{\beta}^{\mathrm{o}})\rightarrow H^1(\Q,I_{\beta})]$ is in bijection with the set $S_{\beta}'$ of such elements in $H^1(\Q,I_{\mathscr{I}_1,\epsilon_1}^{\mathrm{o}})$. So, we just need to show that the inclusion $S_{\beta}\subset S_{\beta}'$ is an equality, i.e. that
if $\beta'\in H^1(\Q,I_{\mathscr{I}_1,\epsilon_1}^{\mathrm{o}})$ has the same image in $H^1(\Q,I_{\mathscr{I}_1,\epsilon_1})$ as $\beta$, then we also have $\beta'\in \widetilde{\Sha}_G(\Q,I_{\mathscr{I}_1,\epsilon_1})^+$. Since $(I_{\mathscr{I}_1})_{\R}$ is a subgroup of the inner form $G'$ of $G_{\R}$ that has compact adjoint group, the map $H^1(\R,I_{\mathscr{I}_1,\epsilon_1}^{\mathrm{o}})\rightarrow H^1(\R,G')$ is injective \cite[4.4.5]{Kisin17}, which implies that $\beta'\in \Sha^{\infty}(\Q,I_{\mathscr{I}_1,\epsilon_1}^{\mathrm{o}})=\ker[ H^1(\Q,I_{\mathscr{I}_1,\epsilon_1}^{\mathrm{o}})\rightarrow H^1(\R,I_{\mathscr{I}_1,\epsilon_1}^{\mathrm{o}})]$. Similarly, the map $\Sha^{\infty}(\Q,I_{\mathscr{I}_1,\epsilon_1}^{\mathrm{o}})\cong \Sha^{\infty}(\Q,G_{\gamma_0}^{\mathrm{o}})\rightarrow \Sha^{\infty}(\Q,G)$ factors through $\Sha^{\infty}(\Q,I_{\mathscr{I}_1})$, which implies (as $\beta$ and $\beta'$ have the same image in $\Sha^{\infty}(\Q,I_{\mathscr{I}_1})$) that $\beta'\in \Sha^{\infty}_G(\Q,I_{\mathscr{I}_1,\epsilon_1}^{\mathrm{o}})$.

For the second equality, we use the following two facts (E1), (E2): 

(E1) There exists an equality of numbers:
\[|\ker[H^1(\Q,I_{\beta}^{\mathrm{o}})\rightarrow H^1(\Q,I_{\beta})]|\cdot [I_{\beta}(\Q):I_{\beta}^{\mathrm{o}}(\Q)] =|\pi_0(G_{\gamma_0})(\Q)|.\] 
Indeed, we recall \cite[5.5]{Serre02} that for any (not-necessarily connected) reductive $\Q$-group $I$ with connected component group $\pi_0(I)$ (especially for $I=I_{\mathscr{I}_1,\epsilon_1}$), the exact sequence $1\rightarrow I^{\mathrm{o}}\rightarrow I\rightarrow \pi_0(I) \rightarrow 1$ gives rise to a natural action $\pi_0(I)(\Q)$ on $H^1(\Q,I^{\mathrm{o}})$ which we normalize to be a left action and write $c\cdot \alpha$ for $c\in \pi_0(I)(\Q)$ and $\alpha \in H^1(\Q,I^{\mathrm{o}})$. 
One easily checks that for $c\in \pi_0(I)(\Q)$ and $\beta \in H^1(\Q,I^{\mathrm{o}})$, $c\cdot \beta$ equals the image of $c$ under the composite map 
\[ \pi_0(I)(\Q) =\pi_0(I_{\beta})(\Q) \stackrel{\partial_a}{\rightarrow}  H^1(\Q,I_{\beta}^{\mathrm{o}}) \isom H^1(\Q,I^{\mathrm{o}}),\] 
where for a choice of a cocycle representative $a$ of $\beta$, $\partial_a$ is the obvious coboundary map attached to the inner twist $I_{\beta}$ of $I$ via $a$ and the last bijection is defined by $[x_{\tau}]\mapsto [x_{\tau}a_{\tau}]$ (and thus sends the distinguished element to $[a_{\tau}]$) \cite[Prop.35bis.]{Serre02}.
So, the stabilizer subgroup of $\beta$ for the action of $\pi_0(I)(\Q)$ on $H^1(\Q,I^{\mathrm{o}})$ is isomorphic to $\ker(\partial_a)=I_{\beta}(\Q)/I_{\beta}^{\mathrm{o}}(\Q)$ and the orbit of $\beta$ is in bijection with $\im(\partial_a)=\ker[H^1(\Q,I_{\beta}^{\mathrm{o}})\rightarrow H^1(\Q,I_{\beta})]$. So, we obtain the equality
\[ |\ker[H^1(\Q,I_{\beta}^{\mathrm{o}})\rightarrow H^1(\Q,I_{\beta})]|\cdot [I_{\beta}(\Q):I_{\beta}^{\mathrm{o}}(\Q)]=|\pi_0(I)(\Q)| \]
(In particular, this product quantity is independent of the inner twist of $I$ by a cocycle in $Z^1(\Q,I^{\mathrm{o}})$.)
Our claim follows since $G_{\gamma_0}=I_{\beta}$ for some class $\beta$ (with $I=I_{\mathscr{I}_1,\epsilon_1}$) and $\pi_0(I_{\mathscr{I}_1,\epsilon_1})=\pi_0(G_{\gamma_0})$.

(E2) Since for $(\mathscr{I},\epsilon):=(\mathscr{I}_1^{\beta},\epsilon_1^{\beta})$, $I_{\mathscr{I},\epsilon}^{\mathrm{o}}$ is an inner twist $I_0$ and the Tamagawa number is invariant under inner twist, by (\ref{eq:Tamagawa_number}) one has the equality:
\[\mathrm{vol}(I_{\mathscr{I},\epsilon}^{\mathrm{o}}(\Q)\backslash I_{\mathscr{I},\epsilon}^{\mathrm{o}}(\A_f))= \tau(I_0)\cdot \mathrm{vol}(A_{G}(\R)^{\mathrm{o}}\backslash I_0(\infty)(\R))^{-1}.\]

Hence, the summand in the first line of (\ref{eq:T(m,f)-2}) indexed by an admissible pair $(\mathscr{I},\epsilon)$ and a class $\beta\in \widetilde{\Sha}_G(\Q,I_{\mathscr{I}_1,\epsilon_1})^+$ depends only on the associated (equivalence class of) stable Kottwitz triple $(\gamma_0;\gamma,\delta)$, and thus the second equality holds. 
In the expression of the second line, the index $(\gamma_0;\gamma,\delta)$ originally should run through a set of representatives of \emph{effective} stable Kottwitz triples (i.e. arising from an admissible pair, so having trivial Kottwitz invariant). But, according to Theorem \ref{cor:LR-Satz5.25}, any stable Kottwitz triple $(\gamma_0;\gamma,\delta)$ with trivial Kottwitz invariant is effective if its twisted-orbital integral $\mathrm{TO}_{\delta}(\phi_p)$ is non-zero. Therefore, in this sum we may as well take simply (a set of representatives of) \emph{all} stable Kottwitz triples with trivial Kottwitz invariant. This finishes the proof of (1).

(2) The first claim follows from (1) plus Deligne's conjecture proved by Fujiwara \cite[Cor.5.4.5]{Fujiwara97} (see also \cite{Var07}) and the equality of the compactly supported etale cohomology groups of generic and special fibers \cite[Thm.6.7]{LanStroh18}. When $\sS_K$ is proper or $f^p$ is the identity, we can simply apply the Grothendieck-Lefschetz fixed point formula for lisse sheaves, which gives the stronger control on $m(f^p)$. The properness holds if $G^{\ad}$ is anisotropic since the valuative criterion holds by \cite{Lee12}.
\end{proof}

%\begin{rem} \label{rem:comments_on_Milne92} Milne \cite[Cor.7.10]{Milne92} claimed to have proved this theorem, in the original setting (i.e when the level group $\tilde{\mathbf{K}}_p$ is hyperspecial and $G^{\der}=G^{\uc}$). His proof is incomplete and flawed, in two respects. First, as was mentioned before, he misquotes the definition of \textit{admissible pair} \cite[p.189]{LR87}  (a pair $(\phi,\epsilon)$ is admissible in his sense if and only if it is admissible in the original sense and also $\mathbf{K}_p$-effective in our sense, cf. Remark \ref{rem:admissible_pair}), so his statements in \textit{loc. cit.} using this terminology/notion require critical reading. Secondly and more seriously, in the proof of his Corollary 7.10, he claims that \textit{if a Frobenius triple does not satisfy the condition of (7.5) then it contributes zero to the sum on the right} (a \textit{Frobenius triple} in Milne's work is the same as a Kottwitz triple with trivial Kottwitz invariant). This non-trivial statement (which is simply an \emph{effectivity criterion of Kottwitz triple}) was never justified in \textit{loc. cit.}, nor elsewhere, until our proof of Theorem \ref{thm:LR-Satz5.21} (which is also valid in a more general setting). Also, we remark that for \emph{general} $g$, one needs extra arguments, more than what Milne outlines based on \cite{Kottwitz84b} which was intended mainly for $g=1$ (or at best for those $g$'s lying in a compact open subgroup of $G(\A_f^p)$). \end{rem}

%The proof of this theorem also shows:

\begin{cor} \label{cor:geom_effectivity_of_K-triple} 
Under the same assumption as in Theorem \ref{thm:Kottwitz_formula:Kisin}, for $[\kappa(\wp):\Fp]|n$, a stable Kottwitz triple $(\gamma_0;\gamma,\delta)$ of level $n$ with trivial Kottwitz invariant is \emph{geometrically effective} in the sense that it arises from a $\F_{p^n}$-valued point of $\sS$ (namely from an admissible pair $(\mathscr{I},\mathrm{Fr}_{\mathcal{A}_x/\F_q}^{-1})$ for the relative $q=p^n$-Frobenius endomorphism $\mathrm{Fr}_{\mathcal{A}_x/\F_q}$ of some point $x\in\mathscr{I}$ defined over $\F_q$) if and only if $\mathrm{O}_{\gamma}(f^p)\cdot \mathrm{TO}_{\delta}(\phi_p)$ is non-zero.

An $\R$-elliptic stable conjugacy class of $\gamma_0\in G(\Q)$ arises from an $\F_{p^n}$-valued point of $\sS_{\mathbf{K}}$ if and only if some $G(\A_f^p)$-conjugate of $\gamma_0$ lies in $\mathbf{K}^p$ and there exists $\delta\in G(L_n)$ such that $\gamma_0$ and $\Nm_n(\delta)$ are $G(\mfk)$-conjugate and $\mathrm{TO}_{\delta}(\phi_p)\neq0$.
\end{cor} 

\begin{proof}
In view of Corollary \ref{cor:LR-Satz5.25} (i.e. the effectivity criterion (\ref{eq:E})), the condition is equivalent to the existence of an admissible pair $(\mathscr{I},\epsilon)$ whose associated stable Kottwitz triple is equivalent to $(\gamma_0;\gamma,\delta)$. By Lemma \ref{lem:fixed-pt_subset_of_Frob-Hecke_corr}, $X_p(\mathscr{I},\epsilon)$ and $X^p(\mathscr{I},\epsilon,g)$ (\ref{eq:fixed_pt_set_of_Heck-corresp2}) are non-empty.
It only remains to see that $\epsilon$ is the relative $q=p^n$-Frobenius endomorphism $\mathrm{Fr}_{\mathcal{A}_x/\F_q}$ of some point $x\in\mathscr{I}$ defined over $\F_q$.
It follows from $\mathrm{O}_{\gamma}(f^p)\neq 0$ that $g^p\epsilon {g^p}^{-1}\in \mathbf{K}^p$ for some $g^p\in G(\A_f^p)$; so by Proposition \ref{prop:phi(delta)=gamma_0_up_to_center} (2), we have $\epsilon^{k/n}=\pi_k$ for sufficiently large $k\in\N$ divisible by $n$. We fix a reference point $x_0\in \mathscr{I}$ whose $\hat{\Z}^p$-Tate module and Dieudonne module are stable under $\mathbf{K}^p$ and $\mathbf{K}_p$ respectively.
Choice of $x_p\in G(\mfk)$ satisfying $i_p(\epsilon)x_p=(b\sigma)^nx_p$ and $g^p$ pin down a point $x\in \mathscr{I}$,
and $\epsilon$ together with the equation gives an isomorphism $\mathcal{A}_x \rightarrow \phi^{\ast n}\mathcal{A}_x$ (where $\phi$ is the absolute Frobenius morphism $\mathcal{A}_x\rightarrow \mathcal{A}_x$), namely $\epsilon$ endows $\mathcal{A}_x$ with a $\F_{p^n}$-rational structure for which $\epsilon$ is the inverse of the relative Frobenius endomorphism (cf. \cite[$\S$10]{Kottwitz92}).

For the second statement, by Theorem \ref{thm:LR-Satz5.21}, the condition $\mathrm{TO}_{\delta}(\phi_p)\neq0$ implies existence of admissible Frobenius pair $(h,\epsilon)$ with $\gamma_0$ being stably conjugate to $\epsilon$, which gives rise to an admissible pair $(\mathscr{I},\epsilon)$ (i.e. $\mathscr{I}$ being the reduction of the special point $[h,1]\in Sh_{\mathbf{K}}(G,X)(\Qb)$). Then by the previous argument, the associated stable Kottwitz triple whose rational component is $\gamma_0$ (but the $p$-component does not need to be $\delta$ in the condition) is geometrically effective.
\end{proof}

%%%%%%%%%%%%%%%%%%%%%%%%%%%%%%%%%%%%%%%% 
%%%%%%%%%%%%%%%%%%%%%%%%%%%%%%%%%%%%%%%%

\begin{appendix}

%%%%%%%%%%%%%%%%%%%%%%%%%%%%%%%%%%%%%%%%
%%%%%%%%%%%%%%%%%%%%%%%%%%%%%%%%%%%%%%%%
\section{A lemma on algebraic groups over $p$-adic fields} \label{sec:existence_of_aniostropic_and_unramified_torus}
%%%%%%%%%%%%%%%%%%%%

The goal of this appendix is to give a proof of the following fact, which is perhaps well-known (for example, it is mentioned without proof in the last paragraph on p. 172 of \cite{LR87}). But, since we could not find its proof in literature, we present a proof here.

\begin{lem} \label{lem:existence_of_aniostropic_and_unramified_torus}
Let $H$ be an unramified semi-simple group defined over a $p$-adic field $F$ with ring of integers $\cO_F$.
Then $H$ admits a maximal $F$-torus which is anisotropic and unramified.
If $K$ is a hyperspecial subgroup of $H(F)$, there exists a such maximal $F$-torus $T'$ (i.e. anisotropic and unramified) such that the (unique) hyperspecial subgroup of $T'(F)$ is contained in $K$. 
\end{lem}

\textsc{Proof.} 
We first review some constructions in the general theory of reductive group schemes over $S=\mathrm{Spec}(\cO_F)$.

Let us use $H$ (by abuse of notation) again to denote the reductive group scheme over $S$ which is an integral model of the given $F$-group in question and such that $H(\cO_F)=K$. 
We fix a maximal $S$-torus $T$ of $H$; by definition \cite[3.2.1]{Conrad14}, it is a $S$-torus $T\subset H$ whose every geometric fiber $T_{\overline{s}}$ is a maximal torus of $G_{\overline{s}}$ in the usual sense.
Let $N=N_H(T)$ be the normalizer scheme of $T$ in $H$ \cite[2.1.1]{Conrad14}: it is a (finitely presented) closed subscheme of $G$, smooth over $S$ \cite[2.1.2]{Conrad14}. Also, let $W:=N_H(T)/Z_H(T)$ be the Weyl group scheme of $T$, where $Z_{H}(T)$ is the centralizer scheme of $T$ in $H$ \cite[2.2]{Conrad14}. It is a finite \'etale scheme over $\cO$ \cite[3.2.8]{Conrad14}. 
Let $X$ be the scheme of maximal tori of $H$: this is the scheme representing the functor on $S$-schemes
\begin{equation*}\underline{\mathrm{Tor}}_{H/S}: S'\mapsto \{\text{ maximal tori in }H_{S'}\}\end{equation*}
(which is a sheaf of sets in the fppf topology, \cite[3.2.6]{Conrad14}), and is a smooth affine scheme (\cite[XII, 5.4]{SGA3}).
Moreover, this also represents the quotient sheaf $H/N$ (\cite{Conrad14}, Theorem 2.3.1 and proof of Theorem 3.2.6). Now, to prove the lemma, it is enough to show the existence of an anisotropic maximal torus of the special fibre $H_{\kappa}=H\otimes_{\cO_F}\kappa$, where $\kappa$ is the residue field of $\cO_F$. Indeed, since $\underline{\mathrm{Tor}}_{H/S}$ is smooth over the henselian base $\cO_F$, there exists a maximal torus $\mathcal{T}'$ over $\cO_F$ whose special fibre is anisotropic. But, since $\mathcal{T}'$ is a torus over $\cO_F$, the $\pi_1(S,\overline{\eta})$-module attached to the \'etale sheaf $X^{\ast}(\mathcal{T}')=\underline{\Hom}_{S_{\et}}(\mathcal{T}',\Gm)$ is unramified ($\overline{\eta}$ is a geometric generic point of $S$), 
hence the $F$-rank of $\mathcal{T}'_F$ equals the $\kappa$-rank of $\mathcal{T}'_{\kappa}$. Also, since $T'$ is a subscheme of $H$, the unique hyperspecial subgroup $\mathcal{T}'(\cO_F)$ of $\mathcal{T}(F)$ is contained in $K=H(\cO_F)$.

For the closed point $s$ of $S$ and a (fixed) algebraic closure $\overline{\kappa}=\overline{\kappa(s)}$ of its residue field $\kappa=\kappa(s)$, there exists a canonical diagram
\begin{equation*}\xymatrix{
X(\kappa) \ar[r]^{\delta} \ar[rd]_{\varphi:=j\circ \delta} & H^1(\kappa,N_{\kappa}) \ar[r] \ar[d]^{j} & H^1(\kappa,H_{\kappa}) \\
& H^1(\kappa,W_{\kappa}) \ar[r] & H^1(\kappa,\mathrm{Aut} T_{\kappa})
}\end{equation*}
Here, the upper sequence of pointed sets is exact and $\delta$ is defined as follows: for a maximal $\kappa$-torus $T'$ of $H_{\kappa}$, the image of $T'\in X(\kappa)$ under $\delta$ is the cohomology class of the cocycle \begin{equation*}\sigma\mapsto a_{\sigma}:=g^{-1}\cdot \sigma(g)\in N(\overline{\kappa}),\end{equation*} where $g\in H(\overline{\kappa})$ satisfies $T'=gTg^{-1}$. The lower horizontal map is induced from the injection $W\hra \mathrm{Aut} T$. 
On the other hand, the pointed set $H^1(\kappa,\mathrm{Aut} T)$ classifies the $\kappa$-isomorphism classes of $\kappa$-forms of $T_{\kappa}$, and for any $T'\in X(\kappa)$, its corresponding class in $H^1(\kappa,\mathrm{Aut} T_{\kappa})$ is the image of $[\varphi(T')]$ under the natural map
$H^1(\kappa,W_{\kappa})\rightarrow H^1(\kappa,\mathrm{Aut} T_{\kappa})$, namely, $T'$ is the twist of $T$ by (a cocycle in $W(\overline{\kappa})$ representing) $\varphi(T')$. 
Therefore, to prove the lemma, it will suffice to find a cohomology class $\xi\in H^1(\kappa,W_{\kappa})$ which lies in the image of $\varphi$ and such that the twist ${}_{\xi}T$ of $T_{\kappa}$ by $\xi$ is anisotropic. 
But, since $H^1(\kappa,H_{\kappa})=H^2(\kappa,T_{\kappa})=\{1\}$ (Lang's theorem),
$\varphi$ is always surjective, hence we only need to find a cohomology class $\xi\in H^1(\kappa,W_{\kappa})$ such that the twist ${}_{\xi}T_{\kappa}$ is anisotropic. 

This is also well-known. More explicitly, since every reductive group over a finite field is quasi-split, choose a Borel pair $(T_1,B_1)$ defined over $\kappa$, and let $\Delta$ be the associated set of simple roots of $(H_{\kappa},T_1)$, and $N(T_1)$ the normalizer of $T_1$; thus, the Frobenius automorphism $\sigma$ of $\kappa$ acts on $X^{\ast}(T_1)$, leaving $\Delta$ stable.
Then, what one needs is an element $w\in N(T_1)$ such that $w \sigma\in H(\overline{\kappa})\rtimes\Gal(\overline{\kappa}/\kappa)$ does not have a fixed vector in $X_{\ast}(T)_{\R}$. 
In this case, it is shown in \cite[Lemma 7.4]{Springer74} that if $\{\alpha_1,\cdots,\alpha_m\}$ is a set of representatives of the $\Gal(\overline{\kappa}/\kappa)$-orbits in $\Delta$, 
the product $\omega$ (in any order) of the corresponding simple reflections $r_i\ (i=1,\cdots,m)$, called ``twisted Coxeter element”, has such property. 
$\qed$

%%%%%%%%%%%%%%%%%%%%%%%%%%%%%%%%%%%%%%%%
%%%%%%%%%%%%%%%%%%%%%%%%%%%%%%%%%%%%%%%%
\section{Complexes of tori attached to connected reductive groups.} \label{sec:abelianization_complex}
Here we collect some general facts on certain complexes of tori attached to a connected reductive group and its Levi subgroups.  If stated otherwise, every two-term complex will be concentrated in degree $-1$ and $0$.

For a connected reductive group $H$ over a field $k$, we let $\rho_H:H^{\uc}\rightarrow H$ denote the canonical homomorphism ($H^{\uc}$ being the simply connected covering of $H^{\der}$) and choosing a maximal $k$-torus $T$ of $H$, define a two-term complex of $k$-tori by
\[H_{\mathbf{ab}}:=\rho_H^{-1}(T)\rightarrow T,\]
where $\rho_H^{-1}(T)$ and $T$ are located in degree $-1$ and $0$ respectively, i.e. $H_{\mathbf{ab}}$ is the mapping cone of the morphism $\rho_H^{-1}(T)\rightarrow T$ in the abelian cateogry $\mathcal{CG}_k$ of commutative algebraic $k$-group schemes;
do not confuse this with the torus $H^{\ab}:=H/H^{\der}$, although they become quasi-isomorphic when $H^{\der}=H^{\uc}$.
This complex of tori is also quasi-isomorphic to the crossed module $H^{\uc}\rightarrow H$ (again with $H$ being placed at degree $0$), and the corresponding object in the (bounded) derived category $\mathbb{D}^b(\mathcal{CG}_k)$ depends only on $H$.

Let $G$ be a connected reductive group over a field $k$ and $I$ be a $k$-subgroup of $G$ which is a $\bar{k}$-Levi subgroup.
For $\tilde{I}=\rho_G^{-1}(I)$ (a connected reductive group), there exists a map $\tilde{i}:I^{\uc}\rightarrow \tilde{I}$:
\[\xymatrix{ I^{\uc} \ar[r]^{\tilde{i}} \ar[rd]_{\rho_I} & \tilde{I} \ar@{^(->}[r] \ar[d]^{\rho_G} & G^{\uc} \ar[d]^{\rho_G} \\ & I \ar[r]^{i} & G. }\]
We choose a maximal $k$-torus $T$ of $I$ and set $T_1:=\rho_G^{-1}(T)$, $T_2:=\rho_I^{-1}(T)$; we have a commutative diagram
\[\xymatrix{ T_2 \ar[r]^{\tilde{i}} \ar[rd]_{\rho_I} & T_1 \ar[d]^{\rho_G} \\ & T }\]
The verification of the following facts are easy and thus are left to readers.
\begin{itemize}
\item[(a)] The complex $T_2\rightarrow T_1$ is quasi-isomorphic to the abelianization complex $\tilde{I}_{\mathbf{ab}}:\rho_{\tilde{I}}^{-1}(T_1)\rightarrow T_1$ attached to $\tilde{I}$, where $T_1$ is located at degree $0$, and thus $i:I\hookrightarrow G$ gives rise, in a canonical manner, to a distinguished triangle (in $\mathbb{D}^b(\mathcal{CG}_k)$):
\begin{equation} \label{eq:DT_of_CX_of_tori}
\tilde{I}_{\mathbf{ab}} \rightarrow I_{\mathbf{ab}}  \rightarrow G_{\mathbf{ab}}  \rightarrow \tilde{I}_{\mathbf{ab}}[1].
\end{equation}
\item[(b)] For a two-term complex $T_{\bullet}=T_{-1}\rightarrow T_0$ of $k$-tori, let $\hat{T}_{\bullet}$ denote the complex of $\C$-tori with $\Gamma$-action:
\[\hat{T}_{\bullet}:=(\hat{T}_0\rightarrow \hat{T}_{-1}),\] 
where $\hat{T}_0$ and $\hat{T}_{-1}$ are located in degree $-1$ and $0$ respectively.
Then, for any connected reductive group $H$ and a maximal $k$-torus $T$ of it, there exists an exact sequence of diagonalizable $\C$-groups with $\Gamma$-action
\[1\rightarrow Z(\hat{H}) \rightarrow \hat{T} \rightarrow \hat{\tilde{T}} \rightarrow 1,\] 
with $\tilde{T}:=\rho_H^{-1}(T)$ and $\Gamma=\Gal(\bar{k}/k)$, i.e. 
\begin{equation} \label{eq:center_of_complex_dual}
Z(\hat{H})[1]=\hat{H}_{\mathbf{ab}}.
\end{equation}
\item[(c)] The distinguished triangle (\ref{eq:DT_of_CX_of_tori}) gives rise, by taking complex dual, to a distinguished triangle in the derived category $\mathbb{D}^b(\operatorname{\Gamma-\mathcal{DG}_{\C}})$ of diagonalizable $\C$-groups with $\Gamma$-action:
\begin{align} \label{eq:DT_of_dualCX_of_tori}
\hat{G}_{\mathbf{ab}} \rightarrow \hat{I}_{\mathbf{ab}} \rightarrow \hat{\tilde{I}}_{\mathbf{ab}} \rightarrow \hat{G}_{\mathbf{ab}}[1]. 
\end{align}
\end{itemize}

\end{appendix}

%\textit{Email:} machhama@gmail.com

%\textit{Address:} Institute for Industrial and Applied Mathematics, Chungbuk National University, \\ (28644), 1, Chungdae-ro, Seowon-gu, Cheongju-si, Chungcheongbuk-do, Korea.

\end{document}